\magnification=1200

\hsize=150truemm  \vsize=222truemm
\advance\hoffset by 10truemm \advance\voffset by 12truemm
\topskip=0pt
\parindent=15mm      
\baselineskip=18pt  

\def\gm{{\smc G{\raise-1mm\hbox{u}}\kern-.1667em y
        M{\raise-1mm\hbox{o}}r{\raise-1mm\hbox{e}}l}}
\def\dli{\ifhmode\par\fi\noindent} 
\def\souli#1{$\underline{\hbox{#1}}$} 
\def\up#1{\raise 1ex\hbox{\sevenrm#1}} 

\def\nom#1#2{  
\line{\hbox to 12mm{[#1]\hfill}\smc#2.\hfill}\nobreak }
\def\refa#1#2#3{  
{\leftskip=1mm{#1.{\ \it#2,}\ #3.}\par}
\goodbreak}
\def\refl#1#2{  
{\leftskip=1mm{{\it#1.}\ #2.}\par}
\goodbreak}

	\def\og{\leavevmode\raise.3ex\hbox{$\scriptscriptstyle\langle\!\langle$}}
	\def\fg{\leavevmode\raise.3ex\hbox{$\scriptscriptstyle\,\rangle\!\rangle$}}


	\def\cc#1{\hfill\kern .7em#1\kern .7em\hfill}
	\def\tvi{\vrule height 12pt depth 5pt width 0pt}
	\def\tv{\tvi\vrule}
	\def\traithori{\noalign{\hrule}}



\font\tensmc=cmcsc10                                
\def\smc{\tensmc}     
\font\couv=cmr10 scaled\magstep5
\font\couvs=cmr10 scaled\magstep4
\font\titre=cmr10 scaled\magstep3
\font\titremoins=cmr10 scaled\magstep2
\font\soustitre=cmr10 scaled\magstep1
\font\defpro=cmsl10


\def\relI{\mathrel I} \def\relR{\mathrel R} \def\IR{\relI\joinrel\relR}
\def\relun{\mathrel {\it 1}} \def\relN{\mathrel N} \def\IN{\relI\joinrel\relN}
\def\II{\relun\joinrel\relI}
\def\relmont#1#2{\mathrel{\mathop{\kern 0pt#1}\limits^{\hbox{
$\scriptstyle #2$}}}}


\catcode`\@=11 

\catcode`\;=\active
\def;{\relax\ifhmode\ifdim\lastskip>\z@
 \unskip\fi\kern.2em\fi\string;}

\catcode`\:=\active
\def:{\relax\ifhmode\ifdim\lastskip>\z@\unskip\fi
 \penalty\@M\ \fi\string:}

\catcode`\!=\active
\def!{\relax\ifhmode\ifdim\lastskip>\z@
 \unskip\fi\kern.2em\fi\string!}

\catcode`\?=\active
\def?{\relax\ifhmode\ifdim\lastskip>\z@
 \unskip\fi\kern.2em\fi\string?}

\def\^#1{\if#1i{\accent"5E\i}\else{\accent"5E #1}\fi}
\def\"#1{\if#1i{\accent"7F\i}\else{\accent"7F #1}\fi}

\frenchspacing 
\catcode`\@=12
\pretolerance=500  \tolerance=1000  \brokenpenalty=5000 

{\nopagenumbers
\null\vglue 2cm
\centerline{\titremoins Guy MOREL\ }
\vglue 2cm
\centerline {\couv EXPERTISES :}
\vglue 1cm
\centerline {\couvs PROC\'EDURES STATISTIQUES}
\vglue 0.5cm
\centerline {\couvs D'AIDE \`A LA D\'ECISION}
\vfill\vfill
\centerline {LAST-UNIVERSIT\'E de TOURS}
\eject
}
\pageno=1
\dli {\titremoins\gm} \hfill UFR Arts et Sciences Humaines
\dli\hfill UNIVERSIT\'E Fran\c cois RABELAIS
\dli\hfill 3 rue des TANNEURS 
\dli\hfill 37041 TOURS Cedex
\dli\hfill e-mail : morel@univ-tours.fr
\vfill

\centerline {\couvs EXPERTISES :}
\vglue 1cm
\centerline {\titre PROC\'EDURES STATISTIQUES}
\vglue 0.5cm
\centerline {\titre D'AIDE \`A LA D\'ECISION}
\vfill\vfill
\centerline {\titremoins Avril 1997}
\eject

\null\vglue 3cm
\centerline{\soustitre SOMMAIRE}
\vglue 2cm
\dli 1--INTRODUCTION.

\dli 2--CHOIX ENTRE DEUX PROBABILIT\'ES.

2.1 D\'efinition des experts. p. 10

2.2 Fonctions de test simples. p. 13

2.3 Ensemble des experts. p. 16

2.4 Votes des experts. p. 19

2.5 Vote pond\'er\'e. p. 23

\dli 3--R\`EGLES DE D\'ECISION DE BOL'SHEV.

3.1 D\'efinitions et propri\'et\'es. p. 25

3.2 Votes des experts et r\`egles de Bol'shev. p. 31

\dli 4--CHOIX ENTRE DEUX HYPOTH\`ESES STABLES.

4.1 D\'efinitions. p. 37

4.2 Ensemble des experts. p. 41

4.3 Votes des experts. p. 56

\dli 5--MOD\`ELES \`A RAPPORT DE VRAISEMBLANCE MONOTONE.

5.1 Hypoth\`eses unilat\'erales. p. 66

5.2 Votes compatibles sur une famille d'hypoth\`eses unilat\'erales. p. 78

5.3 Hypoth\`eses bilat\'erales. p. 106

5.4 Exemples. p. 122
\vfill\eject

\dli 6--HYPOTH\`ESES STABLES ET PARAM\`ETRES FANT\^OMES.

6.1 Introduction. p. 131

6.2 Mod\`eles exponentiels. p. 137

\dli ANNEXE I. p. 143

\dli ANNEXE II. p. 145

\dli ANNEXE III. p. 148

\dli ANNEXE IV. p. 161

\dli BIBLIOGRAPHIE. p. 168

\dli TABLE des d\'efinitions. p. 172

\dli TABLE des propositions. p. 173
\vfill\eject

{\parindent=-10mm\soustitre 1--INTRODUCTION.}
\bigskip

Choisir entre deux hypoth\`eses est un des probl\`emes fondateurs de la statistique math\'ematique. Il n'existe pas de solution miracle, suivant la mod\'e\-lisation choisie pour cette prise de d\'ecision on obtient des solutions diff\'erentes. Le choix entre ces diff\'erents cadres d\'ecisionnels est un probl\`eme externe \`a la statistique math\'ematique, il d\'epend du champ d'application consid\'er\'e. La proc\'edure choisie fournit une aide \`a la d\'ecision qu'il est difficile d'interpr\'eter sans tenir compte des crit\`eres qui ont structur\'e sa production. Pourtant, l'utilisateur interpr\`ete souvent ses r\'esultats sans tenir compte de  la mod\'elisation sous jacente avec laquelle il travaille (cf. [Wan.]). On peut m\^eme se demander si une m\'ethode statistique n'est pas d'autant plus populaire qu'il est possible d'oublier les bases de sa construction dans sa mise en \oe uvre. Les tests statistiques finissent par dire oui ou non \`a des seuils conventionnels sans r\'ef\'erence \`a la dissym\'etrie de traitement entre les deux hypoth\`eses. De plus, quand l'hypoth\`ese $H_0$ traduit la notion d'effet n\'egligeable, elle prend souvent la forme d'un effet parfaitement nul. Il devient alors paradoxalement plus facile de prendre une d\'ecision \`a partir d'un \'echantillon, qui est parfois simplement une sous population, qu'\`a partir d'un travail exhaustif s'il \'etait possible.
Les utilisateurs semblent peu friands de r\'eponses qui laissent place au jugement comme une probabilit\'e a post\'eriori sur l'espace des param\`etres ou des d\'ecisions. Et dans ce cas, il est rare de trouver une probabilit\'e a priori li\'ee au champ d'\'etude, elle est le plus souvent non informative.
Le point de vue que nous allons d\'evelopper n'emp\^eche \'evidemment pas les pratiques magiques. Nous avons simplement essay\'e de rendre difficile l'oubli du caract\`ere al\'eatoire des proc\'edures statistiques, en optant pour une aide sous la forme d'une probabilit\'e sur l'espace des d\'ecisions. Ce type d'inf\'erence a \'et\'e d\'efendu dans de nombreux travaux (cf. [KroM]).

L'\'etude d'un probl\`eme de d\'ecision statistique passe par la donn\'ee de crit\`eres de s\'election entre les diff\'erentes r\`egles de d\'ecision consid\'er\'ees.
Classiquement on commence par se donner une fonction de perte et on compare les proc\'edures de d\'ecision \`a partir des fonctions de risque correspondantes.
Il est rare qu'un choix unique de ce crit\`ere s'impose, bien que l'\'etude de certaines fonctions de perte soit privil\'egi\'ee, par exemple la perte quadratique dans le cadre de l'estimation ou le risque de se tromper pour le choix entre deux hypoth\`eses. M\^eme si ce dernier choix para\^{\i}t assez ``naturel", d'autres pertes sont possibles, en particulier si on regarde ce probl\`eme de choix entre deux hypoth\`eses comme un probl\`eme d'estimation (cf. [HwaC], [Rob.] p. 186). Les fonctions de risque permettent d'introduire un pr\'eordre partiel sur les r\`egles de d\'ecision et ainsi de s\'electionner les r\`egles admissibles. Cette premi\`ere s\'election s'obtient en comparant les r\`egles entre elles, elle ne d\'erive pas de propri\'et\'es intrins\`eques. Nous allons dans cette \'etude partir de telles propri\'et\'es pour d\'efinir les ``bonnes" r\`egles de d\'ecision, celles que nous appellerons ``experts''. Une r\`egle pourra \^etre d\'eclar\'ee expert sans avoir \`a la comparer \`a l'ensemble des autres r\`egles. La propri\'et\'e de base que nous imposons aux experts est simplement que le diagnostic $d$, $\theta$ appartient \`a $\Theta_d$, doit \^etre plus probable quand $\theta$ appartient \`a $\Theta_d$ que lorsque $\theta$ n'y appartient pas. C'est une notion de r\`egle non biais\'ee, mais pour que cela suppose une connaissance fine du mod\`ele nous imposons que cette propri\'et\'e soit aussi v\'erifi\'ee conditionnellement \`a tout \'ev\'enement non n\'egligeable, et pas simplement globalement comme c'est g\'en\'eralement le cas. Quand ils existent les experts ne sont ainsi pas trop d\'ependants de l'espace des r\'ealisations et de l'espace des param\`etres choisis. Cette d\'efinition des experts sera pr\'ecis\'ee par la suite, nous restons dans cette introduction au niveau des id\'ees directrices. L'exemple simple de la famille des lois normales de moyenne inconnue $\theta\in\IR$ et de variance $1$ nous servira \`a imager notre propos. Lorsque les hypoth\`eses sont unilat\'erales les experts sont les r\`egles admissibles pour le risque de se tromper. Mais pour les hypoth\`eses bilat\'erales les experts se r\'eduisent aux r\`egles triviales qui d\'ecident toujours la m\^eme hypoth\`ese. On peut dire que c'est \`a peu pr\`es ce qui se passe dans tout mod\`ele statistique \`a rapport de vraisemblance monotone.
Nous allons \^etre amen\'e \`a traiter deux probl\`emes  diff\'erents : un trop plein d'experts d'un c\^ot\'e et une p\'enurie d'experts de l'autre. 

Commen\c cons par le premier cas : que faire de tous ces experts ? 
Le statisticien se trouve souvent dans une situation semblable, par exemple avec un ensemble de r\`egles admissibles. Il se donne g\'en\'eralement des crit\`eres suppl\'ementaires pour essayer de s\'electionner une r\`egle de d\'ecision : recherche d'un test U.P.P. (uniform\'ement plus puissant), d'un test sans biais U.P.P., d'un test invariant U.P.P., d'une r\`egle de Bayes, d'une r\`egle minimax, etc. Pour nous cela reviendrait \`a chercher un expert qui soit plus expert que les autres. Voil\`a qui sonne \'etrangement car le vocabulaire choisi traduit le fait que nous ne voulons pas s\'electionner un expert. Ceci nous obligerait \`a fournir \`a l'utilisateur une r\'eponse en tout ou rien, ce qui finit par cacher son caract\`ere al\'eatoire. Nous pr\'ef\'erons consid\'erer nos experts comme \'egaux en droit, les faire voter et fournir \`a l'utilisateur le r\'esultat de ce vote, c'est-\`a-dire une probabilit\'e sur l'espace des deux d\'ecisions possibles. Il se trouvera dans une situation semblable \`a celle obtenue quand on transforme une probabilit\'e a post\'eriori en une probabilit\'e sur l'ensemble des d\'ecisions. Pour faire voter nos experts il faut d\'efinir une probabilit\'e sur l'ensemble des experts. Nous en avons associ\'e une \`a chaque probabilit\'e du mod\`ele en essayant d'accorder d'autant plus de poids \`a un ensemble d'experts qu'il est form\'e d'experts donnant des r\'esultats diff\'erents sous cette probabilit\'e. Nous avons alors autant de votes que de valeurs du param\`etre et il y a bien des mani\`eres de les synth\'etiser. 
Reprenons le cas d'hypoth\`eses unilat\'erales dans l'exemple des lois 
$(N(\theta,1))_{\theta\in\IR}$. Pour le choix entre $H : \theta\leq\theta_0$ et $H' : \theta>\theta_0$, on peut dire que le vote associ\'e \`a la probabilit\'e fronti\`ere $N(\theta_0,1)$ est une solution de type minimax, sous chaque hypoth\`ese elle choisit le plus favorable des votes d\'efavorables. C'est un vote neutre au point fronti\`ere, c'est-\`a-dire que lorsque $\theta_0$ est la valeur du param\`etre les votes en faveur de $H$ (resp. $H'$) ont une fr\'equence moyenne \'egale \`a ${1\over 2}$. Pour chaque observation on associe alors \`a la d\'ecision $H : \theta\leq\theta_0$ une fr\'equence de votes qui se trouve \^etre le seuil minimum de rejet du test de $H$ contre $H'$. La fr\'equence des votes en faveur de $H'$ est aussi un seuil minimum de rejet, celui du test de $H'$ contre $H$. Ceci donne un sens nouveau \`a la notion de p-value qui correspond mieux \`a son utilisation courante et peu orthodoxe, comme indice de fiabilit\'e du choix d'une des hypoth\`eses. Cette probabilit\'e sur l'espace des d\'ecisions est obtenue sans l'introduction d'une loi a priori non informative, comme dans le cadre bay\'esien. Il est cependant possible de synth\'etiser l'ensemble des votes en faisant une moyenne \`a partir d'une probabilit\'e sur l'ensemble des param\`etres et donc de prendre en compte des informations a priori. Cette probabilit\'e sert de pond\'eration, nous ne parlerons pas de probabilit\'e a priori car elle ne se comporte pas de fa\c con semblable. Par exemple, en analyse bay\'esienne les p-values pr\'ec\'edentes se trouvent en prenant pour loi a priori la loi impropre non informative d\'efinie par la mesure de Lebesgue, alors que dans notre cas il faut prendre la masse de Dirac en $\theta_0$.

Analysons maintenant le cas o\`u il n'y a que les experts triviaux. C'est ce qui se passe pour les hypoth\`eses bilat\'erales dans l'exemple des lois 
$(N(\theta,1))_{\theta\in\IR}$. On peut penser diminuer les contraintes impos\'ees par la d\'efinition des experts en travaillant sur une sous tribu, ce qui diminue l'ensemble des \'ev\'enements sur lesquels la propri\'et\'e de base des experts doit s'appliquer conditionnellement. Si l'on consid\`ere les hypoth\`eses bilat\'erales $H : \theta\in[-\theta_0,+\theta_0]$ et $H' : \theta\in]-\infty,-\theta_0[\cup]+\theta_0,+\infty[$ on peut sym\'etriser le probl\`eme en se restreignant \`a la tribu engendr\'ee par les intervalles de m\^eme probabilit\'e sous $\theta_0$ et sous $-\theta_0$, c'est-\`a-dire les intervalles sym\'etriques par rapport \`a $0$. Il existe alors des experts, c'est comme si on travaillait sur le mod\`ele image du pr\'ec\'edent par la statistique valeur absolue, $\theta$ et $-\theta$ \'etant confondus car de m\^eme image. Le r\'esultat du vote sous  $\theta_0$ ou $-\theta_0$ fournit le seuil minimum de rejet du test sans biais de $H$ contre $H'$ comme fr\'equence des votes en faveur de $H$. La fr\'equence des votes en faveur de $H'$ est quant \`a elle le seuil minimum de rejet du test de $H'$ contre $H$. Cette solution repose sur une sym\'etrie du probl\`eme de d\'ecision par rapport \`a $\theta=0$, elle est justifi\'ee s'il est \'equivalent pour l'interpr\'etation d'avoir $\theta<-\theta_0$ ou $\theta>+\theta_0$. On casse ainsi la structure d'ordre classique sur $\Theta=\IR$ pour la remplacer par le pr\'eordre induit par la distance \`a $0$. Dans bien des cas la structure de d\'epart garde un sens m\^eme pour des hypoth\`eses bilat\'erales, le choix principal est entre $H$ et $H'$ mais $\theta<-\theta_0$ et $\theta>+\theta_0$ ne signifient pas la m\^eme chose. Il est alors int\'eressant que le vote par rapport \`a ces hypoth\`eses proviennent de votes pour les hypoth\`eses unilat\'erales d\'efinies par le point fronti\`ere $-\theta_0$ d'une part et celles d\'efinies par $\theta_0$ d'autre part. Nous avons pour cela introduit la notion de votes compatibles sur une famille d'hypoth\`eses unilat\'erales. Dans le cas pr\'ec\'edent il suffit que pour toute observation la fr\'equence des votes en faveur de la d\'ecision $\theta<-\theta_0$ soit inf\'erieure \`a la fr\'equence des votes en faveur de la d\'ecision $\theta<+\theta_0$. On peut alors d\'efinir une probabilit\'e sur les trois \'ev\'enements $]-\infty,-\theta_0[$, $[-\theta_0,+\theta_0]$ et $]+\theta_0,+\infty[$ donc sur les hypoth\`eses $H$ et $H'$. Cette mani\`ere de faire a l'avantage d'imposer une coh\'erence entre les solutions de probl\`emes de d\'ecision qui reposent sur une m\^eme structuration de l'espace des param\`etres par rapport aux interpr\'etations possibles.

La notion de votes compatibles prend tout son int\'er\^et  quand l'analyse du param\`etre est structur\'ee par un ordre. L'ensemble des hypoth\`eses unilat\'erales joue alors un r\^ole fondamental. Dans un mod\`ele \`a rapport de vraisemblance monotone, d\'efinir des votes compatibles sur cette famille d'hypoth\`eses revient \`a d\'efinir une probabilit\'e sur l'espace des param\`etres muni de la tribu de l'ordre. Quand $\Theta$ est un intervalle de $\IR$ muni de l'ordre usuel on d\'efinit g\'en\'eralement des votes compatibles en associant \`a chaque probl\`eme de d\'ecision unilat\'eral le vote d\'efini par la probabilit\'e fronti\`ere entre les deux hypoth\`eses. Dans l'exemple des lois 
$(N(\theta,1))_{\theta\in\IR}$ on obtient sur les bor\'eliens de $\Theta=\IR$ la loi normale centr\'ee sur l'observation obtenue et de variance $1$.
On peut aussi obtenir des votes compatibles en utilisant des pond\'erations sur l'espace des param\`etres. Ceci est particuli\`erement int\'eressant dans le cas o\`u on a une information  a priori \`a introduire. La solution du vote fronti\`ere ne suppose, elle, aucune information a priori, elle donne d'ailleurs souvent une loi sur les param\`etres qui peut \^etre consid\'er\'ee, dans le cadre bay\'esien, comme la loi a post\'eriori  d'une loi a priori non informative.
\dli La probabilisation de l'espace des param\`etres muni de la tribu de l'ordre permet de d\'efinir une probabilit\'e sur n'importe quelles hypoth\`eses, un peu comme une probabilit\'e \`a post\'eriori. Il est cependant pr\'ef\'erable que les hypoth\`eses aient une interpr\'etation li\'ee \`a l'ordre qui structure $\Theta$. 
Lorsque l'espace des param\`etres n'est pas ordonn\'e, mais est un produit d'espaces ordonn\'es il est souvent possible d'obtenir une probabilisation de cet espace. Lorsque chaque espace ordonn\'e peut \^etre probabilis\'e de fa\c con ind\'ependante on prendra la loi produit. Dans le cas contraire il faudra hi\'erarchiser les ordres et construire des probabilit\'es de transition. Pour un param\`etre fant\^ome cela revient \`a trouver une solution au probl\`eme de d\'ecision d\'efini en fixant ce param\`etre, et \`a obtenir une probabilisation conditionnelle sur l'espace du param\`etre fant\^ome.

\vfill\eject

{\parindent=-10mm\soustitre 2--CHOIX ENTRE DEUX PROBABILIT\'ES.}
\bigskip
{\parindent=-5mm 2.1 D\'EFINITION DES EXPERTS.}
\medskip
Dans ce cas particulier l'espace des param\`etres se r\'eduit \`a deux \'el\'ements $\Theta = \{0,1\}$. On peut toujours supposer que les deux probabilit\'es
$P_\theta$ admettent une densit\'e $p_\theta$ par rapport \`a une mesure $\mu$
sur $(\Omega ,{\cal A})$ (cf. [Leh.] p. 74). Il y a deux d\'ecisions possibles $d=0$ et $d=1$ qui correspondent respectivement au choix de $\theta =0$ et de $\theta =1$.
Une r\`egle de d\'ecision est alors d\'efinie par une statistique r\'eelle $\phi$ \`a valeurs dans $D=\{0,1\}$. Pour la r\'ealisation $\omega\in\Omega$, on d\'ecide $d=\phi(\omega)$.

Comme nous l'avons expliqu\'e dans l'introduction
nous allons essayer de d\'efinir les experts en imposant des propri\'et\'es qui ne font pas intervenir des comparaisons entre r\`egles de d\'ecision. Ces propri\'et\'es doivent pouvoir se v\'erifier ou s'infirmer en ne consid\'erant que la r\`egle de d\'ecision qui postule au label d'expert.

Nous allons commencer par ce que l'on ne veut pas, donc par une propri\'et\'e qui ne peut pas \^etre l'apanage d'un expert. Face \`a une r\'ealisation $\omega$ un expert doit diagnostiquer entre $\theta=0$ et $\theta=1$, on ne peut pas accepter que la probabilit\'e qu'il d\'ecide $d=0$ quand $\theta=0$ soit plus faible que lorsque $\theta=1$. Il est bien s\^ur \'equivalent de dire que l'on ne veut pas que la probabilit\'e de d\'ecider $d=1$ soit plus faible pour $\theta=1$ que pour $\theta=0$. Ceci revient \`a dire qu'un diagnostic $d$ doit avoir plus de chance d'\^etre produit quand il est bon que quand il est mauvais. Sous cette condition un expert doit v\'erifier :
$$P_1(\{\phi =1\})=E_1(\phi)\geq P_0(\{\phi =1\})=E_0(\phi)$$
($E_\theta$ est l'op\'erateur esp\'erance par rapport \`a $P_\theta$).

C'est la notion classique de sans biais. Cette propri\'et\'e est globale, elle repose sur les moyennes de $\phi$. Elle n'est pas tr\`es contraignante localement. Par exemple si elle est strictement r\'ealis\'ee, $E_1(\phi) - E_0(\phi)=q>0$, on peut inverser la d\'ecision sur tout \'ev\'enement $A$ inclus dans $\{\phi=0\}$ (resp. $\{\phi=1\}$) d\`es que  $P_0(A)\leq q$ (resp. $P_1(A)\leq q)$. Ceci nous montre que la r\'ealisation de cette propri\'et\'e ne suppose pas une connaissance fine de la structure du mod\`ele. Pour cela il faut essayer de l'imposer conditionnellement \`a tout \'ev\'enement $C$. Si $C$ est de probabilit\'e non nulle sous $P_0$ et $P_1$ on peut d\'efinir un mod\`ele statistique conditionnel ; 
\dli $(\Omega,{\cal A})$ est muni des probabilit\'es $P_\theta^C$ d\'efinies par :
$$\forall\theta\in\{0,1\}\qquad\forall B\in {\cal A}\qquad
P_\theta^C(B)=P_\theta(B\cap C)/P_\theta(C).$$
Dans ce nouveau mod\`ele un expert doit alors v\'erifier :
$$P_1^C(\{\phi =1\})=E_1^C(\phi)\geq P_0^C(\{\phi =1\})=E_0^C(\phi)$$ 
($E_\theta^C$ est l'op\'erateur esp\'erance par rapport \`a $P_\theta^C$).
Lorsque $P_0(C)$ ou $P_1(C)$ est nulle on ne peut plus d\'efinir de mod\`ele statistique conditionnel. Mais dans ce cas, si l'une des deux probabilit\'es est non nulle, par exemple $P_\theta(C)\not= 0$, un expert doit d\'ecider $d=\theta$ sur $P_\theta$ presque s\^urement tout $C$. Ces consid\'erations nous conduisent \`a la d\'efinition suivante :
\medskip
{\bf D\'efinition 2.1.1}
\medskip
\medskip
\moveleft 10.4pt\hbox{\vrule\kern 10pt\vbox{\defpro

 Soit le mod\`ele statistique
$\bigl(\Omega ,{\cal A},(P_\theta)_{\theta\in\{0,1\}} \bigr)$. Un expert du choix entre $P_0$ et $P_1$ est une r\`egle de d\'ecision $\phi : (\Omega,{\cal A})\rightarrow\{0,1\}$ qui v\'erifie pour tout \'ev\'enement $C$ :
\dli $P_1(C\cap\{\phi =0\})=0$ si $P_0(C)=0$
\dli $P_0(C\cap\{\phi =1\})=0$ si $P_1(C)=0$
\dli $E_1^C(\phi)\geq E_0^C(\phi)$ si $P_0(C)\not= 0$ et $P_1(C)\not= 0$ .
}}\medskip

Le probl\`eme du choix entre deux probabilit\'es, $P_0$ et $P_1$, est trait\'e de deux mani\`eres dans la th\'eorie des tests, suivant qu'on privil\'egie $P_0$ ou $P_1$. On parle du test de $P_0$ contre $P_1$ ou du test de $P_1$ contre $P_0$. Ce sont les tests de Neyman qui sont s\'electionn\'es.
\dli Rappelons qu'on appelle test de Neyman de $P_1$ contre $P_0$  (cf. [Mon.2] p. 135) une fonction de test pour laquelle il existe k dans $\overline{\IR^+}$ tel que : 
\dli $\phi(\omega)\relmont{=}{p.s.}1$ sur $\{\omega\in\Omega\, ;
\, p_0<kp_1\}$ et
\dli $\phi(\omega)\relmont{=}{p.s.}0$ sur $\{\omega\in\Omega\, ;
\, p_0>kp_1\}$
\dli (avec la convention $\infty\times 0=0$, la notation $\relmont{=}{p.s.}$
signifiant : $P_0$ et $P_1$ presque s\^urement).
\dli Ce sont les tests de Neyman communs aux deux probl\`emes de test :
$P_0$ contre $P_1$ et $P_1$ contre $P_0$, qui vont jouer un r\^ole
fondamental dans notre probl\`eme.
Nous allons commencer par mettre en place cet outil.

\bigskip
\vfill\eject

{\parindent=-5mm 2.2 FONCTIONS DE TEST SIMPLES.}
\medskip
Consid\'erons les \'ev\'enements suivants :
\halign{$#$&#$\ =\ $&$#$\hfill\cr
S_0&&\{\omega\in\Omega\, ;\, p_0(\omega)>0\},\cr
S_1&&\{\omega\in\Omega\, ;\, p_1(\omega)>0\},\cr
\Omega_k&&\{\omega\in\Omega\, ;\, p_0(\omega)=k.p_1(\omega)\}\bigcap
\{S_0\cup S_1\}$ pour $k\in\IR^+,\cr
\Omega_\infty&&\{\omega\in\Omega\, ;\, p_1(\omega)=0\} = S_1^c.\cr}

La famille $\{\Omega_k\}_{k\in\overline{\IR^+}}$ forme une partition de
$\Omega$. On vient ainsi de d\'efinir une statistique $K$ \`a valeurs dans $\overline{\IR^+}$ qui peut \^etre consid\'er\'ee comme le rapport entre les densit\'es $p_0$ et $p_1$ (la forme ind\'etermin\'ee $0/0$ prenant ici la valeur $+\,\infty$). Cette statistique est unique, $P_0$ et $P_1$ presque s\^urement, c'est-\`a-dire que presque s\^urement, elle ne d\'epend pas de la mesure et des densit\'es choisies pour exprimer $P_0$ et $P_1$ (voir annexe I ).

On peut d\'efinir les tests de Neyman de $P_1$ contre $P_0$ (voir 2.1) en utilisant la statistique $K$. $\phi$ est un de ces tests si 
$\phi\relmont{=}{p.s.}\II_{\{K<+\,\infty\}}$ ou si il existe $k\in\IR^+$ tel que $\II_{\{K<k\}}\relmont{\leq}{p.s.}\phi\relmont{\leq}{p.s.}\II_{\{K\leq k\}}$. Pour $k=0$ seul le test $\phi\relmont{=}{p.s.}\II_{\{K\leq 0\}}$ est int\'eressant puisque ceux qui ne valent pas presque s\^urement $1$ sur l'\'ev\'enement $\{K\leq 0\}$ sont de m\^eme puissance ($P_0(\{\phi=0\})=1$) mais de seuil plus grand. Il reste alors les tests de Neyman communs aux deux probl\`emes de test. Ils sont d\'efinis par l'une des propri\'et\'es : $\phi\relmont{=}{p.s.}\II_{\{K\leq 0\}}$,
$\II_{\{K<k\}}\relmont{\leq}{p.s.}\phi\relmont{\leq}{p.s.}\II_{\{K\leq k\}}$ et $k\in\IR_*^+$,
$\phi\relmont{=}{p.s.}\II_{\{K<+\,\infty\}}$.
Les fonctions indicatrices qui permettent de d\'efinir ces tests de Neyman vont jouer un r\^ole essentiel.
\medskip
{\bf D\'efinition 2.2.1}
\medskip
\medskip
\moveleft 10.4pt\hbox{\vrule\kern 10pt\vbox{\defpro

Soit $(k,\beta)\in\Bigl\{(0,1) ; \{(k,\beta)\}_{k\in\IR_*^+,\,\beta\in\{0,1\}}; (\infty ,0)\Bigr\}$. La fonction de test simple associ\'ee \`a $(k,\beta)$ est d\'efinie par : 
$\phi_{(k,\beta)}=\II_{\{K<k\}} + \beta . \II_{\{K=k\}}$
(la statistique $K$ est le rapport des deux densit\'es $p_0$ et $p_1$).
}}\medskip

Consid\'erons les tests de Neyman compris presque s\^urement entre les deux fonctions de test simples $\phi_{(k,0)}$ et $\phi_{(k,1)}$ pour $k\in\IR_*^+$. Ils sont presque s\^urement \'egaux si $P_0(\{K=k\})=0$ et
$P_1(\{K=k\})=0$. Dans le cas contraire, ces deux probabilit\'es sont non nulles puisque $k\in\IR_*^+$. Elles d\'efinissent sur $\{K=k\}$ un mod\`ele conditionnel compos\'e de deux probabilit\'es identiques. Ce qui diff\'erencie principalement les tests de Neyman $\phi$ compris entre $\phi_{(k,0)}$ et $\phi_{(k,1)}$ c'est la probabilit\'e de d\'ecider $d=1$ quand l'\'ev\'enement $\{K=k\}$ est r\'ealis\'e :
$\beta =P_\theta(\{\phi =1\}\cap\{K=k\})/P_\theta(\{K=k\})$. Pour tout $\beta$ dans l'intervalle $[0,1]$ on peut trouver une r\`egle de d\'ecision $\phi : \Omega\rightarrow\{0,1\}$ ayant la propri\'et\'e pr\'ec\'edente, si l'\'ev\'enement $\{K=k\}$ contient des \'ev\'enements  de probabilit\'e conditionnelle $\beta$. Ceci peut \^etre impossible, par exemple si on prend pour mod\`ele le mod\`ele image de la statistique exhaustive $K$.
Afin de rem\'edier \`a cet inconv\'enient, qui peut introduire des discontinuit\'es dans le traitement du probl\`eme, on utilise souvent des fonctions de test al\'eatoires : $\phi$ est \`a valeur dans
$[0,1]$ et $\phi(\omega)$ repr\'esente la probabilit\'e de d\'ecider $d=1$. Pour prendre une d\'ecision il faut alors faire intervenir un al\'ea sur $D=\{0,1\}$ qui donne \`a $d=1$ la probabilit\'e $\phi(\omega)$ d'\^etre obtenu. Si on introduit cet al\'ea dans le mod\`ele on obtient  de nouveau des r\`egles de d\'ecision d\'eterministes, c'est-\`a-dire \`a valeur dans $\{0,1\}$. Stevens [Ste.] par exemple utilise cet artifice pour am\'eliorer les intervalles de confiance dans le cas du param\`etre d'une loi binomiale. 
Nous reviendrons sur cette difficult\'e de la r\'epartition d'une masse en
$k\in\IR_*^+$ entre les d\'ecisions $d=1$ et $d=0$. Elle dispara\^{\i}tra quand nous ferons voter nos experts sans avoir \`a introduire les fonctions de test al\'eatoires ou \`a utiliser un surmod\`ele contenant un al\'ea.
Nous pourrons ne faire voter que les fonctions de test simples.

Consid\'erons l'ensemble des fonctions de test simples :
\dli$\Phi_s=\Bigl\{\phi_{(0,1)}\, ;\,\{\phi_{(k,\beta)}\}_{k\in\IR_*^+,\,\beta\in \{0,1\}}\, ;\,\phi_{(\infty,0)}\Bigr\}$.
\dli Il est totalement ordonn\'e par la relation d'ordre
partiel usuelle sur les fonctions. Aussi nous le noterons parfois
$[\phi_{(0,1)},\phi_{(\infty,0)}]$. L'ordre obtenu co\"{\i}ncide avec celui induit par l'ordre lexicographique sur les couples $(k,\beta)$ :
$$(k,\beta)\leq (k',\beta')\Longleftrightarrow\left\{\matrix{
k<k' &    &\cr
ou   &    &\cr
k=k' & et & \beta\leq\beta' .\cr
}\right.$$
Il permet de d\'efinir les bornes sup\'erieure et inf\'erieure
d'un sous ensemble de fonctions de test simples.

L'op\'erateur $E_\theta$, esp\'erance par rapport \`a $P_\theta$,
est croissant sur $\Phi_s=[\phi_{(0,1)},\phi_{(\infty,0)}]$.
\dli $E_\theta(\phi_{(k,\beta)}) \, = \,
\int_\Omega\phi_{(k,\beta)}\,dP_\theta\,=\,P_\theta(\{\phi_{(k,\beta)}=1\})$.
$$\vbox{
\offinterlineskip
\halign{$\cc{#}$&\tv#&&$\cc{#}$\cr
(k,\beta)&&(0,1)&\nearrow&(\infty,0)\cr
\traithori
E_0(\phi_{(k,\beta)})&&0&\nearrow&1-P_0(\Omega_\infty)\cr
\traithori
E_1(\phi_{(k,\beta)})&&P_1(\Omega_0)&\nearrow&1\cr
}}$$
On a toujours $E_1(\phi_{(k,\beta)}) - E_0(\phi_{(k,\beta)}) \, = \,
\int_\Omega\phi_{(k,\beta)}(p_1-p_0)\,d\mu \,\geq\, 0$ ; en effet
$\phi_{(k,\beta)}(p_1-p_0) \geq 0$ pour $k\leq 1$ et pour $k > 1$
$\int_\Omega(\phi_{(k,\beta)}-\phi_{(1,1)})(p_1-p_0)\,d\mu \,\geq\,
\int_\Omega(1-\phi_{(1,1)})(p_1-p_0)\,d\mu \, =\,
\int_\Omega -\phi_{(1,1)}(p_1-p_0)\,d\mu$.

\bigskip
\vfill\eject

{\parindent=-5mm 2.3 ENSEMBLE DES EXPERTS.}
\medskip
Nous allons d\'emontrer que les experts sont les tests de Neyman d\'eterministes communs aux tests de $P_0$ contre $P_1$ et de $P_1$ contre $P_0$.
Ce que nous allons \'ecrire en utilisant les  fonctions de test simples, comme en 2.2.

\medskip
{\bf Proposition 2.3.1}
\medskip
\medskip
\moveleft 10.4pt\hbox{\vrule\kern 10pt\vbox{\defpro

Les experts sont les r\`egles de d\'ecision $\phi : (\Omega,{\cal A})\rightarrow\{0,1\}$ presque s\^urement ordonnables dans l'ensemble des fonctions de test simples :  $\Phi_s=[\phi_{(0,1)},\phi_{(\infty,0)}]$. C'est-\`a-dire qu'il existe des fonctions de test simples, \'egales ou cons\'ecutives, encadrant presque s\^urement $\phi$ : $\phi_{(k,\beta)}\relmont{\leq}{p.s.}\phi\relmont{\leq}{p.s.}\phi_{(k,\beta')}$. 
}}\medskip

Cette propri\'et\'e est v\'erifi\'ee s'il existe $\phi_{(k,\beta)}$ presque s\^urement \'egale \`a $\phi$ ou s'il existe $k\in\IR_*^+$ tel que $\phi_{(k,0)}\relmont{<}{p.s.}\phi\relmont{<}{p.s.}\phi_{(k,1)}$ (la relation $\relmont{<}{p.s.}$ signifiant $\relmont{\leq}{p.s.}$ sans avoir $\relmont{=}{p.s.}$).
\medskip
{\leftskip=15mm \dli {\bf D\'emonstration}
\medskip
{\parindent=-10mm I ----- Condition n\'ecessaire.}

Soit $\phi$ un expert. On cherche $k\in\overline{\IR^+}$, $\beta$ et $\beta'$ dans $\{0,1\}$ tels que :
$\phi_{(k,\beta)}\relmont{\leq}{p.s.}\phi\relmont{\leq}{p.s.}\phi_{(k,\beta')}$. 

\dli Par construction de $\Omega_0$ et $\Omega_\infty$ (voir 2.2) on a $P_0(\Omega_0)=0$ et $P_1(\Omega_\infty)=0$. La d\'efinition des experts entra\^{\i}ne alors :
$P_1(\Omega_0\cap\{\phi=0\})=0$ et $P_0(\Omega_\infty\cap\{\phi=1\})=0$.
Ce qui peut s'\'ecrire $\phi_{(0,1)}\relmont{\leq}{p.s.}\phi\relmont{\leq}{p.s.}\phi_{(\infty,0)}$. 

Nous pouvons maintenant d\'efinir 
\dli $\phi_{(k_0,\beta_0)}\,=\,sup\bigl\{\phi'\in[\phi_{(0,1)},\phi_
{(\infty,0)}]\, ;\,\phi'\relmont{\leq}{p.s.}\phi\bigr\}$
et
\dli $\phi_{(k_1,\beta_1)}\,=\,inf\bigl\{\phi'\in[\phi_{(0,1)},\phi_
{(\infty,0)}]\, ;\,\phi\relmont{\leq}{p.s.}\phi'\bigr\}$.

Si $k_0 > k_1$, pour $k\in]k_1,k_0[$ on a \'evidemment
\dli $\phi_{(k,0)}\relmont{\leq}{p.s.}\phi\relmont{\leq}{p.s.}\phi_{(k,1)}$.
Ces in\'egalit\'es sont aussi v\'erifi\'ees pour $k=k_0$, lorsque $0<k_0=k_1<
+\,\infty$.
Les cas $k_0=k_1=0$ et $k_0=k_1=+\,\infty$ correspondent \`a
$\phi\relmont{=}{p.s.}\phi_{(0,1)}$ et
$\phi\relmont{=}{p.s.}\phi_{(\infty,0)}$.
\medskip
{\souli{Il nous reste \`a montrer que l'on ne peut pas avoir
$k_0 < k_1$.}}

Pour cela nous allons supposer $k_0 < k_1$ et
trouver un \'ev\'enement $C$ v\'erifiant :
$P_0(C)\not= 0$, $P_1(C)\not= 0$ et $E_1^C(\phi)<E_0^C(\phi)$. Par d\'efinition, ceci contredira le fait que $\phi$ est un expert. 

\dli Soit $k\in]k_0,k_1[$.
\dli Consid\'erons l'\'ev\'enement
$A=\{\phi_{(k,0)}-\phi_{(k_0,\beta_0)}=1\}\cap\{\phi=0\}$.
La d\'efinition de  $\phi_{(k_0,\beta_0)}$ entra\^{\i}ne que $P_0(A)$ ou $P_1(A)$
est non nulle, sinon on aurait $\phi_{(k_0,\beta_0)}<\phi_{(k,0)}\relmont{\leq}{p.s.}\phi$. De plus, sur $A\subseteq\Omega-(\Omega_0\cup\Omega_\infty)$ les densit\'es $p_0$ et $p_1$ sont strictement positives, on a donc : 
$P_0(A)\not= 0$ et $P_1(A)\not= 0$.
\dli De m\^eme, la d\'efinition de $\phi_{(k_1,\beta_1)}$ entra\^{\i}ne que l'\'ev\'enement $B=\{\phi_{(k_1,\beta_1)}-\phi_{(k,1)}=1\}\cap\{\phi=1\}$ v\'erifie : $P_0(B)\not= 0$ et $P_1(B)\not= 0$.
\dli L'\'ev\'enement $C=A\cup B$ est bien s\^ur de probabilit\'e non nulle sous $P_0$ et $P_1$, pour finir nous allons montrer qu'il v\'erifie :
$E_1^C(\phi)<E_0^C(\phi)$. On a $E_0^C(\phi)=P_0(B)/[P_0(A)+P_0(B)]$ et $E_1^C(\phi)=P_1(B)/[P_1(A)+P_1(B)]$, avec
$P_0(A)=\int_Ap_0\,d\mu<\int_Akp_1\,d\mu=k P_1(A)$ et $P_0(B)=\int_Bp_0\,d\mu>\int_Bkp_1\,d\mu=k P_1(B)$, ce qui entra\^{\i}ne 
$P_0(A)/P_0(B)<P_1(A)/P_1(B)$ ; on en d\'eduit alors facilement l'in\'egalit\'e recherch\'ee.

\medskip
{\parindent=-10mm II ----- Condition suffisante.}

Soit $\phi$ une r\`egle de d\'ecision encadr\'ee presque s\^urement par deux fonctions de test simples \'egales ou cons\'ecutives : $\phi_{(k,\beta)}\relmont{\leq}{p.s.}\phi\relmont{\leq}{p.s.}\phi_{(k,\beta')}$. On doit d\'emonter que $\phi$  est un expert. Pour cela on consid\`ere un \'ev\'enement $C$ et on envisage les trois cas de la d\'efinition 2.1.1.

1\up{er} cas : $P_0(C)=0$.

$C$ est donc $P_1$ presque s\^urement inclus dans $\Omega_0$ ; de plus comme
$\phi\relmont{\geq}{p.s.}\phi_{(k,\beta)}\geq\phi_{(0,1)}$, $\phi$ est presque s\^urement \'egale \`a $1$ sur $\Omega_0$ ; l'\'egalit\'e recherch\'ee, $P_1(C\cap\{\phi=0\})=0$ s'en d\'eduit facilement.

2\up{\`eme} cas : $P_1(C)=0$.

Dans ce cas $C$ est $P_0$ presque s\^urement inclus dans $\Omega_\infty$ et 
$\phi$ est presque s\^urement \'egale \`a $0$ sur $\Omega_\infty$ puisque : 
$\phi\relmont{\leq}{p.s.}\phi_{(k,\beta')}\leq\phi_{(\infty,0)}$ ; on obtient alors facilement l'\'egalit\'e recherch\'ee : $P_0(C\cap\{\phi=1\})=0$.

3\up{\`eme} cas : $P_0(C)\not=0$ et $P_1(C)\not=0$.

Consid\'erons le mod\`ele statistique conditionnel :

\dli $\bigl(\Omega ,{\cal A},(P_\theta^C)_{\theta\in\{0,1\}} \bigr)$ avec 
$P^C_\theta$ de densit\'e $(1/P_\theta(C))\II_C$ par rapport \`a $P_\theta$.
\dli Dans ce nouveau mod\`ele $\phi$ est encore un test de Neyman, elle
d\'efinit un test uniform\'ement plus puissant, pour tester $P_0^C$ contre
$P_1^C$ au seuil $(1/P_0(C))E_0(\II_C.\phi)=E_0^C(\phi)$ et pour tester $P_1^C$ contre $P_0^C$ au seuil 
$(1/P_1(C))E_1(\II_C.(1-\phi))=E_1^C(1-\phi)$. La puissance \'etant sup\'erieure au seuil on a bien $E_1^C(\phi)\geq E_0^C(\phi)$ (cf. [Leh.] p. 76).

\medskip\centerline{\hbox to 3cm{\bf \hrulefill}}\par}

\bigskip
\vfill\eject

{\parindent=-5mm 2.4 VOTES DES EXPERTS.}
\medskip
Les experts font partie des proc\'edures de d\'ecision admissibles pour le risque classique : la probabilit\'e de se tromper.
Si on voulait s\'electionner un expert il faudrait imposer des contraintes
suppl\'ementaires.
Dans la th\'eorie des tests on privil\'egie une hypoth\`ese, par exemple
$H_0:\theta=0$, on fixe un seuil $\alpha$ et on s'int\'eresse
aux proc\'edures de d\'ecision qui v\'erifient $E_0(\phi)\leq\alpha$ et maximisent
$E_1(\phi)$. Le lemme fondamental de Neyman et Pearson (cf. [Leh.] p. 74)
nous donne la solution. Cette mani\`ere de faire est bien s\^ur critiquable.
Certains vont pr\'ef\'erer donner un indice d'aide \`a la prise de d\'ecision
comme la {\it p-value}, le seuil minimum de rejet (cf. [Sch.]).
D'autres vont refuser la dissym\'etrie totale
de traitement des deux hypoth\`eses en acceptant une troisi\`eme d\'ecision,
celle de ne pas conclure (cf. [Mor.1] [Ney.] [Nik.1]).
\dli La difficult\'e du choix d'un expert provient du
comportement oppos\'e des deux risques $E_0(\phi)$ et $E_1(1-\phi)$,
quand l'un diminue l'autre augmente.
$$\vbox{
\offinterlineskip
\halign{$\cc{#}$&\tv#&&$\cc{#}$\cr
(k,\beta)&&(0,1)&\nearrow&(\infty,0)\cr
\traithori
E_0(\phi_{(k,\beta)})&&0&\nearrow&1-P_0(\Omega_\infty)\cr
\traithori
E_1(1-\phi_{(k,\beta)})&&1-P_1(\Omega_0)&\searrow&0\cr
}}$$
On peut supprimer cet inconv\'enient en construisant un risque unique, 
synth\`ese des deux pr\'ec\'edents. La difficult\'e est alors transf\'er\'ee au
choix de cette synth\`ese. Dans le cadre Bayesien c'est le choix de
la probabilit\'e a priori sur $\Theta=\{0,1\}$, c'est-\`a-dire les
pond\'erations des deux risques. On s\'electionne ainsi une proc\'edure de d\'ecision qui minimise la moyenne pond\'er\'ee des deux risques. Mais bien souvent on
fournit \`a l'utilisateur la probabilit\'e \`a post\'eriori sur $\Theta=\{0,1\}$
et on le laisse se forger sa
propre opinion \`a partir de cette probabilit\'e sur $\Theta$,
qu'il transformera en probabilit\'e sur $D=\{0,1\}$.

Dans la plupart des pratiques, le choix d'une proc\'edure de d\'ecision est illusoire,
l'utilisateur ne peut pas mettre en avant les yeux ferm\'es la r\'eponse
fournie, m\^eme si elle est d\'eterministe. Il doit d\'efendre
les crit\`eres qui ont pr\'esid\'e au choix de cette proc\'edure. Pour une application donn\'ee ces crit\`eres font rarement l'unanimit\'e des utilisateurs, plusieurs r\`egles sont l\'egitimes.
Comme elles donnent leur avis gratuitement, aucune contrainte \'economique ne nous emp\^eche de les entendre toutes.
C'est ce que nous nous proposons de faire pour les experts.
La difficult\'e sera alors de r\'esumer leurs r\'eponses. Pour cela nous
accorderons \`a chacun d'eux un poids sp\'ecifique en
probabilisant l'ensemble des experts s\'electionn\'es. L'utilisateur
obtiendra le r\'esultat de ce "vote" d'experts sous la forme d'une
probabilit\'e sur $D$. Bien entendu nous ne chercherons pas \`a
probabiliser l'ensemble des experts par une masse de Dirac
car nous retrouverions le probl\`eme de la s\'election d'un expert.
Nous allons plut\^ot essayer de faire voter les experts
``d\'emocratiquement".

Avant d'\'etudier comment r\'epartir les voix des experts, il nous faut mieux d\'efinir les votants. Par exemple, des experts presque s\^urement \'egaux ne feront qu'un, ils n'auront droit qu'\`a une seule carte de vote. Nous voulons aussi que les diff\'erentes fa\c cons de mod\'eliser le probl\`eme de d\'ecision n'influencent pas les r\'esultats. C'est bien s\^ur la valeur de la statistique exhaustive $K$, le rapport des densit\'es, qui joue le r\^ole fondamental. Elle est unique presque s\^urement (voir annexe I). Mais nous avons vu qu'il y a des experts qui ne sont pas uniquement d\'efinis, presque s\^urement, par la valeur de $K$. Ce sont ceux correspondant \`a une fonction de test $\phi$ strictement encadr\'ee par deux fonctions de test simples cons\'ecutives : $\phi_{(k,0)}\relmont{<}{p.s.}\phi\relmont{<}{p.s.}\phi_{(k,1)}$ (voir la proposition 2.3.1). Ces experts ne peuvent exister que pour $k\in\IR_*^+$ et \`a condition qu'il y ait une masse en $k$ : $P_0(\{K=k\})>0$ et $P_1(\{K=k\})>0$. Cela ne suffit pas, il faut pouvoir partager cette masse entre les d\'ecisions $d=0$ et $d=1$. Ce qui n'est pas possible quand $\Omega$ est constitu\'e des valeurs de $K$. L'existence de tels experts suppose donc que l'ensemble des r\'ealisations $\Omega$ contienne des informations superflues, par exemple un al\'ea permettant de prendre une d\'ecision \`a partir d'une r\`egle al\'eatoire (voir 2.2). Quand il existe plusieurs experts de ce type pour un m\^eme $k$, ce qui les diff\'erencie c'est principalement la part de $P_0(\{K=k\})$ et $P_1(\{K=k\})$ qu'ils consacrent \`a la d\'ecision $d=1$. S'ils le font sur des r\'ealisations diff\'erentes dans  $\{K=k\}$ ceci n'a pas d'importance puisque sous $\{K=k\}$, les probabilit\'es conditionnelles sont uniformes. Les experts compris strictement entre $\phi_{(k,0)}$ et $\phi_{(k,1)}$ peuvent donc sur $\{K=k\}$ faire ce qu'ils veulent, le choix entre $P_0$ et $P_1$ n'est plus structurant, seule la limitation des informations inutiles peut les restreindre. Nous allons commencer par refuser toute information superflue et donc nous limiter aux experts d\'ependant de la statistique exhaustive $K$. Ce qui revient \`a faire voter les experts d\'efinis par l'ensemble des fonctions de test simples.

La r\'epartition des voix sur $\Phi_s=[\phi_{(0,1)} , \phi_{(\infty,0)}]$, consiste \`a probabiliser  $\Phi_s$. Pour cela nous le munissons de la tribu bor\'elienne correspondant \`a la topologie de l'ordre.
Un intervalle de $\Phi_s$ aura d'autant plus de poids que
ses extr\'emit\'es $\phi$ et $\phi'$ sont des experts diff\'erents. Il semble naturel d'exprimer ces diff\'erences en comparant les moyennes des experts sous $P_\theta$. Des experts \'egaux presque s\^urement ne compteront alors que pour un, car sur $\Phi_s$, les classes d'\'equivalences form\'ees d'experts de m\^eme moyenne  sont celles d\'efinies par l'\'egalit\'e presque s\^ure. On d\'efinit alors une probabilit\'e sur $\Phi_s$ en posant pour  $\phi\in]\phi_{(0,1)} , \phi_{(\infty,0)}[$ : $m_\theta([\phi_{(0,1)},\phi[)=E_\theta(\phi)$. Il en est de m\^eme si l'on pose $m'_\theta([\phi_{(0,1)},\phi])=E_\theta(\phi)$. Ces deux probabilit\'es sont identiques si la fonction de r\'epartition de $K$ est continue. Dans le cas contraire il existe $k\in\IR_*^+$ poss\'edant une masse : $P_\theta(\{K=k\})>0$. Cette masse est affect\'ee \`a $\phi_{(k,0)}$ par $m_\theta$ et \`a $\phi_{(k,1)}$ par $m'_\theta$.
 En fait, il semble l\'egitime de partager \'equitablement la masse. Pour une r\'ealisation $\omega$, le r\'esultat du vote est alors une
probabilit\'e $Q^\omega_\theta$ sur $D$ d\'efinie par
$Q^\omega_\theta(\{1\})=\int_{\Phi_s}\phi(\omega)
\,d(m_\theta+m'_\theta)/2$.

\medskip
{\bf D\'efinition 2.4.1}
\medskip
\medskip
\moveleft 10.4pt\hbox{\vrule\kern 10pt\vbox{\defpro

Lorsqu'on r\'ealise $\omega$ appartenant \`a $\{K=k\}$, le r\'esultat du vote des experts sous $P_\theta$ est une probabilit\'e $Q^\omega_\theta$ d\'efinie sur l'espace des d\'ecisions par :
$$Q^\omega_\theta(\{1\})=\left\{\matrix{
1\hfill & si & k=0\hfill\cr
(1/2)P_\theta(\{K=k\})+P_\theta(\{K>k\}) = 1-E_\theta(\phi_{(k,1/2)})\hfill 
& si & k\in\IR^+_*\hfill \cr
0\hfill & si & k=\infty\hfill \cr
}\right. $$
(avec $\phi_{(k,1/2)}=\II_{\{K<k\}}+(1/2)\II_{\{K=k\}}=
(1/2)[\phi_{(k,0)}+\phi_{(k,1)}]$).
}}\medskip

On retrouve la notion de {\it p-value}, plus exactement celle de {\it mid-p-value} (cf. [Rou.]). Si on fait voter les experts en utilisant $P_0$ pour les diff\'erencier, la
fr\'equence de la d\'ecision $d=0$ quand
on r\'ealise $\omega$, $Q^\omega_0(\{0\})$, est un seuil minimum de rejet
du test de $H_0:\theta=0$ contre $H_1:\theta=1$. Dans le cas o\`u on
r\'epartit les voix en utilisant $P_1$,
la fr\'equence de la d\'ecision
$d=1$, $Q^\omega_1(\{1\})$, est un seuil minimum de rejet
du test de $H_0:\theta=1$ contre $H_1:\theta=0$.
S'il y a une masse en $k\in\IR^+_*$, ce seuil minimum de rejet correspond \`a une moyenne uniforme entre les seuils des tests al\'eatoires \`a la limite du rejet quand on observe $k$. Ces tests al\'eatoires deviennent des experts si on prend pour mod\`ele le produit du mod\`ele image de $K$ par l'al\'ea sur $[0,1]$ muni de la probabilit\'e uniforme. On obtient l'ensemble $\Phi_e$ des experts d\'efinis par $\phi_{(k,\beta)}(\omega,u)=\II_{\{K<k\}}(\omega)+\II_{\{K=k\}}(\omega).\II_{[0,\beta]}(u)$ avec $k\in\IR^+_*$ et $\beta\in[0,1]$. On peut d\'efinir comme pr\'ec\'edemment des probabilit\'es $m_\theta$ et $m'_\theta$, cette fois elles sont d\'efinies sur $\Phi_e$ muni de la topologie de l'ordre lexicographique, et l'op\'erateur $E_\theta$ d\'esigne la moyenne par rapport au produit de $P_\theta$ par la probabilit\'e uniforme. Dans ce cas on a cr\'e\'e un continuum, les probabilit\'es $m_\theta$ et $m'_\theta$ sont identiques. Il est facile de v\'erifier que pour cette probabilit\'e, le r\'esultat du vote des experts de $\Phi_e$ est celui de la d\'efinition pr\'ec\'edente.

\bigskip
\vfill\eject
{\parindent=-5mm 2.5 VOTE POND\'ER\'E.}
\medskip

Nous allons analyser ce qui se passe quand on fait voter
les experts en utilisant un m\'elange de $P_0$ et $P_1$. Il est d\'efini
par une probabilit\'e $\Lambda$ sur $\Theta=\{0,1\}$ donc par
$\Lambda(\{0\})=\lambda$. Si l'on veut comparer les r\'esultats
obtenus avec ceux de l'analyse Bayesienne \`a partir de la loi a priori
$\Lambda$, il est pr\'ef\'erable que $\lambda$ puisse s'interpr\'eter comme
une prise de position en faveur de l'hypoth\`ese $\theta=0$, donc que la
croissance de $\lambda$ entra\^{\i}ne celle de la fr\'equence des experts
qui d\'ecident $P_0$. 
Comme $E_1(\phi_{(k,\beta)})-E_0(\phi_{(k,\beta)})\geq 0$
(voir fin de 2-2) on a $Q^\omega_0(\{0\})\leq Q^\omega_1(\{0\})$,
le m\'elange associ\'e \`a $\lambda\in [0,1]$ sera :
$(1-\lambda)P_0+\lambda P_1$.

\medskip
{\bf D\'efinition 2.5.1}
\medskip
\medskip
\moveleft 10.4pt\hbox{\vrule\kern 10pt\vbox{\defpro

Soit $\lambda\in[0,1]$. Pour une r\'ealisation $\omega\in\Omega$, on appelle vote des experts pond\'er\'e par
$\lambda$, la probabilit\'e d\'efinie sur l'ensemble des d\'ecisions par :
\dli $Q_\lambda^\omega(\{1\})=(1-\lambda) Q_0^\omega(\{1\})
+\lambda Q_1^\omega(\{1\})$.
}}\medskip

$Q_\lambda^\omega$
est le r\'esultat du vote des experts quand dans la d\'efinition 2.4.1 on remplace la probabilit\'e $P_\theta$ par le m\'elange
$(1-\lambda)P_0+\lambda P_1$.
On retrouve bien les votes $Q_0^\omega$ et $Q_1^\omega$ pour $\lambda=0$ et 
$\lambda=1$.

On maximise la prise de d\'ecision en faveur de $\theta=0$, en prenant $\lambda=1$. Le vote des experts est alors \'etabli \`a partir de $P_1$ et 
la fr\'equence de la d\'ecision
$d=1$, $Q^\omega_1(\{1\})$, est le seuil minimum de rejet
du test de $H_0:\theta=1$ contre $H_1:\theta=0$. C'est le type de test qui
permet de confirmer avec force la faveur accord\'ee \`a $\theta=0$ puisqu'on est alors dans le cas du rejet. On peut faire une remarque semblable dans le cas o\`u on favorise $\theta=1$, en prenant $\lambda=0$. C'est le seuil minimum de rejet du test de $H_0:\theta=0$ contre $H_1:\theta=1$ qui intervient.  

M\^eme si ici $D=\Theta=\{0,1\}$, la probabilit\'e $Q_\lambda^\omega$
ne peut pas \^etre confondue avec la probabilit\'e a posteriori
$\Lambda(\{\theta\}\mid\omega)$ obtenue \`a partir de la probabilit\'e a
priori $\Lambda$ d\'efinie par $\Lambda(\{0\})=\lambda$.
Lorsque $\omega$ appartient \`a $\Omega_k$ on a (cf. [Bor.] p. 283) :
$$\Lambda(\{1\}\mid\omega)\quad=\quad\left\{\matrix{
1\hfill & si & k=0\hfill\cr
(1-\lambda )/[\lambda k+(1-\lambda)]\hfill
 & si & k\in\IR^+_*\hfill \cr
0\hfill & si & k=\infty.\hfill \cr
}\right. $$
Alors que la fr\'equence des experts ayant d\'ecid\'e $P_1$, dans un vote
pond\'er\'e par $\lambda$, est \'egale \`a :
$$Q_\lambda^\omega(\{1\})\quad=\quad\left\{\matrix{
1\hfill & si & k=0\hfill\cr
1-(1-\lambda)E_0(\phi_{(k,1/2)})-\lambda E_1(\phi_{(k,1/2)})\hfill
 & si & k\in\IR^+_*\hfill \cr
0\hfill & si & k=\infty.\hfill \cr
}\right. $$

\medskip
\souli{Exemple} :

Consid\'erons, par rapport \`a la mesure de Lebesgue, les deux densit\'es :
\dli $p_0=(1/6)\II_{[0,1]}+(1/3)\II_{]1,2[}+(1/2)\II_{[2,3]}$ et
\dli $p_1=(1/2)\II_{[0,1]}+(1/3)\II_{]1,2[}+(1/6)\II_{[2,3]}$
\dli Il y a trois valeurs de $k$ utiles : 1/3, 1 et 3. Le tableau
suivant donne le vote des experts et la probabilit\'e a posteriori
en fonction de la pond\'eration $\lambda$.
$$\vbox{
\offinterlineskip
\halign{\tv#&$\cc{#}$&\tv#&&
        $\cc{#}$&\tv#&$\hfill\quad #\,$&$#$&
        ${\,#\,}$&$#$&$\,#\quad\hfill$&\tv#\cr
\omit&&\multispan9\hrulefill&\multispan8\hrulefill \cr
\omit&&&\multispan7$\cc{Q_\lambda^\omega(\{1\})}$&&
        \multispan7$\cc{\Lambda(\{1\}\mid\omega)}$& \cr
\omit&&\multispan9\hrulefill&\multispan8\hrulefill \cr
\omit&&&1-\lambda&&0&\nearrow&1/2&\nearrow&1&
       &1-\lambda&&0&\nearrow&1/2&\nearrow&1&\cr
\multispan3\hrulefill&\multispan8\hrulefill&\multispan8\hrulefill \cr
&\omega\in[0,1]&&{11-2\lambda\over 12}&&3/4&\nearrow&5/6&\nearrow&11/12&
                &{3-3\lambda\over 3-2\lambda}&&0&\nearrow&3/4&\nearrow&1&\cr
\multispan3\hrulefill&\multispan8\hrulefill&\multispan8\hrulefill \cr
&\omega\in]1,2[&&{2-\lambda\over 3}&&1/3&\nearrow&1/2&\nearrow&2/3&
                &1-\lambda&&0&\nearrow&1/2&\nearrow&1&\cr
\multispan3\hrulefill&\multispan8\hrulefill&\multispan8\hrulefill \cr
&\omega\in[2,3]&&{3-2\lambda\over 12}&&1/12&\nearrow&1/6&\nearrow&1/4&
                &{1-\lambda\over 1+2\lambda}&&0&\nearrow&1/4&\nearrow&1&\cr
\multispan3\hrulefill&\multispan8\hrulefill&\multispan8\hrulefill \cr
}}$$
Le choix de $\lambda$ a des cons\'equences moins lourdes sur
$Q_\lambda^\omega(\{1\})$ que sur $\Lambda(\{1\}\mid\omega)$.
Pour une r\'ealisation $\omega$, les r\'eponses ont une amplitude de $1$
dans le cas de la probabilit\'e a posteriori et de $2/12$ ou $1/3$
dans le cas du vote des experts.

\medskip
Cette remarque ne d\'epend pas de l'exemple. Il en est toujours ainsi quand
$\omega$  appartient \`a $\Omega_k$ et $k$ \`a $\IR^+_*$ ;
lorsque $(1-\lambda)$ varie de $0$ \`a $1$, $\Lambda(\{1\}\mid\omega)$ cro\^{\i}t
de $0$ \`a $1$ alors que $Q_\lambda^\omega(\{1\})$ cro\^{\i}t de
$1-E_1(\phi_{(k,1/2)})$
\`a $1-E_0(\phi_{(k,1/2)})$. L'amplitude de ces variations,
$Q_0^\omega(\{1\})-Q_1^\omega(\{1\})$, est maximum lorsque
la r\'ealisation n'apporte aucune information, c'est-\`a-dire que $\omega$
appartient \`a $\Omega_1=\{K=1\}$ ; elle vaut alors
$E_1(\phi_{(1,1/2)})-E_0(\phi_{(1,1/2)})=P_1(B)-P_0(B)$ avec
$B=\{\omega\in\Omega\ ;\ p_0(\omega)<p_1(\omega)\}$.
Le centre de l'intervalle de variation correspond \`a $\lambda=1/2$,
$Q_{1/2}^\omega(\{1\})=1-(1/2)[E_0(\phi_{(k,1/2)})+E_1(\phi_{(k,1/2)}]$,
dans ce cas on ne privil\'egie aucune des deux hypoth\`eses, on les traite
sym\'etriquement.
Ce qualificatif peut s'appliquer \`a d'autres mani\`eres de faire. Dans le paragraphe suivant nous en verrons une qui repose directement sur les deux votes de base sans les m\'elanger.

\vfill\eject

{\parindent=-10mm\soustitre 3--R\`EGLES DE D\'ECISION DE BOL'SHEV.}
\bigskip
{\parindent=-5mm 3.1 D\'EFINITIONS ET PROPRI\'ET\'ES.}
\medskip

Les r\`egles de Bol'shev donnent une solution au probl\`eme du choix
entre deux probabilit\'es quand on accepte la possibilit\'e de ne pas
conclure (cf. [Nik.1] et [Nik.2]).
Bol'shev consid\`ere l'ensemble des d\'ecisions $D'=\{0,1,2\}$, la d\'ecision $2$
voulant dire que les deux probabilit\'es sont plausibles, on refuse de
prendre position en faveur de l'une d'elles. Une r\`egle de d\'ecision $\delta'$
est alors une probabilit\'e de transition de $(\Omega,{\cal A})$ vers
$(D',{\cal D'})$ ;
${\cal D'}$ \'etant l'ensemble des parties de $D'$. Pour toute partie $C$ de $D'$, $\delta'(\omega,C)$ repr\'esente la probabilit\'e de d\'ecider que $d$ appartient \`a $C$, lorsqu'on observe $\omega$.

La perte utilis\'ee par Bol'shev est celle qui donne pour risque la
probabilit\'e de se tromper, on peut l'\'ecrire
$L(\theta,d)=1-\II_{\{\theta,2\}}(d)$. Pour une r\`egle $\delta'$ il
consid\`ere donc les deux erreurs :
\dli $R(\theta,\delta')=\int_\Omega\bigl[\int_{D'} L(\theta,.)\,
d\delta'(\omega,.)\bigr]\,dP_\theta(\omega)=
E_\theta[1-\delta'(.,\{\theta,2\})]$.
\dli Bol'shev se fixe deux seuils $\alpha_0\in[0,1]$ et
$\alpha_1\in[0,1]$ ; il s'int\'eresse
aux r\`egles $\delta'$ v\'erifiant $R(\theta,\delta')\leq\alpha_\theta$
pour $\theta\in\{0,1\}$. Cet ensemble de r\`egles, $\Delta'_\alpha$, est non vide
puisqu'il contient la r\`egle qui d\'ecide toujours $d=2$.
Bien s\^ur, cette r\`egle est inint\'eressante, Bol'shev cherche des r\`egles
qui ne d\'ecident pas trop souvent $d=2$. Pour les r\`egles $\delta'$ de
$\Delta'_\alpha$, il va faire intervenir les deux probabilit\'es de d\'ecider $d=2$ :
$E_\theta(\delta'(.,\{2\}))$ et essayer de les minimiser. Il consid\`ere
donc le pr\'eordre partiel suivant.
\medskip
{\bf D\'efinition 3.1.1}
\medskip
\medskip
\moveleft 10.4pt\hbox{\vrule\kern 10pt\vbox{\defpro

Dans $\Delta'_\alpha$, une r\`egle $\delta'$ est aussi bonne qu'une r\`egle
$\delta''$ si :
\dli $E_0(\delta'(.,\{2\}))\leq E_0(\delta''(.,\{2\}))$ et
$E_1(\delta'(.,\{2\}))\leq E_1(\delta''(.,\{2\}))$.
}}\medskip

Nous allons retrouver les r\'esultats principaux sur les r\`egles de Bol'shev
en utilisant les fonctions de test simples, ce
qui nous permettra de faire plus facilement le lien avec les experts.
Comme les r\`egles de d\'ecision consid\'er\'ees sont al\'eatoires, les fonctions de test $\phi_{(k,\beta)}$ utilis\'ees le seront avec $\beta$ dans $[0,1]$ et pas dans $\{0,1\}$ (voir la d\'efinition 2.2.1). L'ensemble de ces fonctions de test est not\'e $\Phi_a$.
\dli$\Phi_a=\Bigl\{\phi_{(0,1)}\, ;\,\{\phi_{(k,\beta)}\}_{k\in\IR_*^+,\,\beta\in[0,1]}\, ;\,\phi_{(\infty,0)}\Bigr\}$.
\dli Il est muni de l'ordre lexicographique et il a alors les m\^emes propri\'et\'es que $\Phi_s$ (voir la fin de 2.2). L'op\'erateur $E_\theta$ est continu pour la topologie de l'ordre.

\medskip
{\bf D\'efinition 3.1.2}
\medskip
\medskip
\moveleft 10.4pt\hbox{\vrule\kern 10pt\vbox{\defpro

On appelle r\`egle de Bol'shev, les probabilit\'es de transition de
$(\Omega,{\cal A})$ vers $(D',{\cal D'})$ d\'efinies par deux fonctions
de test de $\Phi_a$,
$\phi_{(k',\beta')}\leq\phi_{(k'',\beta'')}$, de la mani\`ere suivante :
$$\delta'(\omega,\{d\})\quad=\quad\left\{\matrix{
1-\phi_{(k'',\beta'')}(\omega)\hfill & si & d=0\hfill\cr
\phi_{(k',\beta')}(\omega)\hfill & si & d=1\hfill \cr
\phi_{(k'',\beta'')}(\omega)-\phi_{(k',\beta')}(\omega)\hfill & si
       & d=2.\hfill \cr
}\right. $$
}}\medskip

Comme le montre la proposition suivante, on peut dans $\Delta'_\alpha$ se
restreindre aux r\`egles de Bol'shev.
\medskip
{\bf Proposition 3.1.1}
\medskip
\medskip
\moveleft 10.4pt\hbox{\vrule\kern 10pt\vbox{\defpro

Soit $\delta'$ un \'el\'ement de $\Delta'_\alpha$. Il existe dans $\Delta'_\alpha$
une r\`egle de Bol'shev $\delta'_B$ aussi bonne que $\delta'$.
}}\medskip

{\leftskip=15mm \dli {\bf D\'emonstration}

$\phi'(\omega)=\delta'(\omega,\{1\})$ d\'efinit un test de $H_0 : \theta=0$
contre $H_1 : \theta=1$ au seuil $\alpha'=E_0(\delta'(.,\{1\}))\leq
\alpha_0$.
Les propri\'et\'es de l'op\'erateur $E_\theta$ sur $\Phi_a$ permettent d'affirmer
qu'il existe dans $\Phi_a$ un plus grand \'el\'ement,
$\phi_{(k',\beta')}$, v\'erifiant :
\dli $E_0(\phi_{(k',\beta')})=inf\{\alpha'\ ,\ (1-P_0(\Omega_\infty))\}$.
\dli D'apr\`es le lemme de Neyman et Pearson (cf. [Leh.] p. 74)
$\phi_{(k',\beta')}$ d\'efinit un test de $H_0$ contre $H_1$
aussi puissant que $\phi'$, ce qui se traduit par :
\dli $E_1(\phi_{(k',\beta')})\geq E_1(\phi')$.

De m\^eme $\phi''(\omega)=\delta'(\omega,\{0\})$ d\'efinit un test de
$H'_0 : \theta=1$ contre $H'_1 : \theta=0$ au seuil
$\alpha''=E_1(\delta'(.,\{0\}))\leq\alpha_1$.
Il existe  dans $\Phi_a$ un plus petit
\'el\'ement $\phi_{(k'',\beta'')}$ qui v\'erifie :
\dli $E_1[1-\phi_{(k'',\beta'')}]=inf\{\alpha''\ ,\ (1-P_1(\Omega_0))\}$.
\dli $1-\phi_{(k'',\beta'')}$ d\'efinit un test de Neyman de $H'_0$
contre $H'_1$ aussi puissant que $\phi''$.
Il v\'erifie donc aussi :
\dli $E_0(1-\phi_{(k'',\beta'')})\geq E_0(\phi'')$.
\medskip
{\parindent=-10mm 1\up{er} cas : $\phi_{(k',\beta')}<\phi_{(k'',\beta'')}$.}

Ces deux fonctions de test d\'efinissent une r\`egle de Bol'shev
$\delta'_B$ qui appartient \`a $\Delta'_\alpha$ car pour $\theta\in\{0,1\}$ on a :
\dli $R(\theta,\delta'_B)=E_\theta[1-\delta'_B(.,\{\theta,2\})]
      \leq\alpha_\theta$.

Il reste \`a d\'emontrer que $\delta'_B$ est aussi bonne que $\delta'$, soit :
\dli $E_\theta[\delta'_B(.,\{2\})]\leq E_\theta[\delta'(.,\{2\})]$
pour $\theta\in\{0,1\}$.

L'in\'egalit\'e $\phi_{(k',\beta')}<\phi_{(k'',\beta'')}$ est
incompatible avec l'\'egalit\'e
\dli $E_0(\phi_{(k',\beta')})=1-P_0(\Omega_\infty)$
car on a alors $\phi_{(k',\beta')}=\phi_{(\infty,0)}$ ;
\dli $E_0(\phi_{(k',\beta')})$ est donc \'egale \`a $\alpha'=E_0(\phi')$.
\dli On a aussi $E_1[1-\phi_{(k'',\beta'')}]=\alpha''$, sinon on aurait
$\phi_{(k'',\beta'')}=\phi_{(0,1)}$.
\dli Pour $\theta\in\{0,1\}$ on obtient alors facilement :
\dli $E_\theta[\delta'_B(.,\{2\})]=
      E_\theta[\phi_{(k'',\beta'')}-\phi_{(k',\beta')}]\leq
      E_\theta(1-\phi''-\phi')=E_\theta(\delta'(.,\{2\}))$.
\medskip
{\parindent=-10mm 2\up{\`eme} cas : $\phi_{(k',\beta')}\geq
   \phi_{(k'',\beta'')}$.}

Soit $\phi_{(k,\beta)}\in[\phi_{(k'',\beta'')}\ ,\ \phi_{(k',\beta')}]$,
nous d\'efinissons une r\`egle de Bol'shev $\delta'_B$ en posant :
\dli $\delta'_B(\omega,\{0\})=1-\phi_{(k,\beta)}(\omega)$ et
     $\delta'_B(\omega,\{1\})=\phi_{(k,\beta)}(\omega)$.

Cette r\`egle $\delta'_B$ appartient \`a $\Delta'_\alpha$ car :
\dli $R(0,\delta'_B)=E_0(\phi_{(k,\beta)})\leq E_0(\phi_{(k',\beta')})
      \leq\alpha'\leq\alpha_0$ et
\dli $R(1,\delta'_B)=E_1(1-\phi_{(k,\beta)})\leq E_1(1-\phi_{(k'',\beta'')})
      \leq\alpha''\leq\alpha_1$.
\dli Elle est bien s\^ur aussi bonne que $\delta'$ puisqu'elle ne d\'ecide
jamais $d=2$ :
\dli $E_\theta[\delta'_B(.,\{2\})]=0\leq E_\theta(\delta'(.,\{2\}))$
      pour $\theta\in\{0,1\}$.
\medskip\centerline{\hbox to 3cm{\bf \hrulefill}}\par}

Nous avons m\^eme d\'emontr\'e que $\delta'_B$ est "aussi bonne" que $\delta'$
au sens suivant :
\dli $E_0[\delta'_B(.,\{1\})]\leq E_0[\delta'(.,\{1\})]$,
$E_1[\delta'_B(.,\{0\})]\leq E_1[\delta'(.,\{0\})]$ et
\dli $E_0[\delta'_B(.,\{2\})]\leq E_0[\delta'(.,\{2\})]$,
$E_1[\delta'_B(.,\{2\})]\leq E_1[\delta'(.,\{2\})]$

\medskip
{\bf Proposition 3.1.2}
\medskip
\medskip
\moveleft 10.4pt\hbox{\vrule\kern 10pt\vbox{\defpro

Soient $\alpha_0\in[0,1]$ et $\alpha_1\in[0,1]$.
\dli Il existe dans $\Delta'_\alpha$ une r\`egle de Bol'shev $\delta'_m$
telle que pour $\theta\in\{0,1\}$ :
\dli $E_\theta[\delta'_m(.,\{2\})]=inf\{E_\theta[\delta'(.,\{2\})]\,;\,
      \delta'\in\Delta'_\alpha\}$.
\dli Cette r\`egle est unique presque s\^urement si et seulement si
$\phi_{(k_0,\beta_0)}\relmont{\leq}{p.s.}\phi_{(k_1,\beta_1)}$ ; les deux fonctions de test $\phi_{(k_0,\beta_0)}$ et $\phi_{(k_1,\beta_1)}$ 
\'etant respectivement le plus grand et le plus petit \'el\'ement de
$\Phi_a$ tels que :
\dli $E_0(\phi_{(k_0,\beta_0)})=inf\{\alpha_0\ ,\ (1-P_0(\Omega_\infty))\}$ ;
\dli $E_1[1-\phi_{(k_1,\beta_1)}]=inf\{\alpha_1\ ,\ (1-P_1(\Omega_0))\}$.
}}\medskip

{\leftskip=15mm \dli {\bf D\'emonstration}

La proposition 3.1.1 nous permet d'affirmer qu'il suffit de consid\'erer
le sous ensemble $\Delta^B_\alpha$ des r\`egles de Bol'shev  de
$\Delta'_\alpha$.
\medskip
{\parindent=-10mm I ----- Construction de $\delta'_m$.}

Les op\'erateurs $E_\theta$ \'etant continus sur $\Phi_a$, il existe un plus grand \'el\'ement $\phi_{(k_0,\beta_0)}$, et un plus petit \'el\'ement $\phi_{(k_1,\beta_1)}$, v\'erifiant respectivement :
\dli $E_0(\phi_{(k_0,\beta_0)})=inf\{\alpha_0\ ,\ (1-P_0(\Omega_\infty))\}$
\dli $E_1[1-\phi_{(k_1,\beta_1)}]=inf\{\alpha_1\ ,\ (1-P_1(\Omega_0))\}$.

$\Delta^B_\alpha$ est constitu\'e des r\`egles de Bol'shev d\'efinies par deux
fonctions de test $\phi_{(k',\beta')}$ et $\phi_{(k'',\beta'')}$
telles que :
\dli $\phi_{(k',\beta')}\leq\phi_{(k'',\beta'')}$ ;
     $\phi_{(k',\beta')}\leq\phi_{(k_0,\beta_0)}$ et
     $\phi_{(k'',\beta'')}\geq\phi_{(k_1,\beta_1)}$.
\medskip
{\parindent=-5mm 1\up{er} cas : $\phi_{(k_0,\beta_0)}\leq
\phi_{(k_1,\beta_1)}$.}

La r\`egle $\delta'_m$ d\'efinie par $\phi_{(k_0,\beta_0)}$ et
$\phi_{(k_1,\beta_1)}$ convient. En effet, pour toute r\`egle $\delta'$
de $\Delta^B_\alpha$ on a pour $\theta\in\{0,1\}$ :
\dli $E_\theta[\delta'(.,\{2\})]=
      E_\theta(\phi_{(k'',\beta'')}-\phi_{(k',\beta')})\geq
      E_\theta(\phi_{(k_1,\beta_1)}-\phi_{(k_0,\beta_0)})=
      E_\theta[\delta'_m(.,\{2\})]$.
\medskip
{\parindent=-5mm 2\up{\`eme} cas : $\phi_{(k_0,\beta_0)}>
\phi_{(k_1,\beta_1)}$.}

Soit $\phi_{(k,\beta)}\in[\phi_{(k_1,\beta_1)}\ ,\ \phi_{(k_0,\beta_0)}]$,
la r\`egle de Bol'shev $\delta'_m$ d\'efinie par $\phi_{(k,\beta)}$ et
$\phi_{(k,\beta)}$ convient. Elle appartient \`a $\Delta^B_\alpha$ et
$E_\theta[\delta'_m(.,\{2\})]=0$ pour $\theta\in\{0,1\}$.
\medskip
{\parindent=-10mm II --- Unicit\'e de $\delta'_m$.}

Dans le premier cas pr\'ec\'edent la r\`egle $\delta'_m$ d\'efinie par
$\phi_{(k_0,\beta_0)}$ et $\phi_{(k_1,\beta_1)}$ est unique presque
s\^urement ; en effet si une autre r\`egle $\delta'$ de $\Delta^B_\alpha$
convient, on a pour $\theta\in\{0,1\}$ :
\dli $E_\theta[\delta'(.,\{2\})]=E_\theta[\delta'_m(.,\{2\})]$ donc 
\dli $E_\theta(\phi_{(k_0,\beta_0)}-\phi_{(k',\beta')})=
      E_\theta(\phi_{(k'',\beta'')}-\phi_{(k_1,\beta_1)})=0$.

Dans le deuxi\`eme cas, $\phi_{(k_0,\beta_0)}>\phi_{(k_1,\beta_1)}$, il
faut montrer qu'il y a unicit\'e presque s\^ure si et seulement si :
$\phi_{(k_0,\beta_0)}\relmont{\leq}{p.s.}\phi_{(k_1,\beta_1)}$, c'est-\`a-dire $E_\theta(\phi_{(k_0,\beta_0)}-\phi_{(k_1,\beta_1)})=0$ pour $\theta\in\{0,1\}$.
\dli C'est une condition n\'ecessaire car l'existence de $\theta\in\{0,1\}$ tel que $E_\theta[\phi_{(k_0,\beta_0)}-\phi_{(k_1,\beta_1)}]>0$ implique que les r\`egles de Bol'shev d\'efinies par $\phi_{(k,\beta)}\in[\phi_{(k_1,\beta_1)}\ ,\ \phi_{(k_0,\beta_0)}]$ ne sont pas presque s\^urement \'egales.
\dli La condition suffisante est \'evidente si l'on montre que toute r\'egle $\delta'$ de $\Delta^B_\alpha$ v\'erifiant $E_\theta[\delta'(.,\{2\})]=0$ pour $\theta\in\{0,1\}$ poss\`ede la propri\'et\'e suivante : 
$\phi_{(k_1,\beta_1)}\leq\phi_{(k',\beta')}\leq\phi_{(k'',\beta'')}\leq\phi_{(k_0,\beta_0)}$. 
$\delta'$ doit bien \^etre de cette forme sinon on aurait 
$\phi_{(k',\beta')}<\phi_{(k_1,\beta_1)}\leq\phi_{(k'',\beta'')}$ et/ou 
$\phi_{(k',\beta')}\leq\phi_{(k_0,\beta_0)}<\phi_{(k'',\beta'')}$, ce qui est impossible puisque la d\'efinition de $\phi_{(k_1,\beta_1)}$ (resp. $\phi_{(k_0,\beta_0)}$) entra\^{\i}nerait 
$0<E_1[\phi_{(k_1,\beta_1)}-\phi_{(k',\beta')}]\leq E_1[\delta'(.,\{2\})]$ 
(resp. $E_0[\delta'(.,\{2\})]>0$).
\medskip\centerline{\hbox to 3cm{\bf \hrulefill}}\par}

La non unicit\'e a lieu lorsqu'il est possible de trouver une r\`egle
$\delta'$ qui d\'ecide toujours $P_0$ ou $P_1$ avec au moins un seuil
strictement inf\'erieur \`a $\alpha_0$ ou $\alpha_1$ ; $\delta'$
v\'erifie pour $\theta\in\{0,1\}$ : 
$E_\theta[\delta'(.,\{2\})]=0$ et 
$R(\theta,\delta')\leq\alpha_\theta$ l'une au moins de ces in\'egalit\'es
\'etant stricte.
\dli On est dans ce cas lorsque les probabilit\'es $P_0$ et $P_1$ sont
suffisamment disjointes pour \^etre s\'epar\'ees aux seuils choisis,
c'est-\`a-dire qu'il existe une fonction de test simple $\phi_{(k,\beta)}$
telle que $E_0(\phi_{(k,\beta)})\leq\alpha_0$ et
$E_1(1-\phi_{(k,\beta)})\leq\alpha_1$ l'une de ces in\'egalit\'es \'etant stricte.
On pourrait donc \^etre plus exigeant en diminuant les
seuils $\alpha_0$ et $\alpha_1$.

Bol'shev se pose le probl\`eme du choix d'une r\`egle de d\'ecision quand il
n'y a pas unicit\'e, ce qui revient \`a choisir un \'el\'ement $\phi_{(k,\beta)}$ de
$[\phi_{(k_1,\beta_1)},\phi_{(k_0,\beta_0)}]$.
Pour cela il consid\`ere une pond\'eration des deux risques :
\dli $R(\lambda,\phi_{(k,\beta)})=(1-\lambda)E_0(\phi_{(k,\beta)})+
\lambda E_1(1-\phi_{(k,\beta)})$.
\dli $\lambda$ peut s'interpr\'eter comme la probabilit\'e a priori accord\'ee
\`a $P_1$. Il se donne un ensemble de pond\'erations possibles,
${\mit\Lambda}\subset [0,1]$ et il cherche $\phi_{(k,\beta)}$ minimisant :
$$sup\{R(\lambda,\phi_{(k,\beta)}) ; \lambda\in{\mit\Lambda}\}=
\left\{\matrix{
R(inf{\mit\Lambda},\phi_{(k,\beta)})\hfill & si &
E_0(\phi_{(k,\beta)})\geq E_1(1-\phi_{(k,\beta)})\hfill \cr
R(sup{\mit\Lambda},\phi_{(k,\beta)})\hfill & si &
E_0(\phi_{(k,\beta)})\leq E_1(1-\phi_{(k,\beta)})\hfill \cr
}\right. $$

\bigskip
\vfill\eject
{\parindent=-5mm 3.2 VOTES DES EXPERTS ET R\`EGLES DE BOL'SHEV.}
\medskip

Nous avons vu que le r\'esultat du vote des experts d\'epend de la probabi\-lit\'e
utilis\'ee pour les faire voter. Les probabilit\'es bas\'ees sur $P_0$ et $P_1$, donnent
les votes les plus diff\'erents. Il est donc tentant d'essayer de prendre
une d\'ecision \`a partir des r\'esultats de ces deux votes, c'est-\`a-dire en
se servant de $Q_0^\omega(\{0\})$, $Q_0^\omega(\{1\})$ et $Q_1^\omega(\{0\})$,
$Q_1^\omega(\{1\})$. Comme
$Q_\theta^\omega(\{0\})+Q_\theta^\omega(\{1\})=1$ et
$E_1(\phi_{(k,\beta)})\geq E_0(\phi_{(k,\beta)})$ (voir fin de 2.2) on a :
$Q_0^\omega(\{1\})\geq Q_1^\omega(\{1\})$ et
$Q_1^\omega(\{0\})\geq Q_0^\omega(\{0\})$. On d\'ecidera $d=1$ (resp. d=0) sans
\'etat d'\^ame lorsque les votes pour $P_1$ (resp. $P_0$) sont jug\'es
pr\'epond\'erants. Pour cela on peut regarder si le vote pour $P_1$
(resp. $P_0$) peut \^etre consid\'er\'e comme un pl\'ebiscite dans le cas le plus
favorable, ce qui revient \`a juger $Q_0^\omega(\{1\})$
(resp. $Q_1^\omega(\{0\})$).
Supposons que $P_1$ (resp. $P_0$) soit pl\'ebiscit\'e lorsque
$Q_0^\omega(\{1\})\geq 1-\alpha_0$ (resp. $Q_1^\omega(\{0\})\geq 1-\alpha_1$),
il y a quatre types de d\'ecision suivant que $P_0$ et $P_1$ sont ou ne sont
pas pl\'ebiscit\'es :
$d=0$ si $P_0$ est pl\'ebiscit\'e alors que $P_1$ ne l'est pas,
$d=1$ si $P_1$ est pl\'ebiscit\'e alors que $P_0$ ne l'est pas,
$d=2$ si ni $P_0$ ni $P_1$ ne sont pl\'ebiscit\'es,
$d=3$ si $P_0$ et $P_1$ sont pl\'ebiscit\'es.
\medskip
{\bf D\'efinition 3.2.1}
\medskip
\medskip
\moveleft 10.4pt\hbox{\vrule\kern 10pt\vbox{\defpro

Soient $\alpha_0<1$ et $\alpha_1<1$.
On appelle r\`egle des pl\'ebiscites, $\delta_p$, la r\`egle de d\'ecision
d\'eterministe \`a valeurs dans $\{0,1,2,3\}$ d\'efinie par :
$$\delta_p(\omega)\quad=\quad\left\{\matrix{
0 & si & Q_0^\omega(\{0\})>\alpha_0\hfill & et
  & Q_1^\omega(\{1\})\leq\alpha_1\hfill\cr
1 & si & Q_0^\omega(\{0\})\leq\alpha_0\hfill & et
  &Q_1^\omega(\{1\})>\alpha_1\hfill\cr
2 & si & Q_0^\omega(\{0\})>\alpha_0\hfill & et
  &Q_1^\omega(\{1\})>\alpha_1\hfill\cr
3 & si & Q_0^\omega(\{0\})\leq\alpha_0\hfill & et
  &Q_1^\omega(\{1\})\leq\alpha_1\hfill\cr
}\right. $$
}}\medskip

La restriction $\alpha_0\not= 1$ et $\alpha_1\not= 1$ est sans cons\'equence
pratique. Elle \'evite la prise en compte de r\`egles $\delta_p$ qui
pl\'ebisciteraient $P_1$ (resp. $P_0$) pour des r\'ealisations $\omega$
dans $\Omega_\infty$ (resp. $\Omega_0$), lorsque $\alpha_0=1$
(resp. $\alpha_1=1$).

$\delta_p$ \'etant d\'eterministe nous allons commencer par la comparer aux
r\`egles $\delta'_m$ de la proposition 3.1.2 quand elles sont presque
s\^urement d\'eterministes. Nous utiliserons les fonctions de test 
$\phi_{(k_0,\beta_0)}$ et $\phi_{(k_1,\beta_1)}$ d\'efinies dans cette
proposition.
\medskip
{\bf Proposition 3.2.1}
\medskip
\medskip
\moveleft 10.4pt\hbox{\vrule\kern 10pt\vbox{\defpro

Soient $\alpha_0<1$ et $\alpha_1<1$.
$\phi_{(k_0,\beta_0)}$ et $\phi_{(k_1,\beta_1)}$ sont
respectivement la plus grande et la plus petite fonction de test de $\Phi_a$ telles que :
\dli $E_0(\phi_{(k_0,\beta_0)})=inf\{\alpha_0\ ,\ (1-P_0(\Omega_\infty))\}$ ;
\dli $E_1[1-\phi_{(k_1,\beta_1)}]=inf\{\alpha_1\ ,\ (1-P_1(\Omega_0))\}$.
\dli  On suppose que $\beta_0$ et $\beta_1$ valent $0$ ou $1$.
\dli a) Si $\phi_{(k_0,\beta_0)}\relmont{\leq}{p.s.}\phi_{(k_1,\beta_1)}$ la r\`egle unique $\delta'_m$ est alors presque s\^urement \'egale \`a la r\`egle des
pl\'ebiscites $\delta_p$.
\dli b) Dans le cas contraire, $\phi_{(k_0,\beta_0)}>\phi_{(k_1,\beta_1)}$ et $\phi_{(k_0,\beta_0)}\relmont{\not=}{p.s.}\phi_{(k_1,\beta_1)}$, la r\`egle des pl\'ebiscites $\delta_p$ d\'ecide $d=0$ (resp. $d=1$) lorsque toutes les r\`egles $\delta'_m$ d\'ecident $d=0$ (resp. $d=1$) et $d=3$ dans le cas contraire.
}}\medskip

{\leftskip=15mm \dli {\bf D\'emonstration}

Les d\'efinitions de $Q^\omega_0(\{0\})$ et $Q^\omega_1(\{1\})$
(voir la d\'efinition 2.4.1) permettent d'obtenir facilement les \'equivalences suivantes :
\dli $Q^\omega_0(\{0\})\leq\alpha_0<1 \Longleftrightarrow
\omega\in\Omega_k$ avec $k<k_0$ ou $k=k_0$ et $\beta_0\geq 1/2$
\dli $Q^\omega_1(\{1\})\leq\alpha_1<1 \Longleftrightarrow
\omega\in\Omega_k$ avec $k>k_1$ ou $k=k_1$ et $\beta_1\leq 1/2$.
\medskip
{\parindent=-5mm D\'emonstration de a) : $\phi_{(k_0,\beta_0)}\relmont{\leq}{p.s.}\phi_{(k_1,\beta_1)}$.}

D'apr\`es la proposition 3.1.2, nous savons que $\delta'_m$ est unique presque s\^urement. 

Si $\phi_{(k_0,\beta_0)}\leq\phi_{(k_1,\beta_1)}$ on peut prendre
pour $\delta'_m$ la r\`egle de Bol'shev d\'efinie par
$\phi_{(k_0,\beta_0)}$ et $\phi_{(k_1,\beta_1)}$. Elle est presque
s\^urement d\'eterministe si et seulement si $\beta_0$ et $\beta_1$
prennent les valeurs $0$ ou $1$. En effet, pour que $\beta_\theta$ appartienne \`a $]0,1[$ il est n\'ecessaire d'avoir : $P_\theta(\Omega_{k_\theta})>0$.
\dli $\delta'_m$ et $\delta_p$ d\'ecident toutes les
deux $d=1$ (resp. $d=0$) sur $\Omega_k$ lorsque $k<k_0$ (resp. $k>k_1$) ou
$k=k_0$ et $\beta_0=1$ (resp. $k=k_1$ et $\beta_1=0$). Sur les autres $\Omega_k$ elles d\'ecident $d=2$.

Il nous reste le cas $\phi_{(k_0,\beta_0)}>\phi_{(k_1,\beta_1)}$ sous la condition $\phi_{(k_0,\beta_0)}\relmont{=}{p.s.}\phi_{(k_1,\beta_1)}$. $\delta'_m$ ne d\'ecide presque jamais $d=2$ (voir le 2\up{\`eme} cas du I de la d\'emonstration de 3.1.2). Comme $\delta_p$ elle d\'ecide $d=1$ (resp. $d=0$) sur $\Omega_k$ lorsque $(k,0)<(k_1,\beta_1)$ (resp. $(k,1)>(k_0,\beta_0)$). $\delta'_m$ et $\delta_p$ sont donc presque s\^urement \'egales. Remarquons que la r\`egle $\delta_p$ ne d\'ecide jamais $d=2$ et presque jamais $d=3$.

\medskip
{\parindent=-5mm D\'emonstration de b) : $\phi_{(k_1,\beta_1)}<
\phi_{(k_0,\beta_0)}$ et $\phi_{(k_0,\beta_0)}\relmont{\not=}{p.s.}\phi_{(k_1,\beta_1)}$.} 

D'apr\`es la fin de la partie II de la d\'emonstration 3.1.2, les r\`egles $\delta'_m$ sont les r\`egles de Bol'shev d\'efinies par
$\phi_{(k',\beta')}$ et $\phi_{(k'',\beta'')}$ tels que : 
$\phi_{(k_1,\beta_1)})\leq\phi_{(k',\beta')}\leq\phi_{(k'',\beta'')}\leq\phi_{(k_0,\beta_0)}]$ et $E_\theta(\phi_{(k'',\beta'')}-\phi_{(k',\beta')})=0$ pour $\theta\in\{0,1\}$.
\dli On a suppos\'e que $\beta_0$ et $\beta_1$ valent $0$ ou $1$. 
Il est facile de voir que les r\`egles $\delta'_m$ d\'ecident toutes $d=1$
(resp. $d=0$) sur $\Omega_k$ lorsque
$k<k_1$ ou $k=k_1$ et $\beta_1=1$
(resp. $k>k_0$ ou $k=k_0$ et $\beta_0=0$). $\delta_p$ fait de m\^eme et
lorsque les r\`egles $\delta'_m$ ne sont pas unanimes elle d\'ecide $d=3$.
\medskip\centerline{\hbox to 3cm{\bf \hrulefill}}\par}

Cette proposition montre que la r\`egle des pl\'ebiscites diff\'erencie les m\^emes
cat\'egories de r\'ealisations que les r\`egles $\delta'_m$ de la proposition
3.1.2 lorsque $\beta_0$ et $\beta_1$ valent $0$ ou $1$.

On a $\beta_0\in]0,1[$ (resp. $\beta_1\in]0,1[$) quand le seuil $\alpha_0$
(resp. $\alpha_1$) ne permet pas de trouver un test de Neyman
d\'eterministe U.P.P. \`a ce seuil lorsque $P_0$ (resp. $P_1$) joue le
r\^ole de $H_0$. La r\`egle des pl\'ebiscites $\delta_p$ \'etant d\'eterministe, on va essayer de la comparer aux r\`egles de Bol'shev d\'eterministes les plus proches des r\`egles $\delta'_m$.
\dli $\delta_m : (\Omega,{\cal A})\rightarrow \{0,1,2\}$ fait partie de ces
r\`egles si il existe $\delta'_m$ qui v\'erifie :
\dli $\forall\omega\in\{\delta_m(\omega)=0\}\cup\{\delta_m(\omega)=1\}
\qquad\delta'_m(\omega,\{\delta_m(\omega)\})\geq 1/2$
\dli $\forall\omega\in\{\delta_m(\omega)=2\}\qquad
\delta'_m(\omega,\{0\})<1/2$ et $\delta'_m(\omega,\{1\})<1/2$.
\dli On peut dire que $\delta_m$ d\'ecide $d\in\{0,1\}$ lorsque la d\'ecision $d$ est majoritaire pour la r\`egle $\delta'_m$. Plusieurs r\`egles $\delta_m$ peuvent correspondre \`a une m\^eme r\`egle $\delta'_m$, il suffit qu'il existe $\omega$ tel que $\delta'_m(\omega,\{0\})=\delta'_m(\omega,\{1\})=1/2$. Bien entendu, si $\delta'_m$ est d\'eterministe la r\`egle $\delta_m$ la plus proche est $\delta_m=\delta'_m$. Il est facile de v\'erifier que la proposition 3.2.1 reste valable lorsque les r\`egles $\delta'_m$ sont remplac\'ees par les r\`egles d\'eterministes $\delta_m$. On peut alors la prolonger en consid\'erant les cas o\`u l'un au moins des param\`etres $\beta_0$ et $\beta_1$ ne  vaut pas $0$ ou $1$.
\medskip
{\bf Proposition 3.2.2}
\medskip
\medskip
\moveleft 10.4pt\hbox{\vrule\kern 10pt\vbox{\defpro

Soient $\alpha_0<1$ et $\alpha_1<1$.
$\phi_{(k_0,\beta_0)}$ et $\phi_{(k_1,\beta_1)}$ sont
respectivement la plus grande et la plus petite fonction de test de $\Phi_a$ telles que :
\dli $E_0(\phi_{(k_0,\beta_0)})=inf\{\alpha_0\ ,\ (1-P_0(\Omega_\infty))\}$ ;
\dli $E_1[1-\phi_{(k_1,\beta_1)}]=inf\{\alpha_1\ ,\ (1-P_1(\Omega_0))\}$.
\dli On suppose avoir $\beta_0\in]0,1[$ ou $\beta_1\in]0,1[$
\dli a) $\phi_{(k_0,\beta_0)}\leq\phi_{(k_1,\beta_1)}$ sans le cas
$(k_0,\beta_0)=(k_1,\beta_1)=(k_0,1/2)$ ;  $\delta_m$
est unique presque s\^urement et \'egale \`a la r\`egle des
pl\'ebiscites $\delta_p$.
\dli b) $\phi_{(k_0,\beta_0)}>\phi_{(k_1,\beta_1)}$ ou
$\phi_{(k_0,\beta_0)}=\phi_{(k_1,\beta_1)}=\phi_{(k_0,1/2)}$ ; 
la r\`egle des pl\'ebiscites $\delta_p$ d\'ecide
$d=0$ (resp. $d=1$) lorsque toutes les r\`egles $\delta_m$ d\'ecident
$d=0$ (resp. $d=1$) et $d=3$ dans le cas contraire.
}}\medskip

La condition  $\phi_{(k_0,\beta_0)}\leq\phi_{(k_1,\beta_1)}$ remplace la condition $\phi_{(k_0,\beta_0)}\relmont{\leq}{p.s.}\phi_{(k_1,\beta_1)}$ de la proposition 3.2.1 car $\beta_0\in]0,1[$ (resp. $\beta_1\in]0,1[$) implique $P_0(\Omega_{k_0})>0$ (resp. $P_1(\Omega_{k_1})>0$) ; on ne peut donc pas avoir $\phi_{(k_0,\beta_0)}>\phi_{(k_1,\beta_1)}$ et $\phi_{(k_0,\beta_0)}\relmont{=}{p.s.}\phi_{(k_1,\beta_1)}$. 

{\leftskip=15mm \dli {\bf D\'emonstration}

{\parindent=-5mm D\'emonstration de a).}

\dli D'apr\`es la proposition 3.1.2, $\delta'_m$ est unique presque s\^urement. 
Elle est d\'efinie par deux fonctions de test $\phi_{(k',\beta')}$ et $\phi_{(k'',\beta'')}$ v\'erifiant : \dli  $\phi_{(k',\beta')}\leq\phi_{(k_0,\beta_0)}\leq\phi_{(k_1,\beta_1)}\leq\phi_{(k'',\beta'')}$, 
\dli $\phi_{(k',\beta')}\relmont{=}{p.s.}\phi_{(k_0,\beta_0)}$ et 
$\phi_{(k'',\beta'')}\relmont{=}{p.s.}\phi_{(k_1,\beta_1)}$.
\dli Lorsque $\beta_0\in]0,1[$ (resp. $\beta_1\in]0,1[$) on a 
$\phi_{(k',\beta')}=\phi_{(k_0,\beta_0)}$ (resp. 
$\phi_{(k'',\beta'')}=\phi_{(k_1,\beta_1)}$). On en d\'eduit facilement que 
$\delta_m$ est aussi presque s\^urement unique si on n'a pas :
$(k_0,\beta_0)=(k_1,\beta_1)=(k_0,1/2)$ ; dans ce cas
$P_\theta(\Omega_{k_0})>0$ et il y a deux $\delta_m$ diff\'erentes
l'une d\'efinie par $\phi_{(k_0,0)}$ l'autre par $\phi_{(k_0,1)}$.

Consid\'erons pour $\delta'_m$ la r\`egle de Bol'shev d\'efinie par
$\phi_{(k_0,\beta_0)}$ et $\phi_{(k_1,\beta_1)}$.
Comme on n'a pas $(k_0,\beta_0)=(k_1,\beta_1)=(k_0,1/2)$, les r\`egles $\delta_p$ et $\delta'_m$ d\'ecident $d=1$ (resp. $d=0$) sur $\Omega_k$ si
et seulement si $k<k_0$ (resp. $k>k_1$) ou
$k=k_0$ et $\beta_0\geq 1/2$ (resp. $k=k_1$ et $\beta_1\leq 1/2$).
La r\`egle $\delta_m$ associ\'ee \`a $\delta'_m$ est alors \'egale \`a $\delta_p$.

\medskip
{\parindent=-5mm D\'emonstration de b).}

Consid\'erons d'abord le cas $\phi_{(k_0,\beta_0)}=\phi_{(k_1,\beta_1)}=\phi_{(k_0,1/2)}$. D'apr\`es la partie a), il existe deux r\`egles $\delta_m$ diff\'erentes l'une d\'efinie par $\phi_{(k_0,0)}$ l'autre par $\phi_{(k_0,1)}$. Elles sont \'egales \`a $\delta_p$ en dehors de $\Omega_{k_0}$ ;
sur $\Omega_{k_0}$ elles sont diff\'erentes et $\delta_p$ d\'ecide $d=3$. La proposition est donc bien v\'erifi\'ee. 

On suppose maintenant $\phi_{(k_1,\beta_1)}<\phi_{(k_0,\beta_0)}$.
\dli Les r\`egles $\delta'_m$ sont les r\`egles de Bol'shev d\'efinies par
$\phi_{(k',\beta')}$ et $\phi_{(k'',\beta'')}$ v\'erifiant : 
$\phi_{(k_1,\beta_1)})\leq\phi_{(k',\beta')}\leq\phi_{(k'',\beta'')}\leq\phi_{(k_0,\beta_0)}$ et $\phi_{(k',\beta')}\relmont{=}{p.s.}\phi_{(k'',\beta'')}$ (voir la fin de la partie II de la d\'emonstration 3.1.2).
\dli $\delta_p$ d\'ecide $d=1$ (resp. $d=0$) sur $\Omega_k$ lorsque
$k<k_1$ ou $k=k_1$ et $\beta_1>1/2$
(resp. $k>k_0$ ou $k=k_0$ et $\beta_0<1/2$) ; dans les autres cas $\delta_p$ d\'ecide $d=3$.
\dli Les $\delta_m$ fournissent bien la m\^eme d\'ecision que $\delta_p$  lorsque $\delta_p$ d\'ecide $d=0$ ou $d=1$. Dans le cas o\`u $\delta_p$ d\'ecide  $d=3$ on trouve des $\delta_m$ qui d\'ecident $d=0$ et d'autres qui d\'ecident $d=1$.
\medskip\centerline{\hbox to 3cm{\bf \hrulefill}}\par}

La r\`egle des pl\'ebiscites $\delta_p$ montre sous une lumi\`ere diff\'erente
les r\`egles $\delta'_m$ s\'electionn\'ees par Bol'shev. Lorsque $\delta'_m$
est unique presque s\^urement ce changement d'\'eclairage est \'equivalent
au passage de la th\'eorie des tests \`a la notion de {\it p-value} vue
sous l'angle d'un vote d'experts.

Dans le cas de non unicit\'e des r\`egles $\delta'_m$ Bol'shev ajoute un crit\`ere
de s\'election pour en garder une seule. Quel que soit le crit\`ere choisi
cette r\`egle ne d\'ecide que $d=0$ ou $d=1$. Nous avons vu que la r\`egle
$\delta_p$ remplace le choix d'un crit\`ere par un nouveau type de d\'ecision,
$d=3$, lorsque les $\delta'_m$ ne sont pas unanimes. $\delta_p$ d\'ecide
$d=3$ sur le m\^eme type de $\Omega_k$ que ceux sur lesquels $\delta'_m$
d\'ecide $d=2$ quand elle est unique. Pour expliquer ceci supposons que
les $\Omega_k$ soient de probabilit\'es nulles, ainsi seules les r\`egles
d\'eterministes sont int\'eressantes, et consid\'erons deux seuils $\alpha_0$
et $\alpha_1$ tels qu'il existe une fonction de test simple $\phi_k$
v\'erifiant : $E_0(\phi_k)=\alpha_0$ et $E_1(\phi_k)=1-\alpha_1$ ;
$\delta'_m$ et $\delta_p$ sont \'egales presque s\^urement et ne d\'ecident
pas $d=2$ ; si l'on diminue les deux seuils elles continueront  \`a \^etre
\'egales presque s\^urement mais d\'ecideront  $d=2$ autour de $\Omega_k$ ;
si l'on augmente les deux seuils il n'y aura plus unicit\'e presque s\^ure
des $\delta'_m$ et autour de $\Omega_k$, $\delta_p$ d\'ecidera  $d=3$.
La d\'ecision $d=3$ intervient d'autant plus facilement que les probabilit\'es
$P_0$ et $P_1$ sont diff\'erentes.
\dli La d\'ecision $d=2$ ou $d=3$ est prise lorsque les fr\'equences des votes
sous $P_1$ et $P_0$, $Q^\omega_0(\{1\})$ et $Q^\omega_1(\{0\})$, sont jug\'ees
de fa\c con semblable : faibles ou fortes.
Il semble int\'eressant de consid\'erer la r\`egle de d\'ecision bas\'ee
sur la diff\'erence $G(\omega)=Q^\omega_1(\{0\})-Q^\omega_0(\{1\})$ qui
varie entre $-1$ et $+1$ ; on ne prendrait aucune d\'ecision autour de $0$ et
on d\'eciderait $d=0$ ou $d=1$ ailleurs suivant que $G(\omega)$ est
positif ou n\'egatif.

\vfill\eject

{\parindent=-10mm\soustitre 4--CHOIX ENTRE DEUX HYPOTH\`ESES STABLES.}
\bigskip
{\parindent=-5mm 4.1 D\'EFINITIONS.}
\medskip
L'ensemble des param\`etres $\Theta$ est partag\'e en deux,
$\Theta=\Theta_0\cup \Theta_1$ avec $\Theta_0\cap \Theta_1=\emptyset$
et $\Theta_i\not=\emptyset$. On consid\`ere le mod\`ele statistique
$(\Omega ,{\cal A},(P_\theta)_{\theta\in\Theta})$ et
l'ensemble des d\'ecisions $D = \{0,1\}$.
Pour $i\in\{0,1\}$, $d=i$ signifie : $\theta$ appartient \`a $\Theta_i$.
Dans la suite pour d\'efinir ce probl\`eme de d\'ecision on \'ecrira simplement le mod\`ele sous la forme :
$(\Omega ,{\cal A},(P_\theta)_{\theta\in\Theta_0\cup \Theta_1})$.
La notation $\relmont{=}{\Theta' p.s.}$ signifiera l'\'egalit\'e presque s\^ure pour la famille de probabilit\'es $\{P_\theta\}_{\theta\in\Theta'}$. Lorsque $\Theta'=\Theta$ on notera simplement $\relmont{=}{p.s.}$.
Une r\`egle de d\'ecision sera d\'efinie par une fonction de test $\phi : (\Omega,{\cal A})\rightarrow\{0,1\}$.

Nous avons d\'ej\`a \'etudi\'e le  cas o\`u $\Theta$ se r\'eduit \`a deux \'el\'ements :
$(\Omega ,{\cal A},(P_\theta)_{\theta\in\{\theta_0\}\cup \{\theta_1\}})$.
Dans ce cas, notre d\'efinition des experts repose sur l'id\'ee qu'un diagnostic $d\in D$ doit avoir plus chance d'\^etre produit quand il est bon que quand il est mauvais.
Cette condition peut s'appliquer \`a chaque couple $(\theta_0,\theta_1)$ de $\Theta_0\times\Theta_1$. Un expert du probl\`eme de d\'ecision $(\Omega ,{\cal A},(P_\theta)_{\theta\in\Theta_0\cup \Theta_1})$ serait alors un expert du choix entre $P_{\theta_0}$ et $P_{\theta_1}$ pour tout $\theta_0\in\Theta_0$ et $\theta_1\in\Theta_1$.
Cette propri\'et\'e est cependant trop forte lorsqu'il existe des \'ev\'enements
de probabilit\'e nulle sous $\theta_0$ (resp.$\theta_1$) alors qu'ils ne
le sont pas pour tous les $\theta$ de $\Theta_0$ (resp. $\Theta_1$).
La restriction au probl\`eme de d\'ecision 
$(\Omega ,{\cal A},(P_\theta)_{\theta\in\{\theta_0\}\cup \{\theta_1\}})$
rend la d\'ecision certaine sur ces \'ev\'enements.
\dli Consid\'erons par exemple le choix
entre $\Theta_1=\{1\}$ et $\Theta_0=\{2,3\}$, $P_\theta$ \'etant la loi uniforme
sur $[\theta -1,\theta+1]$ ; $\phi (x)=\II_{[0,1+\beta]}(x)$, $\beta\in[0,1]$,
d\'efinit un expert du choix entre $\{1\}$ et $\{2\}$ mais pas de celui entre
$\{1\}$ et $\{3\}$ si $\beta <1$ ; en effet les experts de ce choix sont
presque s\^urement \'egaux \`a $\II_{[0,2]}$ ; $\phi$ fait pourtant partie des r\`egles de d\'ecision
int\'eressantes pour le probl\`eme de d\'epart. Ces consid\'erations nous am\`enent \`a
garder des r\`egles de d\'ecision qui ne sont plus des experts pour certains couples $(\theta_0,\theta_1)$,
mais qui dans ce cas d\'ecident $d=1$ (resp. $d=0$) avec la probabilit\'e $1$
sous $\theta_1$ (resp. $\theta_0$).
\medskip
\goodbreak
{\bf D\'efinition 4.1.1}
\nobreak
\medskip
\medskip
\moveleft 10.4pt\hbox{\vrule\kern 10pt\vbox{\defpro

Soient le probl\`eme de d\'ecision $(\Omega ,{\cal A},(P_\theta)_{\theta\in\Theta_0\cup \Theta_1})$ et
une fonction de test $\phi : (\Omega,{\cal A})\rightarrow\{0,1\}$.
\dli $\phi$ d\'efinit un expert du choix entre les deux hypoth\`eses $\theta\in\Theta_0$ et $\theta\in\Theta_1$ s'il poss\`ede les deux propri\'et\'es suivantes : 
\dli i) pour tout $(\theta_0,\theta_1)\in\Theta_0\times\Theta_1$, $\phi$ est un expert du choix entre $P_{\theta_0}$ et $P_{\theta_1}$, sinon il v\'erifie l'une des \'egalit\'es : $E_{\theta_0}(\phi)=0$, $E_{\theta_1}(\phi)=1$ ;
\dli ii) pour tout \'ev\'enement $A$ tel que $\II_A\relmont{=}{\Theta_0 p.s.}0$
(resp. $\II_A\relmont{=}{\Theta_1 p.s.}0$) on a $\phi\II_A\relmont{=}{\Theta_1 p.s.}\II_A$ (resp. $\phi\II_A\relmont{=}{\Theta_0 p.s.}0$).
}}\medskip

La propri\'et\'e i) exprime que pour tout couple $(\theta_0,\theta_1)$, $\phi$ est presque s\^urement ordonnable dans l'ensemble des fonctions de test d\'efinies en 2.2 :  $[\phi_{(0,0)}^{(\theta_0,\theta_1)},\phi_{(\infty,1)}^{(\theta_0,\theta_1)}]$. En effet, si $\phi$ est un expert du choix entre $P_{\theta_0}$ et $P_{\theta_1}$ il est ordonnable dans l'ensemble des fonctions de test simples :  $\Phi_s^{(\theta_0,\theta_1)}=[\phi_{(0,1)}^{(\theta_0,\theta_1)},\phi_{(\infty,0)}^{(\theta_0,\theta_1)}]$ (voir la proposition 2.3.1) et lorsque $E_{\theta_0}(\phi)=0$ (resp. $E_{\theta_1}(\phi)=1$) on a  $\phi\relmont{\leq}{\{\theta_0,\theta_1\} p.s.}\phi_{(0,1)}^{(\theta_0,\theta_1)}$
(resp. $\phi_{(\infty,0)}^{(\theta_0,\theta_1)}\relmont{\leq}{\{\theta_0,\theta_1\} p.s.}\phi$). 

Quant \`a la propri\'et\'e ii) elle traduit l'id\'ee qu'il faut, presque tout le temps, d\'ecider $d=1$ (resp. $d=0$) sur les \'ev\'enements de $P_\theta$ probabilit\'e nulle pour tout $\theta$ de $\Theta_0$ (resp. $\Theta_1$). Cette condition peut \^etre techniquement utile dans le cas o\`u il n'existe pas de support commun \`a toutes les probabilit\'es. 

Dans le cas courant o\`u les probabilit\'es $P_\theta$, $\theta\in\Theta$, admettent des densit\'es strictement positives par rapport \`a une m\^eme mesure, on a pour tous les couples 
$(\theta_0,\theta_1)$ : $\phi_{(0,1)}^{(\theta_0,\theta_1)}=\II_\emptyset$ et $\phi_{(\infty,0)}^{(\theta_0,\theta_1)}=\II_\Omega$.
La d\'efinition pr\'ec\'edente revient \`a dire que $\phi$ est un expert du choix entre $\Theta_0$ et $\Theta_1$ si c'est un expert pour tous les probl\`emes simples : $(\Omega ,{\cal A},(P_\theta)_{\theta\in\{\theta_0\}\cup\{\theta_1\}})$. 
Un expert $\phi$ du mod\`ele $(\Omega ,{\cal A},(P_\theta)_{\theta\in\Theta_0\cup \Theta_1})$ est alors un expert de tout mod\`ele embo\^{\i}t\'e : $(\Omega ,{\cal A},(P_\theta)_{\theta\in\Theta'_0\cup \Theta'_1})$,
$\Theta'_0\subseteq\Theta_0$ et $\Theta'_1\subseteq\Theta_1$.
Le passage d'un mod\`ele embo\^{\i}t\'e au mod\`ele de d\'epart ne peut alors que r\'eduire l'ensemble des experts. Il pourra finir par se restreindre \`a $\II_\emptyset$ et $\II_\Omega$. C'est souvent le cas des probl\`emes bilat\`eres dans un mod\`ele \`a rapport de vraisemblance monotone. Le travail sur les experts sera donc diff\'erent de celui fait sur les r\`egles admissibles dans la th\'eorie de la d\'ecision \`a partir d'une fonction de perte. En effet, les r\`egles admissibles pour un mod\`ele embo\^{\i}t\'e sont g\'en\'eralement encore admissibles dans le mod\`ele de d\'epart. Les r\`egles admissibles sont souvent trop nombreuses, les experts eux sont plut\^ot trop rares.
\dli Le travail sur les r\`egles admissibles est simplifi\'e dans le cas o\`u les probl\`emes de d\'ecision embo\^{\i}t\'es sont sensiblement les m\^emes. Par exemple dans
la th\'eorie des tests, les propri\'et\'es obtenues avec des hypoth\`eses simples
(lemme de Neyman-Pearson principalement) se prolongent
facilement au cas des tests unilat\'eraux dans un mod\`ele \`a rapport de
vraisemblance monotone. Ce qui est important dans ce type de probl\`eme de d\'ecision, c'est que le probl\`eme du choix entre un \'el\'ement $\theta_0$ de
$\Theta_0$ et un \'el\'ement $\theta_1$ de $\Theta_1$ ne change pas
fondamentalement lorsque $\theta_0$ et $\theta_1$ varient.
En fait, pour ces diff\'erents probl\`emes de d\'ecision on peut se res\-treindre
\`a une m\^eme classe de r\`egles de d\'ecision. Nous allons commencer par faire quelque chose de semblable. C'est-\`a-dire analyser le cas  o\`u l'ensemble des experts du mod\`ele global contient les fonctions de test simples des probl\`emes de d\'ecision d\'efinis par les couples $(\theta_0,\theta_1)$.
Pour cela, il faut que les fonctions appartenant aux ensembles  $\Phi_s^{(\theta_0,\theta_1)}$ soient ordonnables. Nous allons le faire en introduisant une notion qui g\'en\'eralise le cas des hypoth\`eses unilat\'erales dans les mod\`eles \`a rapport de vraisemblance monotone.
\medskip
{\bf D\'efinition 4.1.2}
\medskip
\medskip
\moveleft 10.4pt\hbox{\vrule\kern 10pt\vbox{\defpro

Soit $(\Omega ,{\cal A},(P_\theta=p_\theta.\mu)_{\theta\in\Theta_0\cup \Theta_1})$ un probl\`eme de d\'ecision domin\'e par la mesure $\mu$.
Les hypoth\`eses sont stables si il existe une statistique r\'eelle $T$, telle que pour tout couple $(\theta_0,\theta_1)\in\Theta_0\times\Theta_1$ 
il existe une fonction croissante 
$h_{(\theta_0,\theta_1)} :\, \IR \rightarrow \overline{\IR^+}$
v\'erifiant
$p_{\theta_0}/p_{\theta_1}=h_{(\theta_0,\theta_1)}(T)$
sur le domaine de d\'efinition de ce rapport, c'est-\`a-dire en dehors de
$\{\omega\in\Omega\, ;\, p_{\theta_0}(\omega)=p_{\theta_1}(\omega)=0\}$.
}}\medskip

Avant d'\'etudier les experts du choix entre deux hypoth\`eses stables nous allons donner un exemple  de telles hypoth\`eses dans un mod\`ele qui n'est pas \`a rapport de vraisemblance monotone. Il est facile d'en construire puisque notre d\'efinition ne suppose rien sur les rapports $p_{\theta'}/p_{\theta''}$ lorsque $\theta'$ et $\theta''$ appartiennent \`a la m\^eme hypoth\`ese.
Ainsi le probl\`eme de d\'ecision, d\'efini par les deux hypoth\`eses compos\'ees
des probabilit\'es uniformes sur les intervalles de longueur $l\in ]0,2]$
centr\'es sur $1$ pour $\Theta_1$ et sur $2$ pour $\Theta_0$, n'est pas \`a
rapport de vraisemblance monotone, mais les hypoth\`eses sont stables.

\bigskip
\goodbreak
\vfill\eject
{\parindent=-5mm 4.2 ENSEMBLE DES EXPERTS.}
\nobreak
\medskip

Soit $(\Omega ,{\cal A},(P_\theta=p_\theta.\mu)_{\theta\in\Theta_0\cup \Theta_1})$
un probl\`eme de d\'ecision \`a hypoth\`eses stables par rapport \`a la statistique r\'eelle $T$.
\dli Les fonctions de test qui vont jouer un r\^ole fondamental dans ce probl\`eme sont celles d\'efinies \`a partir de l'ordre introduit par la statistique $T$.
C'est-\`a-dire les fonctions de la forme $f_{(t,u)}=\II_{\{T<t\}} + u\II_{\{T=t\}}$. On pose :
$$F=\Bigl\{f_{(t,u)}\, ;\, (t,u)\in\bigl\{(-\infty,1) ; \{(t,u)\}_{t\in\IR, u\in\{0,1\}} ; (+\infty,0)\bigr\}\Bigr\}$$
\dli Les \'el\'ements $(-\infty,1)$ et $(+\infty,0)$ permettent d'avoir dans tous les cas les fonctions de test $\II_\emptyset$ et $\II_\Omega$.
\dli $F$ est totalement ordonn\'e par la relation d'ordre
partiel usuelle sur les fonctions. Nous le noterons parfois
$[f_{(-\infty,1)},f_{(+\infty,0)}]$. Cet ordre co\"{\i}ncide avec l'ordre
lexicographique sur les couples $(t,u)$.
Il permet de d\'efinir les bornes sup\'erieure et inf\'erieure
d'un sous ensemble de $F$ et de munir $F$ de la topologie de l'ordre qui co\"{\i}ncide avec la topologie de la convergence simple.
Sur cet espace l'op\'erateur $E_\theta$, esp\'erance par rapport \`a $P_\theta$,
est croissant et continu. Il peut cependant prendre des valeurs diff\'erentes sur deux \'el\'ements successifs : $f_{(t,0)}$ et $f_{(t,1)}$.

Pour tout couple $(\theta_0,\theta_1)$ de $\Theta_0\times\Theta_1$, on va utiliser l'ensemble des fonctions de test simples bas\'e sur la statistique 
$K=h_{(\theta_0,\theta_1)}(T)$, on le note $\Phi_s^{(\theta_0,\theta_1)}$ (voir la d\'efinition 2.2.1). 
Bien que $K=h_{(\theta_0,\theta_1)}(T)$ soit unique, $P_{\theta_0}$ et $P_{\theta_1}$ presque s\^urement, le choix des $h_{(\theta_0,\theta_1)}$ n'est pas globalement sans cons\'equences. Certaines familles $\{h_{(\theta_0,\theta_1)}\}_{(\theta_0,\theta_1)\in \Theta_0\times\Theta_1}$ facilitent le travail. La d\'efinition suivante donne des propri\'et\'es techniques qui nous seront utiles pour \'ecrire simplement nos propositions.
\medskip
\goodbreak
\vfill\eject
{\bf D\'efinition 4.2.1}
\nobreak
\medskip
\medskip
\moveleft 10.4pt\hbox{\vrule\kern 10pt\vbox{\defpro

Soit $\{h_{(\theta_0,\theta_1)}\}_{(\theta_0,\theta_1)\in \Theta_0\times\Theta_1}$ une famille de fonctions rendant stables les hypoth\`eses $\Theta_0$ et $\Theta_1$ \`a partir de la statistique $T$. 
Notons $D_i=]-\infty,t_i)$ la plus grande demi-droite inf\'erieure ouverte ou ferm\'ee d\'efinissant un \'ev\'enement $T^{-1}(D_i)=\{\omega\in\Omega\, ;\,T(\omega)\in D_i\}$ sur lequel les densit\'es $\{p_{\theta_0}\}_{\theta_0\in\Theta_0}$ sont nulles. De m\^eme 
$D_s=(t_s,+\infty[$ est la plus grande demi-droite sup\'erieure ouverte ou ferm\'ee d\'efinissant un \'ev\'enement $T^{-1}(D_s)$ sur lequel les densit\'es $\{p_{\theta_1}\}_{\theta_1\in\Theta_1}$ sont nulles.

Cette famille est normalis\'ee si chacune des fonctions $h_{(\theta_0,\theta_1)}$ v\'erifie les trois propri\'et\'es suivantes.
\dli i) Si $I$ est un intervalle de $\IR-(D_i\cup D_s)$ d\'efinissant un \'ev\'enement $T^{-1}(I)$ sur lequel $p_{\theta_0}(\omega)=p_{\theta_1}(\omega)=0$,
$h_{(\theta_0,\theta_1)}$ est constante sur $I$.
\dli ii) $h_{(\theta_0,\theta_1)}$ est \'egale \`a $+\infty$ sur $D_s$.
\dli iii) $h_{(\theta_0,\theta_1)}$ est nulle sur $D_i-D_s$.
}}\medskip

Il est facile de montrer qu'\`a partir d'une famille $\{h_{(\theta_0,\theta_1)}\}_{(\theta_0,\theta_1)\in \Theta_0\times\Theta_1}$ on peut toujours construire une famille normalis\'ee (voir annexe II).

Dans la suite nous supposerons toujours travailler avec une famille normalis\'ee.
L'union de tous les $\Phi_s^{(\theta_0,\theta_1)}$, lorsque $(\theta_0,\theta_1)$ parcourt $\Theta_0\times\Theta_1$, est appel\'ee ensemble des fonctions de test simples et not\'e $\Delta_s$. Cet ensemble va jouer le r\^ole de $\Phi_s$ dans le choix entre deux probabilit\'es.
\medskip
{\bf Proposition 4.2.1}
\medskip
\medskip
\moveleft 10.4pt\hbox{\vrule\kern 10pt\vbox{\defpro

Il existe un plus petit intervalle ferm\'e, $\Delta$, de $F$ qui contient l'ensemble $\Delta_s$ des fonctions de test simples.
\dli Les \'el\'ements de $\Delta=[f_{(t_0,u_0)},f_{(t_1,u_1)}]$ d\'efinissent des experts du choix entre $\Theta_0$ et $\Theta_1$. Tout autre expert de ce choix, d\'efini par un \'el\'ement de $F$, est presque s\^urement \'egal \`a $f_{(t_0,u_0)}$ ou $f_{(t_1,u_1)}$.
\dli Les demi-droites $D_i$ et $D_s$ de la d\'efinition 4.2.1 v\'erifient :
$f_{(t_0,u_0)}=\II_{D_i-D_s}\!(T)$ et $f_{(t_1,u_1)}=1-\II_{D_s}\!(T)$.
}}\medskip

\vfill\eject
{\leftskip=15mm \dli {\bf D\'emonstration}
\medskip
{\parindent=-10mm a) --- D\'emonstration de : $\Delta_s\subseteq F$.\par}

Soient $(k,\beta)\in\Bigl\{(0,1) ; \{(k,\beta)\}_{k\in\IR_*^+,\,\beta\in\{0,1\}}; (\infty ,0)\Bigr\}$ et \dli $(\theta_0,\theta_1)\in\Theta_0\times\Theta_1$. La fonction de test simple associ\'ee \`a $(k,\beta)$ dans $\Phi_s^{(\theta_0,\theta_1)}$ est not\'ee : $\phi_{(k,\beta)}$. En utilisant $K=h_{(\theta_0,\theta_1)}(T)$, elle s'\'ecrit $\II_{\{K<k\}}$ ou $\II_{\{K\leq k\}}$ suivant que $\beta=0$ ou $\beta=1$.
\dli Une r\'ealisation $\omega$ appartient \`a $\{K<k\}$ (resp. $\{K\leq k\}$) si et seulement si $T(\omega)$ appartient \`a :
\dli $D_k=\{t\in\IR,\, h_{(\theta_0,\theta_1)}(t)<k\}$ (resp. $D_k=\{t\in\IR,\, h_{(\theta_0,\theta_1)}(t)\leq k\}$).
\dli Lorsque $D_k=\emptyset$ on a $\phi_{(k,\beta)}=\II_\emptyset=f_{(-\infty,1)}$. 
\dli Lorsque $D_k=\IR$ on a $\phi_{(k,\beta)}=\II_\Omega=f_{(+\infty,0)}$. 
\dli Enfin, lorsque $D_k\not=\emptyset$ et $D_k\not=\IR$ il existe $t_k\in\IR$ tel que $D_k=]-\infty,t_k[$ ou $D_k=]-\infty,t_k]$, puisque $h_{(\theta_0,\theta_1)}$ est croissante. On a donc $\phi_{(k,\beta)}=f_{(t_k,0)}$ si $D_k=]-\infty,t_k[$ et $\phi_{(k,\beta)}=f_{(t_k,1)}$ si $D_k=]-\infty,t_k]$. 

\medskip
{\parindent=-10mm b) --- $\Delta=[\II_{D_i-D_s}\!(T),1-\II_{D_s}\!(T)]$.\par}

Posons $f_{(t_0,u_0)}=inf\Delta_s$ et $f_{(t_1,u_1)}=sup\Delta_s$.
\dli $\Delta=[f_{(t_0,u_0)},f_{(t_1,u_1)}]$ est \'evidemment le plus petit intervalle ferm\'e contenant $\Delta_s$.

{\parindent=-5mm i) -- $f_{(t_0,u_0)}=\II_{D_i-D_s}\!(T)$\par}
Lorsque $D_i\cup D_s=\IR$, c'est \'evident puisque $\Delta_s$ est r\'eduit \`a un seul \'el\'ement qui s'\'ecrit $\II_{D_i-D_s}\!(T)$ ou $1-\II_{D_s}\!(T)$. Dans le cas contraire $\II_{D_i-D_s}\!(T)=\II_{D_i}\!(T)$ ; d'apr\`es la d\'efinition 
4.2.1 on a toujours $\II_{D_i}\!(T)\leq\phi_{(0,1)}^{(\theta_0,\theta_1)}$ donc $\II_{D_i}\!(T)\leq f_{(t_0,u_0)}=inf\{\phi_{(0,1)}^{(\theta_0,\theta_1)}\}_{(\theta_0,\theta_1)\in\Theta_0\times\Theta_1}$ ; pour d\'emontrer l'in\'egalit\'e $\II_{D_i}\!(T)\geq f_{(t_0,u_0)}$, d\'efinissons la demi-droite $D_0$ par $\II_{D_0}\!(T)=f_{(t_0,u_0)}$  et montrons qu'elle est incluse dans $D_i$ ; pour tout $(\theta_0,\theta_1)$ la fonction $h_{(\theta_0,\theta_1)}$ est nulle sur $D_0$ puisque $f_{(t_0,u_0)}\leq\phi_{(0,1)}^{(\theta_0,\theta_1)}$, les densit\'es $p_{\theta_0}$ sont donc nulles sur $T^{-1}(D_0)$, ce qui implique 
bien $D_0\subseteq D_i$.

{\parindent=-5mm ii) -- $f_{(t_1,u_1)}=1-\II_{D_s}\!(T)$\par}
La d\'emonstration est semblable \`a celle de i).
D'apr\`es la propri\'et\'e ii) de la d\'efinition 
4.2.1 on a toujours $\phi_{(\infty,0)}^{(\theta_0,\theta_1)}\leq 1-\II_{D_s}\!(T)$ donc $1-\II_{D_s}\!(T)\geq f_{(t_1,u_1)}=sup\{\phi_{(\infty,0)}^{(\theta_0,\theta_1)}\}_{(\theta_0,\theta_1)\in\Theta_0\times\Theta_1}$ ; pour d\'emontrer l'in\'egalit\'e $1-\II_{D_s}\!(T)\leq f_{(t_1,u_1)}$, d\'efinissons $D_1$ par $\II_{D_1}\!(T)=f_{(t_1,u_1)}$  et montrons que la demi-droite $\IR-D_1$ est incluse dans $D_s$ ; pour tout $(\theta_0,\theta_1)$ la fonction $h_{(\theta_0,\theta_1)}$ est infinie sur $\IR-D_1$ puisque $f_{(t_1,u_1)}\geq\phi_{(\infty,0)}^{(\theta_0,\theta_1)}$, les densit\'es $p_{\theta_1}$ sont donc nulles sur $T^{-1}(\IR-D_1)$, ce qui implique 
bien $\IR-D_1\subseteq D_s$.

\medskip
{\parindent=-10mm c) --- Les experts de $F$ sont presque s\^urement \'egaux \`a un \'el\'ement de $\Delta$.\par}

En fait on va d\'emontrer quelque chose de plus fort,
en prenant un \'el\'ement $f_{(t,u)}$ de $F-\Delta$ et en montrant que s'il d\'efinit un expert il v\'erifie : $f_{(t,u)}\relmont{=}{p.s.}f_{(t_0,u_0)}$  ou $f_{(t,u)}\relmont{=}{p.s.}f_{(t_1,u_1)}$.
\dli Pour tout couple $(\theta_0,\theta_1)$ de $\Theta_0\times\Theta_1$, on note $F^{(\theta_0,\theta_1)}$ l'intervalle d\'efini par $[\phi_{(0,1)}^{(\theta_0,\theta_1)},\phi_{(\infty,0)}^{(\theta_0,\theta_1)}]$ dans $F$. On a bien s\^ur : $\Phi_s^{(\theta_0,\theta_1)}\subseteq F^{(\theta_0,\theta_1)}$ et
$f_{(t,u)}$ est en dehors de l'intervalle $F^{(\theta_0,\theta_1)}$ pour tout couple $(\theta_0,\theta_1)$.
\dli Lorsque $f_{(t,u)}\leq f_{(t_0,u_0)}$ on a $f_{(t,u)}\relmont{=}{\Theta_0 p.s.}f_{(t_0,u_0)}$ puisque $E_{\theta_0}(f_{(t_0,u_0)})=0$ pour tout $\theta_0$ de $\Theta_0$.
Si $f_{(t,u)}$ est un expert, le fait que l'\'ev\'enement $A=\{f_{(t_0,u_0)}=1\}$ soit de $P_{\theta_0}$ probabilit\'e nulle pour tout $\theta_0$ de $\Theta_0$ entra\^{\i}ne : $f_{(t,u)}.\II_A\relmont{=}{\Theta_1 p.s.}\II_A$ (voir la propri\'et\'e ii) de la d\'efinition 4.1.1). On en d\'eduit la propri\'et\'e recherch\'ee : $f_{(t,u)}\relmont{=}{p.s.}f_{(t_0,u_0)}$. 
\dli Le cas $f_{(t,u)}\geq f_{(t_1,u_1)}$ se traite de fa\c con semblable. On a $E_{\theta_1}(f_{(t_1,u_1)})=1$ pour tout $\theta_1$ de $\Theta_1$ et donc $f_{(t,u)}\relmont{=}{\Theta_1 p.s.}f_{(t_1,u_1)}$. L'\'ev\'enement $A=\{f_{(t_1,u_1)}=0\}$ \'etant de $P_{\theta_1}$ probabilit\'e nulle pour tout $\theta_1$ de $\Theta_1$, si $f_{(t,u)}$ est un expert, la propri\'et\'e ii) de la d\'efinition 4.1.1 entra\^{\i}ne : $f_{(t,u)}.\II_A\relmont{=}{\Theta_0 p.s.}0$. On a encore la propri\'et\'e recherch\'ee : $f_{(t,u)}\relmont{=}{p.s.}f_{(t_1,u_1)}$. 

\medskip
{\parindent=-10mm d) --- Tout \'el\'ement $f_{(t,u)}$ de $\Delta$ est un expert.\par}

{\parindent=-5mm i) -- D\'emonstration de la propri\'et\'e i) de la d\'efinition 4.1.1.\par}

Soit $(\theta_0,\theta_1)\in\Theta_0\times\Theta_1$. On doit d\'emontrer que 
$f_{(t,u)}$ est un expert du choix entre $P_{\theta_0}$ et $P_{\theta_1}$ ou qu'il v\'erifie : $E_{\theta_0}(f_{(t,u)})=0$ ou $E_{\theta_1}(f_{(t,u)})=1$.

1\up{er} cas : $f_{(t,u)}\in F^{(\theta_0,\theta_1)}=[\phi_{(0,1)}^{(\theta_0,\theta_1)},\phi_{(\infty,0)}^{(\theta_0,\theta_1)}]\subseteq F$.
\dli Posons $k=h_{(\theta_0,\theta_1)}(t)$. Si $k\in\IR_*^+$ on a $\phi_{(k,0)}^{(\theta_0,\theta_1)}\leq f_{(t,u)}\leq\phi_{(k,1)}^{(\theta_0,\theta_1)}$ et d'apr\`es la proposition 2.3.1 $f_{(t,u)}$ est un expert du choix entre $P_{\theta_0}$ et $P_{\theta_1}$. Lorsque $k=0$ (resp. $k=+\infty$) la construction de
$\Phi_s^{(\theta_0,\theta_1)}$ \`a partir de $K=h_{(\theta_0,\theta_1)}(T)$ entra\^{\i}ne : $f_{(t,u)}\leq\phi_{(0,1)}^{(\theta_0,\theta_1)}$ (resp. $f_{(t,u)}\geq\phi_{(\infty,0)}^{(\theta_0,\theta_1)}$) ; comme $f_{(t,u)}$ appartient \`a $F^{(\theta_0,\theta_1)}$ on a
$f_{(t,u)}=\phi_{(0,1)}^{(\theta_0,\theta_1)}$ (resp. $f_{(t,u)}=\phi_{(\infty,0)}^{(\theta_0,\theta_1)}$); $f_{(t,u)}$ est donc encore un expert du choix entre $P_{\theta_0}$ et $P_{\theta_1}$. 

2\up{\`eme} cas : $f_{(t,u)}\notin F^{(\theta_0,\theta_1)}$.
\dli On a alors, soit $f_{(t,u)}<\phi_{(0,1)}^{(\theta_0,\theta_1)}$ , soit 
$f_{(t,u)}>\phi_{(\infty,0)}^{(\theta_0,\theta_1)}$ et donc, soit $E_{\theta_0}(f_{(t,u)})=0$, soit $E_{\theta_1}(f_{(t,u)})=1$. La propri\'et\'e i) est v\'erifi\'ee.

{\parindent=-5mm ii) -- D\'emonstration de la propri\'et\'e ii) de la d\'efinition 4.1.1.\par}

1\up{er} cas : soit un \'ev\'enement $A$ tel que $\II_A\relmont{=}{\Theta_0 p.s.}0$.

\dli On doit obtenir l'\'egalit\'e : $f_{(t,u)}.\II_A\relmont{=}{\Theta_1 p.s.}\II_A$.  D\'ecomposons $A$ en $A'=A\cap\{f_{(t,u)}=0\}$ et $A''=A\cap\{f_{(t,u)}=1\}$. Pour avoir l'\'egalit\'e recherch\'ee on doit d\'emontrer $P_{\theta_1}(A')=0$ pour tout $\theta_1$ dans $\Theta_1$.

Soient $\theta_1\in\Theta_1$ et $B=A'\cap\{p_{\theta_1}>0\}$.
$P_{\theta_1}(A')$ est nulle si $P_{\theta_1}(B)$ l'est. 

\dli Consid\'erons le cas : $f_{(t,u)}>f_{(t_0,u_0)}$.
\dli D'apr\`es le lemme 2 de l'annexe II, il existe $\theta_0\in\Theta_0$  
tel que $\phi_{(0,1)}^{(\theta_0,\theta_1)}\leq f_{(t,u)}$. Posons 
$B'=B\cap\{p_{\theta_0}=0\}$ et $B''=B\cap\{p_{\theta_0}>0\}$. Comme $B''\subseteq A$, on a  $P_{\theta_0}(B'')=0$, donc $\mu (B'')=0$ et $P_{\theta_1}(B'')=0$. Pour finir on va d\'emontrer que $B'$ est vide. 
\dli S'il  
existait $\omega\in B'$ on aurait $h_{(\theta_0,\theta_1)}(T(\omega))=0$ car $B'$ est inclus dans $B$. Ce qui impliquerait $\phi_{(0,1)}^{(\theta_0,\theta_1)}(\omega)=1$ donc $f_{(t,u)}(\omega)=1$ ; 
ceci est impossible puisque $\omega\in B'\subseteq A'$.
\dli Il nous reste le cas : $f_{(t,u)}=f_{(t_0,u_0)}$.
\dli Si $f_{(t_0,u_0)}=f_{(t_1,u_1)}$ on a $P_{\theta_1}(\{f_{(t,u)}=0\})=0$ donc $P_{\theta_1}(A')=0$.
Sinon il existe dans $\Delta -\{f_{(t_0,u_0)}\}$, une suite $\{g_n\}_{n\in\IN}$ d\'ecroissant vers $f_{(t_0,u_0)}$ ;
comme $g_n>f_{(t_0,u_0)}$ on a la propri\'et\'e recherch\'ee sur les $g_n$ : $g_n.\II_A\relmont{=}{\theta_1 p.s.}\II_A$ ; ce qui entra\^{\i}ne bien 
$f_{(t,u)}.\II_A\relmont{=}{\theta_1 p.s.}\II_A$.

2\up{\`eme} cas : soit un \'ev\'enement $A$ tel que $\II_A\relmont{=}{\Theta_1 p.s.}0$. 

\dli On doit obtenir l'\'egalit\'e : $f_{(t,u)}.\II_A\relmont{=}{\Theta_0 p.s.}0$. La d\'emonstration est semblable \`a la pr\'ec\'edente.
On d\'ecompose $A$ en $A'=A\cap\{f_{(t,u)}=1\}$ et $A''=A\cap\{f_{(t,u)}=0\}$. Pour avoir l'\'egalit\'e recherch\'ee on doit d\'emontrer $P_{\theta_0}(A')=0$ pour tout $\theta_0$ dans $\Theta_0$.

Soient $\theta_0\in\Theta_0$ et $B=A'\cap\{p_{\theta_0}>0\}$.
$P_{\theta_0}(A')$ est nulle si $P_{\theta_0}(B)$ l'est. 

\dli Consid\'erons le cas : $f_{(t,u)}<f_{(t_1,u_1)}$.
\dli D'apr\`es le lemme 1 de l'annexe II, il existe $\theta_1\in\Theta_1$  
tel que $\phi_{(\infty,0)}^{(\theta_0,\theta_1)}\geq f_{(t,u)}$. Posons 
$B'=B\cap\{p_{\theta_1}=0\}$ et $B''=B\cap\{p_{\theta_1}>0\}$. Comme $B''\subseteq A$, on a  $P_{\theta_1}(B'')=0$, donc $\mu (B'')=0$ et $P_{\theta_0}(B'')=0$. Pour finir on va d\'emontrer que $B'$ est vide. S'il  
existait $\omega\in B'$ on aurait $h_{(\theta_0,\theta_1)}(T(\omega))=+\infty$ car $B'$ est inclus dans $B$. Ce qui impliquerait $\phi_{(\infty,0)}^{(\theta_0,\theta_1)}(\omega)=0$ donc $f_{(t,u)}(\omega)=0$ ; 
ceci est impossible puisque $\omega\in B'\subseteq A'$.
\dli Il nous reste le cas : $f_{(t,u)}=f_{(t_1,u_1)}$.
\dli Si $f_{(t_0,u_0)}=f_{(t_1,u_1)}$ on a $P_{\theta_0}(\{f_{(t,u)}=1\})=0$ donc $P_{\theta_0}(B)=0$.
Sinon il existe dans $\Delta -\{f_{(t_1,u_1)}\}$, une suite $\{g_n\}_{n\in\IN}$ croissant vers $f_{(t_1,u_1)}$ ;
comme $g_n<f_{(t_1,u_1)}$ on a la propri\'et\'e recherch\'ee sur les $g_n$ : $g_n.\II_A\relmont{=}{\theta_0 p.s.}0$ ; ce qui entra\^{\i}ne bien 
$f_{(t,u)}.\II_A\relmont{=}{\theta_0 p.s.}0$.

\medskip\centerline{\hbox to 3cm{\bf \hrulefill}}\par}

Nous venons de trouver les experts du choix entre $\Theta_0$ et $\Theta_1$
qui sont d\'efinis par un \'el\'ement de $F$. Nous allons maintenant consid\'erer l'ensemble des experts. Dans le choix entre deux probabilit\'es, l'ensemble des fonctions de test simples, $\Phi_s$, a jou\'e un r\^ole fondamental. Les experts n'\'etant rien d'autre que des r\`egles presque s\^urement ordonnables dans $\Phi_s$ (voir proposition 2.3.1). Nous allons g\'en\'eraliser la proposition 2.3.1 au choix entre deux hypoth\`eses stables. C'est l'adh\'erence de $\Delta_s$ dans $F$ qui jouera le r\^ole de $\Phi_s$.

\medskip
{\bf Proposition 4.2.2}
\medskip
\medskip
\moveleft 10.4pt\hbox{\vrule\kern 10pt\vbox{\defpro

Soit $(\Omega ,{\cal A},(p_\theta.\mu)_{\theta\in\Theta_0\cup\Theta_1})$ un 
probl\`eme de d\'ecision \`a hypoth\`eses stables. On note $\Delta_s$ l'ensemble des fonctions de test simples d\'efinies \`a partir d'une famille
$\{h_{(\theta_0,\theta_1)}\}_{(\theta_0,\theta_1)\in \Theta_0\times\Theta_1}$ normalis\'ee.
\dli Les experts du choix entre $\Theta_0$ et $\Theta_1$  sont les r\`egles de d\'ecision $\phi : (\Omega,{\cal A})\rightarrow\{0,1\}$ presque s\^urement ordonnables dans l'adh\'erence $\overline{\Delta_s}$ de $\Delta_s$. C'est-\`a-dire qu'il existe deux \'el\'ements de $\overline{\Delta_s}$, $f$ et $f'$, tels que : $f\relmont{\leq}{p.s.}\phi\relmont{\leq}{p.s.}f'$ et $]f,f'[\cap\overline{\Delta_s}=\emptyset$. 
}}\medskip

{\leftskip=15mm \dli {\bf D\'emonstration}

\medskip
{\parindent=-10mm I ----- Condition suffisante.}

Soit une r\`egle de d\'ecision $\phi$ pour laquelle il existe $f$ et $f'$ dans $\overline{\Delta_s}$ tels que : $f\relmont{\leq}{p.s.}\phi\relmont{\leq}{p.s.}f'$ et $]f,f'[\cap\overline{\Delta_s}=\emptyset$. D\'emontrons que $\phi$ est un expert.
\medskip
{\parindent=-5mm a) -- D\'emonstration de la propri\'et\'e i) de 4.1.1.}

Soit $(\theta_0,\theta_1)\in \Theta_0\times\Theta_1$. 
Si $\phi_{(\infty,0)}^{(\theta_0,\theta_1)}\leq f$ on a $E_{\theta_1}(\phi)=1$.  
Si $\phi_{(0,1)}^{(\theta_0,\theta_1)}\geq f'$ on a $E_{\theta_0}(\phi)=0$. 
Dans le cas contraire nous allons montrer que $\phi$ est un expert du choix entre $P_{\theta_0}$ et $P_{\theta_1}$.
\dli Comme $]f,f'[\cap\Delta_s=\emptyset$, $f<\phi_{(\infty,0)}^{(\theta_0,\theta_1)}$ entra\^{\i}ne $f'\leq\phi_{(\infty,0)}^{(\theta_0,\theta_1)}$ et $\phi_{(0,1)}^{(\theta_0,\theta_1)}<f'$ implique $\phi_{(0,1)}^{(\theta_0,\theta_1)}\leq f$
; nous avons donc :
\dli $\phi_{(0,1)}^{(\theta_0,\theta_1)}\leq f\relmont{\leq}{p.s.}\phi\relmont{\leq}{p.s.}f'\leq
\phi_{(\infty,0)}^{(\theta_0,\theta_1)}$.
\dli Consid\'erons 
$G=\{g\in\Phi_s^{(\theta_0,\theta_1)}\, ;\,g\leq f\}$ et
$G'=\{g\in\Phi_s^{(\theta_0,\theta_1)}\, ;\,g\geq f'\}$ ; $G$ et $G'$ forment une partition de 
$\Phi_s^{(\theta_0,\theta_1)}$ en deux intervalles non vides.

1\up{er} cas : $sup G=\phi_{(k,0)}^{(\theta_0,\theta_1)}$.
\dli $G'$ n'\'etant pas vide on ne peut pas avoir $k=\infty$. On obtient alors l'encadrement suivant :
$\phi_{(k,0)}^{(\theta_0,\theta_1)}\leq f\relmont{\leq}{p.s.}\phi\relmont{\leq}{p.s.}f'\leq
\phi_{(k,1)}^{(\theta_0,\theta_1)}$ ; d'apr\`es
la proposition 2.3.1, $\phi$ est un expert du choix entre $P_{\theta_0}$ et $P_{\theta_1}$.

2\up{\`eme} cas : $sup G=\phi_{(k,1)}^{(\theta_0,\theta_1)}$.
\dli Pour tout $k'>k$ on a :
$\phi_{(k,1)}^{(\theta_0,\theta_1)}\leq f\relmont{\leq}{p.s.}\phi\relmont{\leq}{p.s.}f'\leq
\phi_{(k',0)}^{(\theta_0,\theta_1)}$ ; ceci entra\^{\i}ne 
$f=f'=\phi_{(k,1)}^{(\theta_0,\theta_1)}$ et $\phi\relmont{=}{p.s.}\phi_{(k,1)}^{(\theta_0,\theta_1)}$ 
est un expert du choix entre $P_{\theta_0}$ et $P_{\theta_1}$.

\medskip
{\parindent=-5mm b) -- D\'emonstration de la propri\'et\'e ii) de 4.1.1.}

1\up{er} cas : soit un \'ev\'enement $A$ tel que $\II_A\relmont{=}{\Theta_0 p.s.}0$.

\dli On doit obtenir l'\'egalit\'e : $\phi.\II_A\relmont{=}{\Theta_1 p.s.}\II_A$.
\dli Par hypoth\`ese on a $f\relmont{\leq}{p.s.}\phi$ donc
$f.\II_A\relmont{\leq}{p.s.}\phi.\II_A$. De plus $f$ est un \'el\'ement de $\overline{\Delta_s}\subseteq\Delta$, d'apr\`es la proposition 4.2.1 c'est un expert ; il v\'erifie donc : 
$f.\II_A\relmont{=}{\Theta_1 p.s.}\II_A$ (partie ii) de la d\'efinition 4.1.1). L'\'egalit\'e recherch\'ee est bien v\'erifi\'ee.

2\up{\`eme} cas : soit un \'ev\'enement $A$ tel que $\II_A\relmont{=}{\Theta_1 p.s.}0$.

\dli On doit obtenir l'\'egalit\'e : $\phi.\II_A\relmont{=}{\Theta_0 p.s.}0$. 
\dli Par hypoth\`ese $\phi\relmont{\leq}{p.s.}f'$, donc $\phi.\II_A\relmont{\leq}{p.s.}f'.\II_A$. 
$f'$ appartenant \`a $\Delta$, c'est un expert ; d'apr\`es la partie ii) de la d\'efinition 4.1.1 il v\'erifie : 
$f'.\II_A\relmont{=}{\Theta_0 p.s.}0$. L'\'egalit\'e recherch\'ee est encore v\'erifi\'ee.

\medskip
{\parindent=-10mm II --- Condition n\'ecessaire.}

Soit $\phi$ un expert du choix entre $\Theta_0$ et $\Theta_1$.
On doit trouver deux \'el\'ements $f$ et $f'$ de $\overline{\Delta_s}$, tels que : $f\relmont{\leq}{p.s.}\phi\relmont{\leq}{p.s.}f'$ et $]f,f'[\cap\overline{\Delta_s}=\emptyset$. 
\dli Pour cela nous utiliserons 
$G=\{g\in\overline{\Delta_s}\, ;\,g\relmont{\leq}{p.s.}\phi\}$ et \dli $G'=\{g\in\overline{\Delta_s}\, ;\,g\relmont{\geq}{p.s.}\phi\}$.
Ce sont des intervalles ordonn\'es de $\overline{\Delta_s}$. Nous allons commencer par \'etablir quelques propri\'et\'es de $G$ et $G'$. 
\medskip
{\parindent=-5mm a) -- $G$ est non vide.}

Nous allons montrer qu'il contient $f_{(t_0,u_0)}=inf \Delta$. Pour cela il faut obtenir $f_{(t_0,u_0)}\relmont{\leq}{\theta p.s.}\phi$ quel que soit $\theta$ de $\Theta$.

1\up{er} cas : Soit $\theta\in\Theta_0$.
\dli Par d\'efinition de $f_{(t_0,u_0)}$ on a $E_\theta(f_{(t_0,u_0)})=0$, c'est-\`a-dire $f_{(t_0,u_0)}\relmont{=}{\theta p.s.}0$ et donc :
$\phi\relmont{\geq}{\theta p.s.}f_{(t_0,u_0)}$.    

2\up{\`eme} cas : Soit $\theta\in\Theta_1$.
\dli Consid\'erons $A=\{f_{(t_0,u_0)}=1\}$, on a 
$\II_A\relmont{=}{\Theta_0 p.s.}0$ ; la propri\'et\'e ii) de la d\'efinition 4.1.1 entra\^{\i}ne : $\phi.\II_A\relmont{=}{\Theta_1 p.s.}\II_A$. Comme 
$\II_A=f_{(t_0,u_0)}$, on obtient bien
$\phi\relmont{\geq}{\theta p.s.}f_{(t_0,u_0)}$.    

\medskip
{\parindent=-5mm b) -- $G'$ est non vide.}

Nous allons montrer qu'il contient $f_{(t_1,u_1)}=sup \Delta$. Pour cela il faut obtenir $\phi\relmont{\leq}{\theta p.s.}f_{(t_1,u_1)}$ quel que soit $\theta$ de $\Theta$.

1\up{er} cas : Soit $\theta\in\Theta_1$.
\dli Par d\'efinition de $f_{(t_1,u_1)}$ on a $E_\theta(f_{(t_1,u_1)})=1$, c'est-\`a-dire $f_{(t_1,u_1)}\relmont{=}{\theta p.s.}1$ et donc :
$\phi\relmont{\leq}{\theta p.s.}f_{(t_1,u_1)}$.    

2\up{\`eme} cas : Soit $\theta\in\Theta_0$.
\dli Consid\'erons $A=\{f_{(t_1,u_1)}=0\}$, on a 
$\II_A\relmont{=}{\Theta_1 p.s.}0$ ; la propri\'et\'e ii) de la d\'efinition 4.1.1 entra\^{\i}ne : $\phi.\II_A\relmont{=}{\Theta_0 p.s.}0$. Comme 
$(1-\II_A)=f_{(t_1,u_1)}$, on obtient bien
$\phi\relmont{\leq}{\theta p.s.}f_{(t_1,u_1)}$.

\medskip
{\parindent=-5mm c) -- $G$ et $G'$ sont ferm\'es.}

Pour montrer que $G$ est ferm\'e il suffit de montrer que $g_s=sup G$ appartient \`a $G$. Posons $A=\{g_s>\phi\}$, on doit d\'emontrer que 
$P_\theta(A)=0$ pour tout $\theta\in\Theta$.
\dli On consid\`ere une suite $(g_n)_{n\in\IN}$ de $G$ croissant vers $g_s$.
Les \'ev\'enements $A_n=\{g_n>\phi\}$ croissent vers $A$ ; comme 
$g_n\relmont{\leq}{p.s.}\phi$, on a $P_\theta(A_n)=0$ et donc $P_\theta(A)=0$ pour tout $\theta\in\Theta$.

De m\^eme pour montrer que $G'$ est ferm\'e, on d\'emontre que $g_i=inf G'$ appartient \`a $G'$. On pose $A=\{g_i<\phi\}$
et on consid\`ere une suite $(g_n)_{n\in\IN}$ de $G'$ d\'ecroissant vers $g_i$.
Les \'ev\'enements $A_n=\{g_n<\phi\}$ croissent vers $A$ ; comme 
$g_n\relmont{\geq}{p.s.}\phi$, on a $P_\theta(A_n)=0$ et donc $P_\theta(A)=0$ pour tout $\theta\in\Theta$.

\medskip
{\parindent=-5mm d) -- Existence de $f$ et $f'$ ayant les deux propri\'et\'es requises.}

On utilise les deux bornes $g_s=sup G$ et $g_i=inf G'$. Si $g_i\leq g_s$
on a $\phi\relmont{\leq}{p.s.}g_i\leq g_s\relmont{\leq}{p.s.}\phi$ ; 
le choix $f=f'=g$ avec $g\in [g_i,g_s]$ convient ; $\phi$ est presque s\^urement \'egal \`a un \'el\'ement de $\overline{\Delta_s}$.
\dli Il nous reste \`a trouver $f$ et $f'$ lorsque $g_s<g_i$.
On a, par construction, 
$g_s\relmont{\leq}{p.s.}\phi\relmont{\leq}{p.s.}g_i$. $f=g_s$ et $f'=g_i$
conviennent si $]g_s,g_i[\cap\overline{\Delta_s}=\emptyset$.
Il suffit de v\'erifier : $]g_s,g_i[\cap\Delta_s=\emptyset$. Pour cela on consid\`ere un couple $(\theta_0,\theta_1)$ et on cherche \`a montrer : 
$]g_s,g_i[\cap\Phi_s^{(\theta_0,\theta_1)}=\emptyset$. La propri\'et\'e i) de la d\'efinition 4.1.1 nous conduit \`a distinguer trois cas.

1\up{er} cas : $\phi$ est un expert du choix entre $\theta_0$ et $\theta_1$.
\dli Nous allons raisonner par l'absurde en supposant qu'il existe $\phi_{(k,\beta)}^{(\theta_0,\theta_1)}$ dans $]g_s,g_i[$. $\phi$ \'etant un expert du choix entre $\theta_0$ et $\theta_1$, il est ordonnable par rapport \`a $\phi_{(k,\beta)}^{(\theta_0,\theta_1)}$
(voir proposition 2.3.1).  

\dli 1\up{\`ere} possibilit\'e : $\phi\relmont{<}{\{\theta_0,\theta_1\} p.s.}\phi_{(k,\beta)}^{(\theta_0,\theta_1)}<g_i$.
\dli Rappelons que $\relmont{<}{\{\theta_0,\theta_1\} p.s.}$ signifie $\relmont{\leq}{\{\theta_0,\theta_1\} p.s.}$ sans que l'on ait $\relmont{=}{\{\theta_0,\theta_1\} p.s.}$.
\dli Nous allons commencer par analyser ce qui peut emp\^echer $\phi_{(k,\beta)}^{(\theta_0,\theta_1)}$ d'appartenir \`a $G'$ c'est-\`a-dire de v\'erifier : $\phi\relmont{\leq}{p.s.}\phi_{(k,\beta)}^{(\theta_0,\theta_1)}$.
\dli Posons $A=\{g_i-\phi_{(k,\beta)}^{(\theta_0,\theta_1)}=1\}\cap\{\phi=1\}=\{\phi_{(k,\beta)}^{(\theta_0,\theta_1)}\not= g_i\}\cap\{\phi=1\}$.
Comme $\phi\relmont{\leq}{p.s.}g_i$, on a
$\phi.\II_{A^c}\relmont{\leq}{p.s.}\phi_{(k,\beta)}^{(\theta_0,\theta_1)}.\II_{A^c}$. Par d\'efinition de $A$ on a 
$\phi_{(k,\beta)}^{(\theta_0,\theta_1)}.\II_A=0$, ceci implique
$\phi.\II_A\relmont{=}{\{\theta_0,\theta_1\} p.s.}0$, ce qui peut s'\'ecrire
$P_{\theta_0}(A)=P_{\theta_1}(A)=0$ ;
les \'ev\'enements $A\cap\{p_{\theta_0}>0\}$ et $A\cap\{p_{\theta_1}>0\}$
sont alors $\mu$ n\'egligeables et donc $P_{\theta}$ n\'egligeables pour tout
$\theta$ ; la seule partie de $A$ qui peut ne pas \^etre $P_{\theta}$ n\'egligeable, pour tout $\theta$, est $B=A\cap\{p_{\theta_0}=0\}\cap\{p_{\theta_1}=0\}$ ; on a donc d\'ej\`a :
$\phi.\II_{B^c}\relmont{\leq}{p.s.}\phi_{(k,\beta)}^{(\theta_0,\theta_1)}.\II_{B^c}$. On distingue maintenant deux cas, qui conduiront \`a des contradictions diff\'erentes.
\dli i) Pour tout $\theta\in\Theta_1$ : $P_{\theta}(B)=0$. La partie ii) de la d\'efinition 4.1.1 implique : $\phi.\II_B\relmont{=}{\Theta_0 p.s.}0$.
$\phi$ \'etant \'egal \`a $1$ sur $B$, cet \'ev\'enement est $P_{\theta}$ n\'egligeable pour tout $\theta$ et on a  alors 
$\phi\relmont{\leq}{p.s.}\phi_{(k,\beta)}^{(\theta_0,\theta_1)}$ ; ceci prouve que $\phi_{(k,\beta)}^{(\theta_0,\theta_1)}$ appartient \`a $G'$ et contredit l'hypoth\`ese :
$\phi_{(k,\beta)}^{(\theta_0,\theta_1)}<g_i=inf G'$.
\dli ii) Il existe $\theta\in\Theta_1$ tel que $P_{\theta}(B)>0$. 
L'\'ev\'enement $C=B\cap\{p_{\theta}>0\}$ est non vide et $P_{\theta}(C)>0$. 
Soit $\omega\in C$, comme $p_{\theta_0}(\omega)=0$ on a  
$h_{(\theta_0,\theta)}(T(\omega))=0$ ; $h_{(\theta_0,\theta)}$ \'etant croissante on a, pour $t=T(\omega)$, $E_{\theta_0}(f_{(t,1)})=0$ ; de plus
$\phi_{(k,\beta)}^{(\theta_0,\theta_1)}<f_{(t,1)}\leq g_i$ car $\omega\in A$ ;
$\phi_{(0,1)}^{(\theta_0,\theta_1)}$ et $\phi_{(k,\beta)}^{(\theta_0,\theta_1)}$ sont donc de moyenne nulle sous 
$P_{\theta_0}$ ce qui, par d\'efinition de $\phi_{(0,1)}^{(\theta_0,\theta_1)}$, implique l'\'egalit\'e $P_{\theta_0}$ et
$P_{\theta_1}$ presque s\^urement de ces deux \'el\'ements de
$\Phi_s^{(\theta_0,\theta_1)}$. $\phi$ v\'erifie alors :
$\phi\relmont{<}{\{\theta_0,\theta_1\} p.s.}\phi_{(0,1)}^{(\theta_0,\theta_1)}$ et
$\phi\relmont{=}{\theta_0 p.s.}\phi_{(0,1)}^{(\theta_0,\theta_1)}$ ;
on en d\'eduit $\phi\relmont{<}{\theta_1 p.s.}\phi_{(0,1)}^{(\theta_0,\theta_1)}$, ce qui contredit le fait que $\phi$ est un expert du choix entre $P_{\theta_0}$ et
$P_{\theta_1}$.

\dli 2\up{\`eme} possibilit\'e : $\phi\relmont{>}{\{\theta_0,\theta_1\} p.s.}\phi_{(k,\beta)}^{(\theta_0,\theta_1)}>g_s$.
\dli La d\'emonstration est semblable \`a celle du cas pr\'ec\'edent.
\dli On pose $A=\{\phi_{(k,\beta)}^{(\theta_0,\theta_1)}-g_s=1\}\cap\{\phi=0\}=\{\phi_{(k,\beta)}^{(\theta_0,\theta_1)}\not= g_s\}\cap\{\phi=0\}$.
Comme $g_s\relmont{\leq}{p.s.}\phi$, on a
$\phi_{(k,\beta)}^{(\theta_0,\theta_1)}.\II_{A^c}\relmont{\leq}{p.s.}\phi.\II_{A^c}$. De plus,
$\phi_{(k,\beta)}^{(\theta_0,\theta_1)}.\II_A=\II_A$ entra\^{\i}ne 
$\phi.\II_A\relmont{=}{\{\theta_0,\theta_1\} p.s.}\II_A$, ce qui implique
$P_{\theta_0}(A)=P_{\theta_1}(A)=0$ puisque $A\subseteq\{\phi=0\}$ ;
les \'ev\'enements $A\cap\{p_{\theta_0}>0\}$ et $A\cap\{p_{\theta_1}>0\}$
sont alors $\mu$ n\'egligeables et donc $P_{\theta}$ n\'egligeables pour tout
$\theta$ ; la seule partie de $A$ qui peut ne pas \^etre $P_{\theta}$ n\'egligeable, pour tout $\theta$, est $B=A\cap\{p_{\theta_0}=0\}\cap\{p_{\theta_1}=0\}$ ; on a donc d\'ej\`a :
$\phi_{(k,\beta)}^{(\theta_0,\theta_1)}.\II_{B^c}\relmont{\leq}{p.s.}\phi.\II_{B^c}$. On distingue maintenant deux cas, qui conduiront \`a des contradictions diff\'erentes.
\dli i) Pour tout $\theta\in\Theta_0$ : $P_{\theta}(B)=0$. La partie ii) de la d\'efinition 4.1.1 implique : $\phi.\II_B\relmont{=}{\Theta_1 p.s.}\II_B$.
Comme $\phi$ est \'egal \`a $0$ sur $B$, cet \'ev\'enement est $P_{\theta}$ n\'egligeable pour tout $\theta$ et on a  alors 
$\phi_{(k,\beta)}^{(\theta_0,\theta_1)}\relmont{\leq}{p.s.}\phi$ ; ceci prouve que $\phi_{(k,\beta)}^{(\theta_0,\theta_1)}$ appartient \`a $G$ et contredit l'hypoth\`ese :
$\phi_{(k,\beta)}^{(\theta_0,\theta_1)}>g_s=supG$.
\dli ii) Il existe $\theta\in\Theta_0$ tel que $P_{\theta}(B)>0$. 
L'\'ev\'enement $C=B\cap\{p_{\theta}>0\}$ est non vide et $P_{\theta}(C)>0$. 
Soit $\omega\in C$, comme $p_{\theta_1}(\omega)=0$ on a  
$h_{(\theta,\theta_1)}(T(\omega))=+\infty$ ; $h_{(\theta,\theta_1)}$ \'etant croissante on a, pour $t=T(\omega)$, $E_{\theta_1}(f_{(t,0)})=1$ ; de plus
$g_s\leq f_{(t,0)}<\phi_{(k,\beta)}^{(\theta_0,\theta_1)}$ car $\omega\in A$ ; on a donc
$E_{\theta_1}(\phi_{(k,\beta)}^{(\theta_0,\theta_1)})=1=E_{\theta_1}(\phi_{(\infty,0)}^{(\theta_0,\theta_1)})$,  
 ce qui implique, par d\'efinition de $\phi_{(\infty,0)}^{(\theta_0,\theta_1)}$ :
$\phi_{(k,\beta)}^{(\theta_0,\theta_1)}\relmont{=}{\{\theta_0,\theta_1\} p.s.}\phi_{(\infty,0)}^{(\theta_0,\theta_1)}$.  
$\phi$ v\'erifie alors :
$\phi\relmont{>}{\{\theta_0,\theta_1\} p.s.}\phi_{(\infty,0)}^{(\theta_0,\theta_1)}$ et
$\phi\relmont{=}{\theta_1 p.s.}\phi_{(\infty,0)}^{(\theta_0,\theta_1)}$ ;
on en d\'eduit $\phi\relmont{>}{\theta_0 p.s.}\phi_{(\infty,0)}^{(\theta_0,\theta_1)}$, ce qui contredit le fait que $\phi$ est un expert du choix entre $P_{\theta_0}$ et
$P_{\theta_1}$.

\dli 3\up{\`eme} possibilit\'e : $\phi\relmont{=}{\{\theta_0,\theta_1\} p.s.}\phi_{(k,\beta)}^{(\theta_0,\theta_1)}$.
\dli On a $\phi\relmont{=}{\{\theta_0,\theta_1\} p.s.}\phi_{(k,\beta)}^{(\theta_0,\theta_1)}<g_i$ et 
$g_s<\phi_{(k,\beta)}^{(\theta_0,\theta_1)}\relmont{=}{\{\theta_0,\theta_1\} p.s.}\phi$.
On peut refaire les d\'emonstrations des deux possibilit\'es pr\'ec\'edentes jusqu'au cas i). C'est dans ii) que l'on se sert de l'in\'egalit\'e stricte, 
$P_{\theta_0}$ et $P_{\theta_1}$ presque s\^urement,
entre $\phi$ et $\phi_{(k,\beta)}^{(\theta_0,\theta_1)}$. 
\dli Il nous reste \`a trouver une contradiction lorsque les deux cas ii) sont r\'ealis\'es. Les constructions faites au d\'ebut de ces deux cas sont encore valables. Elles nous conduisent \`a l'existence de $\theta\in\Theta_0$, de
$\theta'\in\Theta_1$ et de deux r\'eels, $t$ et $t'$,  tels que :
\dli $g_s\leq f_{(t,0)}<\phi_{(k,\beta)}^{(\theta_0,\theta_1)}$ et
$h_{(\theta,\theta_1)}(t)=+\infty$ ;
\dli $\phi_{(k,\beta)}^{(\theta_0,\theta_1)}<f_{(t',1)}\leq g_i$ et
$h_{(\theta_0,\theta')}(t')=0$.
\dli Consid\'erons l'intervalle $I=[t,t']$ ; la fonction $h_{(\theta,\theta_1)}$ (resp. $h_{(\theta_0,\theta')}$) est \'egale \`a $+\infty$ (resp. $0$) sur $I$, donc pour tout $\omega$ de $\{T\in I\}$ on a
$p_{\theta_1}(\omega)=0$ et $p_{\theta_0}(\omega)=0$. La famille
$\{h_{(\theta_0,\theta_1)}\}$ \'etant normalis\'ee,
d'apr\`es la d\'efinition 4.2.1, $h_{(\theta_0,\theta_1)}$ est constante sur $I$, puisque $I\subseteq\IR\!-(D_i\cup D_s)$ (voir la proposition 4.2.1) ; ce qui contredit l'existence d'un \'el\'ement
$\phi_{(k,\beta)}^{(\theta_0,\theta_1)}$ strictement entre
$f_{(t,0)}$ et $f_{(t',1)}$ .

2\up{\`eme} cas : $E_{\theta_0}(\phi)=0$.
\dli La condition $E_{\theta_0}(\phi)=0$ est \'equivalente \`a 
$\phi\relmont{=}{\theta_0 p.s.}\phi_{(0,1)}^{(\theta_0,\theta_1)}$.
\dli Montrons que l'on a aussi 
$\phi\relmont{\leq}{\theta_1 p.s.}\phi_{(0,1)}^{(\theta_0,\theta_1)}$.
\dli Par d\'efinition, $\phi_{(0,1)}^{(\theta_0,\theta_1)}$ vaut $1$ sur $\{p_{\theta_0}=0\}\cap\{p_{\theta_1}>0\}$. L'\'ev\'enement $A=\{\phi>\phi_{(0,1)}^{(\theta_0,\theta_1)}\}=\{\phi_{(0,1)}^{(\theta_0,\theta_1)}=0\}\cap\{\phi=1\}$ est bien $\theta_1$ n\'egligeable puisqu'il se d\'ecompose en $A'=A\cap\{p_{\theta_0}>0\}$ qui est $\mu$ n\'egligeable ($E_{\theta_0}(\phi)=0$) et $B=A\cap\{p_{\theta_0}=0\}\cap\{p_{\theta_1}=0\}$.

\dli Le cas $\phi\relmont{=}{\{\theta_0,\theta_1\} p.s.}\phi_{(0,1)}^{(\theta_0,\theta_1)}$ a \'et\'e \'etudi\'e pr\'ec\'edemment puisque $\phi$ est alors un expert du choix entre $P_{\theta_0}$ et $P_{\theta_1}$.

Il reste \`a consid\'erer le cas :
$\phi\relmont{=}{\theta_0 p.s.}\phi_{(0,1)}^{(\theta_0,\theta_1)}$ et
$\phi\relmont{<}{\theta_1 p.s.}\phi_{(0,1)}^{(\theta_0,\theta_1)}$.
Par construction de $B$ on a 
$\phi.\II_{B^c}\relmont{\leq}{p.s.}\phi_{(0,1)}^{(\theta_0,\theta_1)}.\II_{B^c}$.
\dli Nous allons encore raisonner par l'absurde en supposant qu'il existe
$\phi_{(k,\beta)}^{(\theta_0,\theta_1)}$ dans $]g_s,g_i[$ ; on a alors
$\phi_{(0,1)}^{(\theta_0,\theta_1)}\leq\phi_{(k,\beta)}^{(\theta_0,\theta_1)}<g_i$. Si on montre que $B$ est $P_\theta$ n\'egligeable pour tout $\theta$ on aura $\phi\relmont{\leq}{p.s.}\phi_{(0,1)}^{(\theta_0,\theta_1)}$, ce qui se traduit par l'appartenance de $\phi_{(0,1)}^{(\theta_0,\theta_1)}$ \`a $G'$ et est incompatible avec $\phi_{(0,1)}^{(\theta_0,\theta_1)}<g_i=infG'$.

\dli Montrons d'abord que $B$ est $P_\theta$ n\'egligeable pour tout $\theta$ de $\Theta_0$.
S'il existait $\theta\in\Theta_0$ tel que $P_\theta(B)>0$, $\phi$ serait un expert du choix entre $P_\theta$ et $P_{\theta_1}$, d'apr\`es la partie i) de la d\'efinition 4.1.1, puisqu'on aurait $E_\theta(\phi)>0$ et par hypoth\`ese
$\phi\relmont{<}{\theta_1 p.s.}\phi_{(0,1)}^{(\theta_0,\theta_1)}\leq\phi_{(\infty,0)}^{(\theta_0,\theta_1)}$ donc $E_{\theta_1}(\phi)<1$ ; de plus, l'\'ev\'enement $C=B-\{g_i=0\}$ serait non n\'egligeable pour $P_\theta$ puisque $\phi\relmont{\leq}{p.s.}g_i$, il existerait alors $\omega\in C$ tel que
$p_\theta(\omega)>0$ donc $h_{(\theta,\theta_1)}(T(\omega))=+\infty$ et
$\phi_{(\infty,0)}^{(\theta,\theta_1)}\leq f_{(T(\omega),0)}<g_i$ ; on aurait m\^eme $\phi_{(\infty,0)}^{(\theta,\theta_1)}\leq g_s$ puisque le premier cas appliqu\'e \`a $\phi$ expert du choix entre $P_\theta$ et $P_{\theta_1}$ entra\^{\i}ne 
$]g_s,g_i[\cap\Phi_s^{(\theta,\theta_1)}=\emptyset$ ; ce qui impliquerait
$E_{\theta_1}(\phi)=1$ qui est incompatible avec la condition
$\phi\relmont{<}{\theta_1 p.s.}\phi_{(0,1)}^{(\theta_0,\theta_1)}$.
On a donc bien $P_\theta(B)=0$ pour $\theta\in\Theta_0$.
\dli Il reste \`a montrer que l'on a aussi $P_\theta(B)=0$ pour tout $\theta$ dans $\Theta_1$.
Soit $\theta\in\Theta_1$, supposons $P_\theta(B)>0$ ; l'\'ev\'enement $C=B\cap\{p_\theta(\omega)>0\}$ ne serait pas $\mu$ n\'egligeable bien que $P_{\theta'}$ n\'egligeable pour tout $\theta'$ de $\Theta_0$ ; il existerait alors $\omega_{\theta'}\in C$ tel que $p_{\theta'}(\omega_{\theta'})=0$, donc $h_{(\theta',\theta)}(T(\omega_{\theta'}))=0$ et par d\'efinition de
$B\subseteq A$ : $\phi_{(0,1)}^{(\theta_0,\theta_1)}\leq f_{(T(\omega_{\theta'}),1)}\leq \phi_{(0,1)}^{(\theta',\theta)}$ ;
le lemme 2 de l'annexe II impliquerait $f_{(t_0,u_0)}=\phi_{(0,1)}^{(\theta_0,\theta_1)}$, ce qui contredit la condition $\phi\relmont{<}{\theta_1 p.s.}\phi_{(0,1)}^{(\theta_0,\theta_1)}$ car $\phi$ est un expert donc
$\phi\relmont{\geq}{p.s.}f_{(t_0,u_0)}$ (voir II-a)).

3\up{\`eme} cas : $E_{\theta_1}(\phi)=1$.
\dli La d\'emonstration est semblable \`a la pr\'ec\'edente.
\dli On pose $A=\{\phi<\phi_{(\infty,0)}^{(\theta_0,\theta_1)}\}=\{\phi_{(\infty,0)}^{(\theta_0,\theta_1)}=1\}\cap\{\phi=0\}$ qui se d\'ecompose en 
$A'=A\cap\{p_{\theta_1}>0\}$, $A''=A\cap\{p_{\theta_1}=0\}\cap\{p_{\theta_0}>0\}$ et
\dli $B=A\cap\{p_{\theta_1}=0\}\cap\{p_{\theta_0}=0\}$.
\dli $A''$ est vide puisque $\{p_{\theta_1}=0\}\cap\{p_{\theta_0}>0\}$ est inclus dans $\{\phi_{(\infty,0)}^{(\theta_0,\theta_1)}=0\}$ par d\'efinition de 
$\phi_{(\infty,0)}^{(\theta_0,\theta_1)}$ ; $A'$  est $\mu$ n\'egligeable puisque la condition $E_{\theta_1}(\phi)=1$ est \'equivalente \`a 
$\phi\relmont{=}{\theta_1 p.s.}\phi_{(\infty,0)}^{(\theta_0,\theta_1)}$ ; 
on a \'evidemment $P_{\theta_0}(B)=0$ et donc aussi : 
$\phi\relmont{\geq}{\theta_0 p.s.}\phi_{(\infty,0}^{(\theta_0,\theta_1)}$.

\dli Le cas $\phi\relmont{=}{\{\theta_0,\theta_1\} p.s.}\phi_{(\infty,0)}^{(\theta_0,\theta_1)}$ a \'et\'e \'etudi\'e pr\'ec\'edemment puisque $\phi$ est alors un expert du choix entre $P_{\theta_0}$ et $P_{\theta_1}$.

Il reste \`a consid\'erer le cas :
$\phi\relmont{=}{\theta_1 p.s.}\phi_{(\infty,0)}^{(\theta_0,\theta_1)}$ et
$\phi\relmont{>}{\theta_0 p.s.}\phi_{(\infty,0)}^{(\theta_0,\theta_1)}$.
Par construction de $B$ on a 
$\phi.\II_{B^c}\relmont{\geq}{p.s.}\phi_{(\infty,0)}^{(\theta_0,\theta_1)}.\II_{B^c}$.
\dli Nous allons encore raisonner par l'absurde en supposant qu'il existe
$\phi_{(k,\beta)}^{(\theta_0,\theta_1)}$ dans $]g_s,g_i[$ ; on a alors
$g_s<\phi_{(k,\beta)}^{(\theta_0,\theta_1)}\leq\phi_{(\infty,0)}^{(\theta_0,\theta_1)}$. Si on montre que $B$ est $P_\theta$ n\'egligeable pour tout $\theta$ on aura $\phi\relmont{\geq}{p.s.}\phi_{(\infty,0)}^{(\theta_0,\theta_1)}$, ce qui se traduit par l'appartenance de $\phi_{(\infty,0)}^{(\theta_0,\theta_1)}$ \`a $G$ et est incompatible avec $\phi_{(\infty,0)}^{(\theta_0,\theta_1)}>g_s=supG$.

\dli Montrons d'abord que $B$ est $P_\theta$ n\'egligeable pour tout $\theta$ de $\Theta_1$.
S'il existait $\theta\in\Theta_1$ tel que $P_\theta(B)>0$, $\phi$ serait un expert du choix entre $P_{\theta_0}$ et $P_{\theta}$, d'apr\`es la partie i) de la d\'efinition 4.1.1, puisqu'on aurait $E_\theta(\phi)<1$ et par hypoth\`ese
$\phi\relmont{>}{\theta_0 p.s.}\phi_{(\infty,0)}^{(\theta_0,\theta_1)}$ donc $E_{\theta_0}(\phi)>0$ ; de plus, l'\'ev\'enement $C=B-\{g_s=1\}$ serait non n\'egligeable pour $P_\theta$ puisque $\phi\relmont{\geq}{p.s.}g_s$ et $B\subseteq\{\phi=0\}$ ; il existerait alors $\omega\in C$ tel que
$p_\theta(\omega)>0$ donc $h_{(\theta_0,\theta)}(T(\omega))=0$ et
$g_s<f_{(T(\omega),1)}\leq\phi_{(0,1)}^{(\theta_0,\theta)}$ ; on aurait m\^eme $\phi_{(0,1)}^{(\theta_0,\theta)}\geq g_i$ puisque le premier cas appliqu\'e \`a $\phi$ expert du choix entre $P_{\theta_0}$ et $P_{\theta}$ entra\^{\i}ne 
$]g_s,g_i[\cap\Phi_s^{(\theta_0,\theta)}=\emptyset$ ; ce qui impliquerait
$E_{\theta_0}(\phi)=0$ qui est incompatible avec la condition
$\phi\relmont{>}{\theta_0 p.s.}\phi_{(\infty,0)}^{(\theta_0,\theta_1)}$.
On a donc bien $P_\theta(B)=0$ pour $\theta\in\Theta_1$.
\dli Il reste \`a montrer que l'on a aussi $P_\theta(B)=0$ pour tout $\theta$ dans $\Theta_0$.
Soit $\theta\in\Theta_0$, supposons $P_\theta(B)>0$ ; l'\'ev\'enement $C=B\cap\{p_\theta(\omega)>0\}$ ne serait pas $\mu$ n\'egligeable bien que $P_{\theta'}$ n\'egligeable pour tout $\theta'$ de $\Theta_1$ ; il existerait alors $\omega_{\theta'}\in C$ tel que $p_{\theta'}(\omega_{\theta'})=0$, donc $h_{(\theta,\theta')}(T(\omega_{\theta'}))=+\infty$ et par d\'efinition de
$B\subseteq A$ : $\phi_{(\infty,0)}^{(\theta_0,\theta_1)}\geq f_{(T(\omega_{\theta'}),1)}\geq \phi_{(\infty,0)}^{(\theta,\theta')}$ ;
le lemme 1 de l'annexe II impliquerait $f_{(t_1,u_1)}=\phi_{(\infty,0)}^{(\theta_0,\theta_1)}$, ce qui contredit la condition $\phi\relmont{>}{\theta_0 p.s.}\phi_{(\infty,0)}^{(\theta_0,\theta_1)}$ puisque tout expert v\'erifie :
$\phi\relmont{\leq}{p.s.}f_{(t_1,u_1)}$ (voir II-b)).

\medskip\centerline{\hbox to 3cm{\bf \hrulefill}}\par}

\vfill\eject

\bigskip
\goodbreak
{\parindent=-5mm 4.3 VOTES DES EXPERTS.}
\nobreak
\medskip

Consid\'erons un probl\`eme de d\'ecision \`a hypoth\`eses stables :
\dli $(\Omega,{\cal A},(P_\theta=p_\theta.\mu)_{\theta\in\Theta_0\cup \Theta_1})$ (voir la d\'efinition 4.1.2). D'apr\`es l'annexe II on peut lui associer une statistique r\'eelle $T$ et une famille normalis\'ee :
\dli $\{h_{(\theta_0,\theta_1)}\}_{(\theta_0,\theta_1)\in\Theta_0\times\Theta_1}$
(voir la d\'efinition 4.2.1).
$K_{(\theta_0,\theta_1)}=h_{(\theta_0,\theta_1)}(T)$ nous permet de d\'efinir l'ensemble $\Phi_s^{(\theta_0,\theta_1)}=[\phi_{(0,1)}^{(\theta_0,\theta_1)}, \phi_{(\infty,0)}^{(\theta_0,\theta_1)}]$ des fonctions de tests simples bas\'ees sur $K_{(\theta_0,\theta_1)}$ (voir la d\'efinition 2.2.1).
D'apr\`es la proposition 4.2.2, c'est l'adh\'erence $\overline{\Delta_s}$ de $\Delta_s=\cup_{(\theta_0,\theta_1)\in\Theta_0\times\Theta_1} \Phi_s^{(\theta_0,\theta_1)}$ dans $F$ qui contient les experts fondamentaux. Comme nous l'avons fait pour le choix entre deux probabilit\'es (paragraphe 2.4), nous allons commencer par probabiliser cet ensemble d'experts.

$\overline{\Delta_s}\subseteq F$ est muni de la $\sigma$-alg\`ebre trace de la tribu bor\'elienne correspondant \`a la topologie de l'ordre d\'efini sur $F$ (voir le paragraphe 4.2). Pour d\'efinir une probabilit\'e sur $\overline{\Delta_s}$ on peut utiliser les constructions classiques \`a partir d'une semi-alg\`ebre (cf. [Nev.] p. 25). Par exemple, la semi-alg\`ebre engendr\'ee par les intervalles de la forme $]\leftarrow ,f[=\{g\in\overline{\Delta_s}\,;\,g<f\}$ avec $f\in\overline{\Delta_s}$ ou celle engendr\'ee par ceux de la forme $]\leftarrow ,f]=\{g\in\overline{\Delta_s}\,;\,g\leq f\}$. 
 On montre facilement, qu'il suffit alors de d\'efinir sur les intervalles de la forme $]\leftarrow ,f[$ (resp. $]\leftarrow ,f]$) une fonction \`a valeurs dans $[0,1]$, non d\'ecroissante, continue \`a gauche (resp. droite) et valant $0$ (resp. $1$) en $f=inf\overline{\Delta_s}$ (resp.$f=sup\overline{\Delta_s}$).

Comme nous l'avons fait en 2.4, nous choisissons de probabiliser $\overline{\Delta_s}$ en utilisant l'op\'erateur $E_\theta$.
Pour chaque $\theta\in\Theta$, il nous permet de d\'efinir deux probabilit\'es $m_\theta$ et $m'_\theta$ en posant $m_\theta(]\leftarrow ,f[)=E_\theta(f)$ pour $f\not= inf\overline{\Delta_s}$ et $m'_\theta(]\leftarrow ,f])=E_\theta(f)$ pour $f\not= sup\overline{\Delta_s}$.
Ces deux probabilit\'es sont identiques si pour tout couples $(f,g)$  d'\'el\'ements successifs de $\overline{\Delta_s}$, c'est-\`a-dire $f<g$ et $]f,g[\cap\overline{\Delta_s}=\emptyset$, on a : $E_\theta(f)=E_\theta(g)$.
Dans le cas contraire, la masse $E_\theta(g-f)$ est attribu\'ee \`a $f$ par $m_\theta$ et \`a $g$ par $m'_\theta$. Comme nous l'avons fait au paragraphe 2.4, il semble naturel de partager \'equitablement cette masse entre $f$ et $g$, ce qui revient \`a probabiliser $\overline{\Delta_s}$ par $(m_\theta + m'_\theta)/2$. Le r\'esultat du vote des experts \`a partir de cette probabilit\'e d\'efinit pour chaque r\'ealisation $\omega\in\Omega$ une probabilit\'e 	
$Q^\omega_\theta$ sur $D=\{0,1\}$. On a 
$Q^\omega_\theta(\{1\})=\int_{\overline{\Delta_s}}f(\omega)
\,d(m_\theta+m'_\theta)/2$. Afin de pouvoir exprimer simplement ce r\'esultat, on va d\'efinir une statistique qui jouera un r\^ole semblable \`a celui du rapport des densit\'es $K$ dans le choix entre deux probabilit\'es.

\medskip
{\bf D\'efinition 4.3.1}
\medskip
\medskip
\moveleft 10.4pt\hbox{\vrule\kern 10pt\vbox{\defpro

Soit $(\Omega,{\cal A},(P_\theta=p_\theta.\mu)_{\theta\in\Theta_0\cup \Theta_1})$ un probl\`eme de d\'ecision \`a hypoth\`eses stables \`a partir de la 
statistique r\'eelle $T$. Associons \`a $t\in\IR$, deux \'el\'ements de $[f_{(-\infty,0)},f_{(+\infty,1)}]$ d\'efinis par :
\dli $f_{(a_t,u_t)}=sup\{f\in\overline{\Delta_s}\, ;\, f\leq f_{(t,0)}\}\cup\{f_{(-\infty,0)}\}$ 
\dli $f_{(b_t,v_t)}=inf\{f\in\overline{\Delta_s}\, ;\, f\geq f_{(t,1)}\}\cup\{f_{(+\infty,1)}\}$
\dli On appelle statistique essentielle, la statistique $K(T)$ obtenue \`a partir de la fonction croissante $K\, :\, \IR \rightarrow\overline{\IR}$ d\'efinie par :
$$K(t)=\left\{\matrix{
\hfill -\infty\hfill & si & (a_t,u_t)=(-\infty,0)\hfill\cr
\hfill b_t-1\hfill & si & (a_t,u_t)=(-\infty,1)\hfill\cr
\hfill [a_t+b_t]/2\hfill & si & -\infty<a_t\leq b_t<+\infty\hfill \cr
\hfill a_t+1\hfill & si & (b_t,v_t)=(+\infty,0)\ et\ a_t>-\infty\hfill\cr
\hfill +\infty\hfill & si & (b_t,v_t)=(+\infty,1)\hfill \cr
}\right. $$
Sous $P_\theta$, la fonction de r\'epartition moyenne $G_\theta(k)$ de cette statistique est \'egale \`a
$[P_\theta(\{K(T)<k\})+P_\theta(\{K(T)\leq k\})]/2$.
}}\medskip

Cette d\'efinition est coh\'erente car on ne peut pas avoir en m\^eme temps 
$(a_t,u_t)=(-\infty,0)$ et $(b_t,v_t)=(+\infty,1)$. En effet, on a
$\overline{\Delta_s}\subseteq[inf\Delta_s,sup\Delta_s]=[\II_{D_i-D_s}\!(T),1-\II_{D_s}\!(T)]\subseteq[f_{(-\infty,1)},f_{(+\infty,0)}]$ (voir la proposition 4.2.1), le cas $(a_t,u_t)=(-\infty,0)$ (resp. $(b_t,v_t)=(+\infty,1)$) se produit donc uniquement lorsque $t\in (D_i-D_s)$ (resp. $t\in D_s$).
\dli La quatri\`eme possibilit\'e de la d\'efinition de $K$ \'elimine le cas :  $(b_t,v_t)=(+\infty,0)$ et $a_t=-\infty$, il correspond \`a $\Delta_s=\{f_{(-\infty,1)},f_{(+\infty,0)}\}$, donc \`a un choix entre deux 
hypoth\`eses d\'efinies par des densit\'es $p_\theta$ identiques et telles que $D_i=D_s=\emptyset$. Nous aurions pu enlever ce cas sans int\'er\^et, mais comme il ne pose pas de probl\`eme lorsque $D_i\not=\emptyset$ ou $D_s\not=\emptyset$ nous avons pr\'ef\'er\'e choisir entre $K(t)=-\infty$ et $K(t)=+\infty$.
\dli La croissance de $K$ provient du fait que pour $t<t'$ on a soit 
$f_{(a_t,u_t)}=f_{(a'_t,u'_t)}$ et $f_{(b_t,v_t)}=f_{(b'_t,v'_t)}$, 
soit $f_{(b_t,v_t)}\leq f_{(a'_t,u'_t)}$.

Pour $t$ n'appartenant pas  \`a $D_i\cup D_s$, on a $inf\overline{\Delta_s}\leq f_{(t,0)}<f_{(t,1)}\leq sup\overline{\Delta_s}$ ; les deux fonctions de test $f_{(a_t,u_t)}$ et $f_{(b_t,v_t)}$ d\'efinissent alors les deux experts successifs de $\overline{\Delta_s}$ qui ne prennent pas la m\^eme d\'ecision lorsque $T(\omega)=t$.
La statistique essentielle $K(T)$ nous permet de traiter de fa\c con homog\`ene
les \'el\'ements successifs de $\overline{\Delta_s}$ charg\'es par la probabilit\'e
$(m_\theta + m'_\theta)/2$ et ceux qui ne le sont pas.

\medskip
{\bf Proposition 4.3.1}
\medskip
\medskip
\moveleft 10.4pt\hbox{\vrule\kern 10pt\vbox{\defpro

Lorsqu'on r\'ealise $\omega$, le r\'esultat du vote des experts $\overline{\Delta_s}$ sous $P_\theta$ est une probabilit\'e $Q^\omega_\theta$ d\'efinie sur l'espace des d\'ecisions $D=\{0,1\}$ par :
$$Q^\omega_\theta(\{1\})=\left\{\matrix{
1\hfill & si & T(\omega)\in D_i-D_s\hfill\cr
1-G_\theta(K(T(\omega)))\hfill & si & T(\omega)\in\IR\!-(D_i\cup D_s)\hfill \cr
0\hfill & si & T(\omega)\in D_s\hfill \cr
}\right. $$
\dli ($G_\theta$ est la fonction de r\'epartition moyenne de la statistique essentielle $K(T)$, les demi-droites $D_i$ et $D_s$ sont d\'efinies en 4.2.1).

}}\medskip

{\leftskip=15mm \dli {\bf D\'emonstration}
\medskip

Soient $\theta\in\Theta$ et $\omega\in\Omega$.
\dli $Q^\omega_\theta(\{1\})=\int_{\overline{\Delta_s}}f(\omega)
\,d(m_\theta+m'_\theta)/2$ est la fr\'equence des experts qui d\'ecident $d=1$
face \`a la r\'ealisation $\omega$.
\dli Posons $t=T(\omega)$, d'apr\`es la proposition 4.2.1 on a : 
$inf\overline{\Delta_s}=\II_{D_i-D_s}\!(T)$ et $sup\overline{\Delta_s}=1-\II_{D_s}\!(T)$ ; les diff\'erentes conditions de l'expression de $Q^\omega_\theta(\{1\})$ peuvent s'\'ecrire respectivement :
$f_{(t,0)}<inf\overline{\Delta_s}$, $inf\overline{\Delta_s}\leq f_{(t,0)}<f_{(t,1)}\leq sup\overline{\Delta_s}$ et 
$f_{(t,1)}>sup\overline{\Delta_s}$.

1\up{er} cas : $f_{(t,0)}<inf\overline{\Delta_s}$.
\dli Tous les experts de $\overline{\Delta_s}$ d\'ecident $d=1$ en $\omega$ puisque $T(\omega)=t$ ; on a bien $Q^\omega_\theta(\{1\})=1$.

2\up{\`eme} cas : $inf\overline{\Delta_s}\leq f_{(t,0)}<f_{(t,1)}\leq sup\overline{\Delta_s}$.
\dli Consid\'erons les fonctions de test $\delta_t=f_{(a_t,u_t)}$ et
$\delta'_t=f_{(b_t,v_t)}$ de la d\'efinition pr\'ec\'edente. Elles v\'erifient :
$inf\overline{\Delta_s}\leq\delta_t\leq f_{(t,0)}<f_{(t,1)}\leq\delta'_t\leq sup\overline{\Delta_s}$.
\dli $\delta_t$ (resp. $\delta'_t$) est le plus grand (resp. petit) expert de $\overline{\Delta_s}$ d\'ecidant $d=0$ (resp. $d=1$) pour la r\'ealisation $\omega$. On a donc :
\+\kern 15mm$Q^\omega_\theta(\{1\})$&=&
$1-[(m_\theta+m'_\theta)/2](]\leftarrow,\delta_t])$  \cr
\+&=&
$1-(1/2)[m_\theta(]\leftarrow,\delta'_t[)+m'_\theta(]\leftarrow,\delta_t])]$\cr
\+&=&
$1-(1/2)[E_\theta(\delta'_t)+E_\theta(\delta_t)]$\cr
\+&=&
$1-(1/2)[P_\theta(\{\delta'_t=1\})+P_\theta(\{\delta_t=1\})]$.\cr
\dli Ceci d\'emontre le r\'esultat recherch\'e, $Q^\omega_\theta(\{1\})=1-G_\theta(K(t))$, si les deux \'ev\'enements 
$\{\delta_t=1\}$ et $\{\delta'_t=1\}$ s'\'ecrivent respectivement :
$\{K(T)<K(t)\}$ et $\{K(T)\leq K(t)\}$ ; c'est-\`a-dire si 
$\{\delta'_t-\delta_t=1\}=\{K(T)=K(t)\}$.
\dli Par construction de $K$ on a $\{\delta'_t-\delta_t=1\}\subseteq\{K(T)=K(t)\}$. Consid\'erons $\omega_0\notin\{\delta'_t-\delta_t=1\}$, il nous reste \`a montrer que 
$K(x)\not=K(t)$ pour $x=T(\omega_0)$. Nous allons analyser successivement les diff\'erents cas de la d\'efinition de $K(t)$ ; dans notre situation il n'y a que trois cas possibles puisque $\delta_t=f_{(a_t,u_t)}$ et
$\delta'_t=f_{(b_t,v_t)}$ appartiennent \`a $[inf\Delta_s,sup\Delta_s]\subseteq[f_{(-\infty,1)},f_{(+\infty,0)}]$.

i) $K(t)=b_t-1$ pour $(a_t,u_t)=(-\infty,1)$.
\dli Dans ce cas $\{\delta'_t-\delta_t=1\}=\{f_{(b_t,v_t)}=1\}$ et $\omega_0\in\{f_{(b_t,v_t)}=0\}$ ; pour que $\omega_0$ puisse exister il faut $b_t<+\infty$ puisque $x=T(\omega_0)\in\IR$ ; on a alors $f_{(b_t,v_t)}\leq f_{(x,0)}$ donc $f_{(b_t,v_t)}\leq f_{(a_x,u_x)}$ ; comme $f_{(a_t,u_t)}<f_{(b_t,v_t)}$, on a $-\infty<b_t\leq a_x$ ; par d\'efinition de $K(x)$ ceci implique $K(x)\geq a_x\geq b_t>b_t-1=K(t)$ (on a $b_t>b_t-1$ car $b_t$ n'est pas infini).

ii) $K(t)=[a_t+b_t]/2$ pour $-\infty<a_t\leq b_t<+\infty$.
\dli $\omega_0$ appartient \`a $\{f_{(a_t,u_t)}=1\}$ ou $\{f_{(b_t,v_t)}=0\}$, 
dans le premier cas on a $f_{(x,1)}\leq f_{(b_x,v_x)}\leq f_{(a_t,u_t)}$ et
dans le second $f_{(b_t,v_t)}\leq f_{(a_x,u_x)}\leq f_{(x,0)}$, ce qui implique respectivement $b_x\leq a_t$ et $b_t\leq a_x$.
\dli Lorsque $b_x\leq a_t$ la d\'efinition de $K(x)$ entra\^{\i}ne $K(x)<K(t)$ ; c'est \'evident si $a_t<b_t$ ou $a_x=-\infty$ ; dans le cas contraire cela reste vrai car $a_t=b_t$ est \'equivalent \`a  $(a_t,u_t)=(t,0)$ et $(b_t,v_t)=(t,1)$, d'autre part on a  $-\infty<a_x<t$ puisque $x<t$.
\dli De m\^eme, lorsque $b_t\leq a_x$ on a $K(t)<K(x)$ ; c'est \'evident si $a_t<b_t$ ou $b_x=+\infty$ ; dans le cas contraire cela reste vrai car on a $(a_t,u_t)=(t,0)$, $(b_t,v_t)=(t,1)$et  $t<x\leq b_x<+\infty$.

iii) $K(t)=a_t+1$ pour $(b_t,v_t)=(+\infty,0)$ et $a_t>-\infty$.
\dli $\omega_0$ appartient \`a $\{f_{(a_t,u_t)}=1\}$ et comme pr\'ec\'edemment on a $f_{(x,1)}\leq f_{(b_x,v_x)}\leq f_{(a_t,u_t)}$, $b_x\leq a_t$, donc $K(x)<K(t)$. 

3\up{\`eme} cas : $f_{(t,1)}>sup\overline{\Delta_s}$.
\dli Tous les experts de $\overline{\Delta_s}$ d\'ecident $d=0$ en $\omega$ puisque $T(\omega)=t$ ; on a donc bien : $Q^\omega_\theta(\{1\})=0$.

\medskip\centerline{\hbox to 3cm{\bf \hrulefill}}\par}

\bigskip
Les votes $Q_\theta$ que l'on vient de d\'efinir ont \'et\'e construits \`a partir d'une famille normalis\'ee $\{h_{(\theta_0,\theta_1)}\}_{(\theta_0,\theta_1)\in\Theta_0\times\Theta_1}$ particuli\`ere. On peut d\'emontrer (voir la proposition 1 de l'annexe III)qu'en changeant de famille normalis\'ee on obtient les m\^emes votes sauf sur un ensemble de r\'ealisations n\'egligeable pour tout $P_\theta$.

Comparons ce vote avec celui obtenu pour le choix entre $P_{\theta_0}$ et
$P_{\theta_1}$ lorsque $\theta\in\{\theta_0,\theta_1\}$ (voir la d\'efinition 2.4.1). Le r\'esultat du vote des experts 
$\Phi_s^{(\theta_0,\theta_1)}\subseteq\overline{\Delta_s}$ d\'epend de la valeur de la statistique
$K_{(\theta_0,\theta_1)}=h_{(\theta_0,\theta_1)}(T)$.
Posons $K_{(\theta_0,\theta_1)}(\omega)=k'$, la fr\'equence du vote en faveur de $P_{\theta_1}$ est \'egale \`a $1$ si $k'=0$, \`a $0$ si $k'=+\infty$ et \`a
$[P_{\theta}(\{K_{(\theta_0,\theta_1)}>k'\})+P_{\theta}(\{K_{(\theta_0,\theta_1)}\geq k'\})]/2$ sinon. Il est facile de voir que $K_{(\theta_0,\theta_1)}$ se factorise par $K(T)$ :
$K_{(\theta_0,\theta_1)}=h'_{(\theta_0,\theta_1)}(K(T))$, on a alors 
$k'=h'(k)$ avec $k=K(T(\omega))$. Dans le cas o\`u la r\'ealisation $\omega$ pose vraiment un probl\`eme de choix entre $P_{\theta_0}$ et $P_{\theta_1}$, c'est-\`a-dire $0<k'<+\infty$, le vote $Q^\omega_\theta(\{1\})$ des experts $\overline{\Delta_s}$ est compris entre $P_{\theta}(\{K_{(\theta_0,\theta_1)}>k'\})$ et $P_{\theta}(\{K_{(\theta_0,\theta_1)}\geq k'\})$. C'est un vote qui affine celui du choix entre $P_{\theta_0}$ et
$P_{\theta_1}$, en prenant en compte les autres cas possibles.
\medskip

Examinons maintenant le vote de l'ensemble des experts. D'apr\`es la proposition 4.2.2 tout expert $\phi$, non presque s\^urement \'egal \`a un \'el\'ement de $\overline{\Delta_s}$, est presque s\^urement strictement encadr\'e par deux \'el\'ements successifs de $\overline{\Delta_s}$ : $f\relmont{<}{p.s.}\phi\relmont{<}{p.s.}f'$. En dehors de l'\'ev\'enement
$A=\{f\not=f'\}=\{f=0\}\cap\{f'=1\}$, $\phi$ est presque s\^urement \'egal \`a $f$ et $f'$. Sur $A$ la statistique essentielle $K(T)$ est constante, tous les rapports de densit\'es $K_{(\theta_0,\theta_1)}$ le sont donc aussi. Les r\'ealisations $\omega$ appartenant \`a $A$ n'ont aucune raison de se diff\'erencier par rapport \`a la d\'ecision prise. C'est ce que font les experts $f$ et $f'$, qui sont d'ailleurs les r\`egles $f_{(a_t,u_t)}$ et $f_{(b_t,v_t)}$ de la d\'efinition 4.3.1, $t$ \'etant une valeur quelconque prise par la statistique $T$ sur $A$. Les experts compris  entre $f$ et $f'$ diff\`erent principalement par la proportion de $A$ pour laquelle ils d\'ecident $d=1$.
L'ensemble de ces experts d\'epend des \'ev\'enements inclus dans $A$. Si on veut \'eviter l'intervention de ces probl\`emes de mesurabilit\'e on peut consid\'erer sur $A$ les experts al\'eatoires constants. Les r\'ealisations $\omega$ de $A$ sont ainsi trait\'ees de fa\c con semblable et toutes les proportions sur $A$ de la d\'ecision $d=1$ sont permises. Ceci revient \`a consid\'erer le sur-mod\`ele 
$\Omega\times[0,1]$ muni de la probabilit\'e $P'_{\theta}$ produit de $P_{\theta}$ par la loi uniforme sur $[0,1]$ et les experts d\'eterministes de la forme :
$\II_{\{K(T)<k\}}(\omega)\,+\,\II_{\{K(T)=k\}}(\omega).\II_{[0,\beta]}(u)$, $\beta$ variant entre $0$ et $1$ (la restriction de $\beta$ \`a $\{0,1\}$
d\'efinit les experts de $\overline{\Delta_s}$). Dans ce sur-mod\`ele l'ensemble des experts consid\'er\'es peut \^etre probabilis\'e, comme pr\'ec\'edemment,
par l'op\'erateur $E'_\theta$ moyenne par rapport \`a $P'_\theta$.
Les deux probabilit\'es possibles sont ici identiques et on v\'erifie facilement que l'on obtient le m\^eme vote que celui des experts de $\overline{\Delta_s}$.
Dans la suite nous ne consid\`ererons que ce type de votes. Les autres votes que l'on peut construire sur l'ensemble des experts \`a partir de $E_\theta$ 
sont toujours compris entre les votes d\'efinis sur $\overline{\Delta_s}$ \`a partir  de $m_\theta$ et $m'_\theta$. Ils ne peuvent \^etre diff\'erents sur $A$ 
que si $P_\theta(A)>0$.

Nous nous retrouvons avec un vote pour chaque valeur possible du param\`etre $\theta$. Avant d'en choisir un ou d'en construire une synth\`ese nous allons montrer qu'on avantage la d\'ecision $d=1$ (resp. $d=0$) en utilisant $\theta$
appartenant  \`a $\Theta_0$ (resp. $\Theta_1$). Ceci se comprend ais\'ement puisque $P_{\theta_0}$ (resp. $P_{\theta_1}$) a plut\^ot tendance \`a charger les grandes (resp. petites) valeurs de $T$.

\vfill\eject

\medskip
{\bf Proposition 4.3.2}
\medskip
\medskip
\moveleft 10.4pt\hbox{\vrule\kern 10pt\vbox{\defpro

i) $\forall\, \theta_0\in\Theta_0$ et $\forall\, \theta_1\in\Theta_1$ on a :
\dli $Q^\omega_{\theta_0}(\{1\})\geq Q^\omega_{\theta_1}(\{1\})$ et bien s\^ur
$Q^\omega_{\theta_0}(\{0\})\leq Q^\omega_{\theta_1}(\{0\})$.

ii) S'il existe $\theta$ tel que $Q^\omega_{\theta}(\{1\})\relmont{=}{\{\theta\}p.s.}1-G_{\theta}(K(T))$, sous $P_{\theta}$ le vote $Q_{\theta}$ est neutre : 
$E_\theta(Q^\omega_{\theta}(\{1\}))=E_\theta(Q^\omega_{\theta}(\{0\}))={1\over 2}$.
}}\medskip

{\leftskip=15mm \dli {\bf D\'emonstration}
\medskip
i) Soient $\theta_0\in\Theta_0$, $\theta_1\in\Theta_1$ et $\omega\in\Omega$, posons $T(\omega)=t$.
\dli D'apr\`es la proposition 4.3.1, les votes $Q^\omega_{\theta_0}(\{1\})$ et $Q^\omega_{\theta_1}(\{1\})$ sont \'egaux \`a $1$ (resp. $0$) lorsque $t\in(D_i-D_s)$ (resp. $t\in D_s$). Il nous reste \`a consid\'erer le cas : $t\in\IR\!-(D_i\cup D_s)$
\dli On a alors $Q^\omega_{\theta}(\{1\})=1-G_{\theta}(K(t))$, $G_{\theta}$ \'etant la fonction de r\'epartition moyenne de la statistique essentielle $K(T)$ : $G_{\theta}(K(t))=E_{\theta}(f_t)$ avec
$f_t=\II_{\{K(T)<K(t)\}} + {1\over 2}\II_{\{K(T)=K(t)\}}$.
\dli Soient $\theta_0\in\Theta_0$ et $\theta_1\in\Theta_1$, il nous faut d\'emontrer l'in\'egalit\'e :
\dli $E_{\theta_0}(f_t)\leq E_{\theta_1}(f_t)$, c'est-\`a-dire 
$\int_\Omega f_t(p_{\theta_1}-p_{\theta_0})\,d\mu\geq 0$.
\dli Pour cela nous utilisons la stabilit\'e des hypoth\`eses qui suppose l'exis\-tence d'une fonction mesurable croissante 
$h_{(\theta_0,\theta_1)} :\, \IR \rightarrow \overline{\IR^+}$
v\'erifiant
$p_{\theta_0}/p_{\theta_1}=h_{(\theta_0,\theta_1)}(T)$
sur le domaine de d\'efinition de ce rapport, c'est-\`a-dire en dehors de
$\{\omega\in\Omega\, ;\, p_{\theta_0}(\omega)=p_{\theta_1}(\omega)=0\}$.
\dli La valeur de $K(t)$ d\'epend de deux \'el\'ements de $\overline{\Delta_s}$,
$\delta_t=f_{(a_t,u_t)}$ et $\delta'_t=f_{(b_t,v_t)}$, car $t\notin (D_i\cup D_s)$. Ils v\'erifient $\delta_t\leq f_{(t,0)}<f_{(t,1)}\leq\delta'_t$ et 
$]\delta_t,\delta'_t[\cap\overline{\Delta_s}=\emptyset$. Dans le 2\up{\`eme} cas de la d\'emonstration de la proposition 4.3.1 nous avons montr\'e l'\'egalit\'e : $\{K(T)=K(t)\}=\{\delta'_t-\delta_t=1\}$. Quel que soit le cas de figure :
$\phi^{(\theta_0,\theta_1)}_{(\infty,0)}\leq \delta_t$, $\phi^{(\theta_0,\theta_1)}_{(0,1)}\leq \delta_t<\delta'_t\leq\phi^{(\theta_0,\theta_1)}_{(\infty,0)}$ ou $\delta'_t\leq\phi^{(\theta_0,\theta_1)}_{(0,1)}$, on a toujours $\{K(T)=K(t)\}\subseteq\{h_{(\theta_0,\theta_1)}(T)=h_{(\theta_0,\theta_1)}(t)\}$ ($\phi^{(\theta_0,\theta_1)}_{(k,\beta)}$ est une fonction de test simple bas\'ee sur $K_{(\theta_0,\theta_1)}=h_{(\theta_0,\theta_1)}(T)$).
\dli $K$ et $h_{(\theta_0,\theta_1)}$ \'etant croissantes on en d\'eduit :
$\II_A\leq f_t\leq\II_B$ avec \dli $A=\{h_{(\theta_0,\theta_1)}(T)<h_{(\theta_0,\theta_1)}(t)\}$ et $B=\{h_{(\theta_0,\theta_1)}(T)\leq h_{(\theta_0,\theta_1)}(t)\}$.
\dli Si $h_{(\theta_0,\theta_1)}(t)\leq 1$ on a 
$\II_B(p_{\theta_1}-p_{\theta_0})\geq 0$ donc 
$\int_\Omega f_t(p_{\theta_1}-p_{\theta_0})\,d\mu\geq 0$.
\dli Il nous reste le cas $h_{(\theta_0,\theta_1)}(t)> 1$ ; posons
$\Omega_1=\{h_{(\theta_0,\theta_1)}(T)\leq 1\}$, on a alors
$\II_{\Omega_1}\leq f_t$ ; il est facile de v\'erifier l'in\'egalit\'e suivante :
\dli $(f_t-\II_{\Omega_1})(p_{\theta_1}-p_{\theta_0})\geq 
(1-\II_{\Omega_1})(p_{\theta_1}-p_{\theta_0})$ puisque 
$(p_{\theta_1}-p_{\theta_0})$ est n\'egatif sur le compl\'ementaire de 
$\Omega_1$ ; en int\'egrant chacun des membres de cette in\'egalit\'e on obtient le r\'esultat recherch\'e : 
$\int_\Omega f_t(p_{\theta_1}-p_{\theta_0})\,d\mu\geq 0$.
\medskip
ii) Soit $\theta$ tel que $Q^\omega_{\theta}(\{1\})\relmont{=}{\{\theta\}p.s.}1-G_{\theta}(K(T))$.
\dli Il suffit de montrer l'\'egalit\'e $E_\theta(Q^\omega_{\theta}(\{0\}))={1\over 2}$ c'est-\`a-dire \dli $E_\theta(G_{\theta}(K(T)))={1\over 2}$. 
\dli Si $K(T)$ est diffuse $G_\theta$ est la fonction de r\'epartition classique de $K(T)$, elle est de plus continue. Nous savons que $G_{\theta}(K(T))$ suit une loi uniforme sur $[0,1]$ (cf. [Bre.] p. 284) qui est bien de moyenne ${1\over 2}$. 
\dli Si $K(T)$ n'est pas diffuse, on consid\`ere l'ensemble d\'enombrable des points charg\'es par $K(T)$. Le compl\'ementaire de ces points est une union d\'enombrable d'intervalles disjoints sur lesquels $G_\theta$ est continue et \'egale \`a la fonction de r\'epartition classique. La mesure trace de la probabilit\'e image de $P_\theta$ par $G_{\theta}(K(T))$ sur les intervalles $(a_i,b_i)$, images par $G_{\theta}$ des intervalles pr\'ec\'edents, est donc la mesure de Lebesgue. Ceci implique $E_\theta[G_{\theta}(K(T)).\II_{(a_i,b_i)}(G_{\theta}(K(T)))]={b_i^2-a_i^2\over 2}$. Dans l'intervalle $)b_i,a_{i+1}($ la statistique $G_{\theta}(K(T))$ ne prend qu'une valeur ${b_i+a_{i+1}\over 2}$ avec la probabilit\'e $a_{i+1}-b_i$. La moyenne sur ces intervalles est alors ${a_{i+1}^2-b_i^2\over 2}$. En sommant ces r\'esultats on obtient $E_\theta(G_{\theta}(K(T)))=-{0^2\over 2}+{1^2\over 2}={1\over 2}$.
\medskip\centerline{\hbox to 3cm{\bf \hrulefill}}\par}
\bigskip

Cette proposition nous montre que les votes $Q^\omega_{\theta}$ qui avantagent la d\'ecision $d=1$ sont ceux qui se font \`a partir de 
$\theta$ appartenant \`a $\Theta_0$.
Si $\Theta_0$ n'est pas trop limit\'e, en particulier s'il existe des $\theta_0$ qui donnent \`a $T$ des fonctions de r\'epartition presque nulles pour des valeurs r\'eelles aussi grandes que l'on veut, on aura des $\theta_0$ pour lesquels $Q^\omega_{\theta_0}(\{1\})$ sera proche de $1$.
G\'en\'eralement ces $\theta_0$ ne sont pas tr\`es probables pour l'utilisateur qui veut choisir entre $\Theta_0$ et $\Theta_1$. En effet, s'il d\'efinit le principe du partage de $\Theta$ en deux, il est rare qu'il essaie de limiter $\Theta$ \`a partir de ses connaissances empiriques. $\Theta$ contient g\'en\'eralement les valeurs du param\`etre math\'ematiquement possibles.
Le choix d'un vote final prendra alors peu en compte les votes associ\'es  \`a ces $\theta_0$. Cela peut se faire en prenant la moyenne des votes associ\'es \`a $\theta$ appartenant \`a $\Theta_0$, \`a partir d'une probabilit\'e $\Lambda_0$ 
bien choisie. On peut aussi \'eviter les votes $Q^\omega_{\theta_0}$ trop favorables \`a  la d\'ecision $d=1$ en prenant, dans $\Theta_0$, le vote le moins favorable \`a $d=1$. S'il existe c'est un cas particulier du choix pr\'ec\'edent, $\Lambda_0$ \'etant une masse de Dirac, dans le cas contraire on peut tout de m\^eme d\'efinir un vote r\'epondant \`a ce crit\`ere. Nous allons d\'efinir ces deux types de votes ainsi que les deux qui correspondent \`a
$\Theta_1$.

\medskip
{\bf D\'efinition 4.3.2}
\medskip
\medskip
\moveleft 10.4pt\hbox{\vrule\kern 10pt\vbox{\defpro

Soit $(\Omega,{\cal A},(P_\theta=p_\theta.\mu)_{\theta\in\Theta_0\cup \Theta_1})$ un probl\`eme de d\'ecision \`a hypoth\`eses stables. $\Theta=\Theta_0\cup\Theta_1$ est muni d'une tribu ${\cal T}$ et on suppose que $Q^\omega_{\theta}(\{1\})$ est mesurable pour presque tout $\omega$.
\dli $\Lambda_0$ (resp. $\Lambda_1$) \'etant une probabilit\'e sur l'espace mesurable $(\Theta_0,{\cal T}_0)$ (resp. $(\Theta_1,{\cal T}_1)$), on appelle vote pond\'er\'e par $\Lambda_0$ (resp. $\Lambda_1$) la probabilit\'e d\'efinie sur $D=\{0,1\}$, pour presque toute r\'ealisation $\omega$, par :
\dli $Q^\omega_{\Lambda_0}(\{1\})=\int_{\Theta_0}Q^\omega_{\theta}(\{1\})\,d\Lambda_0(\theta)$ (resp. $Q^\omega_{\Lambda_1}(\{1\})=\int_{\Theta_1}Q^\omega_{\theta}(\{1\})\,d\Lambda_1(\theta)$).
\dli On appelle vote le plus favorable sous $\Theta_0$ (resp. $\Theta_1$), la probabilit\'e d\'efinie sur $D=\{0,1\}$, pour toute r\'ealisation $\omega$, par : \dli 
$Q^\omega_{\Theta_0}(\{0\})=sup_{\theta\in\Theta_0}Q^\omega_{\theta}(\{0\})$ (resp. 
$Q^\omega_{\Theta_1}(\{1\})=sup_{\theta\in\Theta_1}Q^\omega_{\theta}(\{1\})$.
\dli Ces quatre types de votes v\'erifient :
$Q^\omega_{\Lambda_1}(\{1\})\leq Q^\omega_{\Theta_1}(\{1\})\leq
Q^\omega_{\Theta_0}(\{1\})\leq Q^\omega_{\Lambda_0}(\{1\})$.
}}\medskip

Cette propri\'et\'e entre les quatre types de votes d\'ecoule directement des d\'efinitions et de la proposition 4.3.2. 
\dli $Q^\omega_{\Theta_0}$ et $Q^\omega_{\Theta_1}$ repr\'esentent les votes les plus proches parmi ceux bas\'es sur $\Theta_0$ d'une part et $\Theta_1$
d'autre part.
Ils peuvent \^etre identiques sur presque tout $\Omega$, nous dirons dans ce cas que les hypoth\`eses sont adjacentes. On peut alors donner le vote 
$Q^\omega_{\Theta_0}\relmont{=}{p.s.}Q^\omega_{\Theta_1}$ comme aide \`a la d\'ecision, quand on observe $\omega$.

Quand les hypoth\`eses ne sont pas adjacentes ou quand les votes pond\'er\'es par $\Lambda_0$ et $\Lambda_1$ ne sont pas identiques, ce qui est g\'en\'eralement le cas, on a deux votes pour nous aider \`a prendre une d\'ecision. L'un 
$Q^\omega_0$ est bas\'e sur $\Theta_0$, l'autre $Q^\omega_1$ est bas\'e sur $\Theta_1$.
Les solutions imagin\'ees dans le cas du choix entre deux probabilit\'es sont encore possibles. Nous allons les passer en revue rapidement. Nous en reparlerons plus en d\'etail au chap\^{\i}tre suivant, dans le cadre classique des mod\`eles \`a rapport de vraisemblance monotone. 
\dli Comme dans la th\'eorie des tests on peut privil\'egier une hypoth\`ese. Ici, ceci revient \`a utiliser un seul des deux votes, $Q^\omega_0$ ou $Q^\omega_1$, comme aide \`a la d\'ecision. 
\dli De fa\c con semblable \`a ce qui a \'et\'e fait au paragraphe 2.5, on peut aussi utiliser un vote pond\'er\'e par
$\lambda\in[0,1]$ : $Q_\lambda^\omega(\{1\})=(1-\lambda) Q_0^\omega(\{1\})
+\lambda Q_1^\omega(\{1\})$. Si $Q^\omega_0$ et $Q^\omega_1$ correspondent \`a deux probabilit\'es $\Lambda_0$ et $\Lambda_1$, $Q_\lambda^\omega$ est le vote pond\'er\'e par la probabilit\'e 
$\Lambda_\lambda =(1-\lambda) \Lambda_0 + \lambda \Lambda_1$. Cette d\'ecomposition de la pond\'eration est int\'eressante car $\lambda$ permet d'exprimer un choix a priori entre les deux hypoth\`eses. Ce choix peut bien s\^ur provenir d'un vote ant\'erieur. Comme nous l'avons vu au paragraphe 2.5, le vote bas\'e sur $\Lambda_\lambda$ n'est pas d\'efini par les probabilit\'es a post\'eriori de $\Theta_0$ et $\Theta_1$ obtenues \`a partir de la loi a priori $\Lambda_\lambda$.
\dli Enfin, on peut prendre une d\'ecision \`a partir de la r\`egle des pl\'ebiscites d\'efinie au paragraphe 3.2.
Ceci fournit une g\'en\'eralisation des r\`egles de Bol'shev.
On peut aussi consid\'erer les r\`egles de d\'ecision bas\'ees
sur la diff\'erence $G(\omega)=Q^\omega_1(\{0\})-Q^\omega_0(\{1\})\in[-1,+1]$,
qui ne prennent aucune d\'ecision autour de $0$ et
d\'ecident $d=0$ ou $d=1$ ailleurs suivant que $G(\omega)$ est positif ou n\'egatif.

\vfill\eject

{\parindent=-10mm\soustitre 5--MOD\`ELES \`A RAPPORT DE VRAISEMBLANCE MONOTONE.}
\bigskip
{\parindent=-5mm 5.1 HYPOTH\`ESES UNILAT\'ERALES.}
\nobreak
\medskip

Les mod\`eles \`a rapport de vraisemblance monotone recouvrent un grand nombre
de cas classiques (cf. [Kar.], [KarR], [Leh.]).
En particulier les mod\`eles exponentiels \`a param\`etre r\'eel, mais aussi les mod\`eles statistiques portant sur le param\`etre de non centralit\'e d'une famille de densit\'es de Student, Fisher ou khi-deux.
\dli Suivant les auteurs, la d\'efinition des mod\`eles \`a rapport de vraisemblance monotone peut prendre diff\'erentes formes, nous utiliserons la suivante :

\medskip
{\bf D\'efinition 5.1.1}
\medskip
\medskip
\moveleft 10.4pt\hbox{\vrule\kern 10pt\vbox{\defpro

Soit $(\Omega ,{\cal A},(p_\theta.\mu)_{\theta\in\Theta})$ un mod\`ele
statistique domin\'e par la mesure $\mu$, $\Theta$ \'etant totalement ordonn\'e par la relation $\preceq$.
Il est \`a rapport de vraisemblance monotone s'il existe une statistique
$T$  \`a valeur dans $\IR$ telle que pour tout $\theta'$ et $\theta''$
de $\Theta$, $\theta'\prec\theta''$, il existe une fonction croissante $h_{(\theta'',\theta')} :\, \IR \rightarrow \overline{\IR^+}$ v\'erifiant
$p_{\theta''}/p_{\theta'}=h_{(\theta'',\theta')}(T)$
sur le domaine de d\'efinition de ce rapport c'est-\`a-dire en dehors de
$\{\omega\in\Omega\, ;\, p_{\theta'}(\omega)=p_{\theta''}(\omega)=0\}$.
}}\medskip

Dans certaines d\'efinitions l'\'egalit\'e
$p_{\theta''}/p_{\theta'}=h_{(\theta'',\theta')}(T)$ est remplac\'ee par une \'egalit\'e presque s\^ure :
$p_{\theta''}/p_{\theta'}\relmont{=}{\mu  p.s.} h_{(\theta'',\theta')}(T)$.
Nous ne l'avons pas fait car
l'\'ev\'enement $\mu$-n\'egligeable sur lequel on n'a pas l'\'egalit\'e d\'epend de
$\theta'$ et $\theta''$. Ceci peut poser probl\`eme quand on utilise cette \'egalit\'e sur une infinit\'e non d\'enombrable de couples $(\theta',\theta'')$.
De plus ce type de d\'efinitions est \'equivalent \`a la d\'efinition choisie ici,
lorsque la mesure $\mu$ est $\sigma$-finie (cf. [Pfa.]).

Comme dans la th\'eorie des tests nous allons commencer par \'etudier les hypoth\`eses unilat\'erales.
Ce sont les hypoth\`eses naturellement stables dans un mod\`ele \`a rapport de vraisemblance monotone, c'est-\`a-dire celles qui le sont \`a partir de la statistique $T$. Nous verrons qu'elles ne sont cependant pas les seules hypoth\`eses stables.
\eject

\medskip
{\bf D\'efinition 5.1.2}
\medskip
\medskip
\moveleft 10.4pt\hbox{\vrule\kern 10pt\vbox{\defpro

Soit $(\Omega ,{\cal A},(p_\theta.\mu)_{\theta\in\Theta})$ un mod\`ele
statistique \`a rapport de vraisem\-blance monotone.
\dli On appelle hypoth\`eses unilat\'erales, les hypoth\`eses d\'efinies par une partition ordonn\'ee $\{\Theta_1 , \Theta_0\}$ de $\Theta$.
\dli $\Theta=\Theta_1 + \Theta_0$, $\Theta_1\not=\emptyset$, 
$\Theta_0\not=\emptyset$ et $\forall \theta_1\in\Theta_1$, $\forall \theta_0\in\Theta_0$ : $\theta_1\prec\theta_0$.
\dli Ces hypoth\`eses sont bien s\^ur stables.
}}\medskip

Les hypoth\`eses unilat\'erales ne sont pas toujours les seules hypoth\`eses stables d'un  mod\`ele \`a rapport de vraisemblance monotone.
Par exemple si $P_\theta$ d\'esigne
la probabilit\'e uniforme sur $[\theta-1 ,\theta+1]$ le mod\`ele
$(\IR ,{\cal B},(P_\theta)_{\theta\in\Theta\subseteq\IR})$ 
est \`a rapport de vraisemblance monotone pour l'ordre ordinaire sur $\Theta\subseteq\IR$ et la statistique identit\'e $T$ ; les deux hypoth\`eses
$\Theta_0 = \{0\}$ et $\Theta_1 = \{-2 , +2\}$ sont stables, les experts
\'etant presque s\^urement \'egaux \`a $1-\II_{[-1,+1]}$ ; 
les deux hypoth\`eses $\Theta_0 = \{-1 , +1\}$ et $\Theta_1 = ]-\infty,-1[\cup ]+1,+\infty[$ le sont aussi.
Remarquons, que dans ces deux cas on peut retrouver des hypoth\`eses unilat\'erales en utilisant la statistique $-\!\mid\! T\!\mid$ et en changeant l'ordre sur $\Theta$.
Quand cet ordre est fondamental pour l'interpr\'etation, les hypoth\`eses unilat\'erales deviennent g\'en\'eralement les seules hypoth\`eses stables int\'eressantes.

Nous allons maintenant comparer nos r\'esultats \`a ceux obtenus pour les hypoth\`eses unilat\'erales classiquement \'etudi\'ees en th\'eorie des tests :
$\Theta\subseteq\IR$ est muni de l'ordre ordinaire et $H_0$ est l'hypoth\`ese $\theta\leq\theta_0$ ou $\theta\geq\theta_0$. Pour cette derni\`ere hypoth\`ese nous savons qu'il existe au seuil $\alpha\in [0,1]$ un test uniform\'ement plus puissant d\'efini par :
$$\phi(\omega)\quad=\quad\left\{\matrix{
1\hfill & si & T(\omega)<c\hfill\cr
\gamma\hfill & si & T(\omega)=c\hfill \cr
0\hfill & si & T(\omega)>c\hfill \cr
}\right. $$
$c\in\overline{\IR}$ et $\gamma\in[0,1]$ v\'erifiant $E_{\theta_0}(\phi)=\alpha$ (cf. par exemple [Leh.] p. 78).
\dli Consid\'erons la notion de p-value, de seuil minimum de rejet, pour ce type de tests. A une r\'ealisation $\omega'$ on peut associer le seuil minimum
$P_{\theta_0}(\{T\leq T(\omega')\})$ ou $P_{\theta_0}(\{T<T(\omega')\})$
suivant que l'on veut une probabilit\'e de rejeter en $\omega'$ \'egale \`a $1$ ou non nulle. Nous prendrons comme p-value la valeur interm\'ediaire :
$P_{\theta_0}(\{T<T(\omega')\}) + (1/2)P_{\theta_0}(\{T=T(\omega')\})$.
Ces trois valeurs sont bien s\^ur identiques lorsque $T$ est diffuse.
\dli Le seuil minimum de rejet ainsi d\'efini est bas\'e sur une famille de tests uniform\'ement plus puissants, construite \`a partir de la statistique $T$. 
Nous pouvons aussi construire une famille semblable \`a partir de la statistique essentielle $K(T)$ li\'ee aux hypoth\`eses stables $H_0\,:\,\theta\geq\theta_0$ et $H_1\,:\,\theta<\theta_0$ (voir la d\'efinition 4.3.1). On obtient alors un seuil minimum de rejet qui peut \^etre interpr\'et\'e comme le r\'esultat d'un vote d'experts.
\dli S'il existe des tests de puissance $1$ sur tout $\Theta_1$ et de niveaux diff\'erents, on peut modifier la d\'efinition pr\'ec\'edente en prenant une p-value de $1$ lorsque la r\'ealisation $\omega'$ appartient \`a la r\'egion de non rejet d'un de ces tests. Ceci est tout \`a fait justifi\'e car dans ce cas aucune probabilit\'e $P_{\theta_1}$ ne charge cette r\'egion, il est donc difficile de d\'ecider $H_1$.

\medskip
{\bf Proposition 5.1.1}
\medskip
\medskip
\moveleft 10.4pt\hbox{\vrule\kern 10pt\vbox{\defpro

Soit $(\Omega ,{\cal A},(p_\theta.\mu)_{\theta\in\Theta\subseteq\IR})$ un mod\`ele statistique \`a rapport de vraisem\-blance monotone fond\'e sur la statistique $T$. Consid\'erons un \'el\'ement $\theta_0$ de $\Theta$ diff\'erent de $inf\Theta$, les hypoth\`eses $\Theta_0=\Theta\cap\{\theta\geq\theta_0\}$ et $\Theta_1=\Theta\cap\{\theta<\theta_0\}$ sont unilat\'erales.
\dli $Q_{\Theta_0}^\omega(\{0\})=Q_{\theta_0}^\omega(\{0\})$ d\'efinit un seuil minimum de rejet de $H_0\,:\,\theta\geq\theta_0$ contre $H_1\,:\,\theta<\theta_0$, pour les tests uniform\'ement plus puissants construits avec la statistique essentielle $K(T)$.
\dli Si les hypoth\`eses sont adjacentes, c'est-\`a-dire
$Q_{\Theta_0}^\omega\relmont{=}{p.s.}Q_{\Theta_1}^\omega$, $Q_{\Theta_1}^\omega(\{1\})\relmont{=}{p.s.}Q_{\theta_0}^\omega(\{1\})$ d\'efinit un seuil minimum de rejet de $H_0\,:\,\theta<\theta_0$ contre $H_1\,:\,\theta\geq\theta_0$.
}}\medskip
\medskip
Dans cette proposition nous faisons r\'ef\'erence aux r\'esultats du paragraphe sur les hypoth\`eses stables. On suppose donc travailler avec une famille 
$\{h_{(\theta_0,\theta_1)}\}_{(\theta_0,\theta_1)\in\Theta_0\times\Theta_1}$ 
normalis\'ee (voir la d\'efinition 4.2.1). 
Il est toujours possible d'en construire une \`a partir des fonctions $h_{(\theta'',\theta')}$, $\theta'\prec\theta''$, d\'efinies avec le mod\`ele \`a rapport de vraisemblance monotone \'etudi\'e (voir l'annexe II).
\dli Il serait int\'eressant de choisir les fonctions
$h_{(\theta'',\theta')}$ d'un mod\`ele \`a rapport de vraisemblance monotone de telle sorte qu'elles d\'efinissent une famille normalis\'ee pour toutes les hypoth\`eses unilat\'erales. Ceci est uniquement possible lorsque les diff\'erentes hypoth\`eses unilat\'erales d\'efinissent les m\^emes demi-droites $D_i$ et $D_s$ (voir l'annexe IV).
Par contre on peut toujours construire des $h_{(\theta'',\theta')}$ v\'erifiant la propri\'et\'e i) de la d\'efinition 4.2.1, il suffit de suivre la premi\`ere \'etape de l'annexe II. On peut alors en d\'eduire facilement une famille normalis\'ee (voir l'annexe IV).

\medskip
{\leftskip=15mm \dli {\bf D\'emonstration}
\medskip
{\parindent=-10mm a) Montrons d'abord que $Q_{\theta}^\omega(\{1\})$ est une fonction croissante en $\theta$.}

Soient $\omega\in\Omega$, $\theta'\in\Theta$ et $\theta''\in\Theta$ tels que
$\theta'<\theta''$, on cherche \`a montrer l'in\'egalit\'e :
$Q_{\theta'}^\omega(\{1\})\leq Q_{\theta''}^\omega(\{1\})$.
\dli Nous avons d\'ej\`a obtenu ce r\'esultat pour $\theta'\in\Theta_1$ et $\theta''\in\Theta_0$ dans le cadre plus g\'en\'eral de la proposition 4.3.2.
Nous devons l'\'etendre ici au cas o\`u $\theta'$ et $\theta''$ appartiennent \`a la m\^eme hypoth\`ese. La d\'emonstration est semblable \`a la pr\'ec\'edente.
\dli D'apr\`es la proposition 4.3.1, $Q_{\theta}^\omega(\{1\})$ d\'epend de la statistique essentielle $K(T)$ li\'ee aux hypoth\`eses stables
$\Theta_0=\Theta\cap\{\theta\geq\theta_0\}$ et \dli $\Theta_1=\Theta\cap\{\theta<\theta_0\}$(voir la d\'efinition 4.3.1).
\dli Si $T(\omega)\in(D_i\cup D_s)$, on a la propri\'et\'e recherch\'ee puisque : $Q_{\theta'}^\omega(\{1\})=Q_{\theta''}^\omega(\{1\})$.
Dans le cas contraire, on a $Q_{\theta}^\omega(\{1\})=1-G_{\theta}(k)$ avec 
$k=K(T(\omega))$ et $G_{\theta}(k)={1\over 2}E_\theta(\II_{\{K(T)<k\}})+{1\over 2}E_\theta(\II_{\{K(T)\leq k\}})$.
L'in\'egalit\'e recherch\'ee est donc \'equivalente \`a 
$G_{\theta'}(k)\geq G_{\theta''}(k)$.
\dli $K$ \'etant croissante, l'ensemble $\{t\in\IR\,;\,K(t)=k\}$ est un intervalle non vide ; il existe donc, dans l'ensemble $F$ des fonctions de test $f_{(t,u)}$ d\'efinies au d\'ebut du paragraphe 4.2, deux \'el\'ements 
$f_{(t',u')}$ et $f_{(t'',u'')}$ tels que : 
$f_{(t',u')}=\II_{\{K(T)<k\}}$ et $f_{(t'',u'')}=\II_{\{K(T)\leq k\}}$.
\dli L'in\'egalit\'e $G_{\theta'}(k)\geq G_{\theta''}(k)$ sera d\'emontr\'ee si pour toute fonction de test $f_{(t,u)}$ de $F$ on montre l'in\'egalit\'e
$E_{\theta'}(f_{(t,u)})\geq E_{\theta''}(f_{(t,u)})$, c'est-\`a-dire
$\int_\Omega f_{(t,u)}(p_{\theta'}-p_{\theta''})\,d\mu\geq 0$. Ceci est une propri\'et\'e classique des mod\`eles \`a rapport de vraisemblance monotone. 
Il suffit de distinguer deux cas :
\dli i) $h_{(\theta'',\theta')}(t)\leq 1$.
\dli Posons $A=\{T\leq t\}$, $h_{(\theta'',\theta')}$ \'etant croissante, le rapport $p_{\theta''}/p_{\theta'}$ est inf\'erieur ou \'egal \`a $1$ sur l'\'ev\'enement $A$, ce qui se traduit par : 
$\II_A . (p_{\theta'}-p_{\theta''})\geq 0$ ; comme $f_{(t,u)}$ est nulle sur le compl\'ementaire de $A$ on a bien :
$\int_\Omega f_{(t,u)}(p_{\theta'}-p_{\theta''})\,d\mu\geq 0$.
\dli ii) $h_{(\theta'',\theta')}(t)>1$.
\dli Posons $\Omega_1=\{h_{(\theta'',\theta')}(T)\leq 1\}$, on a alors
$\II_{\Omega_1}\leq f_{(t,u)}$ ; il est facile de v\'erifier l'in\'egalit\'e suivante :
\dli $(f_{(t,u)}-\II_{\Omega_1})(p_{\theta'}-p_{\theta''})\geq 
(1-\II_{\Omega_1})(p_{\theta'}-p_{\theta''})$ puisque 
$(p_{\theta'}-p_{\theta''})$ est n\'egatif sur le compl\'ementaire de 
$\Omega_1$ ; en int\'egrant chacun des membres de cette in\'egalit\'e on obtient le r\'esultat recherch\'e : 
$\int_\Omega f_{(t,u)}(p_{\theta'}-p_{\theta''})\,d\mu\geq 0$.

{\parindent=-10mm b) $Q_{\Theta_0}^\omega(\{0\})=Q_{\theta_0}^\omega(\{0\})$ est un seuil minimum de rejet.}

L'\'egalit\'e $Q_{\Theta_0}^\omega(\{0\})=Q_{\theta_0}^\omega(\{0\})$ est une cons\'equence directe de la propri\'et\'e d\'emontr\'ee en a), qui est \'equivalente \`a la d\'ecroissance de $Q_{\theta}^\omega(\{0\})$ par rapport $\theta$ ; en effet, d'apr\`es la d\'efinition 4.3.2 on a  $Q^\omega_{\Theta_0}(\{0\})=sup_{\theta\in\Theta_0}Q^\omega_{\theta}(\{0\})$ et par  d\'efinition de $\Theta_0$ : $\forall\, \theta\in\Theta_0\quad\theta\geq\theta_0$.

Soit $K(T)$ la statistique essentielle li\'ee aux hypoth\`eses stables 
$\Theta_0$ et $\Theta_1$ (voir la d\'efinition 4.3.1).
Pour tout $\theta_0$ de $\Theta_0$ et tout $\theta_1$ de $\Theta_1$, le rapport des densit\'es $h_{(\theta_0,\theta_1)}(T)$ peut par construction de $K(T)$ s'\'ecrire $h'_{(\theta_0,\theta_1)}(K(T))$ avec $h'_{(\theta_0,\theta_1)}\, : \,\overline{\IR}\rightarrow\overline{\IR^+}$ croissante. Il suffit de montrer l'existence d'une fonction $h'_{(\theta_0,\theta_1)}$ de $K(\IR)$ dans $\overline{\IR^+}$, v\'erifiant $h'_{(\theta_0,\theta_1)}(K(t))=h_{(\theta_0,\theta_1)}(t)$ pour tout $t$ de $IR$ ; ordonnons les statistiques $f_{(a_t,u_t)}$ et $f_{(b_t,v_t)}$ de la d\'efinition 4.3.1, par rapport aux fonctions de test simples de $\Phi_s^{(\theta_0,\theta_1)}$, il y a trois cas possibles : $h_{(\theta_0,\theta_1)}(t)=0$ on a alors $f_{(b_t,v_t)}\leq \phi_{(0,1)}^{(\theta_0,\theta_1)}$, $h_{(\theta_0,\theta_1)}(t)=k\in]0,+\infty[$ on a alors $\phi_{(k,0)}^{(\theta_0,\theta_1)}\leq f_{(a_t,u_t)}<f_{(b_t,v_t)}\leq \phi_{(k,1)}^{(\theta_0,\theta_1)}$, $h_{(\theta_0,\theta_1)}(t)=+\infty$ on a alors $\phi_{(\infty,0)}^{(\theta_0,\theta_1)}\leq f_{(a_t,u_t)}$ ; l'\'egalit\'e $h'_{(\theta_0,\theta_1)}(K(t))=h_{(\theta_0,\theta_1)}(t)$ d\'efinit bien $h'_{(\theta_0,\theta_1)}$ sur $K(\IR)$ car $\{K(T)=K(t)\}=\{f_{(b_t,v_t)}-f_{(a_t,u_t)}=1\}$ ; nous avons obtenu cette \'egalit\'e pour $t\in\IR\!-(D_i\cup D_s)$ dans le 2\up{\`eme} cas de la d\'emonstration de la proposition 4.3.1, il est facile de la v\'erifier pour $t$ appartenant \`a $D_i-D_s$ ou \`a $D_s$ ($D_i$ et $D_s$ \'etant les demi-droites associ\'ees aux hypoth\`eses unilat\'erales $\{\Theta_1,\Theta_0\}$ dans la d\'efinition 4.2.1).

Nous avons vu qu'il existe, pour tout seuil $\alpha\in [0,1]$, un test uniform\'ement plus puissant (U.P.P.) d\'efini par :
$$\phi(\omega)\quad=\quad\left\{\matrix{
1\hfill & si & T(\omega)<c\hfill\cr
\gamma\hfill & si & T(\omega)=c\hfill \cr
0\hfill & si & T(\omega)>c\hfill \cr
}\right. $$
$c\in\overline{\IR}$ et $\gamma\in[0,1]$ v\'erifiant $E_{\theta_0}(\phi)=\alpha$.

\dli Tous les tests de la forme de $\phi$ sont U.P.P. au seuil $\alpha=E_{\theta_0}(\phi)$ si $\alpha>0$ ; dans le cas $\alpha=0$ il existe au moins un test U.P.P. $\phi_0$, c'est le plus grand test $\phi$ v\'erifiant
$E_{\theta_0}(\phi)=0$ (cf. [Mor.1] p. 36). On a $\phi_0=\II_{D_0}(T)$, $D_0=]-\infty,t_0)$ \'etant la plus grande demi-droite inf\'erieure ouverte ou ferm\'ee telle que $P_{\theta_0}(T^{-1}(D_0))=0$. Par d\'efinition, $D_i\subset D_0$, nous allons montrer que tous les tests $\phi\in[\II_{D_i-D_s}(T),\II_{D_0}(T)]$ sont U.P.P. au seuil $\alpha=0$.
Il suffit de montrer que $\II_{D_i-D_s}(T)$ et $\II_{D_0}(T)$ sont de m\^eme puissance : $\forall\theta_1\in\Theta_1$ $E_{\theta_1}(\II_{D_i-D_s}(T))=E_{\theta_1}(\II_{D_0}(T))$ ; c'est \'evident lorsque $D_i\cup D_s=\IR$, par d\'efinition de $D_s$. Dans le cas contraire
$D_i-D_s=D_i$ et d'apr\`es le lemme 1 de l'annexe IV, $D_i=D_i^{\theta_0}$ (la densit\'e $p_{\theta_0}$ \'etant nulle sur $T^{-1}(D_i^{\theta_0})$) ; posons $A=T^{-1}(D_0-D_i^{\theta_0})$, $A'=A\cap\{p_{\theta_0}=0\}$ et $A''=A\cap\{p_{\theta_0}>0\}$, l'\'ev\'enement $A''$ est $\mu$ n\'egligeable puisque $P_{\theta_0}(A)=0$ ; il nous reste \`a montrer que $A'$ est $P_{\theta_1}$ n\'egligeable pour tout $\theta_1\in\Theta_1$ ; en fait on a m\^eme $p_{\theta_1}$ nulle sur $A'$ car s'il existait $\omega\in A'$ tel que 
$p_{\theta_1}(\omega)>0$ on aurait $h_{(\theta_0,\theta_1)}$  nulle en $t=T(\omega)$ donc nul sur $]-\infty,t]$, $t$ appartiendrait \`a $D_i^{\theta_0}$ ce qui est impossible puisque $A'\subseteq A$.

Nous venons de montrer que les tests de la forme $\phi$ v\'erifiant $\phi\geq\II_{D_i-D_s}(T)$ sont uniform\'ement plus puissants \`a leur niveau, le test $\II_{D_i-D_s}(T)$ \'etant de niveau $0$.

Pour tout test de la forme $\phi$ on peut construire un test de m\^eme puissance de la forme :
$$\phi'(\omega)\quad=\quad\left\{\matrix{
1\hfill & si & K(T(\omega))<c'\hfill\cr
\gamma'\hfill & si & K(T(\omega))=c'\hfill \cr
0\hfill & si & K(T(\omega))>c'\hfill \cr
}\right. $$
$c'\in\overline{\IR}$ et $\gamma'\in[0,1]$ v\'erifiant $E_{\theta_0}(\phi')=E_{\theta_0}(\phi)=\alpha$. 
\dli En effet, nous avons vu que sur l'\'ev\'enement $A_c=\{K(T)=K(c)\}$ le rapport des densit\'es $h_{(\theta_0,\theta_1)}(T)$ est constant, il suffit donc de poser $c'=K(c)$, $\gamma'=[P_{\theta_0}(\{T<c\}\cap A_c)+\gamma P_{\theta_0}(\{T=c\})]/P_{\theta_0}(A_c)$ si $P_{\theta_0}(A_c)>0$ et $\gamma'=1$ sinon.
\dli Par construction de $K$ on a $\II_{D_i-D_s}(T)=\II_{\{K(T)=-\infty\}}$, les tests de la forme $\phi'$ v\'erifiant $\phi'\geq\II_{D_i-D_s}(T)$  constituent une famille de tests U.P.P. \`a leur niveau. Cette famille contient au moins un test de niveau $\alpha$ pour tout $\alpha$ de $[0,1]$.
Elle nous permet de d\'efinir, pour toute r\'ealisation $\omega$, un seuil minimum de rejet $\alpha_m(\omega)$ : \'egal \`a $0$ lorsque $T(\omega)\in D_i-D_s$ et \'egal \`a $P_{\theta_0}(\{K(T)<K(T(\omega))\}) + (1/2)P_{\theta_0}(\{K(T)=K(T(\omega))\})$ sinon.
\dli $\alpha_m(\omega)$ est la valeur de la fonction de r\'epartition moyenne $G_{\theta_0}(k)$ de la statistique essentielle $K(T)$ en $k=K(T(\omega))$. 
\dli Lorsque $T(\omega)\notin D_i\cup D_s$ on obtient donc bien 
$\alpha_m(\omega)=Q_{\theta_0}^\omega(\{0\})$. Il en est de m\^eme pour $T(\omega)\in D_i-D_s$ puisque dans ce cas  
$Q_{\theta_0}^\omega(\{0\})$ et $\alpha_m(\omega)$ valent $0$.
\dli Il nous reste le cas $T(\omega)\in D_s$, on a alors $Q_{\theta_0}^\omega(\{0\})=1$ ; 
les tests construits \`a partir de $K(T)$ avec $c'=+\infty$ ont tous une 
puissance \'egale \`a $1$ sur $\Theta_1$ ; le seul vraiment int\'eressant est celui d\'efini par $\gamma'=0$ ; les autres sont de m\^eme puissance mais ils peuvent avoir un niveau sup\'erieur, quand ce n'est pas le cas on a $P_{\theta_0}(\{K(T)=+\infty\})=0$ donc $\alpha_m(\omega)=1=Q_{\theta_0}^\omega(\{0\})$. Dans tous les cas il est en fait inint\'eressant de consid\'erer les tests $\phi'$ strictement sup\'erieur \`a $1-\II_{D_s}$ ; il n'existe alors plus de test qui puisse rejeter $H_0$ en $\omega$ ; il semble naturel de poser
$\alpha_m(\omega)=1$ qui signifie bien qu'il n'est pas question de rejeter dans un tel cas de figure.

{\parindent=-10mm c) Cas des hypoth\`eses adjacentes.}

On a $Q_{\Theta_1}^\omega\relmont{=}{p.s.}Q_{\Theta_0}^\omega=Q_{\theta_0}^\omega$.
Ce qui implique :  \dli $Q_{\Theta_1}^\omega(\{1\})=sup_{\theta\in\Theta_1}Q^\omega_{\theta}(\{1\})\relmont{=}{p.s.}Q_{\theta_0}^\omega(\{1\})$.
\dli Pour d\'efinir un seuil minimum de rejet de $H_0\,:\,\theta<\theta_0$ contre $H_1\,:\,\theta\geq\theta_0$ on consid\`ere comme en b) les tests uniform\'ement plus puissants de la forme :
$$\phi''(\omega)\quad=\quad\left\{\matrix{
1\hfill & si & K(T(\omega))>c''\hfill\cr
\gamma''\hfill & si & K(T(\omega))=c''\hfill \cr
0\hfill & si & K(T(\omega))<c''\hfill \cr
}\right. $$
$c''\in\overline{\IR}$ et $\gamma''\in[0,1]$ v\'erifiant $sup_{\theta\in\Theta_1}E_{\theta}(\phi'')=\alpha$.
\dli Lorsque $T(\omega)\notin D_i\cup D_s$ le seuil minimum de rejet
$\alpha'_m(\omega)$ est \'egal \`a $sup_{\theta\in\Theta_1}[1-G_{\theta}(k)]=[1-G_{\theta_0}(k)]$, donc
$\alpha'_m(\omega)=Q_{\theta_0}^\omega(\{1\})$. 
\dli Lorsque $T(\omega)\in D_s$, on a $P_{\theta}(T^{-1}(D_s))=0$ pour tout $\theta$ de $\Theta_1$ ($\theta<\theta_0$) ; le seuil minimum de rejet est donc nul, il en est de m\^eme de $Q_{\theta_0}^\omega(\{1\})$.
\dli Il nous reste le cas $T(\omega)\in D_i-D_s$, cette fois c'est l'\'ev\'enement $B=T^{-1}(D_i-D_s)$ qui est de probabilit\'e nulle pour tout
$\theta$ de $\Theta_0$ ; les tests $\phi''$ avec $c''=-\infty$ sont tous de puissance \'egale \`a $1$, le seul vraiment int\'eressant est celui d\'efini par
$\gamma''=0$ ; si pour d\'efinir le seuil minimum de rejet on enl\`eve les tests correspondant \`a $c''=-\infty$ et $\gamma''>0$, il n'existe plus de test qui puisse rejeter $H_0$ en $\omega$ ; il est alors naturel (voir le cas $T(\omega)\in D_s$ de la partie b) de cette d\'emonstration) de poser
$\alpha'_m(\omega)=1$, ce qui est bien la valeur de $Q_{\theta_0}^\omega(\{1\})$.

\medskip\centerline{\hbox to 3cm{\bf \hrulefill}}\par}

L'autre test unilat\'eral classique, correspondant \`a $H_0\,:\,\theta\leq\theta_0$, conduit \`a une proposition semblable. 
Il faut prendre un \'el\'ement $\theta_0$ de $\Theta$ diff\'erent de $sup\Theta$, ceci permet de d\'efinir les hypoth\`eses unilat\'erales $\Theta'_0=\Theta\cap\{\theta\leq\theta_0\}$ et $\Theta'_1=\Theta\cap\{\theta>\theta_0\}$.
$Q_{\theta_0}^\omega(\{1\})$ d\'efinit alors un seuil minimum de rejet de $H_0\,:\,\theta\leq\theta_0$ contre $H_1\,:\,\theta>\theta_0$, pour les tests uniform\'ement plus puissants construits avec la statistique essentielle $K(T)$.

Revenons aux hypoth\`eses unilat\'erales \'etudi\'ees par la proposition 5.1.1 :
\dli $\Theta_0=\Theta\cap\{\theta\geq\theta_0\}$ et $\Theta_1=\Theta\cap\{\theta<\theta_0\}$ avec $\theta_0\not= inf\Theta$.
\dli Nous allons regarder ce que donne la r\`egle des pl\'ebiscites
(voir la d\'efinition 3.2.1) appliqu\'ee aux deux votes : 
$Q_{0}^\omega=Q_{\Theta_0}^\omega$ et $Q_{1}^\omega=Q_{\Theta_1}^\omega$.
Nous le ferons dans le cas o\`u les hypoth\`eses $\Theta_0$ et $\Theta_1$ sont adjacentes. Cette propri\'et\'e signifiant  $Q_{\Theta_0}^\omega=Q_{\Theta_1}^\omega$ presque partout, elle repose sur la famille des  fonctions de r\'epartition moyenne, $(G_\theta)_{\theta\in\Theta}$ de la statistique essentielle $K(T)$ (voir la d\'efinition 4.3.1). Il faut que pour presque tout $\omega$, la fonction de $\theta$, $G_\theta(K(T(\omega)))$ soit continue \`a gauche en $\theta_0$.
Ceci n'est pas une contrainte tr\`es forte, elle est en particulier r\'ealis\'ee lorsque les densit\'es $p_{\theta'_n}$ convergent, dans $L_1$ ou $\mu$ presque s\^urement, vers  $p_{\theta_0}$ quand $\theta'_n$ cro\^{\i}t vers $\theta_0$ (cf. par exemple [Mon.1] p. 138). La plupart des mod\`eles statistiques \`a rapport de vraisem\-blance monotone, classiquement \'etudi\'es, v\'erifient cette continuit\'e en tout point $\theta$ de $\Theta$, aussi bien \`a gauche qu'\`a droite.
\dli Soient $\alpha_0<1$ et $\alpha_1<1$, les hypoth\`eses \'etant adjacentes, la r\`egle des pl\'ebiscites ne d\'epend plus que du vote $Q_{\theta_0}^\omega$.
Elle d\'ecide : 
\dli $d=0$ lorsque $Q_{\theta_0}^\omega(\{0\})>\alpha_0$ et 
$Q_{\theta_0}^\omega(\{1\})\leq\alpha_1$, 
\dli $d=1$ quand $Q_{\theta_0}^\omega(\{0\})\leq\alpha_0$ et 
$Q_{\theta_0}^\omega(\{1\})>\alpha_1$.
\dli Si cette r\`egle ne d\'ecide ni $d=0$ ni $d=1$, on peut dire qu'elle nous conseille de nous abstenir de prendre une d\'ecision. 
\dli Lorsque $\alpha_0=\alpha_1=\alpha<1/2$, la r\`egle des pl\'ebiscites d\'ecide $d=0$, c'est-\`a-dire $\theta\in\Theta_0$, quand $\omega$ appartient \`a la r\'egion de rejet du test de $H_0\,:\,\theta\in\Theta_1$ contre $H_1\,:\,\theta\in\Theta_0$ au seuil $\alpha$ ; elle
d\'ecide $d=1$, c'est-\`a-dire $\theta\in\Theta_1$, quand $\omega$ appartient \`a la r\'egion de rejet du test de $H_0\,:\,\theta\in\Theta_0$ contre $H_1\,:\,\theta\in\Theta_1$ au seuil $\alpha$ ; enfin, cette r\`egle conseille de s'abstenir quand les deux tests pr\'ec\'edents donnent des d\'ecisions contradictoires. Nous avions d\'ej\`a trouv\'e cette r\`egle comme cas particulier des proc\'edures de d\'ecision minimax construites dans une recherche de prise de d\'ecision par rapport \`a une partition ordonn\'ee de $\Theta$ (cf. [Mor.2]).

Nous allons maintenant montrer que les probabilit\'es a post\'eriori de $\Theta_0$ et $\Theta_1$ obtenues \`a partir d'une probabilit\'e a priori $\Lambda$ sur $(\Theta,{\cal T})$ ne sont pas directement comparables avec les r\'esultats du vote pond\'er\'e 
$Q_{\Lambda}^\omega$, qui est la moyenne des votes 
$Q_{\theta}^\omega$ par rapport \`a $\Lambda$ ($Q^\omega_{\theta}(\{1\})$ est suppos\'ee \^etre ${\cal T}$ mesurable presque partout).
$$Q^\omega_{\Lambda}(\{1\})=\int_{\Theta}Q^\omega_{\theta}(\{1\})\,d\Lambda(\theta)$$
\dli Pour d\'efinir la probabilit\'e a post\'eriori de $\Theta_1$, il nous faut supposer que $\Theta_1$ est un \'ev\'enement de ${\cal T}$ et que la famille 
$(P_\theta)_{\theta\in\Theta}$ d\'efinit une probabilit\'e de transition sur 
$\Theta\times{\cal A}$ ($Q^\omega_{\theta}(\{1\})$ est bien ${\cal T}$ mesurable). La donn\'ee de la probabilit\'e a priori $\Lambda$ permet alors de probabiliser l'espace produit $(\Omega,{\cal A})\otimes(\Theta,{\cal T})$. Sous des conditions tr\`es g\'en\'erales, on peut reconstruire cette probabilit\'e \`a partir d'une probabilit\'e de transition sur 
$\Omega\times{\cal T}$ (cf. [DacD] p. 194, [HenT] p. 239). Pour toute r\'ealisation $\omega$, on obtient alors une probabilit\'e conditionnelle, dite a post\'eriori, sur $(\Theta,{\cal T})$.
L'\'ev\'enement $\Theta_1$ \'etant de probabilit\'e :
$$\Lambda(\Theta_1\mid\omega)=\int_{\Theta_1}p_{\theta}(\omega)\,d\Lambda(\theta)\,\biggm/\!\int_{\Theta}p_{\theta}(\omega)\,d\Lambda(\theta)$$
\dli Lorsque $T(\omega)$ appartient \`a $D_i-D_s$ (resp. $D_s$), on a $Q^\omega_{\Lambda}(\{1\})=1$ (resp. $Q^\omega_{\Lambda}(\{1\})=0$)(voir la proposition 4.3.1) et $\Lambda(\Theta_1\mid\omega)$ prend la m\^eme valeur puisque $p_{\theta}(\omega)$ est alors nul pour $\theta$ appartenant \`a $\Theta_0$ (resp. $\Theta_1$).
\dli Le cas important est \'evidemment celui o\`u $T(\omega)$ n'appartient pas \`a $D_i\cup D_s$ . Pour voir que $\Lambda(\Theta_1\mid\omega)$ et $Q^\omega_{\Lambda}(\{1\})$ ne se comportent g\'en\'eralement pas de mani\`ere semblable, on peut \'ecrire $\Lambda$ sous forme du m\'elange de deux probabilit\'es $\Lambda_0$ et $\Lambda_1$ : $\Lambda=\lambda\Lambda_0 + (1-\lambda)\Lambda_1$ avec $\Lambda_0(\Theta_0)=1$, $\Lambda_1(\Theta_1)=1$ et $\lambda=\Lambda(\Theta_0)$. Faisons maintenant varier la probabilit\'e 
$\Lambda$ en changeant uniquement la pond\'eration $\lambda\in[0,1]$. Lorsque $\lambda$ cro\^{\i}t il est facile de montrer que $\Lambda(\Theta_1\mid\omega)$ d\'ecro\^{\i}t alors que $Q^\omega_{\Lambda}(\{1\})$ cro\^{\i}t (voir la proposition 4.3.2). Dans le cadre bay\'esien, plus une probabilit\'e a priori charge $\Theta_1$ plus elle avantage la d\'ecision $\theta\in\Theta_1$, c'est le contraire dans le cadre des votes pond\'er\'es d'experts. Si $\lambda=0$, on a toujours $\Lambda(\Theta_1\mid\omega)=1$, ce qui n'est pas le cas de $Q^\omega_{\Lambda}(\{1\})$, $Q^\omega_{\Lambda}(\{1\})\leq Q^\omega_{\Theta_1}(\{1\})$. Pour le choix entre deux hypoth\`eses simples (voir le paragraphe 2.5) nous avons \'etabli une comparaison mais en utilisant
en fait sur $\Theta$, deux probabilit\'es diff\'erentes, 
$\Lambda=\lambda\Lambda_0 + (1-\lambda)\Lambda_1$ et $\Lambda'=(1-\lambda)\Lambda_0 + \lambda\Lambda_1$. $\lambda$ repr\'esentait 
dans les deux cas un indice de la ``faveur" que l'on accorde \`a l'hypoth\`ese $\Theta_0$. Les conclusions du paragraphe 2.5 restent int\'eressantes mais la mise en oeuvre de la comparaison devient ici plus difficile car il faut encore choisir $\Lambda_0$ et $\Lambda_1$, dans le cas de deux hypoth\`eses simples le seul choix possible \'etait les masses de Dirac en $\theta_0$ et
$\theta_1$. Pour montrer plus en d\'etail cette difficult\'e de comparaison, nous allons \'etudier un exemple en ne privil\'egiant aucune des deux hypoth\`eses, c'est-\`a-dire avec $\lambda=1/2$. La probabilit\'e a priori $\Lambda$ et la pond\'eration $\Lambda'$ sont alors identiques.

Consid\'erons le cas classique d'un n \'echantillon $(X_1,X_2,...,X_n)$ d'une loi normale, $N(\theta,\sigma^2)$, de moyenne inconnue $\theta\in\IR$ et de variance connue $\sigma^2$ ($\sigma>0$). Le mod\`ele statistique est \`a rapport de vraisemblance monotone, pour l'ordre usuel sur $\IR$, par rapport \`a la statistique $\overline{X}=(1/n)\sum_{i=1}^{i=n}X_i$, qui est de loi
$N(\theta,\sigma^2/n)$. $\overline{X}$ est aussi une statistique essentielle pour le choix entre deux hypoth\`eses unilat\'erales. En effet, pour $\theta'<\theta''$, $h_{(\theta'',\theta')}(\overline{x})$ est une fonction strictement croissante de $\overline{x}$.
\dli Etudions le probl\`eme du choix entre $\Theta_0=\{\theta\geq 0\}$ et 
$\Theta_1=\{\theta<0\}$. La partition ordonn\'ee $\{\Theta_1,\Theta_0\}$ d\'efinit des hypoth\`eses stables adjacentes. En effet, la statistique essentielle $\overline{X}$ est diffuse et finie, on a donc 
$Q^\omega_{\theta}(\{1\})=1-F({\sqrt{n}\over\sigma}(\overline{X}(\omega)-\theta))$, F \'etant la fonction de r\'epartition de la loi normale $N(0,1)$, ce qui entra\^{\i}ne bien : 
\dli $Q^\omega_{\Theta_1}(\{1\})=Q^\omega_{\Theta_0}(\{1\})$.
\dli Le seuil minimum de rejet des tests uniform\'ement plus puissants de $\Theta_0$ contre $\Theta_1$ (resp. $\Theta_1$ contre $\Theta_0$) peut s'interpr\'eter comme la valeur du vote $Q^\omega_{\theta_0}(\{0\})$ (resp. $Q^\omega_{\theta_0}(\{1\})$) avec $\theta_0=0$. C'est le r\'esultat de la proposition 5.1.1, mais aussi une cons\'equence directe de l'expression de 
$Q^\omega_{\theta}(\{1\})$. Les seuils minimums de rejet des deux tests s'obtiennent donc en prenant comme pond\'eration la masse de Dirac au point $\theta_0=0$. Dans la th\'eorie bay\'esienne ces deux seuils seront les probabilit\'es a post\'eriori de $\Theta_0$ et $\Theta_1$ en prenant comme loi a priori sur $\Theta=\IR$ une loi impropre, la mesure de Lebesgue (cf. [Ber.] p. 147).
Ceci nous montre bien que probabilit\'e a priori et probabilit\'e de pond\'eration sur $\Theta$ ne peuvent pas s'interpr\'eter de la m\^eme fa\c con.
\dli Regardons tout de m\^eme ce que  donnent ces deux m\'ethodes dans le cas classique o\`u la probabilit\'e $\Lambda$ est la loi normale $N(0,c^2)$ avec $c>0$. On traite ainsi de fa\c con tr\`es sym\'etrique les deux hypoth\`eses 
$\Theta_0$ et $\Theta_1$. Dans le cadre bay\'esien on sait que la loi a post\'eriori sur $\Theta$ est la loi normale $N(m(\omega),v(\omega))$ avec
$m(\omega)={nc^2\over nc^2+\sigma^2}\overline{X}(\omega)$ et $v(\omega)={c^2\sigma^2\over nc^2+\sigma^2}$ (cf. [Ber.] p. 128 ou [Rob.] p. 139).
$\Lambda(\Theta_1\mid\omega)$ est donc la valeur de la fonction de r\'epartition de la loi $N(m(\omega),v(\omega))$ en $\theta_0=0$ : 
\dli $\Lambda(\Theta_1\mid\omega)=1-F({\sqrt{n}\over\sigma}\overline{X}(\omega) \sqrt{nc^2\over nc^2+\sigma^2})$.
\dli Le vote pond\'er\'e correspondant s'\'ecrit :
\cleartabs
\+ $Q^\omega_{\Lambda}(\{1\})$&=&$1-\int_\Theta F({\sqrt{n}\over\sigma}(\overline{X}(\omega)-\theta))\,{1\over\sqrt{2\pi}\,c} exp(-{\theta^2\over 2c^2})\,d\theta$\cr
\+ &=&$1-\int_\Theta\int_{\IR} \II_{]-\infty,\overline{X}(\omega)[}(y) {1\over\sqrt{2\pi}}{\sqrt{n}\over\sigma} exp(-{n\over 2\sigma^2}(y-\theta)^2)\,{1\over\sqrt{2\pi}\,c} exp(-{\theta^2\over 2c^2})\,dy\,d\theta$&\cr
\+ &=&$1-\int_{\IR} \II_{]-\infty,\overline{X}(\omega)[}(y) {1\over\sqrt{2\pi}}{\sqrt{n}\over\sigma} exp(-{ny^2\over 2\sigma^2})\int_\Theta\,{1\over\sqrt{2\pi}\,c} exp({n^2c^2y^2\over 2\sigma^2(nc^2+\sigma^2)})\times$\cr
\+ & &\qquad\hfill $exp(-{nc^2+\sigma^2\over 2c^2\sigma^2}(\theta-{nc^2y\over nc^2+\sigma^2})^2)\,d\theta\,dy$&\cr
\+ &=&$1-\int_{\IR} \II_{]-\infty,\overline{X}(\omega)[}(y) {1\over\sqrt{2\pi}}\sqrt{{n\over nc^2+\sigma^2}} exp(-{ny^2\over 2(nc^2+\sigma^2)})\,dy$\cr
\+ &=&$1-F({\sqrt{n}\over\sigma}\overline{X}(\omega) \sqrt{\sigma^2\over nc^2+\sigma^2})$.\cr
\dli $\Lambda(\Theta_1\mid\omega)$ et $Q^\omega_{\Lambda}(\{1\})$ sont des corrections diff\'erentes apport\'ees au seuil minimum de rejet de $\Theta_1$ contre $\Theta_0$ : $Q^\omega_{\theta_0}(\{1\})= 1-F({\sqrt{n}\over\sigma}\overline{X}(\omega))$. Ces corrections sont n\'egatives (resp. positives) quand $\overline{X}(\omega)$ est n\'egatif (resp. positif) donc quand le seuil minimum est inf\'erieur (resp. sup\'erieur) \`a ${1\over 2}$. Elles accentuent la tendance.
\dli On obtient 
$\Lambda(\Theta_1\mid\omega)=Q^\omega_{\Lambda'}(\{1\})$ en prenant comme pond\'eration $\Lambda'$ la loi $N(0,c'^2)$ avec $c'=\sigma^2/nc$. 

\vfill\eject

{\parindent=-5mm 5.2 VOTES COMPATIBLES SUR UNE FAMILLE D'HYPOTH\`ESES UNILA\-T\'E\-RALES.}
\nobreak
\medskip

Dans un mod\`ele \`a rapport de vraisemblance monotone l'interpr\'etation de la valeur du param\`etre $\theta$ est souvent structur\'ee  par l'ordre d\'efini sur $\Theta$. La prise de d\'ecision peut se voir comme une synth\`ese des r\'esultats d'expertise de plusieurs hypoth\`eses unilat\'erales diff\'erentes : 
$(\{\Theta_1^f,\Theta_0^f\})_{f\in{\cal F}}$. Associons \`a chacun de ces probl\`emes de d\'ecision index\'es par $f\in{\cal F}$, un vote d'experts. Pour toute r\'ealisation $\omega\in\Omega$ on obtient pour chaque $f$ une probabilit\'e sur $D=\{0,1\}$, nous notons $Q^\omega(\Theta_1^f)$ (resp. $Q^\omega(\Theta_0^f)$) la probabilit\'e de la d\'ecision $d=1$ (resp. $d=0$).
$Q^\omega$ est une application de $\{\Theta_i^f\}_{f\in{\cal F},i\in\{0,1\}}\subset{\cal P}(\Theta)$ dans $[0,1]$. Les informations donn\'ees par $Q^\omega$ seront exploitables si les votes choisis sont coh\'erents. Par exemple si $\Theta_1^f\subset\Theta_1^{f'}$, on ne voudrait pas avoir $Q^\omega(\Theta_1^f)>Q^\omega(\Theta_1^{f'})$. De m\^eme si la suite monotone
$(\Theta_1^{f_n})_{n\in\IN}$ a pour limite $\Theta_1^f$, on aimerait bien que $Q^\omega(\Theta_1^{f_n})$ converge vers $Q^\omega(\Theta_1^f)$ pour tout $\omega$.

\medskip
{\bf D\'efinition 5.2.1}
\medskip
\medskip
\moveleft 10.4pt\hbox{\vrule\kern 10pt\vbox{\defpro

Soient $(\{\Theta_1^f,\Theta_0^f\})_{f\in{\cal F}}$ une famille d'hypoth\`eses unilat\'erales dans un mod\`ele \`a rapport de vraisemblance monotone de param\`etre $\Theta$ et $Q$ une application de $\Omega\times\{\Theta_1^f\}_{f\in{\cal F}}$ dans $[0,1]$ telle que $Q^\omega(\Theta_1^{f})$ repr\'esente la valeur en $d=1$ d'un vote d'experts du choix entre $\Theta_1^{f}$ et $\Theta_0^{f}$ lorsqu'on r\'ealise $\omega$.
\dli Nous dirons que $Q$ d\'efinit des votes compatibles lorsqu'elle v\'erifie, pour presque tout $\omega$ de $\Omega$, les trois propri\'et\'es suivantes :
\dli a) si $\Theta_1^{f}\subset\Theta_1^{f'}$ alors
$Q^\omega(\Theta_1^{f})\leq Q^\omega(\Theta_1^{f'})$
\dli b) si la suite $(\Theta_1^{f_n})_{n\in\IN}$ cro\^{\i}t vers $\Theta$ (resp. d\'ecro\^{\i}t vers $\emptyset$), la suite $(Q^\omega(\Theta_1^{f_n}))_{n\in\IN}$ converge vers $1$ (resp. $0$)
\dli c) si $(\Theta_1^{f_n})_{n\in\IN}$ et $(\Theta_1^{f'_n})_{n\in\IN}$ sont deux suites ayant la m\^eme limite, l'une \'etant croissante l'autre d\'ecroissante, les suites $(Q^\omega(\Theta_1^{f_n}))_{n\in\IN}$ et $(Q^\omega(\Theta_1^{f'_n}))_{n\in\IN}$ ont aussi m\^eme limite.
}}\medskip

\medskip
La propri\'et\'e c) implique bien que si la suite monotone $(\Theta_1^{f_n})_{n\in\IN}$ converge  vers $\Theta_1^{f}$ alors $(Q^\omega(\Theta_1^{f_n}))_{n\in\IN}$ converge vers $Q^\omega(\Theta_1^{f})$, il suffit de prendre $\Theta_1^{f'_n}=\Theta_1^{f}$ pour tout $n\in\IN$.
Cette propri\'et\'e est un peu plus forte, afin que $Q^\omega$ puisse se prolonger \`a la tribu engendr\'ee par $\{\Theta_1^f\}_{f\in{\cal F}}$.

\medskip
{\bf Proposition 5.2.1}
\medskip
\medskip
\moveleft 10.4pt\hbox{\vrule\kern 10pt\vbox{\defpro

Si $Q$ d\'efinit des votes compatibles sur la famille d'hypoth\`eses unilat\'erales $(\{\Theta_1^f,\Theta_0^f\})_{f\in{\cal F}}$, pour presque toute r\'ealisation $\omega$ de $\Omega$, $Q^\omega$ se prolonge en une probabilit\'e unique sur la $\sigma$-alg\`ebre engendr\'ee par $\{\Theta_1^f\}_{f\in{\cal F}}$.
}}\medskip

{\leftskip=15mm \dli {\bf D\'emonstration}
\medskip
Notons $N$ le n\'egligeable form\'e des r\'ealisations $\omega$ pour lesquelles une au moins des trois propri\'et\'es de la d\'efinition 5.2.1 n'est pas v\'erifi\'ee.
Nous allons prolonger $Q^\omega$ pour tout $\omega$ de $\Omega-N$.

Consid\'erons la semi-alg\`ebre de Boole, ${\cal S}$, engendr\'ee par : \dli $\{\Theta_1^f\}_{f\in{\cal F}}$. Elle est constitu\'ee du $\emptyset$, de $\Omega$, des demi-droites inf\'erieures $\{\Theta_1^f\}_{f\in{\cal F}}$, des demi-droites sup\'erieures $\{\Theta_0^f\}_{f\in{\cal F}}$ et des intervalles de la forme $\Theta_1^f\cap\Theta_0^{f'}$.
\dli Soit $\omega\in\Omega-N$, la propri\'et\'e a) de la d\'efinition 5.2.1 nous permet de prolonger $Q^\omega$ de fa\c con unique en une fonction additive d'ensembles, de ${\cal S}$ dans $[0,1]$, telle que $Q^\omega(\Omega)=1$.
On a $Q^\omega(\Theta_1^f\cap\Theta_0^{f'})=Q^\omega(\Theta_1^f)-Q^\omega(\Theta_1^{f'})$ si $\Theta_1^f\cap\Theta_0^{f'}\not=\emptyset$.
$Q^\omega(\Theta_0^f)=1-Q^\omega(\Theta_1^f)$ repr\'esente bien la valeur en $d=0$ du vote d'experts d\'efini par $Q$ pour le choix entre $\Theta_1^f$ et $\Theta_0^f$ quand on r\'ealise $\omega$.

Nous aurons d\'emontr\'e la proposition, si nous montrons que $Q^\omega$ est $\sigma$-additive sur ${\cal S}$ (cf. [Nev.] p. 25). Pour cela nous allons montrer que $Q^\omega$ jouit de la propri\'et\'e de continuit\'e monotone s\'equentielle en $\emptyset$.
\dli Soit $(S_n)_{n\in\IN}$ une suite de ${\cal S}$ qui d\'ecro\^{\i}t vers le vide, nous devons d\'emontrer que $(Q^\omega(S_n))_{n\in\IN}$ d\'ecro\^{\i}t vers $0$. C'est \'evident lorsque la suite contient le vide. Dans le cas contraire
($\forall n\, S_n\not=\emptyset$), les suites de ${\cal S}$ qui d\'ecroissent vers le vide sont de trois types : 
\cleartabs
\settabs\+\kern 3cm&\cr
\+ \qquad\hfill i)&\quad$\forall n\quad S_n=\Theta_1^{f_n}$\cr
\+ \qquad\hfill ii)&\quad$\forall n\quad S_n=\Theta_0^{f_n}$\cr
\+ \qquad\hfill iii)&\quad$\forall n>n_0\quad S_n=\Theta_1^{f_n}\cap\Theta_0^{f'_n}\not=\emptyset$\cr

D'apr\`es la propri\'et\'e b) de la d\'efinition 5.2.1 on a \'evidemment
$lim_{n\rightarrow +\infty}Q^\omega(S_n)=0$ dans le cas i), mais aussi dans le cas ii) car alors $Q^\omega(S_n)=1-Q^\omega(\Theta_1^{f_n})$ avec $\Theta_1^{f_n}$ qui cro\^{\i}t vers $\Omega$.
\dli Dans le cas iii), on a $S_n=\Theta_1^{f_n}-\Theta_1^{f'_n}$.
$(\Theta_1^{f_n})_{n\in\IN}$ et $(\Theta_1^{f'_n})_{n\in\IN}$ sont des suites respectivement d\'ecroissante et croissante. Elles ont m\^eme limite car
$\emptyset=\cap_{n\in\IN}S_n=\cap_{n\in\IN}\Theta_1^{f_n}-\cup_{n\in\IN}\Theta_1^{f'_n}$. La propri\'et\'e c) de la d\'efinition 5.2.1 nous dit que les images de ces deux suites par $Q^\omega$ ont aussi m\^eme limite. On a donc bien :
\dli $lim_{n\rightarrow +\infty}Q^\omega(S_n)=lim_{n\rightarrow +\infty}[Q^\omega(\Theta_1^{f_n})-Q^\omega(\Theta_1^{f'_n})]=0$.

\medskip\centerline{\hbox to 3cm{\bf \hrulefill}}\par}

Consid\'erons un probl\`eme de d\'ecision d\'efini par $\{\Theta_0,\Theta_1\}$, 
$\Theta=\Theta_0+\Theta_1$. Si $\Theta_1$ est un \'el\'ement de la tribu engendr\'ee par $(\Theta_1^{f})_{f\in{\cal F}}$ et si ce probl\`eme de d\'ecision ne poss\`ede pas d'expert, il est tentant de prendre comme vote, lorsqu'on r\'ealise $\omega$, les valeurs du prolongement de $Q^\omega$ en 
$\Theta_0$ et $\Theta_1$. On est souvent dans cette situation lorsque $\Theta_0$ et $\Theta_1$ d\'efinissent des hypoth\`eses bilat\'erales. Nous analyserons en d\'etail ce cas classique au paragraphe suivant.
Bien entendu ce proc\'ed\'e de construction d'un vote \`a partir des votes d'une famille d'hypoth\`eses peut se concevoir dans d'autres mod\`eles que les mod\`eles \`a rapport de vraisemblance monotone. Ces mod\`eles s'y pr\^etent particuli\`erement car lorsque $\Theta\subseteq\IR$, les hypoth\`eses unilat\'erales permettent de d\'efinir des votes qui se prolongent \`a l'ensemble des bor\'eliens de $\Theta$. Pour obtenir ce r\'esultat nous avons besoin de relier entre elles les statistiques essentielles d\'efinies \`a partir d'hypoth\`eses unilat\'erales diff\'erentes. Nous allons construire une statistique $K(T)$ dont la fonction de r\'epartition moyenne sera \'egale \`a celle de toute statistique essentielle associ\'ee \`a des hypoth\`eses unilat\'erales $\{\Theta_1,\Theta_0\}$, en dehors de l'image des demi-droites $D_i$ et $D_s$ d\'efinies par $\Theta_0$ et $\Theta_1$.
\vfill\eject
\medskip
{\bf Proposition 5.2.2}
\medskip
\medskip
\moveleft 10.4pt\hbox{\vrule\kern 10pt\vbox{\defpro

Soit $(\Omega ,{\cal A},(p_\theta.\mu)_{\theta\in\Theta})$ un mod\`ele
statistique \`a rapport de vraisemblance monotone pour la statistique r\'eelle $T$, $\Theta$ \'etant muni de l'ordre total : $\preceq$.
\dli A tout probl\`eme de d\'ecision unilat\'eral, $\{\Theta_1,\Theta_0\}$, associons une statistique essentielle $K^{\Theta_1}(T)$ dont la fonction de r\'epartition moyenne sous $P_\theta$ est not\'ee $G_\theta^{\Theta_1}$.
\dli Posons $N=\{t\in\IR\, ;\, \forall\omega\in T^{-1}(t)\quad\forall\theta\in\Theta\quad p_\theta(\omega)=0\}$.

Il existe une fonction croissante $K : \IR\rightarrow\overline{\IR}$, qui d\'efinit une statistique $K(T)$ dont les fonctions de r\'epartition moyenne $G_\theta$ v\'erifient, pour toutes hypoth\`eses unilat\'erales $\{\Theta_1,\Theta_0\}$ : 
\dli 1) $\forall t\in\ \IR-(D_i^{\Theta_0}\cup D_s^{\Theta_1}\cup N)$,\quad 
$\forall \theta\in\Theta$,\quad $G_\theta(K(t))=G_\theta^{\Theta_1}(K^{\Theta_1}(t))$
\dli 2) $\forall t\in(D_i^{\Theta_0}-N)$,\quad $\forall \theta_0\in\Theta_0$,\quad $G_{\theta_0}(K(t))=0$
\dli 3) $\forall t\in(D_s^{\Theta_1}-N)$,\quad $\forall \theta_1\in\Theta_1$,\quad $G_{\theta_1}(K(t))=1$.

$K(T)$ est appel\'ee statistique essentielle globale, pour une r\'ealisation $\omega$ n'appartenant pas \`a $T^{-1}(D_i^{\Theta_0}\cup D_s^{\Theta_1}\cup N)$, $1-G_\theta(K(T(\omega)))$ est le r\'esultat 
$Q^\omega_\theta(\{1\})$ du vote en faveur de $\Theta_1$ sous $P_\theta$.
}}\medskip

\medskip
{\leftskip=15mm \dli {\bf D\'emonstration}
\medskip
Le mod\`ele statistique \'etant \`a rapport de vraisemblance mono\-tone il existe, pour $\theta'\!\prec\theta''$, une fonction croissante $h_{(\theta'',\theta')} :\, \IR \rightarrow \overline{\IR^+}$ v\'erifiant
$p_{\theta''}/p_{\theta'}=h_{(\theta'',\theta')}(T)$
sur le domaine de d\'efinition de ce rapport : 
$\{\omega\in\Omega\, ;\, p_{\theta'}(\omega)\!>\!0\  ou\ p_{\theta''}(\omega)\!>\!0\}$.
Nous allons travailler avec une famille $\{h_{(\theta'',\theta')}\}_{\theta'\prec\theta''}$ de fonctions normalis\'ees ; $h_{(\theta'',\theta')}$ est alors constante sur tout intervalle $I$ ind\'etermin\'e pour $\{\theta',\theta''\}$ : 
$\forall\omega\in T^{-1}(I)\quad p_{\theta'}(\omega)=p_{\theta''}(\omega)=0$
(ceci est toujours possible d'apr\`es l'annexe IV).
\medskip
{\parindent=-10mm I --- D\'efinition de $K : \IR\rightarrow\overline{\IR}$.}

Pour tout $t\in\IR$ et tout couple $(\theta'',\theta')$ d'\'el\'ements de $\Theta$ tels que $\theta'\!\prec\theta''$, consid\'erons les deux fonctions de test de $F=[f_{(-\infty,1)},f_{(+\infty,0)}]$ (voir le paragraphe 4.2) d\'efinies par :
\dli $e_t^{(\theta'',\theta')}=sup\{f\in\Phi_s^{(\theta'',\theta')}\, ;\, f\leq f_{(t,0)}\}\cup\{f_{(-\infty,1)}\}$ 
\dli $g_t^{(\theta'',\theta')}=inf\{f\in\Phi_s^{(\theta'',\theta')}\, ;\, f\geq f_{(t,1)}\}\cup\{f_{(+\infty,0)}\}$ 
\dli $\Phi_s^{(\theta'',\theta')}\subseteq F$ d\'esigne l'ensemble des fonctions de test simples construites \`a partir du rapport des densit\'es $h_{(\theta'',\theta')}(T)$ (voir la d\'efinition 2.2.1).
\dli Par d\'efinition : $e_t^{(\theta'',\theta')}\leq f_{(t,0)}<f_{(t,1)}\leq g_t^{(\theta'',\theta')}$, il est facile de v\'erifier l'\'egalit\'e :  $\{g_t^{(\theta'',\theta')}-e_t^{(\theta'',\theta')}=1\}=\{h_{(\theta'',\theta')}(T)=h_{(\theta'',\theta')}(t)\}$. Posons $h_{(\theta'',\theta')}(t)=k$, on a alors $e_t^{(\theta'',\theta')}=\phi_{(k,0)}^{(\theta'',\theta')}$ lorsque $k>0$ et $g_t^{(\theta'',\theta')}=\phi_{(k,1)}^{(\theta'',\theta')}$ lorsque $k<+\infty$.

La d\'efinition de $K$ repose sur les deux \'el\'ements de $F$ : 
\dli $f_{(a_t,u_t)}=sup\{e_t^{(\theta'',\theta')}\}_{\theta'\!\prec\theta''}$ et $f_{(b_t,v_t)}=inf\{g_t^{(\theta'',\theta')}\}_{\theta'\!\prec\theta''}$.
$$K(t)=\left\{\matrix{
\hfill b_t-1\hfill & si & a_t=-\infty\hfill\cr
\hfill [a_t+b_t]/2\hfill & si & -\infty<a_t\leq b_t<+\infty\hfill \cr
\hfill a_t+1\hfill & si & b_t=+\infty\  et\  a_t>-\infty\hfill\cr
}\right. $$
\dli Ceci d\'efinit bien une application de $\IR$ dans $\overline{\IR}$. V\'erifions qu'elle est croissante. Pour $t'>t$ et $\theta'\!\prec\theta''$ on a $h_{(\theta'',\theta')}(t)\leq h_{(\theta'',\theta')}(t')$ ; nous avons vu que l'\'egalit\'e entra\^{\i}ne :  $e_t^{(\theta'',\theta')}=e_{t'}^{(\theta'',\theta')}$ et 
$g_t^{(\theta'',\theta')}=g_{t'}^{(\theta'',\theta')}$ ; dans le cas contraire on a \'evidemment $g_t^{(\theta'',\theta')}\leq e_{t'}^{(\theta'',\theta')}$.
Soit $t'>t$, si pour tout couple $(\theta'',\theta')$,  $\theta'\!\prec\theta''$, on a $h_{(\theta'',\theta')}(t)=h_{(\theta'',\theta')}(t')$ alors $f_{(a_t,u_t)}=f_{(a_{t'},u_{t'})}$ et $f_{(b_t,v_t)}=f_{(b_{t'},v_{t'})}$, donc $K(t)=K(t')$ ; dans le cas contraire il existe au moins un couple $(\theta'',\theta')$ tel que $h_{(\theta'',\theta')}(t)<h_{(\theta'',\theta')}(t')$, ce qui implique
$g_t^{(\theta'',\theta')}\leq e_{t'}^{(\theta'',\theta')}$ donc $f_{(b_t,v_t)}\leq f_{(a_{t'},u_{t'})}$, on a alors par d\'efinition de $K$ : $K(t)<K(t')$.

En montrant la croissance de $K$, nous avons aussi d\'emontr\'e que $K(t)=K(t')$ si et seulement si : $f_{(a_t,u_t)}=f_{(a_{t'},u_{t'})}$ et $f_{(b_t,v_t)}=f_{(b_{t'},v_{t'})}$. On a donc $\{K(T)<K(t)\}=\{f_{(a_t,u_t)}=1\}$ et $\{K(T)\leq K(t)\}=\{f_{(b_t,v_t)}=1\}$.
La fonction de r\'epartition moyenne de $K(T)$ sous $P_\theta$ v\'erifie 
$G_\theta(K(t))=[E_\theta(f_{(a_t,u_t)})+E_\theta(f_{(b_t,v_t)})]/2$.

\medskip
{\parindent=-10mm II --- $G_\theta(K(t))=G_\theta^{\Theta_1}(K^{\Theta_1}(t))$ lorsque $t\notin D_i^{\Theta_0}\cup D_s^{\Theta_1}\cup N$.}

$\{\Theta_1,\Theta_0\}$ sont des hypoth\`eses unilat\'erales, d'apr\`es l'annexe IV on d\'efinit une famille normalis\'ee $\{h^{\Theta_1}_{(\theta_0,\theta_1)}\}_{(\theta_0,\theta_1)\in\Theta_0\times\Theta_1}$ pour ces hypoth\`eses, en posant :
$$h^{\Theta_1}_{(\theta_0,\theta_1)}(t)\quad=\quad\left\{\matrix{
0\hfill & si & t\in D^{\Theta_0}_i-D^{\Theta_1}_s\hfill\cr
h_{(\theta_0,\theta_1)}(t)\hfill & si & t\in\ \IR-(D^{\Theta_0}_i\cup D^{\Theta_1}_s)\hfill \cr
+\infty\hfill & si & t\in D^{\Theta_1}_s\hfill \cr
}\right. $$
Nous savons que le r\'esultat recherch\'e ne d\'epend pas de la famille normalis\'ee choisie (voir la proposition 1 de l'annexe III).
\dli La statistique essentielle $K^{\Theta_1}(T)$ associ\'ee \`a cette famille est d\'efinie \`a partir des experts $\Delta_s^{\Theta_1}=\cup_{(\theta_0,\theta_1)\in\Theta_0\times\Theta_1} \Phi_{s^{\Theta_1}}^{(\theta_0,\theta_1)}$,  $\Phi_{s^{\Theta_1}}^{(\theta_0,\theta_1)}$ \'etant  l'ensemble des fonctions de test simples construites \`a partir du rapport des densit\'es $h_{(\theta_0,\theta_1)}^{\Theta_1}(T)$ (voir la d\'efinition 2.2.1). $\Phi_s^{(\theta_0,\theta_1)}$ d\'efini pour la partie I de cette d\'emonstration, l'\'etait \`a partir de $h_{(\theta_0,\theta_1)}(T)$.
Par d\'efinition de $h_{(\theta_0,\theta_1)}^{\Theta_1}$ on a 
$\Phi_{s^{\Theta_1}}^{(\theta_0,\theta_1)}\supseteq\Phi_{s}^{(\theta_0,\theta_1)}\cap ]\II_{D^{\Theta_0}_i-D^{\Theta_1}_s}\!(T),1-\II_{D^{\Theta_1}_s}\!(T)[$. $\Phi_{s^{\Theta_1}}^{(\theta_0,\theta_1)}$ contient en plus $\II_{D^{\Theta_0}_i-D^{\Theta_1}_s}\!(T)$ (resp. $1-\II_{D^{\Theta_1}_s}\!(T)$) si pour tout $t\notin D^{\Theta_0}_i-D^{\Theta_1}_s$ (resp. $t\notin D^{\Theta_1}_s$) on a $h_{(\theta_0,\theta_1)}(t)>0$ (resp. $h_{(\theta_0,\theta_1)}(t)<+\infty$).

Soit $t\in\,\IR-(D_i^{\Theta_0}\cup D_s^{\Theta_1})\not=\emptyset$ .
\dli Nous avons vu, au d\'ebut de la d\'emonstration de la proposition 1 de l'annexe III, que les fonctions de test, de la d\'efinition 4.3.1, permettant de construire $K^{\Theta_1}$ sont alors donn\'ees par : 
\dli $f_{(a'_t,u'_t)}^{\Theta_1}=sup\{f\in\Delta_s^{\Theta_1}\, ;\, f\leq f_{(t,0)}\}$ et 
\dli $f_{(b'_t,v'_t)}^{\Theta_1}=inf\{f\in\Delta_s^{\Theta_1}\, ;\, f\geq f_{(t,1)}\}$.
\dli Elles v\'erifient : $inf\Delta_s^{\Theta_1}=\II_{D^{\Theta_0}_i}(T)\leq f_{(a'_t,u'_t)}^{\Theta_1}\leq f_{(t,0)}<f_{(t,1)}\leq f_{(b'_t,v'_t)}^{\Theta_1}\leq 1-\II_{D^{\Theta_1}_s}(T)=sup\Delta_s^{\Theta_1}$.
\dli On peut aussi les d\'efinir par : 
\dli $f_{(a'_t,u'_t)}^{\Theta_1}=sup\{e_t^{(\theta_0,\theta_1)}\}_{(\theta_0,\theta_1)\in\Theta_0\times\Theta_1}\cup\{\II_{D^{\Theta_0}_i}\!(T)\}$,
\dli $f_{(b'_t,v'_t)}^{\Theta_1}=inf\{g_t^{(\theta_0,\theta_1)}\}_{(\theta_0,\theta_1)\in\Theta_0\times\Theta_1}\cup\{1-\II_{D^{\Theta_1}_s}\!(T)\}$.
\dli La premi\`ere (resp. deuxi\`eme) \'egalit\'e est \'evidente dans le cas o\`u $f_{(a'_t,u'_t)}^{\Theta_1}>\II_{D^{\Theta_0}_i}\!(T)$ (resp. $f_{(b'_t,v'_t)}^{\Theta_1}<1-\II_{D^{\Theta_1}_s}\!(T)$), en effet nous venons de montrer que $\Phi_{s^{\Theta_1}}^{(\theta_0,\theta_1)}$ et $\Phi_{s}^{(\theta_0,\theta_1)}$ ont la m\^eme intersection avec \dli $]\II_{D^{\Theta_0}_i}\!(T),1-\II_{D^{\Theta_1}_s}\!(T)[$.
\dli Lorsque $f_{(a'_t,u'_t)}^{\Theta_1}=\II_{D^{\Theta_0}_i}\!(T)$ (resp. $f_{(b'_t,v'_t)}^{\Theta_1}=1-\II_{D^{\Theta_1}_s}\!(T)$), on a $]\II_{D^{\Theta_0}_i}\!(T),f_{(t,0)}]$ (resp. $[f_{(t,1)},1-\II_{D^{\Theta_1}_s}\!(T)[$) d'intersection vide avec $\Delta_s^{\Theta_1}$, il en est de m\^eme avec tous les $\Phi_{s}^{(\theta_0,\theta_1)}$ (toujours d'apr\`es la propri\'et\'e pr\'ec\'edente) ; on a alors, pour tout $(\theta_0,\theta_1)\in\Theta_0\times\Theta_1$, $e_t^{(\theta_0,\theta_1)}\leq\II_{D^{\Theta_0}_i}\!(T)$ (resp. 
$g_t^{(\theta_0,\theta_1)}\geq 1-\II_{D^{\Theta_1}_s}\!(T)$), ce qui implique bien l'\'egalit\'e recherch\'ee.

\dli Comme $t$ n'appartient pas \`a $D_i^{\Theta_0}\cup D_s^{\Theta_1}$, $K^{\Theta_1}(t)$ est d\'efini par : 
$$K^{\Theta_1}(t)=\left\{\matrix{
\hfill b'_t-1\hfill & si & a'_t=-\infty\hfill\cr
\hfill [a'_t+b'_t]/2\hfill & si & -\infty<a'_t\leq b'_t<+\infty\hfill \cr
\hfill a'_t+1\hfill & si & b'_t=+\infty\ et\  a'_t>-\infty\hfill\cr
}\right. $$
On a bien s\^ur $G^{\Theta_1}_\theta(K^{\Theta_1}(t))=[E_\theta(f^{\Theta_1}_{(a'_t,u'_t)})+E_\theta(f^{\Theta_1}_{(b'_t,v'_t)})]/2$.
\dli Nous devons d\'emontrer que cette quantit\'e est \'egale \`a 
$G_\theta(K(t))=[E_\theta(f_{(a_t,u_t)})+E_\theta(f_{(b_t,v_t)})]/2$ lorsque
$t$ appartient \`a $\IR-(D_i^{\Theta_0}\cup D_s^{\Theta_1}\cup N)$.
\dli Nous le ferons en d\'emontrant : \dli $f_{(a_t,u_t)}\relmont{=}{p.s.}f^{\Theta_1}_{(a'_t,u'_t)}$ et 
$f_{(b_t,v_t)}\relmont{=}{p.s.}f^{\Theta_1}_{(b'_t,v'_t)}$ pour $t\notin D_i^{\Theta_0}\cup D_s^{\Theta_1}\cup N$.

\medskip
{\parindent=-5mm a) Pour $t\notin D_i^{\Theta_0}$ on a $f_{(a_t,u_t)}\relmont{\geq}{p.s.}f^{\Theta_1}_{(a'_t,u'_t)}$.}

D'apr\`es la d\'efinition de $f_{(a'_t,u'_t)}^{\Theta_1}$ \`a partir des $e_t^{(\theta_0,\theta_1)}$, nous avons : 
$f_{(a'_t,u'_t)}^{\Theta_1}\leq sup\{f_{(a_t,u_t)},\II_{D^{\Theta_0}_i}\!(T)\}$. La propri\'et\'e recherch\'ee revient \`a d\'emontrer : $f_{(a_t,u_t)}\relmont{\geq}{p.s.}\II_{D^{\Theta_0}_i}\!(T)$.
Nous allons le faire en montrant qu'il existe $\theta_0\in\Theta_0$ tel que 
$sup\{e_t^{(\theta_0,\theta)}\}_{\theta\in\Theta_1}\relmont{\geq}{p.s.}\II_{D^{\Theta_0}_i}\!(T)$.

Comme $t\notin D_i^{\Theta_0}$, il existe $\theta_0\in\Theta_0$, $x_0\leq t$ et $\omega_0\in T^{-1}(x_0)$ tels que $p_{\theta_0}(\omega_0)>0$. Pour tout $\theta$ de $\Theta_1$ on a $h_{(\theta_0,\theta)}(x_0)>0$, donc 
$h_{(\theta_0,\theta)}(t)=k_\theta>0$ et $e_t^{(\theta_0,\theta)}=\{h_{(\theta_0,\theta)}(T)<k_\theta\}$. 
\dli Pour que $\theta_0$ convienne il nous faut montrer que le cas de figure \dli  $sup\{e_t^{(\theta_0,\theta)}\}_{\theta\in\Theta_1}<\II_{D^{\Theta_0}_i}\!(T)$ implique 
$sup\{e_t^{(\theta_0,\theta)}\}_{\theta\in\Theta_1}\relmont{=}{p.s.}\II_{D^{\Theta_0}_i}\!(T)$. 
\dli Il est \'equivalent de supposer $I=\cap_{\theta\in\Theta_1}(D^{\Theta_0}_i-\{h_{(\theta_0,\theta)}(t)<k_\theta\})$ non vide et de montrer que l'on a 
$P_\theta(T^{-1}(I))=0$ pour tout $\theta$ de $\Theta$. C'est \'evident pour $\theta\in\Theta_0$ puisque la densit\'e $p_\theta$ est alors nulle sur  
$T^{-1}(D^{\Theta_0}_i)\supseteq T^{-1}(I)$ (voir la d\'efinition 4.2.1). Lorsque $\theta\in\Theta_1$ la densit\'e $p_\theta$ est encore nulle sur  
$T^{-1}(I)$ ; en effet pour tout $\omega\in T^{-1}(I)$ on a $p_{\theta_0}(\omega)=0$ et $h_{(\theta_0,\theta)}(T(\omega))\geq k_\theta>0$, ceci n'est possible qu'avec $p_{\theta}(\omega)=0$.
\medskip
{\parindent=-5mm b) Pour $t\notin D_i^{\Theta_0}\cup D_s^{\Theta_1}\cup N$ on a $f_{(a_t,u_t)}\relmont{\leq}{p.s.}f^{\Theta_1}_{(a'_t,u'_t)}$.}

Cette in\'egalit\'e s'\'ecrit : 
\dli $sup\{e_t^{(\theta'',\theta')}\}_{\theta'\!\prec\theta''}\relmont{\leq}{p.s.}sup\{e_t^{(\theta_0,\theta_1)}\}_{(\theta_0,\theta_1)\in\Theta_0\times\Theta_1}\cup\{\II_{D^{\Theta_0}_i}\!(T)\}$.
 Nous allons d\'emontrer ceci en consid\'erant les couples $(\theta'',\theta')$ compos\'es de deux \'el\'ements de $\Theta_0$ (resp. $\Theta_1$), $\theta'\!\prec\theta''$, et en montrant qu'il existe $\theta_a\in\Theta_1$ (resp. $\theta_c\in\Theta_0$) tel que :
$e_t^{(\theta'',\theta')}\relmont{\leq}{p.s.}sup\{e_t^{(\theta',\theta_a)},e_t^{(\theta'',\theta_a)},\II_{D^{\Theta_0}_i}\!(T)\}$ (resp. 
$e_t^{(\theta'',\theta')}\relmont{\leq}{p.s.}sup\{e_t^{(\theta_c,\theta')},e_t^{(\theta_c,\theta'')},\II_{D^{\Theta_0}_i}\!(T)\}$). 

1\up{er} cas : $\theta'\in\Theta_0$, $\theta''\in\Theta_0$ et $\theta'\!\prec\theta''$.

Posons $h_{(\theta'',\theta')}(t)=k_3$. On a $k_3<+\infty$ ; en effet si on avait $h_{(\theta'',\theta')}(t)=+\infty$, la densit\'e $p_{\theta'}$ serait nulle sur $T^{-1}([t,+\infty[)$ et $t$ appartiendrait \`a $D^{\theta'}_s$ donc \`a $D^{\Theta_1}_s$ puisque $\theta'\!\succ\Theta_1$ (voir le lemme 1 de l'annexe IV) ; ceci est impossible car $t\in\,\IR-(D_i^{\Theta_0}\cup D_s^{\Theta_1}\cup N)$.
\dli On peut avoir $h_{(\theta'',\theta')}(t)=k_3<+\infty$ dans deux situations : 

-- $\exists\omega_t\in T^{-1}(t)$ tel que $p_{\theta'}(\omega_t)>0$

-- $\forall\omega\in T^{-1}(t)$, $p_{\theta'}(\omega)=0$ et $p_{\theta''}(\omega)=0$.
\dli Nous allons les analyser successivement.

i) $\exists\omega_t\in T^{-1}(t)$ tel que $p_{\theta'}(\omega_t)>0$.

Montrons d'abord qu'il existe $\theta_a\in\Theta_1$ tel que $p_{\theta_a}(\omega_t)>0$, ceci nous permettra d'utiliser la propri\'et\'e 2 du lemme 2 de l'annexe IV ; si ce n'\'etait pas le cas on aurait 
$h_{(\theta',\theta_1)}(t)=+\infty$ pour tout $\theta_1\in\Theta_1$ ; les densit\'es $p_{\theta_1}$ seraient toutes nulles sur $T^{-1}([t,+\infty[)$, ce qui est impossible puisque $t\notin D^{\Theta_1}_s$.
\dli Nous allons appliquer la propri\'et\'e 2 du lemme 2 de l'annexe IV, avec 
$\theta_a\prec\theta_b=\theta'\prec\theta_c=\theta''$. On obtient : 
\dli 
$\phi^{(\theta_c,\theta_b)}_{(k_3,0)}\relmont{\leq}{p.s.}sup\{\phi^{(\theta_b,\theta_a)}_{(k_1,0)},\phi^{(\theta_c,\theta_a)}_{(k_2,0)},\II_{D_i^{\{\theta\succ\theta_a\}}}(T)\}$, (avec $h_{(\theta_b,\theta_a)}(t)=k_1$, 
$h_{(\theta_c,\theta_a)}(t)=k_2$, $h_{(\theta_c,\theta_b)}(t)=k_3$ et la convention $\phi^{(\theta'',\theta')}_{(0,0)}=\II_\emptyset=f_{(-\infty,1)}$).
Par d\'efinition des $e_t^{(\theta'',\theta')}$ on a $e_t^{(\theta_b,\theta_a)}=\phi^{(\theta_b,\theta_a)}_{(k_1,0)}$,  $e_t^{(\theta_c,\theta_a)}=\phi^{(\theta_c,\theta_a)}_{(k_2,0)}$ et  $e_t^{(\theta_c,\theta_b)}=\phi^{(\theta_c,\theta_b)}_{(k_3,0)}$ (voir le d\'ebut de la partie I de cette d\'emonstration), donc  
$e_t^{(\theta_c,\theta_b)}\relmont{\leq}{p.s.}sup\{e_t^{(\theta_b,\theta_a)},e_t^{(\theta_c,\theta_a)},\II_{D_i^{\{\theta\succ\theta_a\}}}(T)\}$. Ceci implique l'in\'egalit\'e recherch\'ee car $\Theta_0\subseteq \{\theta\succ\theta_a\}$ entra\^{\i}ne : 
${D_i^{\{\theta\succ\theta_a\}}}\subseteq {D^{\Theta_0}_i}$ (voir le lemme 1 de l'annexe IV), donc $\II_{D_i^{\{\theta\succ\theta_a\}}}(T)\leq \II_{D^{\Theta_0}_i}\!(T)$.

ii) $\forall\omega\in T^{-1}(t)$, $p_{\theta'}(\omega)=0$ et $p_{\theta''}(\omega)=0$.

Montrons d'abord qu'il existe $\theta_a\in\Theta_1$ et $\omega_a\in T^{-1}(t)$  tels que $p_{\theta_a}(\omega_a)>0$.
\dli Lorsque la densit\'e $p_{\theta_1}$, $\theta_1\in\Theta_1$, est nulle sur 
$T^{-1}(t)$, notons $I_{t,\theta_1}$ le plus grand intervalle de $\IR-(D_i^{\Theta_0}\cup D_s^{\Theta_1})$ contenant $t$ et ind\'etermin\'e pour $(\theta',\theta_1)$. Comme $t\notin N$ il n'est pas totalement ind\'etermin\'e ; d'apr\`es la propri\'et\'e 2 du lemme 1 de l'annexe III on a 
$h_{(\theta',\theta_1)}^{\Theta_1}(I_{t,\theta_1}^+)=+\infty$ ; on a donc $h_{(\theta',\theta_1)}^{\Theta_1}(t')=+\infty$ pour $t'\in I=\{x>I_{t,\theta_1}\}$, ce qui implique que la densit\'e $p_{\theta_1}$ est nulle sur $T^{-1}(I)$ ; comme il en est de m\^eme sur 
$I_{t,\theta_1}$, $t$ appartient \`a $D_s^{\theta_1}$.
\dli On ne peut pas avoir cette situation pour tous les \'el\'ements de ${\Theta_1}$ puisque $t\notin D_s^{\Theta_1}$. Il existe donc bien $\theta_a\in\Theta_1$ pour lequel la densit\'e $p_{\theta_a}$ est non identiquement nulle sur $T^{-1}(t)$.
\dli $t$ \'etant ind\'etermin\'e pour $(\theta'',\theta')$ on a 
$h_{(\theta',\theta_a)}(t)=h_{(\theta'',\theta_a)}(t)=0$, ce qui implique $]-\infty,t]\subseteq D_i^{\theta'}\cap D_i^{\theta''}=D_i^{\theta'}$ ; $h_{(\theta'',\theta')}$ est donc \'egale \`a $k_3$ sur $D_i^{\theta'}\cap D_i^{\theta''}=D_i^{\theta'}$ puisqu'elle est constante sur tout intervalle ind\'etermin\'e pour $(\theta'',\theta')$ ; 
on a donc $e_t^{(\theta'',\theta')}=f_{(-\infty,1)}$ et l'in\'egalit\'e recherch\'ee est \'evidemment v\'erifi\'ee.

2\up{\`eme} cas : $\theta'\in\Theta_1$, $\theta''\in\Theta_1$ et $\theta'\!\prec\theta''$.

La d\'emonstration est semblable \`a celle du 1\up{er} cas.
\dli Posons $h_{(\theta'',\theta')}(t)=k_1$. On a $k_1>0$ ; en effet si on avait $h_{(\theta'',\theta')}(t)=0$, la densit\'e $p_{\theta''}$ serait nulle sur $T^{-1}(]-\infty,t])$ et $t$ appartiendrait \`a $D^{\theta''}_i$ donc \`a $D^{\Theta_0}_i$ puisque $\theta''\!\prec\Theta_0$ (voir le lemme 1 de l'annexe IV) ; ceci est impossible car $t\in\,\IR-(D_i^{\Theta_0}\cup D_s^{\Theta_1}\cup N)$.
\dli On peut avoir $h_{(\theta'',\theta')}(t)=k_1>0$ dans deux situations : 

-- $\exists\omega_t\in T^{-1}(t)$ tel que $p_{\theta''}(\omega_t)>0$

-- $\forall\omega\in T^{-1}(t)$, $p_{\theta'}(\omega)=0$ et $p_{\theta''}(\omega)=0$.
\dli Nous allons les analyser successivement.

i) $\exists\omega_t\in T^{-1}(t)$ tel que $p_{\theta''}(\omega_t)>0$.

Montrons d'abord qu'il existe $\theta_c\in\Theta_0$ tel que $p_{\theta_c}(\omega_t)>0$ ; si ce n'\'etait pas le cas on aurait 
$h_{(\theta_0,\theta'')}(t)=0$ pour tout $\theta_0\in\Theta_0$ ; les densit\'es $p_{\theta_0}$ seraient toutes nulles sur $T^{-1}(]-\infty,t])$, ce qui est impossible puisque $t\notin D^{\Theta_0}_i$.
\dli Nous allons appliquer la propri\'et\'e 3 du lemme 2 de l'annexe IV, avec 
$\theta_a=\theta'\prec\theta_b=\theta''\prec\theta_c$. On obtient : 
\dli $\phi^{(\theta_b,\theta_a)}_{(k_1,0)}\relmont{\leq}{p.s.}sup\{\phi^{(\theta_c,\theta_a)}_{(k_2,0)},\phi^{(\theta_c,\theta_b)}_{(k_3,0)},\II_{D_i^{\{\theta\succ\theta_b\}}}(T)\}$ (avec $h_{(\theta_b,\theta_a)}(t)=k_1$, 
$h_{(\theta_c,\theta_a)}(t)=k_2$, $h_{(\theta_c,\theta_b)}(t)=k_3$ et la convention $\phi^{(\theta'',\theta')}_{(0,0)}=\II_\emptyset=f_{(-\infty,1)}$).
Comme dans le 1\up{er} cas on a : $e_t^{(\theta_b,\theta_a)}=\phi^{(\theta_b,\theta_a)}_{(k_1,0)}$,  $e_t^{(\theta_c,\theta_a)}=\phi^{(\theta_c,\theta_a)}_{(k_2,0)}$ et  $e_t^{(\theta_c,\theta_b)}=\phi^{(\theta_c,\theta_b)}_{(k_3,0)}$, donc  
$e_t^{(\theta_b,\theta_a)}\relmont{\leq}{p.s.}sup\{e_t^{(\theta_c,\theta_a)},e_t^{(\theta_c,\theta_b)},\II_{D_i^{\{\theta\succ\theta_b\}}}(T)\}$. Ceci implique l'in\'egalit\'e recherch\'ee puisque $\Theta_0\subseteq \{\theta\succ\theta_b\}$ entra\^{\i}ne comme pr\'ec\'edemment $\II_{D_i^{\{\theta\succ\theta_b\}}}(T)\leq \II_{D^{\Theta_0}_i}\!(T)$.

ii) $\forall\omega\in T^{-1}(t)$, $p_{\theta'}(\omega)=0$ et $p_{\theta''}(\omega)=0$.

Montrons d'abord qu'il existe $\theta_c\in\Theta_0$ et $\omega_c\in T^{-1}(t)$  tels que $p_{\theta_c}(\omega_c)>0$.
\dli Lorsque la densit\'e $p_{\theta_0}$, $\theta_0\in\Theta_0$, est nulle sur 
$T^{-1}(t)$, notons $I_{t,\theta_0}$ le plus grand intervalle de $\IR-(D_i^{\Theta_0}\cup D_s^{\Theta_1})$ contenant $t$ et ind\'etermin\'e pour $(\theta_0,\theta'')$. Comme $t\notin N$ il n'est pas totalement ind\'etermin\'e ; d'apr\`es la propri\'et\'e 2 du lemme 1 de l'annexe III on a 
$h_{(\theta_0,\theta'')}^{\Theta_1}(I_{t,\theta_0}^-)=0$ ; $h_{(\theta_0,\theta'')}^{\Theta_1}$ est donc nulle sur $I=\{x<I_{t,\theta_0}\}$, ce qui implique que la densit\'e 
$p_{\theta_0}$ est nulle sur $T^{-1}(I)$ ; comme il en est de m\^eme sur 
$I_{t,\theta_0}$, $t$ appartient \`a $D_i^{\theta_0}$.
\dli On ne peut pas avoir cette situation pour tous les \'el\'ements de ${\Theta_0}$ puisque $t\notin D_i^{\Theta_0}$. Il existe donc bien $\theta_c\in\Theta_0$ pour lequel la densit\'e $p_{\theta_c}$ est non identiquement nulle sur $T^{-1}(I)$.
\dli $t$ \'etant ind\'etermin\'e pour $(\theta'',\theta')$ on a 
$h_{(\theta_c,\theta'')}(t)=h_{(\theta_c,\theta')}(t)=+\infty$, ce qui implique $[t,+\infty[\subseteq D_s^{\theta''}\cap D_s^{\theta'}=D_s^{\theta''}$ ; $h_{(\theta'',\theta')}$ est donc \'egale \`a $k_1$ sur $D_s^{\theta''}\cap D_s^{\theta'}=D_s^{\theta''}$ puisqu'elle est constante sur tout intervalle ind\'etermin\'e pour $(\theta'',\theta')$ ; 
on a donc $e_t^{(\theta'',\theta')}\leq 1-\II_{D_s^{\theta''}}(T)$ alors que 
$e_t^{(\theta_c,\theta'')}\geq 1-\II_{D_s^{\theta''}}(T)$ puisque 
$\{h_{(\theta_c,\theta'')}=+\infty\}\subseteq D_s^{\theta''}$ ; l'in\'egalit\'e recherch\'ee est ainsi d\'emontr\'ee.

\medskip
{\parindent=-5mm c) Pour $t\notin D_s^{\Theta_1}$ on a $f_{(b_t,v_t)}\relmont{\leq}{p.s.}f^{\Theta_1}_{(b'_t,v'_t)}$.}

D'apr\`es la d\'efinition de $f_{(b'_t,v'_t)}^{\Theta_1}$ \`a partir des $g_t^{(\theta_0,\theta_1)}$, nous avons : 
$f_{(b'_t,v'_t)}^{\Theta_1}\geq inf\{f_{(b_t,v_t)},1-\II_{D^{\Theta_1}_s}\!(T)\}$. La propri\'et\'e recherch\'ee revient \`a d\'emontrer : $f_{(b_t,v_t)}\relmont{\leq}{p.s.}1-\II_{D^{\Theta_1}_s}\!(T)$.
Nous allons le faire en montrant qu'il existe $\theta_1\in\Theta_1$ tel que 
$inf\{g_t^{(\theta,\theta_1)}\}_{\theta\in\Theta_0}\relmont{\leq}{p.s.}1-\II_{D^{\Theta_1}_s}\!(T)$.

Comme $t\notin D_s^{\Theta_1}$, il existe $\theta_1\in\Theta_1$, $x_1\geq t$ et $\omega_1\in T^{-1}(x_1)$ tels que $p_{\theta_1}(\omega_1)>0$. Pour tout $\theta$ de $\Theta_0$ on a $h_{(\theta,\theta_1)}(x_1)<+\infty$, donc 
$h_{(\theta,\theta_1)}(t)=k_\theta<+\infty$ et $g_t^{(\theta,\theta_1)}=\{h_{(\theta,\theta_1)}(T)\leq k_\theta\}$. 
\dli Pour que $\theta_1$ convienne il nous faut montrer que le cas de figure \dli $inf\{g_t^{(\theta,\theta_1)}\}_{\theta\in\Theta_0}\!\!>\!\!1-\!\!\II_{D^{\Theta_1}_s}\!\!(T)$ implique 
$inf\{g_t^{(\theta,\theta_1)}\}_{\theta\in\Theta_0}\!\!\relmont{=}{p.s.}\!1-\!\!\II_{D^{\Theta_1}_s}\!\!(T)$. 
\dli Il est \'equivalent de supposer $I=\cap_{\theta\in\Theta_0}\{h_{(\theta,\theta_1)}(t)\leq k_\theta\}\cap D^{\Theta_1}_s$ non vide et de montrer que l'on a 
$P_\theta(T^{-1}(I))=0$ pour tout $\theta$ de $\Theta$. C'est \'evident pour $\theta\in\Theta_1$ puisque la densit\'e $p_\theta$ est alors nulle sur  
$T^{-1}(D^{\Theta_1}_s)\supseteq T^{-1}(I)$ (voir la d\'efinition 4.2.1). Lorsque $\theta\in\Theta_0$ la densit\'e $p_\theta$ est encore nulle sur  
$T^{-1}(I)$ ; en effet pour tout $\omega\in T^{-1}(I)$ on a $p_{\theta_1}(\omega)=0$ et $h_{(\theta,\theta_1)}(T(\omega))\leq k_\theta<+\infty$, ceci n'est possible qu'avec $p_{\theta}(\omega)=0$.

\medskip
{\parindent=-5mm d) Pour $t\notin D_i^{\Theta_0}\cup D_s^{\Theta_1}\cup N$ on a $f_{(b_t,v_t)}\relmont{\geq}{p.s.}f^{\Theta_1}_{(b'_t,v'_t)}$.}

Cette in\'egalit\'e s'\'ecrit : 
\dli $inf\{g_t^{(\theta'',\theta')}\}_{\theta'\!\prec\theta''}\!\!\relmont{\geq}{p.s.}\!inf\{g_t^{(\theta_0,\theta_1)}\}_{(\theta_0,\theta_1)\in\Theta_0\times\Theta_1}\!\cup\!\{1-\!\!\II_{D^{\Theta_1}_s}\!\!(T)\}$. Comme en b), 
nous allons consid\'erer les couples $(\theta'',\theta')$ compos\'es de deux \'el\'ements de $\Theta_0$ (resp. $\Theta_1$), $\theta'\!\prec\theta''$, il suffit de montrer qu'il existe $\theta_a\in\Theta_1$ (resp. $\theta_c\in\Theta_0$) tel que :
$g_t^{(\theta'',\theta')}\relmont{\geq}{p.s.}inf\{g_t^{(\theta',\theta_a)},g_t^{(\theta'',\theta_a)},1-\II_{D^{\Theta_1}_s}\!(T)\}$ (resp. 
$g_t^{(\theta'',\theta')}\relmont{\geq}{p.s.}inf\{g_t^{(\theta_c,\theta')},g_t^{(\theta_c,\theta'')},1-\II_{D^{\Theta_1}_s}\!(T)\}$). 

1\up{er} cas : $\theta'\in\Theta_0$, $\theta''\in\Theta_0$ et $\theta'\!\prec\theta''$.

Ces conditions sont celles du 1\up{er} cas de b), elles nous conduisent \`a analyser les deux m\^emes situations.

i) $\exists\omega_t\in T^{-1}(t)$ tel que $p_{\theta'}(\omega_t)>0$.

Comme en i) du 1\up{er} cas de b), il existe $\theta_a\in\Theta_1$ tel que $p_{\theta_a}(\omega_t)>0$.
La propri\'et\'e 2 du lemme 2 de l'annexe IV appliqu\'ee \`a
$\theta_a\prec\theta_b=\theta'\prec\theta_c=\theta''$ nous donne : 
\dli 
$\phi^{(\theta_c,\theta_b)}_{(k_3,1)}\relmont{\geq}{p.s.}inf\{\phi^{(\theta_b,\theta_a)}_{(k_1,1)},\phi^{(\theta_c,\theta_a)}_{(k_2,1)},1-\II_{D_s^{\{\theta\prec\theta_b\}}}(T)\}$, (avec $h_{(\theta_b,\theta_a)}(t)=k_1$, 
$h_{(\theta_c,\theta_a)}(t)=k_2$, $h_{(\theta_c,\theta_b)}(t)=k_3$ et la convention $\phi^{(\theta'',\theta')}_{(\infty,1)}=\II_\Omega=f_{(+\infty,0)}$).
Par d\'efinition des $g_t^{(\theta'',\theta')}$ on a $g_t^{(\theta_b,\theta_a)}=\phi^{(\theta_b,\theta_a)}_{(k_1,1)}$,  $g_t^{(\theta_c,\theta_a)}=\phi^{(\theta_c,\theta_a)}_{(k_2,1)}$ et  $g_t^{(\theta_c,\theta_b)}=\phi^{(\theta_c,\theta_b)}_{(k_3,1)}$ (voir le d\'ebut de la partie I de cette d\'emonstration), donc  
$g_t^{(\theta_c,\theta_b)}\relmont{\geq}{p.s.}inf\{g_t^{(\theta_b,\theta_a)},g_t^{(\theta_c,\theta_a)},1-\II_{D_s^{\{\theta\prec\theta_b\}}}(T)\}$. Ceci implique l'in\'egalit\'e recherch\'ee car $\Theta_1\subseteq \{\theta\prec\theta_b\}$ entra\^{\i}ne : 
\dli ${D_s^{\{\theta\prec\theta_b\}}}\subseteq {D^{\Theta_1}_s}$ (voir le lemme 1 de l'annexe IV), donc 
\dli $1-\II_{D_s^{\{\theta\prec\theta_b\}}}(T)\geq 1-\II_{D^{\Theta_1}_s}\!(T)$.

ii) $\forall\omega\in T^{-1}(t)$, $p_{\theta'}(\omega)=0$ et $p_{\theta''}(\omega)=0$.

Comme en ii) du 1\up{er} cas de b), il existe $\theta_a\in\Theta_1$ et $\omega_a\in T^{-1}(t)$  tels que $p_{\theta_a}(\omega_a)>0$ ; de plus 
$h_{(\theta'',\theta')}$ est \'egale \`a $h_{(\theta'',\theta')}(t)=k_3$ sur $D_i^{\theta'}\cap D_i^{\theta''}=D_i^{\theta'}$.
\dli On a donc $g_t^{(\theta'',\theta')}\geq \II_{D_i^{\theta'}}\!(T)$, alors que $g_t^{(\theta',\theta_a)}\leq \II_{D_i^{\theta'}}\!(T)$ puisque 
$t\in\{h_{(\theta',\theta_a)}=0\}\subseteq D_i^{\theta'}$. L'in\'egalit\'e recherch\'ee est bien v\'erifi\'ee.

2\up{\`eme} cas : $\theta'\in\Theta_1$, $\theta''\in\Theta_1$ et $\theta'\!\prec\theta''$.

On obtient deux situations possibles, celles du 2\up{\`eme} cas de b).

i) $\exists\omega_t\in T^{-1}(t)$ tel que $p_{\theta''}(\omega_t)>0$.

Comme en i) du 2\up{\`eme} cas de b), il existe $\theta_c\in\Theta_0$ tel que $p_{\theta_c}(\omega_t)>0$. La propri\'et\'e 3 du lemme 2 de l'annexe IV appliqu\'ee \`a 
$\theta_a=\theta'\prec\theta_b=\theta''\prec\theta_c$ donne : 
\dli $\phi^{(\theta_b,\theta_a)}_{(k_1,1)}\relmont{\geq}{p.s.}inf\{\phi^{(\theta_c,\theta_a)}_{(k_2,1)},\phi^{(\theta_c,\theta_b)}_{(k_3,1)},1-\II_{D_s^{\{\theta\prec\theta_c\}}}(T)\}$ (avec $h_{(\theta_b,\theta_a)}(t)=k_1$, 
$h_{(\theta_c,\theta_a)}(t)=k_2$, $h_{(\theta_c,\theta_b)}(t)=k_3$ et la convention $\phi^{(\theta'',\theta')}_{(\infty,1)}=\II_\Omega=f_{(+\infty,0)}$).
Comme dans le 1\up{er} cas de cette partie on a : \dli $g_t^{(\theta_b,\theta_a)}=\phi^{(\theta_b,\theta_a)}_{(k_1,1)}$,  $g_t^{(\theta_c,\theta_a)}=\phi^{(\theta_c,\theta_a)}_{(k_2,1)}$ et  $g_t^{(\theta_c,\theta_b)}=\phi^{(\theta_c,\theta_b)}_{(k_3,1)}$, donc  
$g_t^{(\theta_b,\theta_a)}\relmont{\geq}{p.s.}inf\{g_t^{(\theta_c,\theta_a)},g_t^{(\theta_c,\theta_b)},1-\II_{D_s^{\{\theta\prec\theta_c\}}}(T)\}$. Ceci implique l'in\'egalit\'e recherch\'ee puisque $\Theta_1\subseteq \{\theta\prec\theta_c\}$ entra\^{\i}ne comme pr\'ec\'e\-dem\-ment 
$1-\II_{D_s^{\{\theta\prec\theta_c\}}}(T)\geq 1-\II_{D^{\Theta_1}_s}\!(T)$.

ii) $\forall\omega\in T^{-1}(t)$, $p_{\theta'}(\omega)=0$ et $p_{\theta''}(\omega)=0$.

Comme en ii) du 2\up{\`eme} cas de b), $h_{(\theta'',\theta')}$ est \'egale \`a $h_{(\theta'',\theta')}(t)=k_1$ sur $D_s^{\theta''}\cap D_s^{\theta'}=D_s^{\theta''}$. 
On a donc $g_t^{(\theta'',\theta')}=f_{(+\infty,0)}=\II_\Omega$  et l'in\'egalit\'e recherch\'ee est \'evidente.

\medskip
{\parindent=-10mm III -- $\forall t\in(D_i^{\Theta_0}-N)$,\quad $\forall \theta_0\in\Theta_0$,\quad $G_{\theta_0}(K(t))=0$.}

Soient $t\in(D_i^{\Theta_0}-N)$ et $\theta_0\in\Theta_0$.
Par d\'efinition de $D_i^{\Theta_0}$, les densit\'es $p_{\theta}$, $\theta\in{\Theta_0}$, sont nulles sur $T^{-1}(t)$. Comme $t\notin N$, il existe $\theta_1\in{\Theta_1}$ et $\omega_1\in T^{-1}(t)$ tels que 
$p_{\theta_1}(\omega_1)>0$ ; on a alors $h_{(\theta_0,\theta_1)}(t)=0$ et d'apr\`es la partie I de cette d\'emonstration : 
$g_t^{(\theta_0,\theta_1)}=\phi_{(0,1)}^{(\theta_0,\theta_1)}$, 
$f_{(b_t,v_t)}\leq\phi_{(0,1)}^{(\theta_0,\theta_1)}$ ; $G_{\theta_0}(K(t))$ \'etant \'egal \`a $[E_{\theta_0}(f_{(a_t,u_t)})+E_{\theta_0}(f_{(b_t,v_t)})]/2$ avec $f_{(a_t,u_t)}\leq f_{(b_t,v_t)}$ (voir I), on a 
$G_{\theta_0}(K(t))\leq E_{\theta_0}(f_{(b_t,v_t)})\leq E_{\theta_0}(\phi_{(0,1)}^{(\theta_0,\theta_1)})=0$.
\medskip
{\parindent=-10mm IV --- $\forall t\in(D_s^{\Theta_1}-N)$,\quad $\forall \theta_1\in\Theta_1$,\quad $G_{\theta_1}(K(t))=1$.}

La d\'emonstration est semblable \`a celle de III.
\dli Soient $t\in(D_s^{\Theta_1}-N)$ et $\theta_1\in\Theta_1$.
Par d\'efinition de $D_s^{\Theta_1}$, les densit\'es $p_{\theta}$, $\theta\in{\Theta_1}$, sont nulles sur $T^{-1}(t)$. Comme $t\notin N$, il existe $\theta_0\in{\Theta_0}$ et $\omega_0\in T^{-1}(t)$ tels que 
$p_{\theta_0}(\omega_0)>0$ ; on a alors $h_{(\theta_0,\theta_1)}(t)=+\infty$ et d'apr\`es la partie I de cette d\'emonstration : 
$e_t^{(\theta_0,\theta_1)}=\phi_{(\infty,0)}^{(\theta_0,\theta_1)}$, 
$f_{(a_t,u_t)}\geq\phi_{(\infty,0)}^{(\theta_0,\theta_1)}$ ; cette fois on a 
$G_{\theta_1}(K(t))\geq E_{\theta_1}(f_{(a_t,u_t)})\geq E_{\theta_1}(\phi_{(\infty,0)}^{(\theta_0,\theta_1)})=1$.

\medskip\centerline{\hbox to 3cm{\bf \hrulefill}}\par}

\vfill\eject
\medskip
{\bf Proposition 5.2.3}
\medskip
\medskip
\moveleft 10.4pt\hbox{\vrule\kern 10pt\vbox{\defpro

Soit $(\Omega ,{\cal A},(p_\theta.\mu)_{\theta\in\Theta})$ un mod\`ele
statistique \`a rapport de vraisem\-blance monotone pour une statistique r\'eelle $T$, $\Theta$ \'etant un intervalle de $\IR$ muni de l'ordre usuel.
$K(T)$ est une statistique essentielle globale et $G_\theta$ d\'esigne sa fonction de r\'epartition moyenne sous $P_\theta$.
Notons $(\{\Theta_1^f,\Theta_0^f\})_{f\in{\cal F}}$ l'ensemble des hypoth\`eses unilat\'erales et consid\'erons pour chacun de ces probl\`emes de d\'ecision le vote le plus favorable sous $\Theta_1^f$ (resp. $\Theta_0^f$) ; 
pour tout $\omega\in\Omega$ et tout $f\in{\cal F}$ posons 
$Q^\omega_1({\Theta_1^f})=Q^\omega_{\Theta_1^f}(\{1\})$ (resp. $Q^\omega_0({\Theta_1^f})=Q^\omega_{\Theta_0^f}(\{1\})$). 

\dli Les applications $Q_1$ et $Q_0$, de $\Omega\times\{\Theta_1^f\}_{f\in{\cal F}}$ dans $[0,1]$, d\'efinissent des votes compatibles si et seulement si, pour presque tout $t\!\in\IR$, la fonction $G_\theta(K(t))$ v\'erifie : 
\dli 1) $G_\theta(K(t))$ est continue \`a droite en $\theta_f\in\Theta-\{sup\Theta\}$ si $t\in\IR-D_i^{\{\theta>\theta_f\}}$ 
\dli 2) $G_\theta(K(t))$ est continue \`a gauche en $\theta_f\in\Theta-\{inf\Theta\}$ si $t\in\IR-D_s^{\{\theta<\theta_f\}}$ 
\dli 3) si $inf\Theta\notin\Theta$ on a $lim_{\theta\rightarrow inf\Theta}G_\theta(K(t))=1$
\dli 4) si $sup\Theta\notin\Theta$ on a $lim_{\theta\rightarrow sup\Theta}G_\theta(K(t))=0$.
\dli Sous ces conditions les deux types de votes $Q_1$ et $Q_0$ ont presque partout le m\^eme prolongement en une probabilit\'e sur les bor\'eliens de $\Theta$.
}}\medskip

\medskip
{\leftskip=15mm \dli {\bf D\'emonstration}
\medskip
Posons $N=\{t\in\IR\, ;\, \forall\omega\in T^{-1}(t)\quad\forall\theta\in\Theta\quad p_\theta(\omega)=0\}$ et notons simplement $D_i^f$ et $D_s^f$ les demi-droites inf\'erieures et sup\'erieures $D_i^{\Theta_0^f}$ et $D_s^{\Theta_1^f}$. D'apr\`es les propositions 4.3.1 et 5.2.2, pour une r\'ealisation $\omega\notin T^{-1}(N)$ et des hypoth\`eses unilat\'erales $\{\Theta_1^f,\Theta_0^f\}$, le vote   $Q^\omega_{\theta}$ est d\'efini par : 
$$Q^\omega_\theta(\{1_f\})=\left\{\matrix{
1\hfill & si & T(\omega)\in D_i^f-D_s^f\hfill\cr
1-G_\theta(K(T(\omega)))\hfill & si & T(\omega)\in\IR\!-(D_i^f\cup D_s^f)\hfill \cr
0\hfill & si & T(\omega)\in D_s^f\hfill \cr
}\right. $$

\dli Dans cette proposition nous consid\'erons deux familles de votes, la premi\`ere est d\'efinie presque s\^urement par $Q^\omega_1({\Theta_1^f})=sup_{\theta\in\Theta_1^f}Q^\omega_{\theta}(\{1_f\})$, la seconde par $Q^\omega_0({\Theta_1^f})=1-Q^\omega_{\Theta_0^f}(\{0\})=1-sup_{\theta\in\Theta_0^f}Q^\omega_{\theta}(\{0_f\})=inf_{\theta\in\Theta_0^f}Q^\omega_{\theta}(\{1_f\})$ (voir la d\'efinition 4.3.2).
Ces deux familles sont identiques si toutes les hypoth\`eses unilat\'erales du mod\`ele sont adjacentes. 
\dli Nous devons d\'emontrer que $Q_1$ et $Q_0$ v\'erifient les propri\'et\'es a), b) et c) de la d\'efinition 5.2.1 pour presque tout $\omega$  de $\Omega-T^{-1}(N)$ si et seulement si  $G_\theta$ poss\`ede les quatre propri\'et\'es de l'\'enonc\'e.

\medskip
{\parindent=-10mm I --- Condition suffisante.}

Notons $N_c$ l'ensemble des $t\in\IR$ pour lesquels $G_\theta(K(t))$ ne v\'erifie pas les quatre propri\'et\'es de l'\'enonc\'e. $T^{-1}(N_c)$ \'etant n\'egligeable, il suffit de d\'emontrer les propri\'et\'es a),b) et c) sur $\Omega-T^{-1}(N\cup N_c)$.

Soit $\omega$ tel que $T(\omega)=t\notin (N\cup N_c)$. 
\medskip

\medskip
{\parindent=-5mm 1\up{er} cas : $Q^\omega_1({\Theta_1^f})=sup_{\theta\in\Theta_1^f}Q^\omega_{\theta}(\{1_f\})$ pour tout $f\in{\cal F}$.}
\medskip
a) Si $\Theta_1^f\subset\Theta_1^{f'}$ on a $Q^\omega_1({\Theta_1^f})\leq Q^\omega_1(\Theta_1^{f'})$.

Lorsque $t\in D_i^{f'}-D_s^{f'}$ (resp. $t\in D_s^f$) l'in\'egalit\'e recherch\'ee est \'evidente car on a, pour tout $\theta$ de $\Theta$, $Q^\omega_\theta(\{1_{f'}\})=1$ (resp. $Q^\omega_\theta(\{1_f\})=0$) donc
$Q^\omega_1(\Theta_1^{f'})=1$ (resp. $Q^\omega_1(\Theta_1^{f})=0$).
\dli Il nous reste le cas $t\notin(D_i^{f'}-D_s^{f'})\cup D_s^f$. Comme $\Theta_1^f\subset\Theta_1^{f'}$ implique $D_i^f\subseteq D_i^{f'}$ et $D_s^f\supseteq D_s^{f'}$ (voir le lemme 1 de l'annexe IV), $t$ n'appartient pas non plus \`a $D_i^{f}-D_s^{f}\subseteq D_i^{f'}-D_s^{f'}$ et \`a $D_s^{f'}$. On a, pour tout $\theta$ de $\Theta$, 
\dli $Q^\omega_\theta(\{1_f\})=Q^\omega_\theta(\{1_{f'}\})=1-G_\theta(K(t))$ donc 
\dli $Q^\omega_1({\Theta_1^f})=sup_{\theta\in\Theta_1^f}Q^\omega_{\theta}(\{1_f\})\leq sup_{\theta\in\Theta_1^{f'}\supset\Theta_1^f}Q^\omega_{\theta}(\{1_{f'}\})=Q^\omega_1(\Theta_1^{f'})$.
\medskip
b-1) Si $(\Theta_1^{f_n})_{n\in\IN}$ cro\^{\i}t vers $\Theta$, la suite  $(Q^\omega_1(\Theta_1^{f_n}))_{n\in\IN}$ cro\^{\i}t vers $1$.

$\Theta_1^{f_n}=\Theta\cap]-\infty,\theta_n)$ n'\'etant jamais \'egal \`a $\Theta$, l'existence d'une suite $(\Theta_1^{f_n})$ croissant vers $\Theta$
n'est possible que si $sup\Theta\notin\Theta$. C'est dans ce cas qu'il y a quelque chose \`a d\'emontrer.
\dli Posons $D_i^+=\cup_{n\in\IN}D_i^{f_n}$ et $D_s^-=\cap_{n\in\IN}D_s^{f_n}$ ; la suite $(\Theta_1^{f_n})$ \'etant croissante les suites $(D_i^{f_n})$ et $(D_s^{f_n})$ sont respectivement croissante et d\'ecroissante (voir le lemme 1 de l'annexe IV) ; la suite 
$(D_i^{f_n}-D_s^{f_n})$ cro\^{\i}t vers $D_i^+-D_s^-$.
\dli Lorsque $t\in D_i^+-D_s^-$ on a la propri\'et\'e recherch\'ee puisqu'\`a partir d'un certain rang $t\in(D_i^{f_n}-D_s^{f_n})$ donc $Q^\omega_1(\Theta_1^{f_n})=1$.
\dli Le cas $t\in D_s^-$ est impossible car $D_s^-\subseteq N$ ; en effet, la convergence de $(\Theta_1^{f_n})$ vers $\Theta$ et le lemme 1 de l'annexe IV entra\^{\i}nent $D_s^-=\cap_{\theta\in\Theta}D_s^\theta$.
\dli Il reste le cas $t\in\IR\!-(D_i^+\cup D_s^-)\not=\emptyset$. On a toujours :
\dli $lim_{n\rightarrow +\infty}Q^\omega_1(\Theta_1^{f_n})=lim\!\nearrow_{n\rightarrow +\infty}sup_{\theta\in\Theta_1^{f_n}}Q^\omega_{\theta}(\{1_{f_n}\})$. 
Comme $t\notin D_i^+\cup D_s^-$, \`a partir d'un certain rang $t\notin D_i^{f_n}\cup D_s^{f_n}$, on a donc $lim_{n\rightarrow +\infty}Q^\omega_1(\Theta_1^{f_n})=sup_{\theta\in\Theta}(1-G_{\theta}(K(t)))$. \dli $G_{\theta}(K(t))$ \'etant une fonction d\'ecroissante de $\theta$ (voir la partie a) de la d\'emonstration de 5.1.1), $lim_{n\rightarrow +\infty}Q^\omega_1(\Theta_1^{f_n})=lim_{\theta\rightarrow sup\Theta}(1-G_\theta(K(t)))$ qui est \'egal \`a $1$ d'apr\`es la propri\'et\'e 4) de l'\'enonc\'e.
\medskip
b-2) Si $(\Theta_1^{f_n})_{n\in\IN}$ d\'ecro\^{\i}t vers le vide, la suite  $(Q^\omega_1(\Theta_1^{f_n}))_{n\in\IN}$ d\'ecro\^{\i}t vers $0$.

La d\'emonstration est semblable \`a celle de b-1).
\dli $\Theta_1^{f_n}=\Theta\cap]-\infty,\theta_n)$ n'\'etant jamais vide, l'existence d'une suite $(\Theta_1^{f_n})$ d\'ecroissant vers le vide 
n'est possible que si $inf\Theta\notin\Theta$. On pose alors 
$D_i^-=\cap_{n\in\IN}D_i^{f_n}$ et $D_s^+=\cup_{n\in\IN}D_s^{f_n}$. 
\dli Lorsque $t\in D_s^+$, la propri\'et\'e recherch\'ee est \'evidente puisqu'\`a partir d'un certain rang $t\in D_s^{f_n}$ donc $Q^\omega_1(\Theta_1^{f_n})=0$.
\dli Le cas $t\in D_i^-$ est impossible car $D_i^-\subseteq N$ ; en effet, la convergence de $(\Theta_0^{f_n})$ vers $\Theta$ et le lemme 1 de l'annexe IV entra\^{\i}nent $D_i^-=\cap_{\theta\in\Theta}D_i^\theta$.
\dli Il reste le cas $t\in\IR\!-(D_i^-\cup D_s^+)\not=\emptyset$. On a cette fois :
\dli $lim_{n\rightarrow +\infty}Q^\omega_1(\Theta_1^{f_n})=lim\!\searrow_{n\rightarrow +\infty}sup_{\theta\in\Theta_1^{f_n}}Q^\omega_{\theta}(\{1_{f_n}\})=$
\dli $lim_{n\rightarrow +\infty}sup_{\theta\in\Theta_1^{f_n}}(1-G_{\theta}(K(t)))=lim_{\theta\rightarrow inf\Theta}(1-G_\theta(K(t)))$. 
\dli La propri\'et\'e 3) de l'\'enonc\'e implique bien $lim_{n\rightarrow +\infty}Q^\omega_1(\Theta_1^{f_n})=0$.
\medskip
c) Si $(\Theta_1^{f_n})_{n\in\IN}$ et $(\Theta_1^{f'_n})_{n\in\IN}$ sont deux suites ayant m\^eme limite, l'une \'etant croissante l'autre d\'ecroissante, les suites $(Q^\omega_1(\Theta_1^{f_n}))_{n\in\IN}$ et $(Q^\omega_1(\Theta_1^{f'_n}))_{n\in\IN}$ ont m\^eme limite.

La limite commune de $(\Theta_1^{f_n})$ et $(\Theta_1^{f'_n})$ ne peut pas \^etre vide ou \'egale \`a $\Theta$, elle est donc de la forme $\Theta_1^f=\Theta\cap]-\infty,\theta_f)$ avec $\theta_f\in\Theta$. Lorsque $\theta_f$ appartient \`a $\Theta_1^f$ (resp. $\Theta_0^f=\Theta-\Theta_1^f$), les  $\Theta_1^{f_n}$ (resp. $\Theta_1^{f'_n}$) sont \'egaux \`a $\Theta_1^f$ \`a partir d'un certain rang.
Etudions s\'epar\'ement ces deux cas.

i) $\Theta_1^f=\Theta\cap]-\infty,\theta_f]$ ($\theta_f\in\Theta-\{sup\Theta\}$, $\Theta_1^{f'_n}\not=\Theta_1^f$). 
\dli On a $lim_{n\rightarrow +\infty}Q^\omega_1(\Theta_1^{f_n})=Q^\omega_1(\Theta_1^{f})\leq lim_{n\rightarrow +\infty}Q^\omega_1(\Theta_1^{f'_n})$ et 
\dli $D_i^f=D_i^{\{\theta>\theta_f\}}=\cap\!\searrow_{n\in\IN}D_i^{f'_n}$,\quad 
$D_s^f\supseteq\cup\!\nearrow_{n\in\IN}D_s^{f'_n}=D'_s$. 
\dli L'\'egalit\'e des deux limites est \'evidente lorsque $t$ appartient \`a $D_i^f-D_s^f$ (resp. $D'_s$) car on a alors $Q^\omega_1(\Theta_1^{f})=1$ (resp. $lim_{n\rightarrow +\infty}Q^\omega_1(\Theta_1^{f'_n})=0$).
\dli $t$ n'appartenant pas \`a $N$ il n'appartient pas \`a $D_i^f\cap D_s^f$, il nous reste alors le cas : $t\notin D_i^f\cup D'_s$ ; \`a partir d'un certain rang on a $t\notin D_i^{f'_n}\cup D_s^{f'_n}$ donc $Q^\omega_1(\Theta_1^{f'_n})=sup_{\theta\in\Theta_1^{f'_n}}(1-G_{\theta}(K(t)))$ ; ceci implique $lim_{n\rightarrow +\infty}Q^\omega_1(\Theta_1^{f'_n})=lim_{\theta>\theta_f\!\searrow\theta_f}(1-G_\theta(K(t)))$. D'apr\`es la propri\'et\'e 1) de l'\'enonc\'e on a 
$lim_{n\rightarrow +\infty}Q^\omega_1(\Theta_1^{f'_n})=1-G_{\theta_f}(K(t))$ ; l'\'egalit\'e entre $1-G_{\theta_f}(K(t))$ et $Q^\omega_1(\Theta_1^{f})=Q^\omega_{\theta_f}(\{1_f\})$ (voir la partie a) de la d\'emonstration 5.1.1) est triviale lorsque $t\notin D_i^f\cup D_s^f$, dans le cas contraire $t\in D_s^f- D'_s$ on a $Q^\omega_{\theta_f}(\{1_f\})=0$ mais aussi  $G_{\theta_f}(K(t))=1$ d'apr\`es la propri\'et\'e 3) d'une statistique essentielle globale (voir la proposition 5.2.2).

ii) $\Theta_1^f=\Theta\cap]-\infty,\theta_f[$ ($\theta_f\in\Theta-\{inf\Theta\}$).

\dli On a $lim_{n\rightarrow +\infty}Q^\omega_1(\Theta_1^{f'_n})=Q^\omega_1(\Theta_1^{f})\geq lim_{n\rightarrow +\infty}Q^\omega_1(\Theta_1^{f_n})$ et 
\dli $D_i^f\supseteq\cup\!\nearrow_{n\in\IN}D_i^{f_n}=D_i$,\quad 
$D_s^f=D_s^{\{\theta<\theta_f\}}=\cap\!\searrow_{n\in\IN}D_s^{f_n}$. 
\dli L'\'egalit\'e des deux limites est \'evidente lorsque $t$ appartient \`a $D_i-D_s^f$ (resp. $D_s^f$) car on a alors $lim_{n\rightarrow +\infty}Q^\omega_1(\Theta_1^{f_n})=1$ (resp. $Q^\omega_1(\Theta_1^{f})=0$).
\dli Il nous reste le cas : $t\notin D_i\cup D_s^f$ ; \`a partir d'un certain rang $t\notin D_i^{f_n}\cup D_s^{f_n}$ donc $Q^\omega_1(\Theta_1^{f_n})=sup_{\theta\in\Theta_1^{f_n}}(1-G_{\theta}(K(t)))$ ; ceci implique $lim_{n\rightarrow +\infty}Q^\omega_1(\Theta_1^{f_n})=sup_{\theta\in\Theta_1^{f}}(1-G_\theta(K(t)))$ ; cette quantit\'e est \'evidemment \'egale \`a $Q^\omega_1(\Theta_1^{f})$ lorsque $t\notin D_i^f\cup D_s^f$. Dans le cas con\nobreak traire, $t\in(D_i^f-D_s^f)-D_i$, on a $Q^\omega_1(\Theta_1^{f})=1$ et on doit montrer $sup_{\theta\in\Theta_1^{f}}(1-G_\theta(K(t)))=1$, c'est-\`a-dire $G_{\theta_f}(K(t))=0$ d'apr\`es la propri\'et\'e 2) de l'\'enonc\'e et la d\'ecroissance en $\theta$ de $G_{\theta}(K(t))$ (voir la fin de la partie a) de la d\'emonstration 5.1.1) ; cette \'egalit\'e est une cons\'equence directe de la propri\'et\'e 2) d'une statistique essentielle globale (voir la proposition 5.2.2) puisque $\theta_f\in\Theta_0^f$.

\medskip
{\parindent=-5mm 2\up{\`eme} cas : $Q^\omega_0({\Theta_1^f})=inf_{\theta\in\Theta_0^f}Q^\omega_{\theta}(\{1_f\})$ pour tout $f\in{\cal F}$.}
\medskip
a) Si $\Theta_1^f\subset\Theta_1^{f'}$ on a $Q^\omega_0({\Theta_1^f})\leq Q^\omega_0(\Theta_1^{f'})$.
\dli La d\'emonstration du a) du 1\up{er} cas convient en utilisant la d\'efinition de $Q^\omega_0$ \`a la place de celle de $Q^\omega_1$.

\medskip
b-1) $(\Theta_1^{f_n})_{n\in\IN}\nearrow\Theta\quad\Longrightarrow\quad(Q^\omega_0(\Theta_1^{f_n}))_{n\in\IN}\nearrow 1$.

Il suffit de reprendre le d\'ebut du b-1) du 1\up{er} cas et pour 
$t\in\IR\!-(D_i^+\cup D_s^-)$ on a :
\dli $lim_{n\rightarrow +\infty}Q^\omega_0(\Theta_1^{f_n})=lim\!\nearrow_{n\rightarrow +\infty}inf_{\theta\in\Theta_0^{f_n}}Q^\omega_{\theta}(\{1_{f_n}\})=$ 
\dli $lim_{n\rightarrow +\infty}inf_{\theta\in\Theta_0^{f_n}}(1-G_{\theta}(K(t)))=lim_{\theta\rightarrow sup\Theta}(1-G_{\theta}(K(t)))=1$ 
\dli d'apr\`es la d\'ecroissance de $\Theta_0^{f_n}=\Theta\cap(\theta_n,+\infty[$ vers le vide et la propri\'et\'e 4) de l'\'enonc\'e.
\medskip
b-2) $(\Theta_1^{f_n})_{n\in\IN}\searrow\emptyset\quad\Longrightarrow\quad(Q^\omega_0(\Theta_1^{f_n}))_{n\in\IN}\searrow 0$.

On reprend la d\'emonstration du b-2) du 1\up{er} cas et pour 
$t\in\IR\!-(D_i^-\cup D_s^+)$ on a cette fois :
\dli $lim_{n\rightarrow +\infty}Q^\omega_0(\Theta_1^{f_n})=lim\!\searrow_{n\rightarrow +\infty}inf_{\theta\in\Theta_0^{f_n}}Q^\omega_{\theta}(\{1_{f_n}\})=$
\dli $lim_{n\rightarrow +\infty}inf_{\theta\in\Theta_0^{f_n}}(1-G_{\theta}(K(t)))=inf_{\theta\in\Theta}(1-G_\theta(K(t)))=0$ 
\dli d'apr\`es la d\'ecroissance en $\theta$ de $G_{\theta}(K(t))$ et la propri\'et\'e 3) de l'\'enonc\'e. 

\medskip
c) $(\Theta_1^{f_n})_{n\in\IN}\nearrow\Theta_1^f\quad et\quad (\Theta_1^{f'_n})_{n\in\IN}\searrow\Theta_1^f\quad\Longrightarrow$ 
\dli $(Q^\omega_0(\Theta_1^{f_n}))_{n\in\IN}\nearrow Q_0^\omega(\Theta_1^f)\quad et\quad (Q^\omega_0(\Theta_1^{f'_n}))_{n\in\IN}\searrow Q_0^\omega(\Theta_1^f)$.

Comme pour le c) du 1\up{er} cas on distingue deux \'eventualit\'es.

i) $\Theta_1^f=\Theta\cap]-\infty,\theta_f]$ ($\theta_f\in\Theta-\{sup\Theta\}$). 
\dli On a $lim_{n\rightarrow +\infty}Q^\omega_0(\Theta_1^{f_n})=Q^\omega_0(\Theta_1^{f})\leq lim_{n\rightarrow +\infty}Q^\omega_0(\Theta_1^{f'_n})$. Les demi-droites $D_i^f$ et $D_s^f$ ont les m\^emes propri\'et\'es qu'en i) c) du 1\up{er} cas, la d\'emonstration est identique pour $t\in D_i^f\cup D'_s$. 
\dli Dans le cas contraire on a : 
\dli $lim_{n\rightarrow +\infty}Q^\omega_0(\Theta_1^{f'_n})=lim_{n\rightarrow +\infty}inf_{\theta\in\Theta_0^{f'_n}}(1-G_\theta(K(t)))=$
\dli $inf_{\theta\in\Theta_0^{f}}(1-G_\theta(K(t)))$ ; 
cette derni\`ere quantit\'e est \'evidemment \'egale \`a $Q^\omega_0(\Theta_1^{f})$ quand $t\notin D_i^f\cup D_s^f$ ; il reste le cas $t\in (D_s^f-D'_s)-D_i^f$, 
on a alors $Q^\omega_0(\Theta_1^{f})=0$ et on doit montrer $inf_{\theta\in\Theta_0^{f}}(1-G_\theta(K(t)))=0$, c'est-\`a-dire $G_{\theta_f}(K(t))=1$ d'apr\`es la propri\'et\'e 1) de l'\'enonc\'e et la d\'ecroissance en $\theta$ de $G_\theta(K(t))$ ; c'est une cons\'equence directe de la propri\'et\'e 3) d'une statistique essentielle globale (voir la proposition 5.2.2) puisque $\theta_f\in\Theta_1^{f}$.

ii) $\Theta_1^f=\Theta\cap]-\infty,\theta_f[$ ($\theta_f\in\Theta-\{inf\Theta\}$, $\Theta_1^{f_n}\not=\Theta_1^f$).
\dli On reprend le d\'ebut de la d\'emonstration du ii) c) du 1\up{er} cas en rempla\c cant $Q^\omega_1$ par $Q^\omega_0$. Lorsque $t\notin D_i\cup D_s^f$ on a maintenant : 
\dli $lim_{n\rightarrow +\infty}Q^\omega_0(\Theta_1^{f_n})=lim\!\nearrow_{n\rightarrow +\infty}inf_{\theta\in\Theta_0^{f_n}}(1-G_\theta(K(t)))=$
\dli $lim_{\theta<\theta_f\nearrow\theta_{f}}(1-G_\theta(K(t)))=1-G_{\theta_f}(K(t))$ d'apr\`es la d\'ecroissance en $\theta$ de $G_\theta(K(t))$ et la propri\'et\'e 2) de l'\'enonc\'e ($D_s^f=D_s^{\{\theta<\theta_f\}}$).
On obtient bien $Q^\omega_0(\Theta_1^{f})$ quand $t\notin D_i^f\cup D_s^f$ car $Q^\omega_0(\Theta_1^{f})=inf_{\theta\geq\theta_{f}}(1-G_\theta(K(t)))=1-G_{\theta_f}(K(t))$ ; lorsque $t\in(D_i^f-D_s^f)-D_i$, on a $Q^\omega_0(\Theta_1^{f})=1$ et il en est de m\^eme de $1-G_{\theta_f}(K(t))$ d'apr\`es la propri\'et\'e 2) d'une statistique essentielle globale (voir la proposition 5.2.2) puisque $\theta_f\in\Theta_0^f$.

\medskip
{\parindent=-10mm II -- Condition n\'ecessaire.}

$Q_1$ et $Q_0$ sont suppos\'es d\'efinir des votes compatibles.
\dli Sur $\Omega-T^{-1}(N)$ on a $Q^\omega_1({\Theta_1^f})=sup_{\theta\in\Theta_1^f}Q^\omega_{\theta}(\{1_f\})$ et $Q^\omega_0({\Theta_1^f})=inf_{\theta\in\Theta_0^f}Q^\omega_{\theta}(\{1_f\})$. Par d\'efinition les votes $Q^\omega_{\theta}(\{1_f\})$ ne d\'ependent que de $t$ ; les propri\'et\'es a), b) et c) de la d\'efinition 5.2.1 sont donc v\'erifi\'ees en dehors d'un n\'egligeable de la forme $T^{-1}(M\cup N)$. Il suffit de d\'emontrer que $G_\theta(K(t))$ v\'erifie les quatre propri\'et\'es de l'\'enonc\'e pour $t\notin M\cup N$.

Soit $\omega$ tel que $T(\omega)=t\notin M\cup N$.
\medskip
1) Soit $\theta_f\in\Theta-\{sup\Theta\}$, $G_\theta(K(t))$ est continue \`a droite en $\theta_f$ si $t\notin D_i^{\{\theta>\theta_f\}}$.

Consid\'erons une suite $\{\theta_n\}_{n\in\IN}$ d\'ecroissant vers $\theta_f$, 
$\theta_n\not=\theta_f$ ; posons $\Theta_1^{f'_n}=\Theta\cap]-\infty,\theta_n]$ et $\Theta_1^{f_n}=\Theta\cap]-\infty,\theta_f]=\Theta_1^f$. La propri\'et\'e c) de la d\'efinition 5.2.1 \'etant v\'erifi\'ee pour la r\'ealisation $\omega$ nous avons : 
\dli $lim_{n\rightarrow +\infty}Q^\omega_1(\Theta_1^{f'_n})=Q^\omega_1(\Theta_1^{f})$.
\dli Reprenons les notations de la partie c)-i) du 1\up{er} cas de I.
\dli Lorsque $t\notin(D_i^f\cup D'_s)=(D_i^{\{\theta>\theta_f\}}\cup D'_s)$, 
nous avons obtenu : $lim_{n\rightarrow +\infty}Q^\omega_1(\Theta_1^{f'_n})=lim_{\theta>\theta_f\!\searrow\theta_f}(1-G_\theta(K(t)))$. La continuit\'e \`a droite en $\theta_f$ de $G_\theta(K(t))$
s'en d\'eduit imm\'ediatement pour $t\notin(D_i^f\cup D_s^f)$, car dans ce cas : $Q^\omega_1(\Theta_1^{f})=sup_{\theta\in\Theta_1^f}Q^\omega_{\theta}(\{1_f\})=sup_{\theta\in\Theta_1^f}(1-G_\theta(K(t)))=1-G_{\theta_f}(K(t))$ ; 
cette continuit\'e est encore valable pour $t\in(D_s^f-D'_s)$ car on a alors 
$Q^\omega_1(\Theta_1^{f})=0$ mais aussi $1-G_{\theta_f}(K(t))=0$ d'apr\`es la propri\'et\'e 3) d'une statistique essentielle globale (voir la proposition 5.2.2).
\dli Il reste le cas $t\in D'_s-D_i^f$ (ce qui suit est valable pour $t\in D'_s$). A partir d'un certain rang $t$ appartient aux  $D_s^{f'_n}$ et d'apr\`es la propri\'et\'e 3) d'une statistique essentielle globale on a $G_{\theta_n}(K(t))=1$ car $\theta_n\in\Theta_1^{f'_n}$ ; ceci d\'emontre la continuit\'e \`a droite en $\theta_f$ puisque $D'_s\subseteq D_s^f$ et la propri\'et\'e 3) pr\'ec\'edente impliquent :  $G_{\theta_f}(K(t))=1$.

\medskip
2) Soit $\theta_f\in\Theta-\{inf\Theta\}$, $G_\theta(K(t))$ est continue \`a gauche en $\theta_f$ si $t\notin D_s^{\{\theta<\theta_f\}}$ 

La d\'emonstration est semblable \`a la pr\'ec\'edente. On consid\`ere cette fois une suite $\{\theta_n\}_{n\in\IN}$ croissant vers $\theta_f$, 
$\theta_n\not=\theta_f$ ; on pose $\Theta_1^{f_n}=\Theta\cap]-\infty,\theta_n[$ et $\Theta_1^{f'_n}=\Theta\cap]-\infty,\theta_f[=\Theta_1^f$. Les votes $Q_0$ \'etant compatibles en $\omega$ on a : $lim_{n\rightarrow +\infty}Q^\omega_0(\Theta_1^{f_n})=Q^\omega_0(\Theta_1^{f})$.
\dli Reprenons les notations de la partie c)-ii) du 1\up{er} cas de I.
\dli Lorsque $t\notin(D_i\cup D_s^f)=(D_i\cup D_s^{\{\theta<\theta_f\}})$, 
la partie c)-ii) du 2\up{\`eme} cas implique : $lim_{n\rightarrow +\infty}Q^\omega_0(\Theta_1^{f_n})=lim_{\theta<\theta_f\!\nearrow\theta_f}(1-G_\theta(K(t)))$. La continuit\'e \`a gauche en $\theta_f$ de $G_\theta(K(t))$
s'en d\'eduit imm\'ediatement pour $t\notin(D_i^f\cup D_s^f)$, car dans ce cas : $Q^\omega_0(\Theta_1^{f})=inf_{\theta\in\Theta_0^f}Q^\omega_{\theta}(\{1_f\})=inf_{\theta\geq\theta_f}(1-G_\theta(K(t)))=1-G_{\theta_f}(K(t))$ ; 
cette continuit\'e est encore valable pour $t\in(D_i^f-D_i)-D_s^f$ car on a alors $Q^\omega_0(\Theta_1^{f})=1$ mais aussi $1-G_{\theta_f}(K(t))=1$ d'apr\`es la propri\'et\'e 2) d'une statistique essentielle globale (voir la proposition 5.2.2).
\dli Il reste le cas $t\in D_i-D_s^f$ (ce qui suit est valable pour $t\in D_i$). A partir d'un certain rang $t$ appartient aux  $D_i^{f_n}$ et d'apr\`es la propri\'et\'e 2) d'une statistique essentielle globale on a $G_{\theta_n}(K(t))=0$ car $\theta_n\in\Theta_0^{f_n}$ ; ceci d\'emontre la continuit\'e \`a gauche en $\theta_f$ puisque $D_i\subseteq D_i^f$ et la propri\'et\'e 2) pr\'ec\'edente impliquent :  $G_{\theta_f}(K(t))=0$.

\medskip
3) si $inf\Theta\notin\Theta$ on a $lim_{\theta\rightarrow inf\Theta}G_\theta(K(t))=1$

Consid\'erons une suite $\{\theta_n\}_{n\in\IN}$ d\'ecroissant vers $inf\Theta$ et posons $\Theta_1^{f_n}=\Theta\cap]-\infty,\theta_n]$. On a  $\cap_{n\in\IN}\Theta_1^{f_n}=\emptyset$ et d'apr\`es la propri\'et\'e b) de la d\'efinition 5.2.1 : 
$lim_{n\rightarrow +\infty}Q^\omega_1(\Theta_1^{f_n})=0$.
\dli Reprenons les notations de la partie b-2) du 1\up{er} cas.
\dli Lorsque $t\notin(D_i^-\cup D_s^+)$ nous avons obtenu $lim_{n\rightarrow +\infty}Q^\omega_1(\Theta_1^{f_n})=lim_{\theta\rightarrow inf\Theta}(1-G_\theta(K(t)))$, ce qui implique bien l'\'egalit\'e recherch\'ee.
\dli De plus nous avons vu que $t$ n'appartient pas \`a $D_i^-\subseteq N$. Quant au cas $t\in D_s^+$, il est \'evident puisqu'\`a partir d'un certain rang  $t\in D_s^{f_n}$ et $G_{\theta_n}(K(t))=1$ d'apr\'es la propri\'et\'e 3) d'une statistique essentielle globale(voir la proposition5.2.2).

\medskip
4) si $sup\Theta\notin\Theta$ on a $lim_{\theta\rightarrow sup\Theta}G_\theta(K(t))=0$.

Consid\'erons une suite $\{\theta_n\}_{n\in\IN}$ croissant vers $sup\Theta$ et posons $\Theta_1^{f_n}=\Theta\cap]-\infty,\theta_n[$. On a  $\cup_{n\in\IN}\Theta_1^{f_n}=\Theta$ et d'apr\`es la propri\'et\'e b) de la d\'efinition 5.2.1 : 
$lim_{n\rightarrow +\infty}Q^\omega_1(\Theta_1^{f_n})=1$.
\dli Reprenons les notations de la partie b-1) du 1\up{er} cas.
\dli Lorsque $t\notin(D_i^+\cup D_s^-)$ nous avons obtenu $lim_{n\rightarrow +\infty}Q^\omega_1(\Theta_1^{f_n})=lim_{\theta\rightarrow sup\Theta}(1-G_\theta(K(t)))$, ce qui implique bien l'\'egalit\'e recherch\'ee.
\dli De plus, $t$ n'appartient pas \`a $D_s^-\subseteq N$. Quant au cas $t\in D_i^+$, il est \'evident puisqu'\`a partir d'un certain rang  $t\in D_i^{f_n}$ et $G_{\theta_n}(K(t))=0$ d'apr\`es la propri\'et\'e 2) d'une statistique essentielle globale ($\theta_n\in\Theta_0^{f_n}$).

\medskip
{\parindent=-10mm III - $Q_1$ et $Q_0$ ont presque partout le m\^eme prolongement.}

Sous les conditions de l'\'enonc\'e, $Q_1$ et $Q_0$ d\'efinissent des votes compatibles ; ils sont donc presque partout prolongeables \`a la $\sigma$-alg\`ebre engendr\'ee par $\{\Theta_1^f\}_{f\in{\cal F}}$ qui est la tribu bor\'elienne de $\Theta$. 
\dli Consid\'erons une r\'ealisation $\omega$ v\'erifiant les quatre propri\'et\'es de l'\'enonc\'e et n'appartenant pas a $T^{-1}(N\cup N_c)$, c'est-\`a-dire 
$T(\omega)=t\notin (N\cup N_c)$ (voir le d\'ebut de I). Il suffit de montrer que $Q_1^\omega(\Theta_1^f)=Q_0^\omega(\Theta_1^f)$ pour 
$\Theta_1^f=\Theta\cap]-\infty,\theta_f]$, avec 
$\theta_f\in\Theta-\{sup\Theta\}$.
\dli Lorsque $t\notin (D_i^f\cup D_s^f)$, 
$Q^\omega_{\theta}(\{1_f\})=1-G_\theta(K(t))$ donc 
$Q^\omega_1({\Theta_1^f})=sup_{\theta\in\Theta_1^f}Q^\omega_{\theta}(\{1_f\})=1-G_{\theta_f}(K(t))$ et 
$Q^\omega_0({\Theta_1^f})=inf_{\theta\in\Theta_0^f}Q^\omega_{\theta}(\{1_f\})=lim_{\theta>\theta_f\searrow\theta_f}(1-G_{\theta_f}(K(t)))$ ; 
la propri\'et\'e 1) de l'\'enonc\'e entra\^{\i}ne : 
\dli $Q_1^\omega(\Theta_1^f)=Q_0^\omega(\Theta_1^f)$.
\dli On a aussi cette \'egalit\'e lorsque $t\in (D_i^f-D_s^f)$ (resp. $t\in D_s^f$) puisque pour tout $\theta$ de $\Theta$, $Q^\omega_{\theta}(\{1_f\})$ est \'egal \`a $1$ (resp. $0$).
\medskip\centerline{\hbox to 3cm{\bf \hrulefill}}\par}

La plupart des mod\`eles statistiques \`a rapport de vraisemblance mono\-tone classiquement \'etudi\'es v\'erifient les conditions de cette proposition.
Pour que les conditions 1) et 2) soient r\'ealis\'ees il suffit que l'application, qui associe \`a $\theta$ l'\'el\'ement de $L_1(\Omega;{\cal A},\mu)$ correspondant \`a la densit\'e $p_\theta$, soit continue pour la topologie faible de $L_1$ (en effet $G_\theta(K(t))=\int_\Omega f_t.p_\theta d\mu$ avec $f_t$ born\'ee). Dans ces mod\`eles, la probabilit\'e ainsi obtenue sur $\Theta$ peut aider \`a choisir entre des hypoth\`eses non expertisables. C'est ce que nous \'etudierons dans les paragraphes suivants.

Il y a bien d'autres mani\`eres d'obtenir des votes compatibles sur l'ensemble des hypoth\`eses unilat\'erales $(\{\Theta_1^f,\Theta_0^f\})_{f\in{\cal F}}$. Consid\'erons pour chacun de ces probl\`emes de d\'ecision un vote pond\'er\'e par la probabilit\'e $\Lambda^f$ d\'efinie sur l'intervalle $\Theta$ muni de la tribu des bor\'eliens. Posons $Q^\omega(\Theta_1^f)=Q^\omega_{\Lambda^f}(\{1_f\})=\int_\Theta Q^\omega_{\theta}(\{1_f\})d\Lambda^f(\theta)$, 
$Q^\omega_{\theta}(\{1_f\})$ \'etant la fr\'equence des votes en faveur de $\Theta_1^f$ sous $P_\theta$, pour la r\'ealisation $\omega$ (cette int\'egrale a un sens car $Q^\omega_{\theta}(\{1_f\})$ est une fonction croissante de $\theta$ \`a valeurs dans $[0,1]$ (voir le a) de la d\'emonstration 5.1.1)).
Si l'on veut que la famille des pond\'erations $(\Lambda^f)_{f\in{\cal F}}$  d\'efinisse des votes compatibles, il faut que l'application $Q$, de $\Omega\times\{\Theta_1^f\}_{f\in{\cal F}}$ dans $[0,1]$, v\'erifie les propri\'et\'es a), b) et c) de la d\'efinition 5.2.1 pour presque toutes les r\'ealisations $\omega$ . 
$Q$ a un comportement moins simple que les votes $Q_1$ et $Q_0$ pr\'ec\'edents, il d\'epend \`a la fois du mod\`ele et des pond\'erations. Nous allons d\'efinir des conditions peu contraignantes qui donnent \`a $Q$ les bonnes propri\'et\'es.

\medskip
{\bf Proposition 5.2.4}
\medskip
\medskip
\moveleft 10.4pt\hbox{\vrule\kern 10pt\vbox{\defpro

Soit $(\Omega ,{\cal A},(p_\theta.\mu)_{\theta\in\Theta})$ un mod\`ele
statistique \`a rapport de vraisem\-blance monotone pour une statistique r\'eelle $T$, $\Theta$ \'etant un intervalle de $\IR$ muni de l'ordre usuel et de la tribu des bor\'eliens ${\cal B}$.
$K(T)$ est une statistique essentielle globale et $G_\theta$ d\'esigne sa fonction de r\'epartition moyenne sous $P_\theta$. On suppose avoir les deux  propri\'et\'es suivantes :\quad  i) si $\theta_f\in\Theta-\{sup\Theta\}$,\quad  $D_s^{\theta_f}=\cup_{\theta>\theta_f}D_s^{\theta}$, 
\dli ii) si $\theta_f\in\Theta-\{inf\Theta\}$,\quad   $D_i^{\theta_f}=\cup_{\theta<\theta_f}D_i^{\theta}$. 
Notons $(\{\Theta_1^f,\Theta_0^f\})_{f\in{\cal F}}$ l'ensemble des hypoth\`eses unilat\'erales. ${\cal F}$ est totalement ordonn\'e par la relation :
$f\leq f'$ $\Leftrightarrow$ $\Theta_1^f\subseteq\Theta_1^{f'}$, on le munit de la topologie d\'efinie par cet ordre.  
\dli Soit $(\Lambda^f)_{f\in{\cal F}}$ une famille de probabilit\'es sur $(\Theta,{\cal B})$, on consid\`ere, pour tout $t\in\IR$, la fonction d\'efinie sur ${\cal F}$ par :
$H_t(f)=\int_\Theta(1-G_\theta(K(t)))d\Lambda^f(\theta)$ ; pour presque tout $t$, on suppose que $H_t$ est croissante, continue et v\'erifie : 
\dli 1) si ${\cal F}$ n'a pas de borne inf\'erieure ($inf\Theta\notin\Theta$) on a $inf_{f\in{\cal F}}H_t(f)=0$
\dli 2) si ${\cal F}$ n'a pas de borne sup\'erieure ($sup\Theta\notin\Theta$) on a $sup_{f\in{\cal F}}H_t(f)=1$.

Si on associe \`a chaque hypoth\`ese unilat\'erale $\Theta_1^f$ le vote pond\'er\'e $Q_{\Lambda^f}$, on obtient alors des votes compatibles.
}}\medskip

\medskip
{\leftskip=15mm \dli {\bf D\'emonstration}
\medskip
Posons $N=\{t\in\IR\, ;\, \forall\omega\in T^{-1}(t)\quad\forall\theta\in\Theta\quad p_\theta(\omega)=0\}$ et notons simplement $D_i^f$ et $D_s^f$ les demi-droites inf\'erieures et sup\'erieures $D_i^{\Theta_0^f}$ et $D_s^{\Theta_1^f}$. D'apr\`es les propositions 4.3.1 et 5.2.2, pour une r\'ealisation $\omega\notin T^{-1}(N)$ et des hypoth\`eses unilat\'erales $\{\Theta_1^f,\Theta_0^f\}$, le vote   $Q^\omega_{\theta}$ est d\'efini par : 
$$Q^\omega_\theta(\{1_f\})=\left\{\matrix{
1\hfill & si & T(\omega)\in D_i^f-D_s^f\hfill\cr
1-G_\theta(K(T(\omega)))\hfill & si & T(\omega)\in\IR\!-(D_i^f\cup D_s^f)\hfill \cr
0\hfill & si & T(\omega)\in D_s^f\hfill \cr
}\right. $$

Nous consid\'erons la famille de votes d\'efinie par : 
$Q^\omega({\Theta_1^f})=Q^\omega_{\Lambda^f}(\{1_f\})=\int_\Theta Q^\omega_\theta(\{1_f\})d\Lambda^f(\theta)$ ; cette int\'egrale (resp. $H_t(f)$) est bien d\'efinie puisque $Q^\omega_\theta(\{1_f\})$ (resp. $1-G_\theta(K(t))$) est une fonction croissante en $\theta$, \`a valeurs dans $[0,1]$ (voir le a) de la d\'emonstration 5.1.1).
\dli Nous devons d\'emontrer que $Q$ v\'erifie les propri\'et\'es a), b) et c) de la d\'efinition 5.2.1 pour presque tout $\omega$  de $\Omega-T^{-1}(N)$. On va suivre les \'etapes de la partie I de la d\'emonstration 5.2.3.

Notons $N_c$ l'ensemble des $t\in\IR$ pour lesquels $H_t$ ne v\'erifie pas les propri\'et\'es suppos\'ees dans l'\'enonc\'e. $T^{-1}(N_c)$ \'etant n\'egligeable, il suffit de d\'emontrer les propri\'et\'es a),b) et c) sur $\Omega-T^{-1}(N\cup N_c)$.

Soit $\omega$ tel que $T(\omega)=t\notin (N\cup N_c)$. 

\medskip
{\parindent=-10mm a) --- Si $\Theta_1^f\subset\Theta_1^{f'}$ on a $Q^\omega({\Theta_1^f})\leq Q^\omega(\Theta_1^{f'})$.}

Lorsque $t\in D_i^{f'}-D_s^{f'}$ (resp. $t\in D_s^f$) l'in\'egalit\'e recherch\'ee est \'evidente car on a, pour tout $\theta$ de $\Theta$, $Q^\omega_\theta(\{1_{f'}\})=1$ (resp. $Q^\omega_\theta(\{1_f\})=0$) donc
$Q^\omega(\Theta_1^{f'})=1$ (resp. $Q^\omega(\Theta_1^{f})=0$).
\dli Il nous reste le cas $t\notin(D_i^{f'}-D_s^{f'})\cup D_s^f\not=\emptyset$. Comme $\Theta_1^f\subset\Theta_1^{f'}$ implique $D_i^f\subseteq D_i^{f'}$ et $D_s^f\supseteq D_s^{f'}$ (voir le lemme 1 de l'annexe IV), $t$ n'appartient pas non plus \`a $D_i^{f}-D_s^{f}\subseteq D_i^{f'}-D_s^{f'}$ et \`a $D_s^{f'}$. On a, pour tout $\theta$ de $\Theta$, 
$Q^\omega_\theta(\{1_f\})=Q^\omega_\theta(\{1_{f'}\})=1-G_\theta(K(t))$ donc 
$Q^\omega({\Theta_1^f})=H_t(f)$ et $Q^\omega(\Theta_1^{f'})=H_t(f')$ ; la croissance de $H_t$ entra\^{\i}ne l'in\'egalit\'e recherch\'ee.
\medskip
{\parindent=-10mm b-1) - Si $(\Theta_1^{f_n})_{n\in\IN}$ cro\^{\i}t vers $\Theta$, la suite  $(Q^\omega(\Theta_1^{f_n}))_{n\in\IN}$ cro\^{\i}t vers $1$.}

$\Theta_1^{f_n}=\Theta\cap]-\infty,\theta_n)$ n'\'etant jamais \'egal \`a $\Theta$, l'existence d'une suite $(\Theta_1^{f_n})$ croissant vers $\Theta$
n'est possible que si $sup\Theta\notin\Theta$. Il est \'equivalent de dire que 
${\cal F}$ n'a pas de borne sup\'erieure.
\dli Posons $D_i^+=\cup_{n\in\IN}D_i^{f_n}$ et $D_s^-=\cap_{n\in\IN}D_s^{f_n}$ ; la suite $(\Theta_1^{f_n})$ \'etant croissante les suites $(D_i^{f_n})$ et $(D_s^{f_n})$ sont respectivement croissante et d\'ecroissante (voir le lemme 1 de l'annexe IV) ; la suite 
$(D_i^{f_n}-D_s^{f_n})$ cro\^{\i}t vers $D_i^+-D_s^-$.
\dli Lorsque $t\in D_i^+-D_s^-$ on a la propri\'et\'e recherch\'ee puisqu'\`a partir d'un certain rang $t\in(D_i^{f_n}-D_s^{f_n})$ donc $Q^\omega(\Theta_1^{f_n})=1$.
\dli Le cas $t\in D_s^-$ est impossible car $D_s^-\subseteq N$ ; en effet, la convergence de $(\Theta_1^{f_n})$ vers $\Theta$ et le lemme 1 de l'annexe IV entra\^{\i}nent $D_s^-=\cap_{\theta\in\Theta}D_s^\theta$.
\dli Il reste le cas $t\notin D_i^+\cup D_s^-$. A partir d'un certain rang $t\notin D_i^{f_n}\cup D_s^{f_n}$, on a alors $Q^\omega(\Theta_1^{f_n})=H_t(f_n)$. La propri\'et\'e 2) de l'\'enonc\'e et la croissance de $H_t$ inpliquent bien $lim_{n\rightarrow +\infty}H_t(f_n)=1$.
\medskip
{\parindent=-10mm b-2) - Si $(\Theta_1^{f_n})_{n\in\IN}$ d\'ecro\^{\i}t vers le vide, la suite  $(Q^\omega(\Theta_1^{f_n}))_{n\in\IN}$ d\'ecro\^{\i}t vers $0$.}

La d\'emonstration est semblable \`a celle de b-1).
\dli $\Theta_1^{f_n}=\Theta\cap]-\infty,\theta_n)$ n'\'etant jamais vide, l'existence d'une suite $(\Theta_1^{f_n})$ d\'ecroissant vers le vide 
n'est possible que si $inf\Theta\notin\Theta$, c'est-\`a-dire que ${\cal F}$ n'a pas de borne inf\'erieure. On pose alors 
$D_i^-=\cap\searrow_{n\in\IN}D_i^{f_n}$ et $D_s^+=\cup\nearrow_{n\in\IN}D_s^{f_n}$. 
\dli Lorsque $t\in D_s^+$, la propri\'et\'e recherch\'ee est \'evidente puisqu'\`a partir d'un certain rang $t\in D_s^{f_n}$ donc $Q^\omega(\Theta_1^{f_n})=0$.
\dli Le cas $t\in D_i^-$ est impossible car $D_i^-\subseteq N$ ; en effet, la convergence de $(\Theta_0^{f_n})$ vers $\Theta$ et le lemme 1 de l'annexe IV entra\^{\i}nent $D_i^-=\cap_{\theta\in\Theta}D_i^\theta$.
\dli Il reste le cas $t\in\IR\!-(D_i^-\cup D_s^+)\not=\emptyset$. A partir d'un certain rang $t\notin D_i^{f_n}\cup D_s^{f_n}$, on a alors $Q^\omega(\Theta_1^{f_n})=H_t(f_n)$. La propri\'et\'e 1) de l'\'enonc\'e et la croissance de $H_t$ impliquent bien $lim_{n\rightarrow +\infty}H_t(f_n)=0$.
\medskip
{\parindent=-10mm c) --- Si $(\Theta_1^{f_n})_{n\in\IN}$ et $(\Theta_1^{f'_n})_{n\in\IN}$ sont deux suites ayant m\^eme limite, res\-pectivement croissante et d\'ecroissante, les suites $(Q^\omega(\Theta_1^{f_n}))_{n\in\IN}$ et $(Q^\omega(\Theta_1^{f'_n}))_{n\in\IN}$ ont m\^eme limite.}

La limite commune de $(\Theta_1^{f_n})$ et $(\Theta_1^{f'_n})$ ne peut pas \^etre vide ou \'egale \`a $\Theta$, elle est donc de la forme $\Theta_1^f=\Theta\cap]-\infty,\theta_f)$ avec $\theta_f\in\Theta$. Lorsque $\theta_f$ appartient \`a $\Theta_1^f$ (resp. $\Theta_0^f=\Theta-\Theta_1^f$), les  $\Theta_1^{f_n}$ (resp. $\Theta_1^{f'_n}$) sont \'egaux \`a $\Theta_1^f$ \`a partir d'un certain rang.
Etudions s\'epar\'ement ces deux cas.

{\parindent=-5mm i) - $\Theta_1^f=\Theta\cap]-\infty,\theta_f]$ ($\theta_f\in\Theta-\{sup\Theta\}$).} 

On a   $lim_{n\rightarrow +\infty}Q^\omega(\Theta_1^{f_n})=Q^\omega(\Theta_1^{f})$,\quad $D_i^f=D_i^{\Theta_0^f}=\cap\!\searrow_{n\in\IN}D_i^{f'_n}$,\quad 
$D_s^f=D_s^{\theta^f}\supseteq\cup\!\nearrow_{n\in\IN}D_s^{f'_n}$
et la suite $({f'_n})_{n\in\IN}$ d\'ecro\^{\i}t vers $f$. D'apr\`es la propri\'et\'e i) de l'\'enonc\'e on a m\^eme $D_s^f=\cup\!\nearrow_{n\in\IN}D_s^{f'_n}$. 
\dli L'\'egalit\'e des deux limites est \'evidente lorsque $t\in D_i^f-D_s^f$ ; $t$ appartenant aussi \`a $D_i^{f'_n}-D_s^{f'_n}$ on a \`a la fois  $Q^\omega(\Theta_1^{f})=1$ et $Q^\omega(\Theta_1^{f'_n})=1$.
\dli Lorsque $t\in D_s^f$ on a $Q^\omega(\Theta_1^{f})=0$ mais aussi $lim_{n\rightarrow +\infty}Q^\omega(\Theta_1^{f'_n})=0$ car $t$ appartient \`a $D_s^{f'_n}$ \`a partir d'un certain rang.
\dli Il nous reste le cas : $t\notin D_i^f\cup D_s^f$ ; \`a partir d'un certain rang $t$ n'appartient pas \`a $D_i^{f'_n}\cup D_s^{f'_n}$, ce qui entra\^{\i}ne  $Q^\omega(\Theta_1^{f'_n})=H_t(f'_n)$ ; $H_t$ \'etant continue on a 
$lim_{n\rightarrow +\infty}Q^\omega(\Theta_1^{f'_n})=H_t(f)$ ; l'\'egalit\'e entre $H_t(f)$  et $Q^\omega(\Theta_1^{f})$ est triviale puisque $t\notin D_i^f\cup D_s^f$. 

{\parindent=-5mm ii) - $\Theta_1^f=\Theta\cap]-\infty,\theta_f[$ ($\theta_f\in\Theta-\{inf\Theta\}$).}

On a $lim_{n\rightarrow +\infty}Q^\omega(\Theta_1^{f'_n})=Q^\omega(\Theta_1^{f})$,\quad 
$D_i^f=D_i^{\theta^f}\supseteq\cup\!\nearrow_{n\in\IN}D_i^{f_n}$,\quad 
$D_s^f=D_s^{\Theta_1^f}=\cap\!\searrow_{n\in\IN}D_s^{f_n}$ et la suite $({f_n})_{n\in\IN}$ cro\^{\i}t vers $f$. D'apr\`es la propri\'et\'e ii) de l'\'enonc\'e on a m\^eme $D_i^f=\cup\!\nearrow_{n\in\IN}D_i^{f_n}$.
\dli L'\'egalit\'e des deux limites est \'evidente lorsque $t\in D_s^f$ ; 
il appartient aussi aux $D_s^{f_n}$, on a donc $Q^\omega(\Theta_1^{f})=0$ et $Q^\omega(\Theta_1^{f_n})=0$. 
\dli Lorsque $t\in D_i^f-D_s^f$, on a $Q^\omega(\Theta_1^{f})=1$ ; comme \`a partir d'un certain rang, $t\in D_i^{f_n}-D_s^{f_n}$, on a aussi 
$lim_{n\rightarrow +\infty}Q^\omega(\Theta_1^{f_n})=1$. 
\dli Il nous reste le cas : $t\notin D_i^f\cup D_s^f$ ; \`a partir d'un certain rang $t$ n'appartient pas \`a $D_i^{f_n}\cup D_s^{f_n}$, on a alors  $Q^\omega(\Theta_1^{f_n})=H_t(f_n)$ et la continuit\'e de $H_t$ implique $lim_{n\rightarrow +\infty}Q^\omega(\Theta_1^{f_n})=H_t(f)$ ; cette quantit\'e est \'evidemment \'egale \`a $Q^\omega(\Theta_1^{f})$ puisque $t\notin D_i^f\cup D_s^f$. 

\medskip\centerline{\hbox to 3cm{\bf \hrulefill}}\par}

\dli Les conditions d'application de cette proposition ne sont pas tr\`es contraignantes. 
\dli Les propri\'et\'es i) et ii) interdisent certaines discontinuit\'es des demi-droites $D_s^\theta$ et $D_i^\theta$. La plupart du temps, elles ne varient pas ou varient contin\^ument : $\cap_{\theta<\theta_f}D_s^\theta=D_s^{\theta_f}=\cup_{\theta>\theta_f}D_s^\theta$ et $\cup_{\theta<\theta_f}D_i^\theta=D_i^{\theta_f}=\cap_{\theta>\theta_f}D_i^\theta$. Dans cette situation les votes de la proposition 5.2.3 sont d'ailleurs neutres pour toutes les hypoth\`eses unilat\'erales (voir la proposition 4.3.2).
\dli On a suppos\'e les fonctions $H_t(f)=\int_\Theta(1-G_\theta(K(t)))d\Lambda^f(\theta)$ croissantes sur ${\cal F}$ pour presque tout $t$. Comme $1-G_\theta(K(t))$ est une fonction croissante en $\theta$, pour avoir la croissance de $H_t$ il suffit que la probabilit\'e $\Lambda^f$ soit d'autant plus grande que $f$ est grand (voir [DacD] p. 79-80), c'est-\`a-dire que les fonctions de r\'epartitions $F_f$ d\'ecroissent lorsque $f$ cro\^{\i}t. Les pond\'erations $\Lambda^f$ chargent alors d'autant plus les grandes valeurs de $\theta$ que $\Theta_1^f$ est grand, ceci semble assez ``naturel". Si $G_\theta(K(t))$ est continue en $\theta$, ce qui est souvent le cas, il suffit d'ajouter la continuit\'e  des $\Lambda^f$ au sens de la convergence \'etroite ($f_n\rightarrow f\ \Rightarrow\ \Lambda^{f_n}\relmont{\rightarrow}{et.}\Lambda^f$) pour avoir la continuit\'e de $H_t$. 
\dli En ce qui concerne les propri\'et\'es 1) et 2) de $H_t$, elles n'imposent aucune contrainte lorsque $\Theta$ est ferm\'e. Lorsque $\Theta$ est ouvert \`a gauche (resp. droite) la condition 3) (resp. 4)) de la proposition 5.2.3 entra\^{\i}ne la propri\'et\'e 1) (resp. 2)) si les fonctions de r\'epartition $F_f$ tendent vers $1$ (resp. $0$) quand $\Theta_1^f$ d\'ecro\^{\i}t vers le vide (resp. cro\^{\i}t vers $\Theta$). 

Consid\'erons le cas classique d'un n-\'echantillon $(X_1,X_2,...,X_n)$ d'une loi normale, $N(\theta,\sigma^2)$, de moyenne inconnue $\theta\in\IR$ et de variance connue $\sigma^2$ ($\sigma>0$). Nous avons vu \`a la fin du paragraphe 5.1 que ce mod\`ele d'\'echantillonnage est \`a rapport de vraisemblance monotone par rapport \`a $\overline{X}$, qui est une statistique essentielle pour le choix entre deux hypoth\`eses unilat\'erales quelconques. $\overline{X}$ est donc une statistique essentielle globale, sa fonction de r\'epartition moyenne est donn\'ee par : 
$G_{\theta}(\overline{x})=F({\sqrt{n}\over\sigma}(\overline{x}-\theta))$ (F est la fonction de r\'epartition de la loi normale $N(0,1)$).
\dli $G_\theta(\overline{x})$ v\'erifie bien les quatre conditions de la proposition 5.2.3 puisque $F({\sqrt{n}\over\sigma}(\overline{x}-\theta))$ est continu en $\theta$, $lim_{\theta\rightarrow -\infty}F({\sqrt{n}\over\sigma}(\overline{x}-\theta))=1$ et 
$lim_{\theta\rightarrow +\infty}F({\sqrt{n}\over\sigma}(\overline{x}-\theta))=0$. Quelles que soient les hypoth\`eses unilat\'erales $\{\Theta_1^f=]-\infty,\theta_f),\Theta_0^f\}$ les deux votes les plus favorables sous $\Theta_1^f$ et sous $\Theta_0^f$ sont identiques ; pour toute r\'ealisation $\omega$ on a : 
$Q^\omega_{\Theta_1^f}(\{1_f\})=Q^\omega_{\Theta_0^f}(\{1_f\})=Q^\omega_{\theta_f}(\{1_f\})= 1-F({\sqrt{n}\over\sigma}(\overline{X}(\omega)-\theta_f))$. Ils d\'efinissent des votes compatibles sur l'ensemble des hypoth\`eses unilat\'erales : 
$Q^\omega(\Theta_1^f)= 1-F({\sqrt{n}\over\sigma}(\overline{X}(\omega)-\theta_f))$. 
Le prolongement en une probabilit\'e unique sur les bor\'eliens de $\Theta=\IR$ v\'erifie $Q^\omega(]-\infty,\theta_f[)= 1-F({\sqrt{n}\over\sigma}(\overline{X}(\omega)-\theta_f))=F({\sqrt{n}\over\sigma}(\theta_f-\overline{X}(\omega)))$. C'est donc la loi normale 
$N(\overline{x},{\sigma^2\over n})$. En th\'eorie bay\'esienne, ce r\'esultat est la loi a post\'eriori associ\'ee \`a la mesure de Lebesgue sur $\IR$ (cf. [Ber.] p. 132).
\dli Consid\'erons maintenant des votes pond\'er\'es. Comme dans le paragraphe pr\'ec\'edent nous allons prendre pour chaque probl\`eme unilat\'eral 
$\{\Theta_1^f=]-\infty,\theta_f),\Theta_0^f\}$ une pond\'eration 
$\Lambda^f=N(\theta_f,c^2)$ avec $c>0$. Les deux hypoth\`eses $\Theta_1^f$ et $\Theta_0^f$  sont ainsi trait\'ees de fa\c con sym\'etrique. On a :
\dli $Q^\omega(\Theta_1^f)=Q^\omega_{\Lambda^f}(\{1_f\})=1-\int_\Theta F({\sqrt{n}\over\sigma}(\overline{X}(\omega)-\theta))\,{1\over\sqrt{2\pi}\,c} exp(-{(\theta-\theta_f)^2\over 2c^2})\,d\theta$.
\dli En posant $\theta'=(\theta-\theta_f)$, les calculs de la fin de 5.1 donnent :
\dli $Q^\omega(\Theta_1^f)=1-F({\sqrt{n}\over\sigma}(\overline{X}(\omega)-\theta_f) \sqrt{\sigma^2\over nc^2+\sigma^2})=F(\sqrt{n\over nc^2+\sigma^2}(\theta_f-\overline{X}(\omega)))$.
\dli D'apr\`es la proposition 5.2.4, $Q^\omega$ se prolonge en une probabilit\'e sur $\Theta=\IR$ qui est donc la loi $N(\overline{x},{\sigma^2\over n}+c^2)$. Lorsque $c$ tend vers $0$ on obtient le prolongement pr\'ec\'edent, celui des votes les plus favorables. Plus l'\'ecart type $c$ est grand, plus les votes extr\^emes prennent de l'importance et la probabilit\'e sur $\Theta$ a une dispersion plus grande autour de $\overline{x}$, ce qui pr\'esente peu d'int\'er\^et. Pour avoir une famille de pond\'erations int\'eressante il faut partir d'une information a priori. Nous le ferons dans le paragraphe 5.4. Les votes que nous venons de d\'efinir sont neutres, on a $E_{\theta_f}(Q^\omega(\Theta_1^f))={1\over 2}$ pour tout $\theta_f$. Dans le premier cas c'est une cons\'equence directe de la proposition 4.3.2. Pour les votes pond\'er\'es on a : 
\cleartabs
\+ $E_{\theta_f}(Q^\omega(\Theta_1^f))$&=&$\int_{\IR}F(\sqrt{n\over nc^2+\sigma^2}(\theta_f-y))\,{\sqrt{n}\over\sqrt{2\pi}\sigma}exp(-{n\over 2\sigma^2}(y-\theta_f)^2)\,dy$\cr
\+ &=&$\int_{\IR}F(\sqrt{n\over nc^2+\sigma^2}z)\,{\sqrt{n}\over\sqrt{2\pi}\sigma}exp(-{n\over 2\sigma^2}z^2)\,dz$\cr
\dli il suffit d'int\'egrer s\'epar\'ement sur $]-\infty,0[$ et $]0,+\infty[$ en utilisant l'\'egalit\'e $F(\sqrt{n\over nc^2+\sigma^2}(-z))=1-F(\sqrt{n\over nc^2+\sigma^2}z)$.

\vfill\eject

{\parindent=-5mm 5.3 HYPOTH\`ESES BILAT\'ERALES.}
\nobreak
\medskip
$(\Omega ,{\cal A},(p_\theta.\mu)_{\theta\in\Theta\subseteq\IR})$ d\'esigne un mod\`ele statistique \`a rapport de vraisem\-blance monotone par rapport \`a la statistique $T$. Consid\'erons deux couples d'hypoth\`eses unilat\'erales : $\{\Theta'_1 , \Theta'_0\}$ et $\{\Theta''_1 , \Theta''_0\}$, $\Theta'_1\subset\Theta''_1$ (voir la d\'efinition 5.1.2). Ils d\'efinissent une partition ordonn\'ee de $\Theta$ en trois intervalles non vides :  $\Theta'_1$, $\Theta''_1-\Theta'_1$ et $\Theta''_0$. On \'etudie des hypoth\`eses bilat\'erales lorsque l'appartenance de $\theta$ \`a $\Theta'_1$ et celle de $\theta$ \`a $\Theta''_0$  sont consid\'er\'ees comme \'equivalentes pour l'interpr\'etation. On casse ainsi la structure d'ordre sur $\Theta$ et on \'etudie les hypoth\`eses : $\Theta_1=\Theta'_1\cup\Theta''_0$ et $\Theta_0=\Theta'_0-\Theta''_0=\Theta''_1-\Theta'_1$. Cette cassure entra\^{\i}ne des cons\'equences lourdes sur les possibilit\'es d'expertise lorsque les deux parties de $\Theta_1$ ne peuvent pas s'analyser s\'epar\'ement.

\medskip
{\bf D\'efinition 5.3.1}
\medskip
\medskip
\moveleft 10.4pt\hbox{\vrule\kern 10pt\vbox{\defpro

$(\Omega ,{\cal A},(p_\theta.\mu)_{\theta\in\Theta\subseteq\IR})$ est un mod\`ele statistique \`a rapport de vraisem\-blance monotone par rapport \`a la statistique $T$. 

Soient $\{\Theta_1 , \Theta_0\}$ des hypoth\`eses bilat\'erales bas\'ees sur les hypoth\`eses unilat\'erales $\{\Theta'_1 , \Theta'_0\}$ et $\{\Theta''_1 , \Theta''_0\}$, $\Theta'_1\subset\Theta''_1$. 
\dli Ces hypoth\`eses sont dites bilat\'erales impropres si pour tout $\theta_0$ de $\Theta_0$ elles v\'erifient $P_{\theta_0}(T^{-1}(\IR-D_s^{\Theta'_1}))=0$ ou $P_{\theta_0}(T^{-1}(\IR-D_i^{\Theta''_0}))=0$. Elles sont dites bilat\'erales pures si aucune sous hypoth\`ese de la forme $\{\Theta_1 , \Theta_0^h\subseteq\Theta_0\}$ n'est impropre ($\forall \theta_0\in\Theta_0$ \quad $P_{\theta_0}(T^{-1}(\IR-D_s^{\Theta'_1}))>0$ et $P_{\theta_0}(T^{-1}(\IR-D_i^{\Theta''_0}))>0$).
}}\medskip

Le cas des hypoth\`eses bilat\'erales impropres revient \`a traiter deux probl\`emes de d\'ecision unilat\'eraux distincts. 
\dli Posons $\Theta_0^i=\{\theta\in\Theta_0\ ,\ P_\theta(T^{-1}(\IR-D_s^{\Theta'_1}))>0\}$ et \dli $\Theta_0^s=\{\theta\in\Theta_0\ ,\ P_\theta(T^{-1}(\IR-D_i^{\Theta''_0}))>0\}$, le mod\`ele \'etant \`a rapport de vraisemblance monotone on a : $\Theta'_1<\Theta_0^i<\Theta_0-\Theta_0^i-\Theta_0^s<\Theta_0^s<\Theta''_0$ (parmi ces trois nouveaux sous espaces de param\`etres certains peuvent \^etre vides). 
Lorsque l'\'ev\'enement $(\IR-D_s^{\Theta_0^i})\cap(\IR-D_i^{\Theta''_0})$ (resp. $(\IR-D_s^{\Theta'_1})\cap(\IR-D_i^{\Theta_0^s})$) n'est pas vide, il est $\Theta-(\Theta_0^s\cup\Theta''_0)$ (resp. $\Theta-(\Theta'_1\cup\Theta_0^i)$) n\'egligeable ; ceci permet une analyse s\'epar\'ee des deux parties de l'hypoth\`ese $\Theta_1$. Les experts du choix entre $\Theta_1$ et $\Theta_0$ s'obtiennent presque s\^urement en ajoutant un expert du choix entre $\Theta'_1$ et $\Theta_0$ \`a un expert du choix entre $\Theta''_0$ et $\Theta_0$. C'est une cons\'equence directe de la d\'efinition 4.1.1 et de la propri\'et\'e suivante.

\medskip
{\bf Proposition 5.3.1}
\medskip
\medskip
\moveleft 10.4pt\hbox{\vrule\kern 10pt\vbox{\defpro

Soit $\delta$ un expert du choix entre $\Theta_1$ et $\Theta_0$. Si l'\'ev\'enement $C$ v\'erifie : 
$\Theta_1^C=\{\theta\in\Theta_1\ ,\ P_\theta(C)>0\}\not=\emptyset$ et 
$\Theta_0^C=\{\theta\in\Theta_0\ ,\ P_\theta(C)>0\}\not=\emptyset$, la restriction de $\delta$ \`a $C$ est un expert du choix entre $\Theta_1^C$ et $\Theta_0^C$ sur le sous mod\`ele statistique conditionn\'e par $C$.
}}\medskip
\medskip
{\leftskip=15mm \dli {\bf D\'emonstration}
\medskip

$\delta$ est un expert du probl\`eme de d\'ecision 
$(\Omega ,{\cal A},(P_\theta)_{\theta\in\Theta_0\cup\Theta_1})$.
$C\in{\cal A}$ est muni de la tribu ${\cal C}$ trace de ${\cal A}$ sur $C$. 
Le mod\`ele statistique conditionn\'e par $C$, $(C ,{\cal C},(P_\theta^C)_{\theta\in\Theta_0^C\cup\Theta_1^C})$ est d\'efini par 
$P_\theta^C(C')=P_\theta(C')/P_\theta(C)$ ($C'\in{\cal C}$). On doit d\'emontrer que $\delta^C$, la restriction de $\delta$ \`a $C$ est un expert du choix entre $\Theta_1^C$ et $\Theta_0^C$.

La condition ii) de la d\'efinition 4.1.1 est \'evidente car pour $A\in{\cal C}\subseteq{\cal A}$ tel que $\II_A\relmont{=}{\Theta_0^C p.s.}0$ 
(resp. $\II_A\relmont{=}{\Theta_1^C p.s.}0$) on a aussi $\II_A\relmont{=}{\Theta_0 p.s.}0$ (resp. $\II_A\relmont{=}{\Theta_1 p.s.}0$).

Soient $\theta_0\in\Theta_0^C\subseteq\Theta_0$ et 
$\theta_1\in\Theta_1^C\subseteq\Theta_1$, la condition i) est aussi triviale lorsque $E_{\theta_0}(\delta)=0$ ou $E_{\theta_1}(\delta)=1$ car pour 
$\theta\in\Theta_0^C\cup\Theta_1^C$, 
$E_{\theta}^C(\delta)=\int_C\delta^C dP_\theta^C={1\over P_\theta(C)}\int_\Omega\delta\II_C dP_\theta$. Il nous reste le cas o\`u $\delta$ est un expert du choix entre $P_{\theta_1}$ et $P_{\theta_0}$ ; on peut facilement v\'erifier que $\delta^C$ est aussi un expert du choix entre $P_{\theta_1}^C$ et $P_{\theta_0}^C$ (voir la d\'efinition 2.1.1). 
\medskip\centerline{\hbox to 3cm{\bf \hrulefill}}\par}

Faire voter les experts d'hypoth\`eses bilat\'erales impropres est un probl\`e\-me \'equivalent \`a celui du vote des experts sous des hypoth\`eses unilat\'erales. Pour une r\'ealisation $\omega$, le vote $Q_\theta^\omega$ sous $P_\theta$ ne d\'epend que d'un seul des deux experts unilat\'eraux qui d\'efinissent l'expert bilat\'eral.

Les hypoth\`eses bilat\'erales le plus souvent trait\'ees dans le cadre de la th\'eorie des tests sont en fait des hypoth\`eses bilat\'erales pures. C'est en particulier le cas dans les mod\`eles exponentiels \`a un param\`etre,  
$D_i^{\Theta''_0}$ et $D_s^{\Theta'_1}$ sont alors n\'egligeables. Plus g\'en\'eralement, les hypoth\`eses bilat\'erales sont parfois \'etudi\'ees sur les familles de densit\'es strictement totalement positives \`a l'ordre trois ou plus (cf. [Leh.] p. 120 et 140). Les densit\'es \'etant strictement positives on a $D_i^{\Theta''_0}=D_s^{\Theta'_1}=\emptyset$ et les  hypoth\`eses  bilat\'erales sont pures. Ces hypoth\`eses ne poss\`edent g\'en\'eralement que des experts triviaux, comme le montre la proposition suivante.
\medskip
{\bf Proposition 5.3.2}
\medskip
\medskip
\moveleft 10.4pt\hbox{\vrule\kern 10pt\vbox{\defpro

$(\Omega ,{\cal A},(p_\theta.\mu)_{\theta\in\Theta\subseteq\IR})$ est un mod\`ele statistique \`a rapport de vraisem\-blance monotone par rapport \`a la statistique $T$. 

Soient $\{\Theta_1 , \Theta_0\}$ des hypoth\`eses bilat\'erales pures bas\'ees sur les hypoth\`eses unilat\'erales $\{\Theta'_1 , \Theta'_0\}$ et $\{\Theta''_1 , \Theta''_0\}$, $\Theta'_1\subset\Theta''_1$. 
\dli L'ensemble des experts de ce probl\`eme de d\'ecision se r\'eduit aux experts triviaux, presque s\^urement \'egaux \`a $\II_{D_i^{\Theta_0}\cup D_s^{\Theta_0}}(T)$ ou $1-\II_{D_i^{\Theta''_0}\cap D_s^{\Theta'_1}}(T)$, 
lorsque $B=T^{-1}(\IR-D_i^{\Theta''_0}-D_s^{\Theta'_1})$ est n\'egligeable et dans le cas contraire d\`es qu'il existe $\theta_a\in\Theta'_1$, $\theta_b\in\Theta_0$ et $\theta_c\in\Theta''_0$ pour lesquels on ne peut pas trouver un recouvrement de $B$, $B_a\cup B_c\relmont{\supseteq}{p.s.}B$, tel que le rapport des densit\'es $p_{\theta_b}/p_{\theta_a}$ (resp. $p_{\theta_c}/p_{\theta_b}$) soit constant sur $B_a$ (resp. $B_c$) quand il est d\'efini.
}}\medskip
\medskip
{\leftskip=15mm \dli {\bf D\'emonstration}

Nous allons commencer par \'etablir une propri\'et\'e  des \'ev\'ene\-ments n\'egligeables dans un mod\`ele \`a rapport de vraisemblance monotone.
\medskip
{\parindent=-10mm I --- Lemme : Soient deux sous espaces ordonn\'es de $\Theta$ : $\Theta'<\Theta''$.}
\dli Si $N\subset T^{-1}(\IR-D_s^{\Theta'})$ (resp. $N\subset T^{-1}(\IR-D_i^{\Theta''})$) est un \'ev\'enement $\Theta'$ (resp. $\Theta''$) n\'egligeable, il est $\Theta'\cup\Theta''$ n\'egligeable.
\medskip
{\parindent=-5mm 1\up{er}cas : $N\subset T^{-1}(\IR-D_s^{\Theta'})$ et $N$ est $\Theta'$ n\'egligeable.}
\dli Soit $\theta''\in\Theta''$, on doit montrer que l'\'ev\'enement $N$ est $P_{\theta''}$ n\'egligeable. Nous distinguerons deux \'eventualit\'es.

i) $\theta_s=sup\Theta'\in\Theta'$. Dans ce cas $P_{\theta_s}(N)=0$ et 
$D_s^{\Theta'}=D_s^{\theta_s}$ (voir le lemme 1 de l'annexe IV).
Posons $N=N_1+N_2+N_3$ avec $N_1=N\cap\{p_{\theta''}=0\}$, $N_2=N\cap\{p_{\theta''}>0\}\cap\{p_{\theta_s}>0\}$ et 
$N_3=N\cap\{p_{\theta''}>0\}\cap\{p_{\theta_s}=0\}$. On a \'evidemment 
$P_{\theta''}(N_1)=0$ mais aussi $P_{\theta''}(N_2)=0$ car $P_{\theta_s}(N_2)=0$ implique que $N_2$ est $\mu$ n\'egligeable. Pour finir nous allons montrer que $N_3$ est vide. S'il existait $\omega\in N_3$ on aurait $h_{(\theta'',\theta_s)}(T(\omega))=+\infty$ donc 
$h_{(\theta'',\theta_s)}(t)=+\infty$ pour $t\in D=[T(\omega),+\infty[$ ; 
la densit\'e $p_{\theta_s}$ serait nulle sur $T^{-1}(D)$ ; $T(\omega)$ appartiendrait \`a $D_s^{\theta_s}$, ce qui est impossible puisque   $T(\omega)$ appartient \`a $T(N)$ qui est d'intersection vide avec $D_s^{\Theta'}=D_s^{\theta_s}$.

ii) $\theta_s=sup\Theta'\notin\Theta'$. On consid\`ere dans $\Theta'$ une suite $(\theta'_n)_{n\in\IN}$ croissant vers $\theta_s$. On a 
$D_s^{\Theta'}=\cap_{n\in\IN}D_s^{\theta'_n}$. Pour tout $n\in\IN$,  l'\'ev\'enement $N$  est $P_{\theta'_n}$ n\'egligeable et comme pr\'ec\'edemment on le d\'ecompose en $N_1+N_2^n+N_3^n$ avec $N_1=N\cap\{p_{\theta''}=0\}$, $N_2^n=N\cap\{p_{\theta''}>0\}\cap\{p_{\theta'_n}>0\}$,  
$N_3^n=N\cap\{p_{\theta''}>0\}\cap\{p_{\theta'_n}=0\}$ ; on a encore 
$P_{\theta''}(N_1)=0$ et $P_{\theta''}(N_2^n)=0$. Il nous suffit maintenant de montrer que $N'=\cap_{n\in\IN}N_3^n$ est vide puisque $N=N_1+(\cup_{n\in\IN}N_2^n)+N'$. S'il existait $\omega\in N'$ on aurait pour tout $\theta'_n$, comme pr\'ec\'edemment pour $\theta_s$ : 
$T(\omega)\in D_s^{\theta'_n}$ ; ceci est impossible puisque   $T(\omega)$ appartient \`a $T(N)$ qui est d'intersection vide avec $D_s^{\Theta'}=\cap_{n\in\IN}D_s^{\theta'_n}$.

{\parindent=-5mm 2\up{\`eme}cas : $N\subset T^{-1}(\IR-D_i^{\Theta''})$ et $N$  est $\Theta''$ n\'egligeable.}
\dli Soit $\theta'\in\Theta'$, on doit montrer que l'\'ev\'enement $N$ est $P_{\theta'}$ n\'egligeable. La d\'emonstration est semblable \`a celle du premier cas. 

Pour tout $\theta''\in\Theta''$ on pose 
$N=N_1+N_2^{\theta''}+N_3^{\theta''}$ avec $N_1=N\cap\{p_{\theta'}=0\}$, $N_2^{\theta''}=N\cap\{p_{\theta'}>0\}\cap\{p_{\theta''}>0\}$ et 
$N_3^{\theta''}=N\cap\{p_{\theta'}>0\}\cap\{p_{\theta''}=0\}$. On a toujours $P_{\theta'}(N_1)=0$ et $N_2^{\theta''}$ $\mu$ n\'egligeable. Quant \`a $N_3^{\theta''}$ il est inclus dans $T^{-1}(D_i^{\theta''})$. En effet pour  $\omega\in N_3^{\theta''}$ on a $h_{(\theta'',\theta')}(T(\omega))=0$ ; la fonction $h_{(\theta'',\theta')}$ est alors nulle sur $D=]-\infty,T(\omega)]$ et il en est donc de m\^eme de la densit\'e $p_{\theta''}$ sur $T^{-1}(D)$ ; $T(\omega)$ appartient bien \`a $D_i^{\theta''}$.

Si $\theta_i=inf\Theta''$ appartient \`a $\Theta''$ on a bien $P_{\theta'}(N)=0$ car l'\'ev\'enement $N_3^{\theta_i}$ est vide. En effet $N_3^{\theta_i}$ est alors inclus \`a la fois dans  $T^{-1}(D_i^{\theta_i})$ et dans $N\subset T^{-1}(\IR-D_i^{\Theta''})$ avec $D_i^{\Theta''}=D_i^{\theta_i}$. 

Dans le cas contraire, $\theta_i=inf\Theta''\notin\Theta''$, on consid\`ere  une suite $(\theta''_n)_{n\in\IN}$ de $\Theta''$ d\'ecroissant vers $\theta_i$ et la d\'ecomposition : 
\dli $N=N_1+(\cup_{n\in\IN}N_2^{\theta''_n})+N'$ avec 
$N'=\cap_{n\in\IN}N_3^{\theta''_n}$. L'\'ev\'enement $N'$ est encore vide  car 
$N'\subset\cap_{n\in\IN}T^{-1}(D_i^{\theta''_n})=T^{-1}(\cap_{n\in\IN}D_i^{\theta''_n})$ et $\cap_{n\in\IN}D_i^{\theta''_n}=D_i^{\Theta''}$. L'\'ev\'enement $N$  est donc bien $P_{\theta'}$ n\'egligeable.
\medskip
{\parindent=-10mm II -- $\phi_0=\II_{D_i^{\Theta_0}\cup D_s^{\Theta_0}}(T)$ et $\phi_1=1-\II_{D_i^{\Theta''_0}\cap D_s^{\Theta'_1}}(T)$ sont des experts.}

$\phi_0$ (resp. $\phi_1$) v\'erifient bien la propri\'et\'e i) de la d\'efinition 4.1.1 puisque $E_\theta(\phi_0)=0$ (resp. $E_\theta(\phi_1)=1$) pour tout $\theta$ de $\Theta_0$ (resp. $\Theta_1$).

La propri\'et\'e ii) de cette m\^eme d\'efinition contient deux propri\'et\'es \`a d\'emontrer.

1) $\II_A\relmont{=}{\Theta_0 p.s.}0\ \Longrightarrow\ \phi_0\II_A\relmont{=}{\Theta_1 p.s.}\II_A$ (resp. $\phi_1\II_A\relmont{=}{\Theta_1 p.s.}\II_A$). 
\dli C'est \'evident pour $\phi_1$ qui est nul sur un \'ev\'enement $\Theta_1$ n\'egligeable : $T^{-1}(D_i^{\Theta''_0}\cap D_s^{\Theta'_1})$. On a la propri\'et\'e recherch\'ee sur tout \'ev\'enement $A$.
\dli Dans le cas de $\phi_0$ on doit d\'emontrer que $N=A- T^{-1}(D_i^{\Theta_0}\cup D_s^{\Theta_0})$ est $\Theta_1$ n\'egligeable. $N$ \'etant $\Theta_0$ n\'egligeable, il suffit d'appliquer deux fois le lemme pr\'ec\'edent avec $\Theta'_1<\Theta_0$ puis $\Theta_0<\Theta''_0$.

2) $\II_A\relmont{=}{\Theta_1 p.s.}0\ \Longrightarrow\ \phi_0\II_A\relmont{=}{\Theta_0 p.s.}0$ (resp. $\phi_1\II_A\relmont{=}{\Theta_0 p.s.}0$).
\dli C'est \'evident pour $\phi_0$ qui est \'egal \`a $1$ sur un \'ev\'enement $\Theta_0$ n\'egligeable : $T^{-1}(D_i^{\Theta_0}\cup D_s^{\Theta_0})$. On a  la propri\'et\'e recherch\'ee sur tout \'ev\'enement $A$.
\dli Consid\'erons le cas de $\phi_1$. Le lemme pr\'ec\'edent appliqu\'e \`a  $\Theta'_1<\Theta_0$ et \`a l'\'ev\'enement $\Theta'_1$ n\'egligeable :   $N'=A- T^{-1}(D_s^{\Theta'_1})$, nous dit que $N'$ est aussi $\Theta_0$ n\'egligeable. En utilisant $\Theta_0<\Theta''_0$ on obtient que $N''=A- T^{-1}(D_i^{\Theta''_0})$, est $\Theta_0$ n\'egligeable .
Posons $A'=A-(N'\cup N'')$, il nous reste \`a montrer l'\'egalit\'e : 
$\phi_1\II_{A'}\relmont{=}{\Theta_0 p.s.}0$, c'est \'evident puisque 
$\phi_1$ est nul sur $T^{-1}(D_i^{\Theta''_0}\cap D_s^{\Theta'_1})$ et $A'\subseteq T^{-1}(D_i^{\Theta''_0}\cap D_s^{\Theta'_1})$. 

\medskip
{\parindent=-10mm III - $B=T^{-1}(\IR-D_i^{\Theta''_0}-D_s^{\Theta'_1})$ est n\'egligeable.}

Cette situation est \'equivalente \`a :  $\II_{\IR-D_s^{\Theta'_1}}(T)\relmont{\leq}{p.s.}\II_{D_i^{\Theta''_0}}(T)$. Les hypoth\`eses \'etant bilat\'erales pures on a m\^eme : 
\dli $\II_{D_i^{\Theta_0}}(T)<\II_{\IR-D_s^{\Theta'_1}}(T)\relmont{\leq}{p.s.}\II_{D_i^{\Theta''_0}}(T)<\II_{\IR-D_s^{\Theta_0}}(T)$. 

Soit $\phi$ un expert du choix entre $\Theta_1$ et $\Theta_0$, on doit d\'emontrer qu'il est presque s\^urement \'egal \`a $\phi_0$ ou $\phi_1$. 
\dli Les \'ev\'enements $F=T^{-1}(D_i^{\Theta_0}\cup D_s^{\Theta_0})$ et 
$G=T^{-1}(D_i^{\Theta''_0}\cap D_s^{\Theta'_1})$ \'etant respectivement $\Theta_0$ et $\Theta_1$ n\'egligeables, d'apr\`es la propri\'et\'e ii) de la d\'efinition 4.1.1, l'expert $\phi$ v\'erifie : $\phi\II_F\relmont{=}{\Theta_1 p.s.}\II_F$ et $\phi\II_G\relmont{=}{\Theta_0 p.s.}0$ ; 
il est, comme $\phi_0$ et $\phi_1$, presque s\^urement \'egal \`a $1$ sur $F$ et \`a $0$ sur $G$.
\dli Pour obtenir $\phi\relmont{=}{p.s.}\phi_0$ ou  $\phi\relmont{=}{p.s.}\phi_1$ il faut montrer que $\phi$ est presque s\^urement constant sur $A+C$ avec  
$A=T^{-1}(D_i^{\Theta''_0}-D_i^{\Theta_0}-D_s^{\Theta'_1})$ et  
$C=T^{-1}(D_s^{\Theta'_1}-D_s^{\Theta_0}-D_i^{\Theta''_0})$ (on a 
$F^c=A+G+C$ ou $F^c=A+B+C$ si $G=\emptyset$), c'est-\`a-dire que 
$\{\phi=1\}\cap (A+C)=A_1+C_1$ ou 
\dli $\{\phi=0\}\cap (A+C)=A_0+C_0$ est n\'egligeable.
\dli $A+C$ \'etant inclus dans $F^c$, d'apr\`es le lemme de la partie I, il suffit de d\'emontrer que $A_1+C_1$ ou $A_0+C_0$ est $\Theta_0$  n\'egligeable. 

Supposons que $A_1+C_1$ ne soit pas $\Theta_0$  n\'egligeable, c'est-\`a-dire qu'il existe $\theta_0\in\Theta_0$ tel que $P_{\theta_0}(A_1+C_1)>0$, nous  devons montrer que l'\'ev\'enement $A_0+C_0$ est $\Theta_0$  n\'egligeable.

{\parindent=-5mm 1\up{\`ere} \'etape : $P_{\theta_0}(A_1)>0\ \Longrightarrow\ C_0$ est $\Theta_0$  n\'egligeable.}
\dli $C_0$ \'etant inclus dans $T^{-1}(\IR-D_i^{\Theta''_0})$, il suffit de d\'emontrer que $C_0$ est $\Theta''_0$  n\'egligeable (voir le lemme de I). 
Si ce n'\'etait pas le cas il existerait $\theta_1\in\Theta''_0$ tel que 
$P_{\theta_1}(C_0)>0$ et on aurait $E_{\theta_1}(\phi)<1$ ; comme $P_{\theta_0}(A_1)>0$ implique $E_{\theta_0}(\phi)>0$, l'expert $\phi$ serait un expert du choix entre $P_{\theta_1}$ et $P_{\theta_0}$ (voir le i) de la d\'efinition 4.1.1) ; ceci est impossible, en effet $A_1\subseteq T^{-1}(D_i^{\Theta''_0})$ \'etant $\theta_1$ n\'egligeable on devrait avoir 
$P_{\theta_0}(A_1\cap\{\phi=1\})=P_{\theta_0}(A_1)=0$ (voir la d\'efinition 2.1.1).

{\parindent=-5mm 2\up{\`eme} \'etape : $P_{\theta_0}(C_1)>0\ \Longrightarrow\ A_0$ est $\Theta_0$  n\'egligeable.}
\dli La d\'emonstration est semblable \`a celle de l'\'etape pr\'ec\'edente. 
Comme $A_0\subseteq T^{-1}(\IR-D_s^{\Theta'_1})$, il suffit de d\'emontrer que $A_0$ est $\Theta'_1$  n\'egligeable (voir le lemme de I). 
Si ce n'\'etait pas le cas il existerait $\theta_1\in\Theta'_1$ tel que 
$P_{\theta_1}(A_0)>0$ et on aurait $E_{\theta_1}(\phi)<1$ ; $P_{\theta_0}(C_1)>0$ impliquant $E_{\theta_0}(\phi)>0$, l'expert $\phi$ serait un expert du choix entre $P_{\theta_1}$ et $P_{\theta_0}$ (voir le i) de la d\'efinition 4.1.1) ; ceci est impossible, en effet $C_1\subseteq T^{-1}(D_s^{\Theta'_1})$ est $\theta_1$ n\'egligeable, on devrait avoir 
$P_{\theta_0}(C_1\cap\{\phi=1\})=P_{\theta_0}(C_1)=0$ (voir la d\'efinition 2.1.1).

A partir des deux propri\'et\'es pr\'ec\'edentes, il nous suffit maintenant de montrer que l'on ne peut pas avoir $P_{\theta_0}(A_1)=0$ ou  $P_{\theta_0}(C_1)=0$.

Si on avait $P_{\theta_0}(A_1)=0$, on aurait $P_{\theta_0}(C_1)>0$ puisque 
$P_{\theta_0}(A_1+C_1)>0$ et d'apr\`es la deuxi\`eme \'etape $P_{\theta_0}(A_0)=0$, donc $P_{\theta_0}(A+T^{-1}(D_i^{\Theta_0}))=0$ avec 
$A+T^{-1}(D_i^{\Theta_0})=T^{-1}(D_i^{\Theta''_0}-D_s^{\Theta'_1})$ puisque 
$D_i^{\Theta_0}\cap D_s^{\Theta'_1}=\emptyset$ ; $T^{-1}(D_i^{\Theta''_0}-D_s^{\Theta'_1})$ \'etant \'egal \`a $T^{-1}(\IR-D_s^{\Theta'_1})-B$ avec $B$ vide ou n\'egligeable, on aurait $P_{\theta_0}(T^{-1}(\IR-D_s^{\Theta'_1}))=0$ ce qui est impossible pour des hypoth\`eses bilat\'erales pures (voir la d\'efinition 5.3.1). 

De m\^eme $P_{\theta_0}(C_1)=0$ et $P_{\theta_0}(A_1)>0$ est impossible pour des hypoth\`eses bilat\'erales pures. Dans ce cas on aurait :
\dli $0=P_{\theta_0}(C_0)=P_{\theta_0}(C)=P_{\theta_0}(C+T^{-1}(D_s^{\Theta_0}))=P_{\theta_0}(T^{-1}(D_s^{\Theta'_1}-D_i^{\Theta''_0}))=P_{\theta_0}(T^{-1}(\IR-D_i^{\Theta''_0})-B)=P_{\theta_0}(T^{-1}(\IR-D_i^{\Theta''_0}))$.

\medskip
{\parindent=-10mm IV -- $B=T^{-1}(\IR-D_i^{\Theta''_0}-D_s^{\Theta'_1})$ n'est pas n\'egligeable et il n'existe pas $B_a\cup B_c\relmont{\supseteq}{p.s.}B$ avec $p_{\theta_b}/p_{\theta_a}$ (resp. $p_{\theta_c}/p_{\theta_b}$) constant sur $B_a$ (resp. $B_c$).}

Dans ce cas on a : 
\dli $\II_{D_i^{\Theta_0}}(T)\leq\II_{D_i^{\Theta''_0}}(T)<\II_{\IR-D_s^{\Theta'_1}}(T)\leq\II_{\IR-D_s^{\Theta_0}}(T)$ et $\phi_1=\II_\Omega$.

Soit $\phi$ un expert du choix entre $\Theta_1$ et $\Theta_0$. Comme en III il est presque s\^urement \'egal \`a $1$, donc a $\phi_0$ et $\phi_1$, sur 
$F=T^{-1}(D_i^{\Theta_0}\cup D_s^{\Theta_0})$. On doit d\'emontrer qu'il est presque s\^urement \'egal \`a $0$ ou $1$ sur $F^c=A+B+C$ avec 
$A=T^{-1}(D_i^{\Theta''_0}-D_i^{\Theta_0})$ et  
$C=T^{-1}(D_s^{\Theta'_1}-D_s^{\Theta_0})$. Ceci revient \`a montrer que 
$\{\phi=1\}\cap (A+B+C)=A_1+B_1+C_1$ ou $\{\phi=0\}\cap (A+B+C)=A_0+B_0+C_0$ est n\'egligeable. Comme pr\'ec\'edemment il suffit de d\'emontrer que l'un de ces deux \'ev\'enements est $\Theta_0$  n\'egligeable. 
\dli Nous distinguerons deux cas suivant que $B_0$ est ou n'est pas n\'egligeable. Nous aurons besoin des propri\'et\'es suivantes : 

(1) si $\theta_0\in\Theta_0$ et $P_{\theta_0}(A_1)>0$, $B_0+C_0$ est $\Theta_0$  n\'egligeable

(2) si $\theta_0\in\Theta_0$ et $P_{\theta_0}(C_1)>0$, $A_0+B_0$ est $\Theta_0$  n\'egligeable

La propri\'et\'e (1) (resp. (2)) se d\'emontre en suivant le raisonnement de la 1\up{\'ere} (resp. 2\up{\'eme}) \'etape de III et en rempla\c cant $C_0$ par $B_0+C_0$ (resp. $A_0$ par $A_0+B_0$).

\medskip
{\parindent=-5mm 1\up{er} cas : $B_0$ est n\'egligeable.}

$B$ n'\'etant pas n\'egligeable, $B_1$ ne l'est pas. On doit donc d\'emontrer que $A_0$ et $C_0$ sont $\Theta_0$ n\'egligeables.

a) Montrons que $A_0$ est $\Theta_0$ n\'egligeable.
\dli Si $C_1$ n'est pas $\Theta_0$ n\'egligeable c'est une cons\'equence de la propri\'et\'e (2).
\dli Lorsque $C_1$ est $\Theta_0$ n\'egligeable, montrons d'abord que $\phi$ est un expert du choix entre $\Theta'_1$ et $\Theta_0$. $\phi$ \'etant un expert du probl\`eme de d\'ecision  $\{\Theta_1,\Theta_0\}$, il suffit de v\'erifier que pour tout \'ev\'enement $N$, $\Theta'_1$ n\'egligeable, on a 
$\phi\II_N\relmont{=}{\Theta_0 p.s.}0$ c'est-\`a- dire $(\{\phi=1\}\cap N)$ $\Theta_0$ n\'egligeable (voir la d\'efinition 4.1.1). Pour cela d\'ecomposons $N$ en $N_1=N\cap T^{-1}(\IR-D_s^{\Theta'_1})$, $N_2=N\cap C$ et 
$N_3=N\cap T^{-1}(D_s^{\Theta_0})$ ; $N_3\subseteq T^{-1}(D_s^{\Theta_0})$ et $N_2\cap\{\phi=1\}\subseteq C_1$ sont \'evidemment $\Theta_0$ n\'egligeables, il en est de m\^eme de $N_1\subseteq N$, d'apr\`es le lemme I, puisqu'il est $\Theta'_1$ n\'egligeable.
\dli Les hypoth\`eses unilat\'erales $\{\Theta'_1,\Theta_0\}$ \'etant stables, d'apr\`es la proposition 4.2.2 il existe deux \'el\'ements successifs $f$ et $f'$ de $\overline{\Delta_s}\subseteq F$ qui encadrent $\phi$, 
$\Theta'_1\cup\Theta_0=\Theta''_1$ presque s\^urement :  
\dli $\II_{D_i^{\Theta_0}}(T)\leq f\relmont{\leq}{\Theta''_1 p.s.}\phi\relmont{\leq}{\Theta''_1 p.s.}f'\leq\II_{\IR-D_s^{\Theta'_1}}(T)$ et $]f,f'[\cap\overline{\Delta_s}=\emptyset$. 
\dli $B$ n'\'etant pas n\'egligeable, il n'est pas $\Theta_0$ n\'egligeable (voir le lemme I) ; comme $\phi$ vaut presque s\^urement $1$ sur $B$ on a : 
$f'\relmont{=}{\Theta_0 p.s.}\II_{\IR-D_s^{\Theta'_1}}(T)$.
\dli La condition de l'\'enonc\'e implique en particulier que le rapport $p_{\theta_b}/p_{\theta_a}$ n'est pas constant sur 
$B_a\relmont{\supseteq}{p.s.}B$. Il existe donc au moins une fonction de test simple $\phi_{(k,\beta)}^{(\theta_b,\theta_a)}\in\Phi_s^{(\theta_b,\theta_a)}\subseteq\Delta_s$ qui v\'erifie : 
\dli  $\II_{D_i^{\Theta''_0}}(T)\relmont{<}{p.s.}\phi_{(k,\beta)}^{(\theta_b,\theta_a)}\relmont{<}{p.s.}\II_{\IR-D_s^{\Theta'_1}}(T)$ ; d'apr\`es le lemme I on a aussi : $\II_{D_i^{\Theta''_0}}(T)\relmont{<}{\Theta_0 p.s.}\phi_{(k,\beta)}^{(\theta_b,\theta_a)}\relmont{<}{\Theta_0 p.s.}\II_{\IR-D_s^{\Theta'_1}}(T)$.
\dli On en d\'eduit $\phi_{(k,\beta)}^{(\theta_b,\theta_a)}<f'$ et comme 
$]f,f'[\cap\Phi_s^{(\theta_b,\theta_a)}=\emptyset$ on a : 
\dli $f\geq\phi_{(k,\beta)}^{(\theta_b,\theta_a)}>\II_{D_i^{\Theta''_0}}(T)$. $\phi$ est donc $\Theta_0$ presque s\^urement \'egal \`a $1$ sur $T^{-1}(D_i^{\Theta''_0})\supseteq A$, ce qui implique bien que $A_0$ est  $\Theta_0$ n\'egligeable.

b) Montrons que $C_0$ est $\Theta_0$ n\'egligeable.
\dli La d\'emonstration est semblable \`a la pr\'ec\'edente.
\dli Pour $A_1$ non $\Theta_0$ n\'egligeable c'est la propri\'et\'e (1) qui donne le r\'esultat.
\dli Lorsque $A_1$ est $\Theta_0$ n\'egligeable, $\phi$ est cette fois un expert du choix entre $\Theta''_0$ et $\Theta_0$ car 
$\II_N\relmont{=}{\Theta''_0 p.s.}0$ implique bien
$\phi\II_N\relmont{=}{\Theta_0 p.s.}0$ ; il suffit de d\'ecomposer $N$ en : 
$N_1=N\cap T^{-1}(\IR-D_i^{\Theta''_0})$, $N_2=N\cap A$ et 
$N_3=N\cap T^{-1}(D_i^{\Theta_0})$ ;  $N_3\subseteq T^{-1}(D_i^{\Theta_0})$ et $N_2\cap\{\phi=1\}\subseteq A_1$ sont bien $\Theta_0$ n\'egligeables, il en est de m\^eme de $N_1$ d'apr\`es le lemme de I.
\dli L'application de la proposition 4.2.2 aux hypoth\`eses $\{\Theta_0,\Theta''_0\}$ nous conduit \`a deux \'el\'ements $g$ et $g'$ de $F$ v\'erifiant : 
\dli $\II_{D_i^{\Theta''_0}}(T)\leq g\relmont{\leq}{\Theta'_0 p.s.}(1-\phi)\relmont{\leq}{\Theta'_0 p.s.}g'\leq\II_{\IR-D_s^{\Theta_0}}(T)$ et $]g,g'[\cap\overline{\Delta'_s}=\emptyset$ avec 
$\Delta'_s=\cup_{(\theta,\theta')\in\Theta_0\times\Theta''_0}\Phi_s^{(\theta',\theta)}$. 
\dli $(1-\phi)$ valant presque s\^urement $0$ sur $B$ non $\Theta_0$ n\'egligeable on a : 
\dli $g\relmont{=}{\Theta_0 p.s.}\II_{D_i^{\Theta''_0}}(T)$.
\dli La condition de l'\'enonc\'e implique que le rapport $p_{\theta_c}/p_{\theta_b}$ n'est pas constant sur $B_c\relmont{\supseteq}{p.s.}B$. Il existe donc au moins une fonction de test simple $\phi_{(k,\beta)}^{(\theta_c,\theta_b)}\in\Phi_s^{(\theta_c,\theta_b)}$ qui v\'erifie : 
\dli  $\II_{D_i^{\Theta''_0}}(T)\relmont{<}{\Theta_0 p.s.}\phi_{(k,\beta)}^{(\theta_c,\theta_b)}\relmont{<}{\Theta_0  p.s.}\II_{\IR-D_s^{\Theta'_1}}(T)$.
\dli On en d\'eduit $g<\phi_{(k,\beta)}^{(\theta_c,\theta_b)}$ et comme 
$]g,g'[\cap\Phi_s^{(\theta_c,\theta_b)}=\emptyset$ on a : 
\dli $g'\leq\phi_{(k,\beta)}^{(\theta_c,\theta_b)}<\II_{\IR-D_s^{\Theta'_1}}(T)$. $(1-\phi)$ est donc $\Theta_0$ presque s\^urement \'egal \`a $0$ sur $T^{-1}(\IR-D_s^{\Theta'_1})\supseteq C$, ce qui implique bien que $C_0=C\cap\{(1-\phi)=1\}$ est  $\Theta_0$ n\'egligeable.

\medskip
{\parindent=-5mm 2\up{\`eme} cas : $B_0$ n'est pas n\'egligeable.}

$B_0$ \'etant inclus dans $F^c$ il n'est pas $\Theta_0$ n\'egligeable (voir le lemme I). D'apr\`es les propri\'et\'es (1) et (2) on a donc : 
$\forall\theta_0\in\Theta_0$ $P_{\theta_0}(A_1)=0$ et $P_{\theta_0}(C_1)=0$.

Nous devons d\'emontrer que $B_1$ est $\Theta_0$ n\'egligeable.
\dli Le fait que $C_1$ (resp. $A_1$) soit $\Theta_0$ n\'egligeable entra\^{\i}ne que $\phi$ est un expert du choix entre $\Theta'_1$ et $\Theta_0$ (resp. $\Theta''_0$ et $\Theta_0$)(voir les parties a) et b) du 1\up{er}cas).
Comme pr\'ec\'edemment on obtient des \'el\'ements de $F$ v\'erifiant : 
\dli $\II_{D_i^{\Theta_0}}(T)\leq f\relmont{\leq}{\Theta''_1 p.s.}\phi\relmont{\leq}{\Theta''_1 p.s.}f'\leq\II_{\IR-D_s^{\Theta'_1}}(T)$,  $]f,f'[\cap\overline{\Delta_s}=\emptyset$ et 
\dli $\II_{D_i^{\Theta''_0}}(T)\leq g\relmont{\leq}{\Theta'_0 p.s.}(1-\phi)\relmont{\leq}{\Theta'_0 p.s.}g'\leq\II_{\IR-D_s^{\Theta_0}}(T)$,  $]g,g'[\cap\overline{\Delta'_s}=\emptyset$.
\dli Ces deux encadrements impliquent que les \'ev\'enements 
$\{f=1\}\cap\{g=1\}$ et $\{f'=0\}\cap\{g'=0\}$ sont $\Theta_0$ n\'egligeables car sur ces \'ev\'enements $\phi$ est $\Theta_0$ presque s\^urement \'egal \`a $0$ et \`a $1$. On a donc : 
\dli $\II_{D_i^{\Theta_0}}(T)\relmont{=}{\Theta_0 p.s.}inf\{f,g\}\relmont{\leq}{\Theta_0 p.s.}\II_{D_i^{\Theta''_0}}(T)<\II_{\IR-D_s^{\Theta'_1}}(T)\relmont{\leq}{\Theta_0 p.s.}sup\{f',g'\}\relmont{=}{\Theta_0 p.s.}\II_{\IR-D_s^{\Theta_0}}(T)$.
\dli Si $f'\relmont{\leq}{\Theta_0 p.s.}\II_{D_i^{\Theta''_0}}(T)$ ou 
$g\relmont{\geq}{\Theta_0 p.s.}\II_{\IR-D_s^{\Theta'_1}}(T)$ l'expert $\phi$ 
est $\Theta_0$ presque s\^urement \'egal \`a $0$ sur $B$ et $B_1$ est bien $\Theta_0$ n\'egligeable.

Consid\'erons l'autre possibilit\'e : 
\dli $f'\relmont{>}{\Theta_0 p.s.}\II_{D_i^{\Theta''_0}}(T)$ et  
$g\relmont{<}{\Theta_0 p.s.}\II_{\IR-D_s^{\Theta'_1}}(T)$. Nous allons distinguer deux cas suivant que $inf\{f,g\}$ est \'egal \`a $f$ ou $g$.

a) $inf\{f,g\}=f$ donc $f\relmont{\leq}{\Theta_0 p.s.}\II_{D_i^{\Theta''_0}}(T)$.
\dli Nous allons utiliser la fonction de test simple 
$\phi_{(k,\beta)}^{(\theta_b,\theta_a)}$ d\'efinie dans la partie a) du 1\up{er} cas. Comme elle n'appartient pas \`a $]f,f'[$ on a : 
\dli  $f'\leq\phi_{(k,\beta)}^{(\theta_b,\theta_a)}\relmont{<}{\Theta_0 p.s.}\II_{\IR-D_s^{\Theta'_1}}(T)$.
\dli Ceci implique : 
$sup\{f',g'\}=g'\relmont{\geq}{\Theta_0 p.s.}\II_{\IR-D_s^{\Theta'_1}}(T)$. 
\dli La fonction de test simple $\phi_{(k,\beta)}^{(\theta_c,\theta_b)}$ d\'efinie dans la partie b) du 1\up{er} cas n'appartenant pas \`a $]g,g'[$, on a : $\II_{D_i^{\Theta''_0}}(T)\relmont{<}{\Theta_0 p.s.}\phi_{(k,\beta)}^{(\theta_c,\theta_b)}\leq g$.
\dli Lorsque $f'\relmont{\leq}{\Theta_0 p.s.}g$ on a \'evidemment 
$\phi\relmont{=}{\Theta_0 p.s.}0$ et $B_1$ est bien $\Theta_0$ n\'egligeable.
Il nous reste le cas : 
\dli $f\relmont{\leq}{\Theta_0 p.s.}\II_{D_i^{\Theta''_0}}(T)\relmont{<}{\Theta_0 p.s.}g\relmont{<}{\Theta_0 p.s.}f'\relmont{<}{\Theta_0 p.s.}\II_{\IR-D_s^{\Theta'_1}}(T)\relmont{\leq}{\Theta_0 p.s.}g'$. Il est impossible de l'avoir sous la condition de l'\'enonc\'e. En effet sur l'\'ev\'enement $\{f'-f=1\}$ (resp. $\{g'-g=1\}$) le rapport des densit\'es 
$p_{\theta_b}/p_{\theta_a}$ (resp. $p_{\theta_c}/p_{\theta_b}$) est constant et $\{f'-f=1\}\cup\{g'-g=1\}$ formerait un recouvrement, $\Theta_0$ donc $\Theta$ presque s\^ur, de $B$.

b) $inf\{f,g\}=g$ donc $g\relmont{=}{\Theta_0 p.s.}\II_{D_i^{\Theta''_0}}(T)\relmont{=}{\Theta_0 p.s.}\II_{D_i^{\Theta_0}}(T)$.
\dli Nous allons utiliser la fonction de test simple 
$\phi_{(k,\beta)}^{(\theta_c,\theta_b)}$ d\'efinie dans la partie b) du 1\up{er} cas. Comme elle n'appartient pas \`a $]g,g'[$ on a : 
\dli  $g'\leq\phi_{(k,\beta)}^{(\theta_c,\theta_b)}\relmont{<}{\Theta_0 p.s.}\II_{\IR-D_s^{\Theta'_1}}(T)$.
\dli Ceci implique : 
$sup\{f',g'\}=f'\relmont{=}{\Theta_0 p.s.}\II_{\IR-D_s^{\Theta'_1}}(T)\relmont{=}{\Theta_0 p.s.}\II_{\IR-D_s^{\Theta_0}}(T)$. 
\dli La fonction de test simple $\phi_{(k,\beta)}^{(\theta_b,\theta_a)}$ d\'efinie dans la partie a) du 1\up{er} cas n'appartenant pas \`a $]f,f'[$, on a : $\II_{D_i^{\Theta''_0}}(T)\relmont{<}{\Theta_0 p.s.}\phi_{(k,\beta)}^{(\theta_b,\theta_a)}\leq f$.
\dli Le cas $g'\relmont{\leq}{\Theta_0 p.s.}f$ est impossible car on aurait 
$\phi\relmont{=}{\Theta_0 p.s.}1$ et $B_0$ serait $\Theta_0$ n\'egligeable donc n\'egligeable.
Il nous reste le cas : 
\dli $g\relmont{=}{\Theta_0 p.s.}\II_{D_i^{\Theta''_0}}(T)\relmont{<}{\Theta_0 p.s.}f\relmont{<}{\Theta_0 p.s.}g'\relmont{<}{\Theta_0 p.s.}\II_{\IR-D_s^{\Theta'_1}}(T)\relmont{=}{\Theta_0 p.s.}f'$. Il est aussi impossible du fait de la condition impos\'ee par l'\'enonc\'e. Comme nous l'avons vu a la fin du paragraphe a) pr\'ec\'edent, $p_{\theta_b}/p_{\theta_a}$ (resp. $p_{\theta_c}/p_{\theta_b}$) est constant sur $\{f'-f=1\}$ (resp. $\{g'-g=1\}$) et on a encore 
$B\relmont{\subseteq}{p.s.}\{f'-f=1\}\cup\{g'-g=1\}$.
\nobreak
\medskip\centerline{\hbox to 3cm{\bf \hrulefill}}\par}

L'existence de $\theta_a$, $\theta_b$ et $\theta_c$, v\'erifiant les conditions de l'\'enonc\'e, est r\'ealis\'ee dans les mod\`eles \`a rapport de vraisemblance monotone classiquement \'etudi\'es. Nous allons cependant donner un exemple simple o\`u des hypoth\`eses bilat\'erales pures poss\`edent des experts non triviaux.

Soit $(\IR ,{\cal B},(p_\theta.\lambda)_{\theta\in\Theta_0\cup\Theta_1})$, $\lambda$ \'etant la mesure de Lebesgue, le probl\`eme de d\'ecision d\'efini par : 
$\Theta_0=\{\theta_0\}$, $\Theta_1=\{\theta'_1,\theta''_0\}$, $p_{\theta_0}={1\over 2}\II_{[-1,+1]}$, $p_{\theta'_1}={1\over 2}\II_{[-{3\over 2},+{1\over 4}]}+{1\over 4}\II_{]+{1\over 4},+{3\over 4}]}$ 
et $p_{\theta''_0}={1\over 4}\II_{[-{3\over 4},-{1\over 4}[}+{1\over 2}\II_{[-{1\over 4},+{3\over 2}]}$. Ce mod\`ele est \`a rapport de vraisemblance monotone pour la statistique identit\'e et les hypoth\`eses sont bien bilat\'erales pures. On a $B=[-{3\over 4},+{3\over 4}]$ et la condition de la proposition 5.3.2 n'est pas r\'ealis\'ee puisque : $p_{\theta_0}/p_{\theta'_1}=1$ sur $B_a=[-1,+{1\over 4}]$, $p_{\theta''_0}/p_{\theta_0}=1$ sur $B_c=[-{1\over 4},+1]$ et $B\subset B_a\cup B_c$. Consid\'erons les r\`egles de d\'ecision valant $1$ sur 
$]-\infty,-1[\cup]+1,+\infty[$, $0$ sur $[-1,-{1\over 4}[\cup]+{1\over 4},+1]$, quelconques sur $[-{1\over 4},+{1\over 4}]$ mais non triviales, c'est-\`a-dire non presque s\^urement \'egales \`a $0$ sur $[-{1\over 4},+{1\over 4}]$. Ces r\`egles sont des experts du choix entre $\Theta_1$ et $\Theta_0$ ; la propri\'et\'e ii) de la d\'efinition 4.1.1 est \'evidente et il est facile de v\'erifier que ces r\`egles de d\'ecision sont des experts du choix entre $\theta'_1$ (resp. $\theta''_0$) et $\theta_0$ ; en effet, elles sont 
$\{\theta'_1,\theta_0\}$ (resp. $\{\theta''_0,\theta_0\}$) presque s\^urement comprises entre $\II_{]-\infty,-1[}$ et $\II_{]-\infty,+{1\over 4}]}$ 
(resp. $1-\II_{]-\infty,+1]}$ et $1-\II_{]-\infty,-{1\over 4}[}$) (voir la proposition 2.3.1).

Revenons au cas des hypoth\`eses bilat\'erales pures non expertisables. Deux traitements de cette p\'enurie d'experts nous semblent possibles suivant que l'ordre stucturant les param\`etres est ou n'est pas important pour l'interpr\'etation. 
\dli Dans bien des cas les deux parties $\Theta'_1$ et $\Theta''_0$ de l'hypoth\`ese $\Theta_1$ ne sont pas vraiment \'equivalentes pour l'utilisateur. Son choix principal est entre  $\Theta_0$ et $\Theta_1$ mais pour autant $\Theta'_1$ et $\Theta''_0$ s'interpr\`etent diff\'eremment, m\^eme si cette diff\'erence est mise au second plan. Les hypoth\`eses 
$\{\Theta'_1,\Theta'_0\}$, $\{\Theta''_1,\Theta''_0\}$, voire le choix entre les trois \'eventualit\'es $\{\Theta'_1,\Theta_0,\Theta''_0\}$, sont envisageables. Dans un tel cadre une bonne solution nous para\^{\i}t \^etre l'utilisation de votes compatibles pour les hypoth\`eses unilat\'erales d\'efinissant le probl\`eme bilat\'eral, et m\^eme pour l'ensemble des hypoth\`eses unilat\'erales (voir les propositions 5.2.3 et 5.2.4). Nous \'etudierons des exemples ult\'erieurement. 
\dli Consid\'erons maintenant le cas o\`u l'ordre sur les param\`etres n'a pas d'int\'er\^et pour diff\'erencier $\Theta'_1$ et $\Theta''_0$, 
ces deux \'eventualit\'es conduisent \`a la m\^eme interpr\'etation alors que l'ordre 
les oppose. Cette \'equivalence entre $\Theta'_1$ et $\Theta''_0$ se traduit sur les r\'ealisations par une \'equivalence entre les petites et les grandes valeurs de la statistique $T$ rendant le rapport de vraisemblance monotone. 
Reprenons par exemple le cas d'un n-\'echantillon de la loi normale de moyenne inconnue $\theta\in\IR$ et de variance connue $\sigma^2$. La statistique $T$ est alors la moyenne empirique $\overline{X}$. Int\'eressons nous  
aux hypoth\`eses $\{\Theta_0=[\theta_1,\theta_2],\Theta_1=\IR-\Theta_0\}$, ce probl\`eme de d\'ecision poss\`ede des sym\'etries fortes par rapport \`a 
$\theta_0=(\theta_1+\theta_2)/2$. Si les deux demi-droites $]-\infty,\theta_1[$ et $]\theta_2,+\infty[$ sont \'equivalentes pour l'interpr\'etation il est tentant de dire qu'il est \'equivalent de r\'ealiser 
$\overline{x}$ ou son sym\'etrique par rapport \`a $\theta_0$ : 
$2\theta_0-\overline{x}$. C'est la statistique $T'=\mid\overline{X}-\theta_0\mid$ qui devient alors primordiale et on est amen\'e \`a travailler sur son mod\`ele statistique image.
Dans ce nouveau mod\`ele les hypoth\`eses $\{\Theta_0,\Theta_1\}$ sont expertisables sans probl\`eme puisque $({\sqrt n\over\sigma}T')^2$ suit une loi de khi-deux d\'ecentr\'ee \`a un degr\'e de libert\'e et de param\`etre d'excentricit\'e ${n\over\sigma^2}(\theta-\theta_0)^2$. Sur le mod\`ele de d\'epart ceci revient \`a remplacer, dans la d\'efinition des experts, la tribu $\cal A$ par la sous tribu $\cal C$ engendr\'ee par $T'$. On impose ainsi des contraintes moins fortes pour le label expert et les hypoth\`eses bilat\'erales pures non expertisables peuvent devenir stables. C'est un proc\'ed\'e \'equivalent \`a celui qui consiste, dans la th\'eorie de la d\'ecision \`a partir d'une fonction de perte, \`a diminuer l'ensemble des r\`egles admissibles en imposant une contrainte suppl\'ementaire. Ici il faut enlever des contraintes car il y a p\'enurie et non pas trop-plein. Il reste le probl\`eme important du choix de la sous tribu $\cal C$. Dans un cadre g\'en\'eral il est impossible de faire intervenir les sym\'etries du mod\`ele, comme dans l'exemple des lois normales, pour proposer une statistique $T'$ d\'efinissant les petites et grandes valeurs de $T$ \'equivalentes. Disons que l'on doit choisir 
les \'el\'ements de $\cal C$ de la forme $C'=T^{-1}(]-\infty,t)\cup(t',+\infty[)$ de telle sorte que l'on puisse dire qu'ils se comportent de fa\c con semblable sous $\Theta'_1$ et $\Theta''_0$. 
Les probabilit\'es de $C'$ sous $\Theta'_1$ et celles sous $\Theta''_0$ doivent se ressembler. Lorsque $\theta$ tend vers $inf\Theta$ (resp. $sup\Theta$) c'est la partie $]-\infty,t)$ (resp. $(t',+\infty[$) qui prend de l'importance. Une mani\`ere d'exprimer la ressemblance, sans probabiliser $\Theta'_1$ et $\Theta''_0$, est d'imposer l'\'egalit\'e des probabilit\'es les plus proches : 
\dli $lim_{\theta\in\Theta'_1\nearrow sup\Theta'_1}P_\theta(C')=lim_{\theta\in\Theta''_0\searrow inf\Theta''_0}P_\theta(C')$. 
Si $\Theta_0$ est un intervalle d'extr\'e\-mi\-t\'es $\theta_1<\theta_2$, sous des conditions de continuit\'e courantes (par exemple la continuit\'e en $\theta$ de $p_\theta(t)$ pour presque tout $t$), l'\'egalit\'e pr\'ec\'edente s'\'ecrit 
$P_{\theta_1}(C')=P_{\theta_2}(C')$. Ceci fait penser aux tests sans biais. La r\'egion de rejet $C'$ doit v\'erifier $P_{\theta}(C')\geq\alpha$ pour tout $\theta$ de $\Theta_1$. Sous la condition de continuit\'e pr\'ec\'edente et si la famille $\{p_\theta\}_{\theta\in\Theta}$ est totalement positive \`a l'ordre trois, cette propri\'et\'e de la r\'egion de rejet est \'equivalente \`a 
$P_{\theta_1}(C')=P_{\theta_2}(C')$ (cf. [Mor.1] p.63). La notion de sans biais est utile pour tester $H_0 : \theta\in\Theta_0$ contre $H_1 : \theta\in\Theta_1$ mais pas pour tester $H_0 : \theta\in\Theta_1$ contre 
$H_1 : \theta\in\Theta_0$. En fait dans ce dernier cas c'est la notion de seuil sous $\Theta_1=\Theta'_1\cup\Theta''_0$ qui impose une r\'egion de rejet de probabilit\'e $\alpha$ sous $\theta_1$ et $\theta_2$. Une fois d\'efinie la sous tribu $\cal C$ \`a partir des $C'$ pr\'ec\'edents il faut r\'e\'ecrire le mod\`ele avec des densit\'es $\cal C$ mesurables. Dans ce nouveau mod\`ele on pourra souvent trouver une statistique r\'eelle $T'$ rendant les hypoth\`eses stables en r\'e\'ecrivant les $C'$ sous la forme $T'\in(x,+\infty[$. Les r\'esultats du paragraphe 4 peuvent alors s'appliquer. Pour l'exemple pr\'ec\'edent le mod\`ele peut m\^eme \^etre param\'etr\'e par $[\theta_0,+\infty[$, les hypoth\`eses sont alors unilat\'erales et le mod\`ele \`a rapport de vraisemblance monotone. Nous n'insisterons pas plus sur ce type de solutions car nous pensons que dans la majorit\'e des applications, l'ordre sur l'espace des param\`etres influence les interpr\'etations de $\theta\in\Theta'_1$ et $\theta\in\Theta''_0$. 

L'\'etude des tests bilat\'eraux est classiquement faite dans les mod\`eles exponentiels \`a un param\`etre : 
$(\Omega ,{\cal A},(P_\theta=c(\theta)e^{\theta T(\omega)}.\mu)_{\theta\in\Theta\subseteq\IR})$. $\Theta$ \'etant un intervalle de $\IR$ muni de l'ordre usuel, ces mod\`eles sont \`a rapport de vraisem\-blance monotone par rapport \`a $T$. Pour $\theta'<\theta''$, le rapport des densit\'es $h_{(\theta'',\theta')}(t)=p_{\theta''}(t)/p_{\theta'}(t)=[c(\theta'')/c(\theta')]e^{(\theta''-\theta')t}$ est toujours d\'efini et strictement croissant \`a valeur dans $]0,+\infty[$.
$T$ est \'evidemment une statistique essentielle globale de fonction de r\'epartition moyenne : $G_\theta(t)=F_\theta(t)+{1\over 2}p_\theta(t)\mu(T^{-1}[t])$ ($F_\theta$ \'etant la fonction de r\'epartition de l'image de $P_\theta$ par $T$). Comme fonction de $\theta$, $G_\theta(t)$ est m\^eme continue pour tout $t$ puisque $p_\theta(t)$ l'est, ce qui implique la continuit\'e des $p_\theta$ dans $L_1(\mu)$ (cf. [Mon.1] p. 138). Toutes les hypoth\`eses unilat\'erales $\{\Theta^f_1,\Theta^f_0\}$ sont donc adjacentes : 
$Q^\omega_{\Theta^f_1}(1_f)=Q^\omega_{\Theta^f_0}(1_f)=Q^\omega_{\theta^f}(1_f)=1-G_{\theta^f}(T(\omega))$ avec $\theta^f=sup\Theta^f_1=inf\Theta^f_0$ (voir la d\'efinition 4.3.2 et les propositions 4.3.1-2). Ces votes sont \'evidemment compatibles sur les deux hypoth\`eses unilat\'erales $\{\Theta'_1,\Theta'_0\}$ et $\{\Theta''_1\supset\Theta'_1,\Theta''_0\}$ qui d\'efinissent les hypoth\`eses bilat\'erales $\{\Theta_1,\Theta_0\}$. Posons $\theta_1=inf\Theta_0$ et $\theta_2=sup\Theta_0$, le prolongement de ces votes nous donne pour les hypoth\`eses bilat\'erales les votes suivants, lorsqu'on r\'ealise $\omega$ ou $t=T(\omega)$ : 
\dli $Q^\omega(\Theta_0)=Q^\omega_{\theta_2}(\Theta''_1)-Q^\omega_{\theta_1}(\Theta'_1)=[F_{\theta_1}(t)-F_{\theta_2}(t)]+{1\over 2}[p_{\theta_1}(t)-p_{\theta_2}(t)]\mu(T^{-1}[t])$ et $Q^\omega(\Theta_1)=1-Q^\omega(\Theta_0)$.
\dli On obtiendrait les m\^emes r\'esultats \`a partir de la probabilit\'e d\'efinie sur $\Theta$ par la proposition 5.2.3 lorsqu'elle s'applique. Sur les exemples du paragraphe suivant nous verrons que c'est g\'en\'eralement le cas pour les mod\`eles exponentiels et nous donnerons des applications de la proposition 5.2.4. 

Les votes que nous venons de d\'efinir ne s'interpr\`etent pas comme seuil minimum des deux tests bilat\'eraux. Pour obtenir ceci il faudrait choisir la solution pr\'ec\'edente du changement de mod\`ele. Le vote $Q^\omega(\Theta_0)$ tend \'evidemment vers $0$ lorsque $t=T(\omega)$ tend vers $-\infty$ ou $+\infty$ ; si le rapport des densit\'es $h_{(\theta_2,\theta_1)}(T)$ est \'egal \`a $1$ en $\omega$, le vote $Q^\omega(\Theta_0)$ cro\^{\i}t avant $t=T(\omega)$ et d\'ecro\^{\i}t ensuite. Le maximum de $Q^\omega(\Theta_0)$ est d'autant plus grand que la probabilit\'e $P_{\theta_2}$ est grande par rapport \`a $P_{\theta_1}$ donc que l'intervalle $\Theta_0$ est grand. Lorsque $\Theta_0$ se r\'eduit \`a un point $\theta_0$ on a $Q^\omega(\Theta_0)=0$, on ne d\'ecide donc jamais $\theta=\theta_0$. Ceci nous semble totalement l\'egitime car la continuit\'e autour de $\theta_0$ fait que dans $\Theta_1$ il y a des probabilit\'es aussi proches que l'on veut de 
$P_{\theta_0}$. Ce type d'hypoth\`eses, bien que souvent utilis\'e, correspond rarement \`a une bonne traduction de ce que veut v\'erifier l'utilisateur : il cherche g\'en\'eralement \`a savoir si le param\`etre est autour de $\theta_0$. 
Cet ``autour'' d\'epend du contexte : ordre de grandeur des erreurs de mesure, ordre de grandeur des diff\'erences qui entra\^{\i}nent une interpr\'etation diff\'erente, etc... R\'eduire cet intervalle autour de $\theta_0$ \`a $[\theta_0]$ est tr\`es populaire car ceci \'evite bien des questions \`a l'utilisateur. Malgr\'e des critiques renouvel\'ees (cf. par exemple [Reu.],[Wan.]) ces hypoth\`eses restent tr\`es utilis\'ees, elles procurent un confort attirant si l'on ne se pose pas de question sur la mani\`ere dont fonctionne la solution statistique employ\'ee. Nous n'en proposerons pas de nouvelle, il faudrait prendre en compte des votes autres que celui sous $\theta_0$. Une telle d\'emarche est envisageable si on veut traduire un a priori positif par rapport \`a la d\'ecision $\theta=\theta_0$. Par exemple, lorsque cette hypoth\`ese est l'approximation d'une hypoth\`ese qui se concentre autour de $\theta_0$ et qui a priori a des chances d'\^etre vraie (cf. [BerD]).

\vfill\eject

\bigskip
{\parindent=-5mm 5.4 EXEMPLES.}
\nobreak
\medskip
La notion de votes compatibles nous a permis de proposer une aide \`a la d\'ecision pour des hypoth\`eses bilat\'erales. Bien souvent, cette aide peut \^etre d\'efinie \`a partir d'une probabilisation de l'espace des param\`etres. Il est alors possible de traiter toutes les hypoth\`eses dont la structure repose sur l'ordre de l'espace des param\`etres. Cette probabilit\'e peut aussi aider \`a prendre une d\'ecision dans le cas de plus de deux hypoth\`eses, en particulier lorsque $\Theta$ est partag\'e en une partition d'intervalles ordonn\'es. L'explicitation de la probabilit\'e d\'efinie sur $\Theta$ par les propositions 5.2.3 et 5.2.4 est inutile dans les probl\`emes ne faisant intervenir que des intervalles, ce qui est le cadre normal d'utilisation de ces propositions. Les fonctions de r\'epartition moyenne, $(G_\theta)_{\theta\in\Theta}$, de la statistique essentielle globale suffisent. On n'a pas besoin d'inverser les r\^oles de $\theta$ et de la r\'ealisation $\omega$ afin de d\'efinir la probabilisation de $\Theta$ par les fonctions de r\'epartition $(Q^\omega(]\leftarrow,\theta[))_{\omega\in\Omega}$.
Il est toutefois int\'eressant d'obtenir ces probabilit\'es sur quelques exemples classiques. Nous devrons en particulier v\'erifier les conditions d'application de la proposition 5.2.3. Comme nous l'avons vu pr\'ec\'edemment, les deux premi\`eres le sont d\`es que $p_\theta(\omega)$ est continu en $\theta$ pour presque tout $\omega$ (cf. [Mon.1] p. 138).

\medskip
EXEMPLE 1 : param\`etre de position.
\medskip

Consid\'erons un mod\`ele statistique \`a rapport de vraisemblance monotone et une statistique essentielle globale $K(T)$ dont les fonctions de r\'epartition moyenne $G_\theta$ v\'erifient les conditions de la proposition 5.2.3.
Si $\theta$ est un param\`etre de position pour $K(T)$, c'est-\`a-dire si $K(T)$ admet des densit\'es de la forme $g(x-\theta)$ par rapport \`a la mesure de Lebesgue, il est facile d'expliciter la probabilit\'e induite sur $\Theta$ par les votes les plus favorables. En effet lorsque l'\'egalit\'e 
$Q^\omega(]\leftarrow,\theta[)=1-G_\theta(K(T(\omega)))$ s'applique on peut l'\'ecrire : 
\dli $Q^\omega(]\leftarrow,\theta[)=\int\II_{]K(T(\omega)),+\infty[}\!(x)g(x-\theta)\,dx=\int\II_{]-\infty,\theta[}\!(\lambda)g(K(T(\omega))-\lambda)\,d\lambda$ avec $\lambda=[\theta+K(T(\omega))]-x$.

La famille des lois uniformes sur $[\theta-1,\theta+1]$ entre dans ce cadre. Les densit\'es $({1\over 2}\II_{[\theta-1,\theta+1]}\!(x))_{\theta\in\IR}$ d\'efinissent un mod\`ele \`a rapport de vraisemblance monotone par rapport \`a la statistique identit\'e, qui est aussi une statistique essentielle globale. 
On a $G_\theta(t)=\int^t_{-\infty}{1\over 2}\II_{[-1,+1]}\!(x-\theta)\,dx$ qui est continue en $\theta$ pour tout $t$ et tend vers $0$ (resp. $1$) lorsque $\theta$ tend vers $+\infty$ (resp. $-\infty$). Les quatre conditions de la proposition 5.2.3 sont donc v\'erifi\'ees, $Q^t(]-\infty,\theta[)$ est le vote $Q^t_\theta(1)$ pour les hypoth\`eses unilat\'erales $\{\Theta_1=]-\infty,\theta[,\Theta_0=[\theta,+\infty[\}$. \dli $Q^t(]-\infty,\theta[)$ est \'egal \`a $1$ si $\theta\geq t+1$, \`a $0$ si 
$\theta<t-1$ et \`a $1-G_\theta(t)=\int^\theta_{-\infty}{1\over 2}\II_{[-1,+1]}\!(t-\lambda)\,d\lambda$ sinon. Quand on r\'ealise $t$, on obtient sur l'espace des param\`etres la loi uniforme sur $[t-1,t+1]$. C'est la loi a posteriori de la mesure de Lebesgue, loi a priori impropre et non informative.

L'exemple le plus classique de ce type de param\`etre de position est celui de la famille des lois normales $(N(\theta,a^2))_{\theta\in\IR}$, l'\'ecart-type  $a>0$ \'etant connu. Nous avons d\'ej\`a trait\'e ce cas \`a la fin du paragraphe 5.2, avec l'exemple d'un n-\'echantillon d'une loi normale de moyenne $\theta$ inconnue et de variance $\sigma^2$ connue ($a^2=\sigma^2/n$). La proposition 5.2.3 s'applique et on trouve sur $\Theta=\IR$ la loi $N(t,a^2)$ lorsqu'on r\'ealise $t$. Si pour chaque probl\`eme unilat\'eral $\{\Theta^f_1=]-\infty,\theta_f),\Theta^f_0\}$ on consid\`ere la pond\'eration 
$N(\theta_f,c^2)$, la proposition 5.2.4 s'applique et nous avons obtenu sur 
$\Theta$ la loi $N(t,a^2+c^2)$. Nous avons aussi vu que ces votes \'etaient neutres. Ceci n'est pas \'etonnant puisque cette pond\'eration traite les hypoth\`eses de chaque probl\`eme unilat\'eral de fa\c con semblable et sym\'etriquement. Pour chacun de ces probl\`emes elle r\'eduit l'\'ecart entre le vote pour $\Theta^f_1$ et le vote pour $\Theta^f_0$, ce qui ne facilite pas la conclusion. \dli Essayons maintenant de construire une pond\'eration qui tienne compte d'une information a priori.
Par exemple, d'une valeur m\'ediane $\theta_0$, pour laquelle l'utilisateur consid\'ererait qu'il y a autant de chance que la vraie valeur de $\theta$ soit en dessous qu'au dessus. On est alors amen\'e \`a privil\'egier l'hypoth\`ese 
$\Theta^f_1=]-\infty,\theta_f)$ (resp. $\Theta^f_0$) lorsque $\theta_f$ est sup\'erieur (resp. inf\'erieur) \`a $\theta_0$. Ceci ne peut se faire qu'en consi\-d\'erant d'autres possibilit\'es de votes que celle sous $\theta_f$. Une pond\'eration gaus\-sienne des votes du type $N(\mu_f,c^2)$, $c^2>0$, donne des calculs simples, nous reviendrons sur ce choix \`a la fin. Pour avantager $\Theta^f_1$ (resp. $\Theta^f_0$) lorsque $\theta_f>\theta_0$ (resp. $\theta_f<\theta_0$) il faut que $\mu_f$ soit sup\'erieur (resp. inf\'erieur) \`a $\theta_f$. Un $\mu_f$ du type $\theta_f+\lambda(\theta_f-\theta_0)$, avec $\lambda>0$, semble une solution int\'eressante. Plus $\theta_f$ est loin de $\theta_0$ plus on avantage l'une des hypoth\`eses. Le param\`etre $\lambda$ sera d'autant plus grand que l'utilisateur consid\`ere qu'\`a une faible distance de $\theta_0$ il n'y a pratiquement plus qu'une hypoth\`ese de plausible. Dans un tel cadre, les hypoth\`eses $\{\Theta^f_1,\Theta^f_0\}$ sont trait\'ees avec un vote pond\'er\'e par $\Lambda^f$ qui est la loi normale : $N(\mu_f=\theta_f+\lambda(\theta_f-\theta_0),c^2)$. Nous devons regarder si la propri\'et\'e 5.2.4 s'applique. Les propri\'et\'es i) et ii) sont \'evidentes puisque $D^\theta_i=D^\theta_s=\emptyset$. Quant \`a $H_t(f)$, il v\'erifie :
\cleartabs
\+ $H_t(f)$&=&$\int_\Theta(1-G_\theta(t))d\Lambda^f(\theta)$\cr
\+ &=&$1-\int_{\IR} F({t-\theta\over a})\,{1\over\sqrt{2\pi}\,c} exp(-{(\theta-\mu_f)^2\over 2c^2})\,d\theta$\cr
\+ &=&$1-\int F({t-\mu_f-\nu\over a})\,{1\over\sqrt{2\pi}\,c} exp(-{\nu^2\over 2c^2})\,d\nu$\cr
\+ &=&$1-\int[\int_{-\infty}^{t-\mu_f} {1\over\sqrt{2\pi}\,a}exp(-{(x-\nu)^2\over 2a^2}\,dx]{1\over\sqrt{2\pi}\,c} exp(-{\nu^2\over 2c^2})\,d\nu$\cr
\+ &=&$1-\int_{-\infty}^{t-\mu_f}\,{1\over\sqrt{2\pi}\,a}exp(-{x^2\over 2a^2})[\int\,{1\over\sqrt{2\pi}\,c} exp(-{1\over 2}({\nu^2\over a^2}+{\nu^2\over c^2}-{2x\nu\over a^2}))\,d\nu]\,dx$\cr
\+ &=&$1-\int_{-\infty}^{t-\mu_f}\,{1\over\sqrt{2\pi}\,a}exp(-{x^2\over 2(a^2+c^2)})[\int\,{1\over\sqrt{2\pi}\,c} exp(-{a^2+c^2\over 2a^2c^2}(\nu-{xc^2\over a^2+c^2})^2)\,d\nu]\,dx$\cr
\+ &=&$1-F({t-\mu_f\over\sqrt{a^2+c^2}})=F({\mu_f-t\over\sqrt{a^2+c^2}})=F({1+\lambda\over\sqrt{a^2+c^2}}(\theta_f-{t+\lambda\theta_0\over 1+\lambda}))$.\cr

\dli Il est \'evident que les votes $Q_{\Lambda^f}$ sont compatibles, d'ailleurs $H_t(f)$ est bien une fonction croissante et continue de $\theta_f$, donc de $f$, et elle tend vers $0$ (resp. $1$) lorsque $\theta_f$ tend vers $-\infty$ (resp. $+\infty$). On obtient sur $\Theta$ la loi normale : $N({t+\lambda\theta_0\over 1+\lambda},{a^2+c^2\over(1+\lambda)^2})$. 
\dli Plus $\lambda$ est grand plus le point $\theta_0$ choisi prend le pas sur la r\'ealisation $t$. $c^2$, lui, n'intervient que dans la variance. Au d\'epart il a \'et\'e introduit pour prendre en compte des votes autres que celui de $\theta_f$, limite entre les deux hypoth\`eses $\Theta^f_1$ et $\Theta^f_0$. Plus $c^2$ est petit, plus la loi sur $\Theta$ est concentr\'ee. On peut d'ailleurs prendre $c^2=0$, c'est-\`a-dire consid\'erer le vote sous $\theta=\mu_f$ pour les hypoth\`eses $\{\Theta^f_1,\Theta^f_0\}$. $H_t(f)$ est alors \'egal \`a $1-F({t-\mu_f\over a})$ et on obtient encore la loi 
$N({t+\lambda\theta_0\over 1+\lambda},{a^2\over(1+\lambda)^2})$. 
Cependant, si l'information a priori demand\'ee \`a l'utilisateur avait surtout pour but d'influencer le param\`etre $t$ de la distribution $N(t,a^2)$ on peut 
prendre $c^2=\lambda(\lambda+2)a^2$. Remarquons enfin que la loi 
$N({t+\lambda\theta_0\over 1+\lambda},{a^2+c^2\over(1+\lambda)^2})$ est la loi a posteriori que l'on obtient \`a partir de la loi a priori 
$N(\theta_0,\tau^2)$ si $\lambda={a^2\over\tau^2}$ et $c^2=\lambda a^2$ (cf. [Ber.] p. 127).
\vfill\eject

\medskip
EXEMPLE 2 : param\`etre d'\'echelle.
\medskip

Soient un mod\`ele statistique \`a rapport de vraisemblance monotone et une statistique essentielle globale $K(T)$ positive dont les fonctions de r\'epartition moyenne $G_\theta$ v\'erifient les conditions de la proposition 5.2.3. Nous dirons que
$\theta\in\Theta\subseteq]0,+\infty[$ est un param\`etre d'\'echelle pour $K(T)$ si $K(T)$ admet des densit\'es de la forme ${1\over\theta}g({x\over\theta})$ par rapport \`a la mesure de Lebesgue sur $\IR^+$ (cf. [Leh.] p. 510). Lorsque l'\'egalit\'e $Q^\omega(]0,\theta[)=1-G_\theta(K(T(\omega)))$ s'applique on peut l'\'ecrire : 
\dli $Q^\omega(]0,\theta[)=\int\II_{]K(T(\omega)),+\infty[}\!(x){1\over\theta}g({x\over\theta})\,dx=\int\II_{]0,\theta[}\!(\lambda){K(T(\omega))\over\lambda^2}g({K(T(\omega))\over\lambda})\,d\lambda$ avec $\lambda={\theta K(T(\omega))\over x}$. Il est alors facile d'expliciter la probabilisation de $\Theta$ induite par $Q^\omega$.

La famille des lois uniformes sur $[0,\theta]$, $\theta\!>\!0$, entre dans ce cadre. Les densit\'es $({1\over \theta}\II_{[0,\theta]}\!(x))_{\theta\in\IR^+_*}$ d\'efinissent sur $\IR^+$ un mod\`ele \`a rapport de vraisemblance monotone pour la statistique identit\'e, qui est aussi une statistique essentielle globale. 
Les quatre conditions de la proposition 5.2.3 sont bien v\'erifi\'ees puisque  $G_\theta(t)=\int^t_{0}{1\over\theta}\II_{[0,1]}\!({x\over\theta})\,dx$. 
$Q^t(]0,\theta[)$ est le vote $Q^t_\theta(1)$ pour les hypoth\`eses unilat\'erales $\{\Theta_1=]0,\theta[,\Theta_0=[\theta,+\infty[\}$. 
$Q^t(]0,\theta[)$ est \'egal \`a $0$ si $\theta\leq t$ et \`a $1-G_\theta(t)=\int^\theta_{0}{t\over\lambda^2}\II_{[0,1]}\!({t\over\lambda})\,d\lambda$ sinon. Quand on r\'ealise $t$, on obtient sur l'espace des param\`etres $\IR^+_*$ la loi de densit\'e ${t\over\theta^2}\II_{[t,+\infty[}\!(\theta)$. C'est la loi a posteriori associ\'ee \`a la mesure de densit\'e ${1\over\theta}$ par rapport \`a la mesure de  Lebesgue sur $\IR^+$, cette mesure est consid\'er\'ee comme une loi a priori impropre, non informative (cf. [Rob.] p. 108).

 Un autre exemple classique de ce type de param\`etre d'\'echelle est celui de la famille des lois normales $(N(m,\theta^2))_{\theta\in\IR^+_*}$, la moyenne $m$ \'etant connue. 
Les densit\'es ${1\over\sqrt{2\pi}\,\theta} exp(-{1\over 2}({x-m\over\theta})^2)$ d\'efinissent sur $\IR$ un mod\`ele \`a rapport de vraisemblance monotone pour la statistique $T(x)=\mid\!x-m\!\mid$, qui est aussi une statistique essentielle globale. 
Les quatre conditions de la proposition 5.2.3 sont bien v\'erifi\'ees puisque  $G_\theta(t)=\int^t_{0}{2\over\sqrt{2\pi}\,\theta} exp(-{1\over 2}({x\over\theta})^2)\,dx=2F({t\over\theta})-1$. On obtient $Q^t(]0,\theta[)=\int_0^\theta{2\over\sqrt{2\pi}}{t\over\lambda^2}exp(-{1\over 2}({t\over\lambda})^2)\,d\lambda$. Cette loi sur $\Theta=\IR^+_*$ est encore la loi a posteriori associ\'ee \`a la loi a priori impropre et non informative de densit\'e ${1\over\theta}$ puisque 
$\int_{\IR^+_*}{1\over\theta}g({t\over\theta})\times{1\over\theta}\,d\theta=\int_{\IR^+_*}{1\over t^2}g({\nu\over t})\,d\nu={1\over t}$. 

Consid\'erons enfin la famille des lois gamma $(\gamma(p,\theta))_{\theta\in\IR^+_*}$, le param\`etre $p>0$ est fixe et $\theta>0$ est le param\`etre d'\'echelle. Elles d\'efinissent sur $\IR^+_*$ muni de la mesure de Lebesgue une famille de densit\'e $({1\over\Gamma(p)\theta}({x\over\theta})^{p-1}exp(-{x\over\theta}))_{\theta\in\IR^+_*}$ \`a rapport de vraisemblance monotone pour la statistique identit\'e qui est aussi une statistique essentielle globale. $G_\theta(t)$ est alors \'egal \`a $\Gamma(p,{t\over\theta})=\int_0^{t\over\theta}{1\over\Gamma(p)}y^{p-1}e^{-y}\,dy$ et l'application de la proposition 5.2.3 nous donne les votes compatibles suivant :
$Q^t(]0,\theta[)=\int_0^\theta{1\over\Gamma(p)t}({t\over\lambda})^{p+1}exp(-{t\over\lambda})\,d\lambda$ ($t>0$). Leur prolongement est la loi inverse d'une loi $\gamma(p,{1\over t})$ qui est bien s\^ur la loi a posteriori pour la loi a priori impropre et non informative de densit\'e ${1\over\theta}$ (cf. [Ber.]p. 255). Comme cas particulier consid\'erons celui d'un n-\'echantillon d'une loi normale de moyenne $m$ connue et de variance inconnue $\sigma^2$. Nous avons sur $\IR^n$ un mod\`ele statistique \`a rapport de vraisemblance monotone pour la statistique $T(x_1,...,x_n)=\sum_{i=1}^n(x_i-m)^2$ qui est aussi une statistique essentielle globale. $T/\sigma^2$ suivant une loi de khi-deux \`a n degr\'es de libert\'e, $T$ est de loi $\gamma({n\over 2},\theta=2\sigma^2)$. La loi obtenue pour le param\`etre $\theta$ nous donne pour $v=\sigma^2$ une loi sur $\IR^+_*$ de densit\'e $f_t(v)={1\over 2^{n/2}\Gamma({n\over 2})t}({t\over v})^{{n\over 2}+1}exp(-{t\over 2v})$. C'est la loi inverse d'une loi $\gamma({n\over 2},{2\over t})$ (cf. [Ber.] p. 561). Lorsque $n=1$, la loi du param\`etre $\theta=\sigma$ est bien s\^ur celle obtenue dans l'exemple pr\'ec\'edent. 
\dli Essayons maintenant de construire des votes qui tiennent compte d'une information a priori. Comme pr\'ec\'edemment l'utilisateur est suppos\'e fournir une valeur m\'ediane $\theta_0$. Les hypoth\`eses $\{]0,\theta_0),)\theta_0,+\infty\}$ \'etant pour lui aussi probables l'une que l'autre, on choisit un vote neutre pour ces hypoth\`eses. Le plus simple est de prendre le vote $Q_{\theta_0}$. Les autres probl\`emes unilat\'eraux $\{\Theta^f_1=]0,\theta_f),\Theta^f_0\}$ doivent \^etre trait\'es dissym\'etriquement, lorsque $\theta_f$ est sup\'erieur (resp. inf\'erieur) \`a $\theta_0$ l'avantage se porte sur $\Theta^f_1$ (resp. $\Theta^f_0$). Ceci peut se faire en choisissant un vote $Q_{\mu_f}$ v\'erifiant $\mu_f>\theta_f$ (resp. $\mu_f<\theta_f$) lorsque $\theta_f>\theta_0$ (resp. $\theta_f<\theta_0$). Bien entendu, plus $\theta_f$ est loin de $\theta_0$ plus l'\'ecart entre $\mu_f$ et $\theta_f$ doit augmenter. $\theta$ \'etant un param\`etre d'\'echelle, nous exprimerons ces distances par la valeur du rapport, donc par la diff\'erence des logarithmes.
La solution utilis\'ee dans l'exemple 1 devient :
$ln(\mu_f)=ln(\theta_f)+\lambda(ln(\theta_f)-ln(\theta_0))$ ou 
$\mu_f=\theta_f({\theta_f\over\theta_0})^\lambda$ ($\lambda>0$). Le param\`etre $\lambda$ sera d'autant plus grand que l'utilisateur consid\`ere qu'\`a une faible distance de $\theta_0$ il n'y a pratiquement plus qu'une hypoth\`ese plausible. Nous devons appliquer la proposition 5.2.4 avec des $\Lambda^f$ qui sont les masses de Dirac en $\mu_f$. Les propri\'et\'es i) et ii) sont \'evidentes puisque $D^\theta_i=D^\theta_s=\emptyset$. Quant \`a $H_t(f)$, il v\'erifie bien les conditions de l'\'enonc\'e :
\cleartabs
\+ $H_t(f)$&=&$\int_\Theta(1-G_\theta(t))d\Lambda^f(\theta)
=1-G_{\mu_f}(t)=\int_0^{\mu_f}{1\over\Gamma(p)} {t^p\over\nu^{p+1}}exp(-{t\over\nu})\,d\nu$\cr
\+ &=&$\int_0^{\theta_f}{(\lambda+1)\over\Gamma(p)} {t^p(\theta_0)^{p\lambda}\over\theta^{p(\lambda+1)+1}}exp(-{t(\theta_0)^{\lambda}\over\theta^{\lambda+1}})\,d\theta$\quad ($\nu={\theta^{\lambda+1}\over (\theta_0)^{\lambda}}$)\cr
\+ &=&$Q^t(]0,\theta_f))$\cr
\dli La loi sur $\Theta=\IR^+_*$ est la loi de $({1\over X})^{1\over \lambda+1}$ avec $X$ de loi $\gamma(p,{1\over t(\theta_0)^{\lambda}})$.
En prenant comme nouveau param\`etre $\beta=ln({\theta^{\lambda+1}\over\theta_0^\lambda})$ donc $\beta_f=ln(\mu_f)$ et en posant $u=ln(t)$, le changement de variable  $x=ln(\nu)-u$ nous donne : \dli $Q^t(]0,\theta_f))=\int_{-\infty}^{\beta_f-u}{1\over\Gamma(p)} exp[-px-e^{-x}]\,dx$. Dans le cas d'une loi exponentielle de param\`etre $\theta$ on a $p=1$ et  $Q^t(]0,\theta_f))=exp(-exp(u-\beta_f))$.

\medskip
EXEMPLE 3 : param\`etre d'une loi de Poisson.
\medskip

Le mod\`ele statistique, $(\IN,{\cal P}(\IN),(p_\theta.\mu)_{\theta\in\IR^+})$ est d\'efini par 
$p_\theta(x)=e^{-\theta}.{\theta^x\over x!}$, $\mu$ \'etant la mesure de comptage sur $\IN$. Il est \`a rapport de vraisemblance monotone pour la statistique identit\'e qui est une statistique essentielle globale. La fonction de r\'epartition moyenne $G_\theta(t)=\sum_{k=0}^te^{-\theta}.{\theta^k\over k!}\,-{1\over 2}e^{-\theta}.{\theta^t\over t!}$ est \'evidemment continue en $\theta$ et elle tend bien vers $0$ lorsque $\theta$ tend vers $+\infty$. Par contre $G_\theta(0)$ ne tend pas vers $1$ mais vers ${1\over 2}$ lorsque $\theta$ tend vers $0$. C'est pour cela que nous avons pris $\IR^+$ comme espace des param\`etres, contrairement \`a l'usage de le restreindre \`a $\IR^+_*$. Ainsi la proposition 5.2.3 s'applique, elle nous donne les votes compatibles 
$Q^t([0,\theta])=1-G_\theta(t)$. Nous allons utiliser l'expression de la fonction de r\'epartition d'une loi de Poisson sous la forme d'une int\'egrale de la densit\'e de la loi $\gamma(t+1,1)$ : 
$\sum_{k=0}^te^{-\theta}.{\theta^k\over k!}=1-\int_0^\theta{x^t\over\Gamma(t+1)}e^{-x}\,dx$ (cf. [R\'en.] p. 112).
Si $t>0$, les votes  $Q^t([0,\theta])=1-[\sum_{k=0}^{t-1}e^{-\theta}.{\theta^k\over k!}\,+{1\over 2}e^{-\theta}.{\theta^t\over t!}]$ d\'efinissent une probabilit\'e sur $\Theta$ qui est le m\'elange \'equipond\'er\'e de deux lois gamma : 
${1\over 2}\gamma(t,1)+{1\over 2}\gamma(t+1,1)$. Lorsque $t=0$, $Q^0([0,\theta])=1-{1\over 2}e^{-\theta}$ et la premi\`ere gamma du m\'elange est remplac\'ee par la masse de Dirac en $0$. On obtient sur $\Theta$ la loi ${1\over 2}\delta_0+{1\over 2}\gamma(1,1)$, $\gamma(1,1)$ \'etant la loi exponentielle de param\`etre $1$.
\dli Les m\'elanges obtenus sont \'equipond\'er\'es parce que nous avons pris en compte \`a part \'egale les deux votes possibles sous $P_\theta$ lorsque la r\'ealisation n'est pas de probabilit\'e nulle (voir le paragraphe 2.4). Ce choix a l'avantage de donner des votes neutres pour toutes les hypoth\`eses unilat\'erales (voir la proposition 4.3.2). 
\dli Les lois gamma et leurs m\'elanges sont les lois a priori conjugu\'ees privil\'egi\'ees par l'analyse bay\'esienne des lois de Poisson (cf. [Rob.] p. 98 et 100). Mais la loi que nous venons d'obtenir n'est pas une loi a posteriori pour une loi a priori m\'elange de deux gamma. Elle a plut\^ot \`a voir avec les lois  a posteriori des lois a priori impropres et non informatives souvent utilis\'ees : les densit\'es ${1\over\theta}$ et ${1\over\sqrt{\theta}}$ sur $\IR^+_*$ (cf. [Ber.] p. 114).
\dli Consid\'erons un n-\'echantillon $(X_1,...,X_n)$ d'une loi de Poisson. Il est facile de voir que $T=\sum_{i=1}^nX_i$ est une statistique essentielle globale qui suit une loi de Poisson de param\`etre $n\theta$. Ce qui pr\'ec\`ede nous donne, pour la r\'ealisation $t=x_1+...+x_n$ une loi sur $\Theta=\IR^+$ \'egale \`a ${1\over 2}\gamma(t,{1\over n})+{1\over 2}\gamma(t+1,{1\over n})$ lorsque $t>0$ et \`a ${1\over 2}\delta_0+{1\over 2}\gamma(1,{1\over n})$ lorsque $t=0$.

\medskip
EXEMPLE 4 : param\`etre d'une loi binomiale.
\medskip
Soit un mod\`ele binomial $(\Omega=\{0,1,...,n\},{\cal P}(\Omega),(B(n,\theta))_{\theta\in[0,1]})$. Ce mod\`ele est \`a rapport de vraisemblance monotone par rapport \`a la mesure $\mu$ de masse $C_n^\omega$ en $\omega$ et \`a la statistique identit\'e : $p_\theta(\omega)=\theta^\omega(1-\theta)^{n-\omega}$ avec la convention $0^0=1$. La fonction de r\'epartition moyenne de cette statistique essentielle globale est \'egale \`a : $G_\theta(\omega)=\sum_{i=0}^\omega C_n^i\theta^i(1-\theta)^{n-i}\,-{1\over 2}C_n^\omega\theta^\omega(1-\theta)^{n-\omega}$. Elle est continue en $\theta$, ce qui permet l'application de la proposition 5.2.3. On obtient pour chaque r\'ealisation $\omega$ une probabilit\'e sur les bor\'eliens de $\Theta=[0,1]$ dont la fonction de r\'epartition est d\'efinie pour $\theta>0$ par : 
\dli $Q^\omega([0,\theta[)=1-G_\theta(\omega)=\sum_{i=\omega}^n C_n^i\theta^i(1-\theta)^{n-i}\,-{1\over 2}C_n^\omega\theta^\omega(1-\theta)^{n-\omega}$. 
\dli En utilisant la formule $\sum_{i=k}^n C_n^i\theta^i(1-\theta)^{n-i}=\int_0^\theta{\Gamma(n+1)\over\Gamma(k)\Gamma(n+1-k)}x^{k-1}(1-x)^{n-k}\,dx$ lorsque $k>0$ (cf. [R\'en.] p. 88), on a pour $0<\omega<n$ : $Q^\omega([0,\theta[)={1\over 2}\int_0^\theta{\Gamma(n+1)\over\Gamma(\omega)\Gamma(n+1-\omega)}x^{\omega-1}(1-x)^{n-\omega}\,dx\,+\,{1\over 2}\int_0^\theta{\Gamma(n+1)\over\Gamma(\omega+1)\Gamma(n-\omega)}x^{\omega}(1-x)^{n-\omega-1}\,dx$. 
\dli La loi sur $\Theta$ est alors un m\'elange \'equipond\'er\'e de deux loi b\^eta \`a densit\'e sur $[0,1]$, ${1\over 2}\beta(\omega,n+1-\omega)+{1\over 2}\beta(\omega+1,n-\omega)$.
Lorsque $\omega=0$ (resp. $\omega=n$) on obtient sur $\Theta$ la loi 
${1\over 2}\delta_0+{1\over 2}\beta(1,n)$ (resp. ${1\over 2}\beta(n,1)+{1\over 2}\delta_1$), $\delta_0$ et $\delta_1$ \'etant les masses de Dirac en $0$ et $1$. 
\dli Les lois b\^eta et leurs m\'elanges sont les lois a priori conjugu\'ees privil\'egi\'ees par l'analyse bay\'esienne des lois binomiales (cf. [Rob.] p. 98 et 100). La loi que nous avons obtenue n'est cependant pas la loi a posteriori d'un m\'elange de lois b\^eta. La neutralit\'e des votes qui la d\'efinissent la rapproche plut\^ot de lois a posteriori obtenues \`a partir de lois a priori consid\'er\'ees comme non informatives (cf. [Gei.]).
Ce m\'elange de lois b\^eta a d\'ej\`a \'et\'e trouv\'e et d\'efendu, dans une \'etude s'appuyant sur les fonctions de perte ``propres'' (cf. [KroM]).

\medskip
EXEMPLE 5 : param\`etres d'une analyse de variance.
\medskip
En analyse de variance \`a effets fixes on utilise des statistiques qui suivent des lois de Fisher d\'ecentr\'ees (cf. [Sch.] p. 38). Le param\`etre de non centralit\'e $\lambda$ \'etant la quantit\'e qui int\'eresse l'utilisateur exprim\'ee par rapport \`a l'\'ecart-type commun des variables. Par exemple, dans une analyse de variance classique \`a deux facteurs, $\lambda$ est \'egal \`a $\sqrt{\sum(\theta_{ij})^2}/\sigma$ pour le test sur l'additivit\'e des facteurs, et \'egal \`a $\sqrt{\sum(\alpha_{i})^2}/\sigma$ pour le test sur la nullit\'e des effets additifs du premier facteur. Si une statistique $W$ suit une loi de Fisher d\'ecentr\'ee de param\`etre $\lambda$ et de degr\'es de libert\'e $(k,l)$, nous savons que $(k/l)W$ suit une loi b\^eta sur $\IR^+$ : 
$\beta({k\over 2},{l\over 2},{\lambda^2\over 2})$ (cf. [Bar.] p. 84).Nous allons \'etudier ces lois.
\dli Consid\'erons sur $\IR^+$ la famille des lois b\^eta d\'ecentr\'ees :  $\beta(p,q,\theta)$, les param\`etres $p>0$ et $q>0$ sont connus, le param\`etre de non centralit\'e $\theta\geq 0$, lui, est inconnu. Comme pour la famille des lois de Fisher d\'ecentr\'ees, on montre qu'elle est \`a rapport de vraisemblance strictement monotone pour la statistique identit\'e (cf. [Kar.]). Cette statistique est donc une statistique essentielle globale. Sa fonction de r\'epartition moyenne s'\'ecrit : 
\dli $G_\theta(t)=\int_0^t[\sum_{m=0}^{+\infty}e^{-\theta}.{\theta^m\over m!}{\Gamma(p+m+q)\over\Gamma(p+m)\Gamma(q)}{x^{p+m-1}\over (1+x)^{p+m+q}}]\,dx\,=\,\int_{\IN} F(p+m,q,t)\,d{\cal P}_\theta(m)$, ${\cal P}_\theta$ \'etant une loi de Poisson de param\`etre $\theta\geq 0$ et 
$F(p+m,q,t)$ la valeur en $t$ de la fonction de r\'epartition d'une loi $\beta(p+m,q)$. Cette loi est celle du quotient $Z=X/Y$ de deux variables ind\'ependantes $X$ et $Y$ de loi $\gamma(p+m,1)$ et $\gamma(q,1)$. La densit\'e de $\beta(p+m,q)$ peut s'\'ecrire comme int\'egrale en $y$ de la loi du couple $(Y,Z)$, on obtient alors :  
\cleartabs
\+ $F(p+m,q,t)$&=&$\int_0^t\II_{\IR^+_*}\!(z)[\int_0^{+\infty}{1\over\Gamma(p+m)\Gamma(q)}z^{p+m-1}exp(-y(z+1)).y^{p+m+q-1}\,dy]\,dz$\cr
\+ &=&$\int_0^{+\infty}{1\over\Gamma(q)}exp(-y).y^{q-1}[\int_0^{ty}{1\over\Gamma(p+m)}exp(-x).x^{p+m-1}\,dx]\,dy$.\cr 
\dli La loi $\gamma(p+m)$ est la loi de la somme de $m+1$ variables ind\'ependantes, l'une de loi $\gamma(p)$ et les autres de loi $\gamma(1)$. Sa fonction de r\'epartition tend donc vers $0$ lorsque $m$ tend vers $+\infty$. On en d\'eduit facilement qu'il en est de m\^eme de $F(p+m,q,t)$. Comme fonction de $\theta$, $G_\theta(t)$ v\'erifie donc bien les conditions de la proposition 5.2.3 sur $\Theta=[0,+\infty[$. Les votes compatibles $Q^t([0,\theta])=\int_{\IN}[1-F(p+m,q,t)]\,d{\cal P}_\theta(m)$ d\'efinissent une loi de probabilit\'e sur $\Theta$ muni de la tribu des bor\'eliens.Cette probabilit\'e poss\`ede une masse en $0$ \'egale \`a $1-F(p,q,t)$, c'est le seuil minimum de rejet du test de $H_0 : \theta=0$ contre $H_1 : \theta>0$. On peut exprimer $Q^t(]0,\theta])$ sous la forme $\int_0^\theta f_t(\lambda)\,d\lambda$ lorsque $p\geq{1\over 2}$, en utilisant l'expression de $F(p+m,q,t)$ pr\'ec\'edente et quand $p>{1\over 2}$, le fait que la loi $\gamma(p+m,1)$ d\'ecentr\'ee de $\theta$ est la convolution d'une loi $\gamma({1\over 2},1)$ d\'ecentr\'ee de $\theta$ par une loi $\gamma(p+m-{1\over 2},1)$ (cf. [Bar.] p. 82). Ceci permet de faire intervenir les propri\'et\'es de sym\'etrie des lois normales, puisque la loi $\gamma({1\over 2},1)$ d\'ecentr\'ee de $\theta$ est la loi du carr\'e d'une variable $N(\sqrt{\theta},{1\over 2})$. On trouve une densit\'e sous la forme d'une int\'egrale double.  Ne faisant pas partie des densit\'es classiques elle pr\'esente peu d'int\'er\^et, il vaut mieux calculer la probabilit\'e d'un intervalle de $\Theta$ pour la r\'ealisation $t$ \`a partir de la premi\`ere expression de $G_\theta(t)$.

\vfill\eject

{\parindent=-10mm\soustitre 6--HYPOTH\`ESES STABLES ET PARAM\`ETRES FANT\^OMES.}
\bigskip
{\parindent=-5mm 6.1 INTRODUCTION.}
\medskip
Dans le mod\`ele statistique $(\Omega ,{\cal A},(P_{(\theta,\upsilon)})_{(\theta,\upsilon)\in\Theta\times\Upsilon})$ consid\'erons des hypoth\`eses de la forme $\{\Theta_1\times\Upsilon,\Theta_0\times\Upsilon\}$, $\Theta_0$ et $\Theta_1$ d\'efinissant une partition de  $\Theta$ ne contenant pas le vide.
Le param\`etre $\upsilon$ est un param\`etre fant\^ome pour ce probl\`eme de d\'ecision (cf. [Bar.] p. 51). Nous noterons ce type de probl\`eme de d\'ecision : 
$(\Omega ,{\cal A},(P_{(\theta,\upsilon)})_{(\theta,\upsilon)\in(\Theta_0\cup\Theta_1)\times\Upsilon})$. Sauf cas exceptionnel, les hypoth\`eses $\{\Theta_1\times\Upsilon,\Theta_0\times\Upsilon\}$ ne sont pas stables. Il est cependant absurde de consid\'erer le probl\`eme du choix entre $P_{(\theta_1,\upsilon_1)}$ et $P_{(\theta_0,\upsilon_0)}$ lorsque $\upsilon_1\not=\upsilon_0$, puisqu'on ne cherche aucune information sur $\upsilon$. Ce qui est int\'eressant c'est la stabilit\'e des hypoth\`eses 
$\{\Theta_1\times\{\upsilon\},\Theta_0\times\{\upsilon\}\}$ pour chaque $\upsilon\in\Upsilon$.
\medskip
{\bf D\'efinition 6.1.1}
\medskip
\medskip
\moveleft 10.4pt\hbox{\vrule\kern 10pt\vbox{\defpro
Soit $(\Omega ,{\cal A},(p_{(\theta,\upsilon)}.\mu)_{(\theta,\upsilon)\in(\Theta_0\cup\Theta_1)\times\Upsilon})$ un probl\`eme de d\'ecision domin\'e par la mesure $\mu$. Les hypoth\`eses $\{\Theta_0,\Theta_1\}$ \`a param\`etre fant\^ome $\upsilon\in\Upsilon$ sont stables si il existe une statistique r\'eelle $T$ rendant stables les hypoth\`eses des sous probl\`emes de d\'ecision $(\Omega ,{\cal A},(p_{(\theta,\upsilon)}.\mu)_{\theta\in\Theta_0\cup\Theta_1})$ pour tout 
$\upsilon\in\Upsilon$.
}}\medskip
Par exemple dans un mod\`ele exponentiel de la forme $p_{(\theta,\upsilon)}(\omega)=exp(\theta T(\omega)+\langle\upsilon,U(\omega)\rangle-\psi(\theta,\upsilon))$, $\Theta\times\Upsilon\subseteq\IR\times\IR^k$, les hypoth\`eses unilat\'erales sur $\Theta$ \`a param\`etre fant\^ome $\upsilon$ sont stables par rapport \`a la statistique r\'eelle $T$. Dans ce cas la famille 
$(p_{(\theta,\upsilon)})_{\theta\in\Theta}$ est m\^eme \`a rapport de vraisemblance monotone.

Pour toute r\'ealisation $t=T(\omega)$ notons $Q^t_{(\theta,\upsilon)}$ le vote des experts sous $P_{(\theta,\upsilon)}$. Nous devons construire un vote \`a partir de la famille $(Q^t_{(\theta,\upsilon)})_{(\theta,\upsilon)\in\Theta\times\Upsilon}$. Nous avons d\'ej\`a fait ce type de travail pour les sous mod\`eles param\`etr\'es par $\Theta\times\{\upsilon\}$ (voir le paragraphe 4.3). On peut utiliser une pond\'eration $\Lambda_\upsilon$ sur $\Theta$ pour traduire une information a priori ou prendre le vote $Q^t_{\Theta_0\times\{\upsilon\}}$ (resp. $Q^t_{\Theta_1\times\{\upsilon\}}$) le plus favorable sous $\Theta_0\times\{\upsilon\}$ (resp. $\Theta_1\times\{\upsilon\}$).
A chaque valeur $\upsilon$ du param\`etre fant\^ome correspond alors un vote 
$Q^t(d,\upsilon)$ construit \`a partir des votes $(Q^t_{(\theta,\upsilon)}(d))_{\theta\in\Theta}$. Il peut arriver que les votes $Q^t(.,\upsilon)$ ne d\'ependent pas du param\`etre fant\^ome $\upsilon$. Consid\'erons par exemple la famille de lois uniformes par rapport \`a la mesure de Lebesgue : $({1\over\theta}\II_{[0,\upsilon\theta]})_{(\theta,\upsilon)\in\IR^+_*\times\{-1,+1\}}$. Nous avons \'etudi\'e ce probl\`eme dans l'exemple 2 de 5.4 avec $\upsilon=+1$. Ici nous ajoutons une direction miroir mais seule la distance $\theta$ \`a $0$ int\'eresse l'utilisateur. Toutes les hypoth\`eses unilat\'erales sur $\Theta$ \`a param\`etre fant\^ome $\upsilon\in\{-1,+1\}$ sont stables par rapport \`a la statistique valeur absolue. Si l'on choisit pour les deux sous probl\`emes de d\'ecision les votes compatibles de la proposition 5.2.3 on obtient la m\^eme probabilit\'e sur $\Theta=\IR^+_*$ lorsque $\upsilon=+1$ et $\upsilon=-1$. 
\dli Le plus souvent la famille $(Q^t(.,\upsilon))_{\upsilon\in\Upsilon}$ est constitu\'ee de votes diff\'erents, il faut construire un r\'esum\'e de ces votes. G\'en\'eralement il existe dans $\Upsilon$ des valeurs extr\^emes qui avantagent outrageusement la d\'ecision $d=1$ (resp. $d=0$), on ne peut alors pas s'appuyer sur ces deux types de votes extr\^emes pour essayer de prendre une d\'ecision. Il est par contre difficile de d\'efinir un vote ``neutre" comme nous l'avons fait pour les hypoth\`eses stables avec les votes les plus favorables sous $\Theta_0$ ou sous $\Theta_1$. Pour cela il faudrait faire intervenir un a priori sur une valeur m\'ediane dans $\Upsilon$. 
Il nous reste la possibilit\'e de d\'efinir un vote moyen, $Q^t(.,\upsilon)$ prenant d'autant plus d'importance que $\upsilon$ est probable. Si la r\'ealisation $\omega$ ne donne aucune information sur $\upsilon$ on ne peut que faire une moyenne des votes \`a partir d'une loi a priori $\tau$ sur $\Upsilon$. Dans ce cas il est aussi possible de travailler avec les densit\'es $q_\theta(\omega)=\int p_{(\theta,\upsilon)}(\omega)\,d\tau(\upsilon)$ si elles existent et forment un probl\`eme de d\'ecision expertisable. 
\dli En fait, bien souvent la r\'ealisation $\omega$ contient des informations sur le param\`etre $\upsilon$. M\^eme si l'utilisateur ne veut rien savoir sur $\upsilon$ il est int\'eressant, pour r\'esumer les votes $Q^t(.,\upsilon)$, de 
tenir compte des r\'esultats d'expertises sur $\Upsilon$. On peut m\^eme parfois obtenir, \`a partir de votes compatibles, une probabilisation de $\Upsilon$ muni d'une tribu rendant $Q^t(.,\upsilon)$ mesurable en $\upsilon$. Pour que ces expertises sur $\Upsilon$ soient int\'eressantes il faut qu'elles donnent des informations ne d\'ependant pas de $\theta$. Ceci nous conduit \`a la d\'efinition suivante.
\vfill\eject
\medskip
{\bf D\'efinition 6.1.2}
\medskip
\medskip
\moveleft 10.4pt\hbox{\vrule\kern 10pt\vbox{\defpro
Soit $(\Omega ,{\cal A},(P_{(\theta,\upsilon)}=p_{(\theta,\upsilon)}.\mu)_{(\theta,\upsilon)\in(\Theta_0\cup\Theta_1)\times\Upsilon})$ un probl\`eme de d\'ecision \`a hypoth\`eses stables par rapport \`a la statistique r\'eelle $T$ de loi image $P^T_{(\theta,\upsilon)}$. Le param\`etre fant\^ome est dit expertisable s'il existe une statistique r\'eelle $U$ et, pour tout $\upsilon\in\Upsilon$, une transition $\Pi_\upsilon$ de $(\IR,{\cal B})$ dans $(\IR,{\cal B})$ de la forme $\Pi_\upsilon(t,du)=q_\upsilon(t,u).\nu_t(du)$ telle que : 

i) $P^{(T,U)}_{(\theta,\upsilon)}=P^T_{(\theta,\upsilon)}\Pi_\upsilon$ 
c'est-\`a-dire : \dli $\forall\, B_1\!\times\!B_2\in{\cal B}^2\quad P^{(T,U)}_{(\theta,\upsilon)}(B_1\!\times\!B_2)=\int_{B_1}[\int_{B_2}q_\upsilon(t,u)\,d\nu_t(u)]\,dP^T_{(\theta,\upsilon)}(t)$

ii) pour presque tout $t$, la famille des densit\'es $(q_\upsilon(t,u))_{\upsilon\in\Upsilon}$ est \`a rapport de vraisemblance monotone par rapport \`a $u$ et \`a un ordre $\preceq$ sur $\Upsilon$.
}}\medskip

Dans ce cas de figure la r\'ealisation de $t=T(\omega)$ d\'efinit les votes des experts et permet de choisir un r\'esum\'e $Q^t(.,\upsilon)$ de ces votes dans chacun des sous mod\`eles index\'es par $\upsilon\in\Upsilon$. L'information contenue dans le mod\`ele image de $U$ conditionnellement \`a $T=t$ peut servir \`a probabiliser $\Upsilon$ pour prendre la moyenne des $Q^t(.,\upsilon)$.
En effet ce mod\`ele $(\IR,{\cal B},(\Pi_\upsilon(t,.))_{\upsilon\in\Upsilon})$ est \`a rapport de vraisemblance monotone et on peut choisir des votes compatibles pour les hypoth\`eses unilat\'erales. Bien entendu il est pr\'ef\'erable que la statistique $(T,U)$ soit exhaustive et que $T$ soit une statistique essentielle pour chacune des hypoth\`eses $\{\Theta_1\times\upsilon,\Theta_0\times\upsilon\}$.
Ceci est en particulier v\'erifi\'e dans certains mod\`eles exponentiels pour les hypoth\`eses unilat\'erales. Nous les \'etudierons au paragraphe suivant.
\medskip
{\bf Exemple :} analyse de variance \`a effets fixes.

Dans de nombreux cas la statistique $U$ de la d\'efinition 6.1.2 est ind\'ependante de $T$, la transition $\Pi_\upsilon$ est alors la probabilit\'e image de $U$. Reprenons l'exemple 5 de 5.4. Nous y avons \'etudi\'e le param\`etre de non centralit\'e $\lambda$ (ou $\lambda^2$) des lois de Fisher. Dans le cas de l'analyse de variance \`a effets fixes ce param\`etre exprime la quantit\'e qui int\'eresse l'utilisateur dans une unit\'e qui est la valeur inconnue $\sigma$ de l'\'ecart-type. Nous allons supprimer cette r\'ef\'erence \`a $\sigma$ en travaillant sur le param\`etre $\theta=\sigma^2\lambda^2$. La statistique de Fisher d\'ecentr\'ee est construite \`a partir de deux statistiques ind\'ependantes $T$ et $U$ ; $T/\sigma^2$ suit une loi de khi-deux d\'ecentr\'e \`a $k$ degr\'es de libert\'e et de param\`etre de non centralit\'e $\theta$ ;  $U/\sigma^2$ suit une loi de khi-deux \`a $l$ degr\'es de libert\'e (cf. [Sch.] p. 38). On retrouve aussi cette situation dans les mod\`eles de r\'egression multiple pour les hypoth\`eses lin\'eaires (cf. [Mon2] p. 261, [Leh.] p. 370). La r\'eduction du probl\`eme de base \`a l'\'etude du mod\`ele engendr\'e par $(T,U)$ se fait en imposant des propri\'et\'es d'invariance. Le mod\`ele image s'\'ecrit 
$((\IR^+)^2,{\cal B}^2,(P^T_{(\theta,\sigma^2)}\otimes P^U_{\sigma^2})_{(\theta,\sigma^2)\in\IR^+\times\IR^+_*})$, $P^T_{(\theta,\sigma^2)}$ est une loi $\gamma({k\over 2},2\sigma^2)$ d\'ecentr\'ee de $\theta/(2\sigma^2)$ et $P^U_{\sigma^2}$ une loi 
$\gamma({l\over 2},2\sigma^2)$.
\dli Plus g\'en\'eralement nous allons \'etudier l'exemple des probl\`emes de d\'ecision de la forme : 
$((\IR^+)^2,{\cal B}^2,(P^T_{(\theta,\upsilon)}\otimes P^U_{\upsilon})_{(\theta,\upsilon)\in\{\Theta_1=[0,\theta_1],\Theta_0=]\theta_1,+\infty[\}\times\IR^+_*})$ o\`u $P^T_{(\theta,\upsilon)}$ est une loi $\gamma(p,\upsilon)$ d\'ecentr\'ee de $\theta/\upsilon$ et $P^U_{\upsilon}$ une loi $\gamma(q,\upsilon)$. Les hypoth\`eses $\{\Theta_1,\Theta_0\}$ \`a param\`etre fant\^ome $\upsilon\in\IR^+_*$ sont stables puisque pour $\upsilon$ fix\'e le probl\`eme de d\'ecision se r\'eduit \`a 
$(\IR^+,{\cal B},(P^T_{(\theta,\upsilon)})_{\theta\in\{\Theta_1=[0,\theta_1],\Theta_0=]\theta_1,+\infty[\}})$ donc au choix entre deux hypoth\`eses unilat\'erales dans un mod\`ele \`a rapport de vraisemblance strictement monotone par rapport \`a l'identit\'e (cf. [Kar.]). La statistique identit\'e est une statistique essentielle globale puisque le rapport de vraisemblance est strictement monotone. Notons $G_{(\theta,\upsilon)}$ sa fonction de r\'epartition moyenne, c'est la fonction de r\'epartition de $P^T_{(\theta,\upsilon)}$ qui a pour densit\'e : 
\dli $f_{(\theta,\upsilon)}(t)=\int_{\IN}{1\over\Gamma(p+m) \upsilon}\,({t\over\upsilon})^{p+m-1}\, exp(-{t\over\upsilon})\,d{\cal P}_{\theta\over\upsilon}(m)$ o\`u ${\cal P}_{\theta\over\upsilon}$ est une loi de poisson de param\`etre ${\theta\over\upsilon}\geq 0$ (cf. [Bar.] p. 82).
\dli $G_{(\theta,\upsilon)}\!(t)=\int_{\IN}\Gamma(p+m,\upsilon,t)\,d{\cal P}_{\theta\over\upsilon}(m)$, le terme $\Gamma(p+m,\upsilon,t)$ d\'esignant la valeur de la fonction de r\'epartition d'une loi $\gamma(p+m,\upsilon)$ en $t$ : \dli $\Gamma(p+m,\upsilon,t)=\int_0^t{1\over\Gamma(p+m) \upsilon}\,({x\over\upsilon})^{p+m-1}\, exp(-{x\over\upsilon})\,dx=\Gamma(p+m,1,{t\over\upsilon})$.
\dli Les hypoth\`eses $\Theta_1\times\{\upsilon\}$ et $\Theta_0\times\{\upsilon\}$ sont \'evidemment adjacentes, les votes les plus favorables sous $\Theta_1\times\{\upsilon\}$ et $\Theta_0\times\{\upsilon\}$ sont donc \'egaux au vote des experts sous $P_{(\theta_1,\upsilon)}$. Si l'on ne veut ou peut pas faire intervenir une information a priori ce vote jouera le r\^ole de $Q^t(.,\upsilon)$ qui sera alors d\'efini par : 
\dli $Q^t(1,\upsilon)=Q^t([0,\theta_1],\upsilon)=1-G_{(\theta_1,\upsilon)}\!(t)$.
\dli Nous allons maintenant probabiliser $\Upsilon=\IR^+_*$ en utilisant la statistique $U$, ceci nous permettra d'obtenir une moyenne $Q^{(t,u)}([0,\theta_1])$ des votes $Q^t([0,\theta_1],\upsilon)$. Le mod\`ele image de $U$ conditionnellement \`a $T=t$ est d\'efini par la transition constante $\Pi_\upsilon(t,.)=P^U_\upsilon=\gamma(q,\upsilon)$. Nous avons \'etudi\'e ce mod\`ele dans les exemples 2 de 5.4. Si on n'a pas d'information a priori sur le param\`etre $\upsilon$, la proposition 5.2.3 nous fournit des votes compatibles sur les hypoth\`eses unilat\'erales pour chaque r\'ealisation $U=u>0$. Ils se prolongent en une probabilit\'e sur les bor\'eliens de $\Upsilon=\IR^+_*$ qui est l'inverse d'une loi $\gamma(q,{1\over u})$.
\dli La moyenne des votes $Q^t([0,\theta_1],\upsilon)$ est alors \'egale \`a : 

\cleartabs
\+ $Q^{(t,u)}([0,\theta_1])$&=&$1-\int_{\IR^+_*}G_{(\theta_1,\upsilon)}\!(t)
{1\over\Gamma(q) u}\,({u\over\upsilon})^{q+1}\, exp(-{u\over\upsilon})\,d\upsilon$ \cr
\+ &=&$1-\int_{\IR^+_*}\int_{\IN}[\int_0^t{1\over\Gamma(p+m) \upsilon}\,({x\over\upsilon})^{p+m-1}\, exp(-{x\over\upsilon})\,dx]\,d{\cal P}_{\theta_1\over\upsilon}(m)\times$& \cr
\+ & &\qquad\hfill ${1\over\Gamma(q) u}\,({u\over\upsilon})^{q+1}\, exp(-{u\over\upsilon})\,d\upsilon$& \cr
\+ &=&$1-\sum_{m=0}^{+\infty}\int_0^t[\int_{\IR^+_*}{exp(-{\theta_1\over\upsilon})\,({\theta_1\over\upsilon})^m\over m!\Gamma(p+m) \upsilon}\,({x\over\upsilon})^{p+m-1}\, exp(-{x\over\upsilon})\,{1\over\Gamma(q) u}\,\times$& \cr
\+ & &\qquad\hfill $({u\over\upsilon})^{q+1}\, exp(-{u\over\upsilon})\,d\upsilon]\,dx$& \cr
\+ &=&$1-\sum_{m=0}^{+\infty}\int_0^t[\int_{\IR^+_*}{(\theta_1)^m\,u^q\,x^{p+m-1}\over m!\Gamma(p+m) \Gamma(q)}\,({1\over\upsilon})^{p+q+2m+1}\,exp(-{\theta_1+u+x\over\upsilon}) \,d\upsilon]\,dx$& \cr
\+ &=&$1-\sum_{m=0}^{+\infty}\int_0^t{\Gamma(p+q+2m)\over m!\Gamma(p+m)\Gamma(q)}\,{(\theta_1)^m\,u^q\,x^{p+m-1}\over (\theta_1+u+x)^{p+q+2m}}\,dx$& \cr
\+ &=&$1-\sum_{m=0}^{+\infty}\int_0^{t\over \theta_1+u}{\Gamma(p+q+2m)\over m!\Gamma(p+m)\Gamma(q)}\,{(\theta_1)^m\,u^q\over (\theta_1+u)^{q+m}}{y^{p+m-1}\over (1+y)^{p+q+2m}}\,dy$& \cr
\+ &=&$1-\sum_{m=0}^{+\infty}{\Gamma(q+m)\over m!\Gamma(q)}\,{(\theta_1)^m\,u^q\over (\theta_1+u)^{q+m}}\,F(p+m,q+m,{t\over \theta_1+u})$& \cr
\dli $F(p+m,q+m,{t\over \theta_1+u})$ \'etant la valeur en ${t\over \theta_1+u}$ de la fonction de r\'epartition d'une loi $\beta(p+m,q+m)$ sur $\IR^+$.
Pour $\theta_1=0$ on obtient une masse $Q^{(t,u)}([0])=1-F(p,q,{t\over u})$.
Dans le cas de l'analyse de variance, $p={k\over 2}$ et $q={l\over 2}$, cette masse est \'egale \`a la probabilit\'e qu'une loi de Fisher, de degr\'es de libert\'e $k$ et $l$, soit sup\'erieure \`a ${t/k\over u/l}$. C'est le seuil minimum de rejet du test de $H_0 : \theta=0$ contre $H_1 : \theta>0$. Nous avions d\'ej\`a obtenu ce type de r\'esultat pour le param\`etre ${\theta\over\sigma^2}$ dans l'exemple 5 de 5.4.
\dli Ce qui pr\'ec\`ede nous permet de d\'efinir une probabilit\'e $Q^{(t,u)}$ sur $\Theta=\IR^+$ comme m\'elange des probabilit\'es $Q^t(.,\upsilon)$ d\'efinies par les votes compatibles $(Q^t([0,\theta_f],\upsilon))_{\theta_f\in\IR^+}$. Pour tout $\upsilon\in\Upsilon$, la compatibilit\'e de ces votes d\'ecoule de la proposition 5.2.3 ; la fonction de r\'epartition $G_{(\theta,\upsilon)}(t)$ est \'evidemment continue en $\theta$ et elle tend bien vers $0$ lorsque $\theta$ tend vers $+\infty$ puisque $\Gamma(p+m,\upsilon,t)=\Gamma(p+m,1,{t\over\upsilon})$ tend vers $0$ lorsque $m$ tend vers $+\infty$ (voir l'exemple 5 de 5.4).
Cette probabilit\'e sur $\Theta$ compl\`ete l'information donn\'ee par le r\'esultat du test classique de $H_0 : \theta=0$ contre $H_1 : \theta>0$. Le seuil minimum de rejet $\alpha_m({t\over u})$ de ce test est vu comme la fr\'equence $Q^{(t,u)}([0])$ des experts qui votent pour l'hypoth\`ese $H_0 : \theta=0$. 
Dans le cas du non rejet de $H_0$, c'est-\`a-dire lorsque $\alpha_m({t\over u})$ est sup\'erieur au seuil choisi, on peut pr\'eciser cette r\'eponse en regardant si $Q^{(t,u)}([0,\theta])$ se rapproche rapidement de $1$ lorsque $\theta$ cro\^{\i}t. C'est plus simple que d'analyser la fonction puissance du test. Le cas du rejet de $H_0$ pose probl\`eme lorsque l'hypoth\`ese $H_0 : \theta=0$ est une id\'ealisation de l'hypoth\`ese r\'eelle \`a tester. Bien souvent l'utilisateur se demande si $\theta$ est petit et non pas si $\theta$ est nul. Une interpr\'etation trop rapide du rejet peut conduire \`a consid\'erer que $\theta$ est notable alors qu'il est n\'egligeable. Une analyse de la croissance de $Q^{(t,u)}([0,\theta])$ lorsque $\theta$ s'\'eloigne de $0$ permet d'\'eviter facilement ce pi\`ege. On peut tout simplement porter un jugement sur la valeur $\theta_1$ (resp. $\theta_2$) qui donne une fr\'equence de votes $Q^{(t,u)}([0,\theta_1])$ (resp. $Q^{(t,u)}([0,\theta_2])$) consid\'er\'ee comme petite (resp. grande). Il est cependant plus satisfaisant d'essayer de traduire l'hypoth\`ese ``$\theta$ est petit'' par $\theta\in[0,\theta_0]$ et de porter un jugement \`a partir de $Q^{(t,u)}([0,\theta_0])$.
\vfill\eject

{\parindent=-5mm 6.2 MOD\`ELES EXPONENTIELS.}
\medskip
Consid\'erons un mod\`ele exponentiel $(\Omega ,{\cal A},(p_{(\theta,\upsilon)}.\mu)_{(\theta,\upsilon)\in\Theta\times\Upsilon})$ de la forme $p_{(\theta,\upsilon)}(\omega)=exp(f(\theta,\upsilon)T(\omega)+g(\upsilon)U(\omega)-\psi(\theta,\upsilon))$. 
\dli La fonction $f(\theta,\upsilon)$ (resp. $g(\upsilon)$) est suppos\'ee strictement croissante en $\theta$ (resp. $\upsilon$) sur l'intervalle $\Theta$ (resp. $\Upsilon$) de $\IR$.
\dli Les hypoth\`eses unilat\'erales $\{\Theta_1,\Theta_0\}$ \`a param\`etre fant\^ome 
$\upsilon\in\Upsilon$ sont stables par rapport \`a la statistique r\'eelle $T$, qui est m\^eme une statistique essentielle globale pour chacun des mod\`eles \`a rapport de vraisemblance monotone : 
$(\Omega ,{\cal A},(p_{(\theta,\upsilon)}.\mu)_{\theta\in\Theta})$. En effet, si $\theta_1<\theta_2$ on a pour tout $\upsilon\in\Upsilon$ : 
\dli $p_{(\theta_2,\upsilon)}(\omega)/p_{(\theta_1,\upsilon)}(\omega)=exp[(f(\theta_2,\upsilon)-f(\theta_1,\upsilon))T(\omega)].exp[\psi(\theta_1,\upsilon)-\psi(\theta_2,\upsilon)]$, qui est une fonction strictement croissante de $T$.
\dli Les lois de $T$ forment une famille exponentielle, ainsi que les lois conditionnelles de $U$ quand $T=t$, cette deuxi\`eme famille ne d\'ependant que de $\upsilon$ (cf. [Mon.2] p. 60 et 62). Plus pr\'ecis\'ement, il existe sur $\IR$ des mesures $\mu_\upsilon$ et $\nu_t$ telles que : 
\dli $P_{(\theta,\upsilon)}^T(dt)=exp[f(\theta,\upsilon)t-\psi(\theta,\upsilon)].\mu_\upsilon(dt)$ et 
\dli $P_{(\theta,\upsilon)}^{U\mid t}(du)=exp[g(\upsilon)u-\psi_t(\upsilon)].\nu_t(du)$. 
\dli Les hypoth\`eses unilat\'erales $\{\Theta_1,\Theta_0\}$ d\'efinissent donc un probl\`eme de d\'ecision \`a param\`etre fant\^ome expertisable, par rapport \`a la statistique $U$ et pour l'ordre ordinaire sur $\Upsilon\subseteq\IR$ (voir la d\'efinition 6.1.2). Il suffit de poser 
$q_\upsilon(t,u)=exp[g(\upsilon)u-\psi_t(\upsilon)]$.
Dans le sous probl\`eme de d\'ecision correspondant \`a la valeur $\upsilon$ du param\`etre fant\^ome, $Q_{(\theta,\upsilon)}^t(1)$ est le vote en faveur de $\Theta_1$ sous $P_{(\theta,\upsilon)}$ lorsqu'on r\'ealise $t=T(\omega)$. Notons $G_{(\theta,\upsilon)}$ la fonction de r\'epartition moyenne de la statistique essentielle $T$, elle permet de d\'efinir le vote pr\'ec\'edent puisque $D_i=D_s=\emptyset$ (voir la proposition 4.3.1) et on a : 
\cleartabs 
\+ $Q_{(\theta,\upsilon)}^t(1)=1-G_{(\theta,\upsilon)}(t)$&=&${1\over 2}exp[f(\theta,\upsilon)t-\psi(\theta,\upsilon)]\mu_\upsilon(\{t\})+$\cr
\+ & & $\int\II_{]t,+\infty[}\!(x)exp[f(\theta,\upsilon)x-\psi(\theta,\upsilon)]\,d\mu_\upsilon(x)$.\cr
\dli Les votes $(Q_{(\theta,\upsilon)}^t(1))_{\theta\in\Theta}$ sont alors r\'esum\'es en un vote $Q^t(1,\upsilon)$. Il y a plusieurs choix possibles (voir le paragraphe 4.3), on peut utiliser une pond\'eration de ces votes ou choisir le vote le plus favorable sous $\Theta_0$ (resp. $\Theta_1$). Ces deux derniers votes sont g\'en\'eralement \'egaux, il suffit par exemple que les densit\'es $p_{(\theta,\upsilon)}(\omega)$ soient continues en $\theta$ (voir les commentaires sur la proposition 5.1.1). Dans ce cas on a : 
$Q^t(1,\upsilon)=Q_{(\theta_1,\upsilon)}^t(1)$ avec $\theta_1=sup\Theta_1=inf\Theta_0$.
\dli Le param\`etre fant\^ome \'etant expertisable nous allons utiliser le mod\`ele image de $U$ conditionnellement \`a $T=t$, $(\IR,{\cal B},(exp[g(\upsilon)u-\psi_t(\upsilon)].\nu_t(du))_{\upsilon\in\Upsilon})$, pour  pro\-ba\-bi\-liser $\Upsilon$ \`a partir d'un ensemble de votes compatibles sur l'ensemble des hypoth\`eses unilat\'erales. La proposition 5.2.3 permet souvent de d\'efinir une probabilit\'e $\Lambda^{(t,u)}$, sur les bor\'eliens de $\Upsilon$, \`a partir des votes les plus favorables. Cette probabilit\'e est particuli\`erement int\'eressante car elle ne repose sur aucune information suppl\'ementaire concernant le param\`etre fant\^ome. Les conditions d'application de la proposition 5.2.3 portent ici sur la fonction de r\'epartition de la statistique identit\'e 
\dli $G_\upsilon^{U\mid t}(u)={1\over 2}exp[g(\upsilon)u-\psi_t(\upsilon)]\nu_t(\{u\})+\int\II_{]-\infty,u[}\!(y)exp[g(\upsilon)y-\psi_t(\upsilon)]\,d\nu_t(y)$. 
\dli Le point important est qu'elle soit continue en $\upsilon$ pour presque tout $t$, c'est en particulier le cas lorsque $g(\upsilon)$ est continue.
\dli Le vote final est alors : $Q^{(t,u)}(1)=\int_\Upsilon Q^t(1,\upsilon)\,d\Lambda^{(t,u)}(\upsilon)$.
\dli Nous allons expliciter cette solution sur quelques exemples classiques. 
\medskip
{\bf Exemple 1} : moyenne d'une loi normale.
\dli Soit $(X_1,X_2,...,X_n)$ un n-\'echantillon de la loi $N(\mu,\sigma^2)$, $\mu$ et $\sigma^2$ inconnus ($n>1$). Consid\'erons les hypoth\`eses $\{\mu\leq\mu_0\}$ et $\{\mu>\mu_0\}$.
\dli Posons $\theta=\mu-\mu_0$ et $\upsilon=\sigma^2$, en travaillant sur les variables $Y_i=X_i-\mu_0$ on doit traiter le probl\`eme de d\'ecision suivant : $(\IR^n,{\cal B},(P_{(\theta,\upsilon)})_{(\theta,\upsilon)\in\{\Theta_1\cup\Theta_0\}\times\IR^+_*})$ 
\dli avec $\Theta_1=\{\theta\leq 0\}$, $\Theta_0=\{\theta>0\}$ et 
$P_{(\theta,\upsilon)}$ de densit\'e 
\cleartabs
\+ 
$p_{(\theta,\upsilon)}(y_1,...,y_n)$&=&$({1\over\sqrt{2\pi\upsilon}})^n exp[-{1\over 2\upsilon}\sum_{i=1}^n(y_i-\theta)^2]$\cr 
\+ &=&$exp[{\theta\over\upsilon}(\sum_{i=1}^n y_i)-{1\over 2\upsilon}(\sum_{i=1}^ny_i^2)-{n\theta^2\over 2\upsilon}-{n\over2}ln(2\pi\upsilon)]$\cr
\dli par rapport \`a la mesure de Lebesgue sur $\IR^n$. C'est bien un cas particulier du mod\`ele exponentiel \'etudi\'e dans ce paragraphe, $T(y_1,...,y_n)=\sum_{i=1}^n y_i$, $U(y_1,...,y_n)=\sum_{i=1}^n y_i^2$, $f(\theta,\upsilon)={\theta\over\upsilon}$, $g(\upsilon)=-{1\over2\upsilon}$ et $\psi(\theta,\upsilon)={n\over2}[{\theta^2\over \upsilon}+ln(2\pi\upsilon)]$.
$T$ est de loi $N(n\theta,n\upsilon)$ de densit\'e $p_{(\theta,\upsilon)}^{T}$, quant \`a la statistique ${U\over\upsilon}$ elle suit une loi de khi-deux \`a $n$ degr\'es de libert\'e et \`a param\`etre de non centralit\'e ${n\theta^2}$. Pour obtenir la loi conditionnelle de $U$ quand $T=t$, not\'ee $P_{(\theta,\upsilon)}^{U\mid t}$, nous allons nous servir de l'ind\'ependance entre $T$ et $Z(y_1,...,y_n)={1\over\upsilon}\sum_{i=1}^n(y_i-{\overline y})^2={U\over\upsilon}-{T^2\over n\upsilon}$ qui suit une loi de khi-deux \`a $n-1$ degr\'es de libert\'e. $(T,Z)$ admet donc, par rapport \`a la mesure de Lebesgue sur $\IR^2$, la densit\'e : 
$p_{(\theta,\upsilon)}^{T}(t)\times\II_{\IR^+}\!\!(z){1\over 2^{(n-1)/2}\Gamma({n-1\over 2})}z^{{n-1\over 2}-1}\,e^{-{z\over 2}}$.
Le changement de variable qui associe $(t,u=\upsilon z+{t^2\over n})\!\in\!\IR^2$ \`a $(t,z)\!\in\!\IR^2$ conduit \`a l'expression suivante de la densit\'e de $(T,U)$ : \dli  $p_{(\theta,\upsilon)}^{T}(t)\times\II_{\IR^+}\!\!({u\over\upsilon}-{t^2\over n\upsilon}){1\over 2^{(n-1)/2}\Gamma({n-1\over 2})\upsilon}({u\over\upsilon}-{t^2\over n\upsilon})^{{n-1\over 2}-1}\,e^{-{u\over2\upsilon}+{t^2\over 2n\upsilon}}$.
\dli La densit\'e de la loi $P_{(\theta,\upsilon)}^{U\mid t}$ est alors donn\'ee par : 
\cleartabs 
\+ $p_{(\theta,\upsilon)}^{U\mid t}(u)$&=&$[\II_{\IR^+}\!\!({u\over\upsilon}-{t^2\over n\upsilon})({u\over\upsilon}-{t^2\over n\upsilon})^{{n-3\over 2}}\,e^{-{u\over2\upsilon}}]/\int\II_{\IR^+}\!\!({u\over\upsilon}-{t^2\over n\upsilon})({u\over\upsilon}-{t^2\over n\upsilon})^{{n-3\over 2}}\,e^{-{u\over2\upsilon}}\,du$\cr
\+ &=&$[\II_{[{t^2\over n},+\infty[}\!\!(u)(u-{t^2\over n})^{{n-3\over 2}}\,e^{-{u\over2\upsilon}}]/\int^{+\infty}_0x^{{n-1\over 2}-1}\,e^{-{x\over2\upsilon}}\,e^{-{t^2\over2n\upsilon}}\,dx$\cr
\+ &=&$\II_{[{t^2\over n},+\infty[}\!\!(u){1\over \Gamma({n-1\over 2})(2\upsilon)^{n-1\over 2}}(u-{t^2\over n})^{{n-1\over 2}-1}\,e^{-{1\over2\upsilon}(u-{t^2\over n})}$\cr 
\dli C'est une loi $\gamma({n-1\over 2},2\upsilon)$ translat\'ee de ${t^2\over n}$. Comme nous l'avons fait dans l'exemple 2 de 5.4, on peut obtenir une probabilit\'e $\Lambda^{(t,u)}$ sur $(\Upsilon=\IR^+,{\cal B})$ en appliquant la proposition 5.2.3, puisque la fonction de r\'epartition moyenne de 
$P_{(\theta,\upsilon)}^{U\mid t}$ s'\'ecrit : 
\+ $G_\upsilon^{U\mid t}(u)$&=&$\int_{-\infty}^u\II_{[{t^2\over n},+\infty[}\!\!(x){1\over \Gamma({n-1\over 2})(2\upsilon)^{n-1\over 2}}(x-{t^2\over n})^{{n-1\over 2}-1}\,e^{-{1\over2\upsilon}(x-{t^2\over n})}\,dx$ \cr
\+ &=&$\Gamma({n-1\over 2},(u-{t^2\over n})/2\upsilon)=\Gamma({n-1\over 2},{(n-1)s^2\over 2\upsilon})$\cr
\dli Lorsque $s^2>0$, $\Lambda^{(t,u)}$ est la loi inverse de la loi $\gamma({n-1\over 2},{2\over (n-1)s^2})$, c'est-\`a-dire que le param\`etre ${1\over\upsilon}$ suit une loi de densit\'e $\II_{\IR^+}\!\!(\lambda){(n-1)s^2\over2\Gamma({n-1\over 2})}({(n-1)s^2\over2}\lambda)^{{n-1\over 2}-1}\,e^{-{(n-1)s^2\over2}\lambda}$. On retrouve la loi a post\'eriori correspondant a la loi a priori impropre et non informative de densit\'e : 
${1\over\upsilon}\II_{\IR^+_*}\!(\upsilon)$ (cf. [Ber.] p. 289). Nous aurions obtenu le m\^eme r\'esultat en travaillant avec la statistique exhaustive $(T,(n-1)S^2=\upsilon Z)$ dont les composantes sont ind\'ependantes (voir l'exemple de 6.1). 
\dli Dans le sous probl\`eme de d\'ecision correspondant \`a la valeur $\upsilon$ du param\`etre fant\^ome, le vote le plus favorable sous $\Theta_0$ est \'egal au vote le plus favorable sous $\Theta_1$. C'est le vote  $Q_{(0,\upsilon)}^t$ correspondant \`a la probabilit\'e fronti\`ere $P_{(0,\upsilon)}$. Si l'on choisit ce vote dans chacun des sous probl\`emes de d\'ecision et qu'on r\'ealise $t=T(y_1,...,y_n)$, la famille des votes en faveur de $\Theta_1=\{\theta\leq 0\}$, donc en faveur de $\{\mu\leq\mu_0\}$, est d\'efinie par : \dli $Q^t(1,\upsilon)=Q_{(0,\upsilon)}^t(1)=1-P_{(0,\upsilon)}^{T}(]-\infty,t[)=1-F({t\over\sqrt{n\upsilon}})$, 
\dli $F$ \'etant la fonction de r\'epartition de la loi $N(0,1)$.
\dli La moyenne de ces votes par rapport \`a la probabilit\'e $\Lambda^{(t,u)}$ est alors presque partout ($s^2>0$) \'egale \`a : 
\cleartabs
\+ $Q^{(t,u)}(1)$&=&$\int[1-F({t\over\sqrt{n\upsilon}})]\,d\Lambda^{(t,u)}(\upsilon)$\cr
\+ &=&$\int_{\IR^+}[1-F({t\over\sqrt{n}}\sqrt{\lambda})]{(n-1)s^2\over2\Gamma({n-1\over 2})}({(n-1)s^2\over2}\lambda)^{{n-1\over 2}-1}\,e^{-{(n-1)s^2\over2}\lambda}\,d\lambda$\cr
\+ &=&$\int_{\IR^+}[\int^{+\infty}_{t\over\sqrt{n}}{\sqrt{\lambda}\over\sqrt{2\pi}}e^{-{x^2\over2}\lambda}\,dx]{(n-1)s^2\over2\Gamma({n-1\over 2})}({(n-1)s^2\over2}\lambda)^{{n-1\over 2}-1}\,e^{-{(n-1)s^2\over2}\lambda}\,d\lambda$\cr
\+ &=&$\int_{\IR^+}[\int^{+\infty}_{t\over\sqrt{ns^2}}{\sqrt{s^2\lambda}\over\sqrt{2\pi}}e^{-{y^2\over2}s^2\lambda}\,dy]{(n-1)s^2\over2\Gamma({n-1\over 2})}({(n-1)s^2\over2}\lambda)^{{n-1\over 2}-1}\,e^{-{(n-1)s^2\over2}\lambda}\,d\lambda$\cr
\+ &=&$\int^{+\infty}_{t\over\sqrt{ns^2}}[\int_0^{+\infty}{\sqrt{\nu}\over\sqrt{2\pi}}e^{-{y^2\over2}\nu}\times{({n-1\over 2})^{n-1\over 2}\over\Gamma({n-1\over 2})}\nu^{{n-1\over 2}-1}\,e^{-{n-1\over2}\nu}\,d\nu]\,dy$\cr
\dli L'int\'egrale entre crochets donne la densit\'e d'une loi de Student \`a $(n-1)$ degr\'es de libert\'e comme m\'elange des lois $N(0,{1\over\nu})$ par la loi $\gamma({n-1\over2},{2\over n-1})$ (cf. [Dic.]). Le vote $Q^{(t,u)}(0)$ est donc \'egal \`a la valeur de la fonction de r\'epartition d'un Student $(n-1)$ en ${t\over\sqrt{ns^2}}=\sqrt{n}{\overline{y}\over\sqrt{s^2}}=\sqrt{n}{(\overline{x}-\mu_0)\over\sqrt{s^2}}$. C'est le seuil minimum de rejet du test de Student de $H_0 : \{\mu>\mu_0\}$ contre $H_1 : \{\mu\leq\mu_0\}$. De m\^eme 
$Q^{(t,u)}(1)$ est le seuil minimum de rejet du test de Student de $H_0 : \{\mu\leq\mu_0\}$ contre $H_1 : \{\mu>\mu_0\}$.
\dli Lorsque $\mu_0$ parcourt $\IR$, les votes pr\'ec\'edents sont \'evidemment compatibles. Ils d\'efinissent, sur l'espace $\IR$ du param\`etre $\mu$, une probabilit\'e qui est une loi de Student \`a $(n-1)$ degr\'es de libert\'e, de moyenne $\overline{x}$ et de param\`etre d'\'echelle ${s^2\over n}$ (cf. [Dic.] ou [Ber.] p. 561). La loi de $\mu^2$ n'est pas celle trouv\'ee dans l'exemple d'analyse de variance trait\'e au paragraphe 6.1. Il n'y a rien d'\'etonnant puisque dans ce dernier cas on traite des hypoth\`eses unilat\'erales pour $\mu^2$ donc bilat\'erales pour $\mu$. L'ordre qui int\'eresse l'utilisateur est diff\'erent.

\medskip
{\bf Exemple 2} : comparaison de deux fr\'equences \`a partir d'\'echantillons ind\'ependants.
\dli Soient $X_1$ et $X_2$ deux variables ind\'ependantes de lois binomiales 
${\cal B}(n_1,p_1)$ et ${\cal B}(n_2,p_2)$, $0<p_1<1$ et $0<p_2<1$.
Consid\'erons les hypoth\`eses $\{p_1\leq p_2\}$ et $\{p_1>p_2\}$. Le mod\`ele statistique image de $(X_1,X_2)$ est d\'efini sur $\Omega=\{0,...,n_1\}\times\{0,...,n_2\}$, les param\`etres $(p_1,p_2)$ appartenant \`a $]0,1[\times]0,1[$. Par rapport \`a la mesure de masse 
$C_{n_1}^{x_1}C_{n_2}^{x_2}$ en $(x_1,x_2)\in\Omega$, il admet des densit\'es de forme exponentielle : $exp[x_1ln({p_1\over1-p_1})+x_2ln({p_2\over1-p_2})+n_1ln(1-p_1)+n_2ln(1-p_2)]$.
\dli Posons $\theta=ln({p_1\over1-p_1})-ln({p_2\over1-p_2})=ln({p_1(1-p_2)\over(1-p_1)p_2})\in\IR$ et $\upsilon=ln({p_2\over1-p_2})\in\IR$. 
\dli La fonction $ln({p\over1-p})$ \'etant strictement croissante en $p$, les hypoth\`eses pr\'ec\'edentes s'\'ecrivent : $\{\theta\leq0\}$ et $\{\theta>0\}$. La densit\'e devient : 
\dli $exp[\theta x_1+\upsilon(x_1+x_2)-n_1ln(1+e^{\theta+\upsilon})-n_2ln(1+e^{\upsilon})]$, elle est de la forme \'etudi\'ee avec $T(x_1,x_2)=x_1$ et $U(x_1,x_2)=x_1+x_2$.
La loi $P_{(\theta,\upsilon)}^{T}$ est une binomiale ${\cal B}(n_1,p_1={exp(\theta+\upsilon)\over1+exp(\theta+\upsilon)})$. D'apr\`es l'exemple 4 de 5.4 les votes les plus favorables dans $\Theta_1=\{\theta\leq 0\}\times\{\upsilon\}$ ou dans $\Theta_0=\{\theta>0\}\times\{\upsilon\}$ sont identiques au vote sous $P_{(0,\upsilon)}^T={\cal B}(n_1,p_2={exp(\upsilon)\over1+exp(\upsilon)})$. Si l'on fait ce choix on a : 
\dli $Q^t(1,\upsilon)=Q_{(0,\upsilon)}^t(1)={1\over2}F(t,n_1+1-t,{exp(\upsilon)\over1+exp(\upsilon)})+{1\over2}F(t+1,n_1-t,{exp(\upsilon)\over1+exp(\upsilon)})$, 
\dli $F(p,q,x)$ est la valeur en $x$ de la fonction de r\'epartition de la loi $\beta(p,q)$ sur $[0,1]$ lorsque $p>0$ et $q>0$, pour $p=0$ (resp. $q=0$), c'est-\`a-dire $t=0$ (resp. $t=n_1$), c'est la fonction de r\'epartition de la masse de Dirac en $0$ (resp. $1$).
\dli Nous allons maintenant nous servir de la loi $P_{(\theta,\upsilon)}^{U\mid t}$ qui ne d\'epend que de $\upsilon$, pour probabiliser $\Upsilon=\IR$ et consid\'erer la moyenne des votes pr\'ec\'edents.
$P_{(\theta,\upsilon)}^{U\mid t}$ est une loi binomiale
${\cal B}(n_2,p_2={exp(\upsilon)\over1+exp(\upsilon)})$ translat\'ee de $t$. Pour $u\in\{t,...,t+n_2\}$ sa fonction de r\'epartition moyenne est \'egale \`a : 
\dli $G_{\upsilon}^{U\mid t}(u)=\sum_{i=0}^{u-t}C_{n_2}^ip_2^i(1-p_2)^{n_2-i}-{1\over2}C_{n_2}^{u-t}p_2^{u-t}(1-p_2)^{n_2-u+t}$. 
\dli Comme dans l'exemple 4 de 5.4 on obtient pour le param\`etre $p_2\in]0,1[$ la loi ${1\over2}\beta(u-t,n_2+1-u+t)+{1\over2}\beta(u-t+1,n_2-u+t)$ avec la convention $\beta(0,n_2+1)=\delta_0$ et $\beta(n_2+1,0)=\delta_1$. L'expression de la loi de $\upsilon\in\IR$ est inutile pour faire la moyenne des votes $Q^t(1,\upsilon)$ puisqu'ils ne d\'ependent que de $p_2={exp(\upsilon)\over1+exp(\upsilon)})$. 
La moyenne des votes en faveur de $\{p_1\leq p_2\}$, $Q^{(t,u)}(1)$, correspond aux observations $x_1=t$ et $x_2=u-t$. Pour les observations $n_1-x_1$ et $n_2-x_2$ les votes sont invers\'es, on trouve  $Q^{(t,u)}(1)$ comme moyenne des votes en faveur de $\{p_1>p_2\}$, car $F(q,p,x)=1-F(p,q,1-x)$. C'est une propri\'et\'e classique des proc\'edures de s\'election entre deux binomiales (cf. [DhaM]).
\dli Remarquons enfin que le vote moyen $Q^{(t,u)}(1)$ trouv\'e est \'egal \`a la probabilit\'e de l'\'ev\'enement $\{p_1\leq p_2\}$ lorsque l'espace des param\`etres $[0,1]^2$ est muni de la probabilit\'e produit obtenue en probabilisant s\'epar\'ement les param\`etres $p_1$ et $p_2$ \`a partir des r\'ealisations ind\'ependantes $x_1$ et $x_2$ (voir l'exemple 4 de 5.4). Ceci nous donne une solution pour d'autres types d'hypoth\`eses construites \`a partir de $p_1$ et $p_2$. D\`es que ces hypoth\`eses d\'ependent de $p_1$ et $p_2$ elles ne sont d'ailleurs pas stables puisque le choix entre $(p_1,p_2)$ et $(p'_1,p_2)$ d\'epend de $X_1$ alors que le choix entre $(p_1,p_2)$ et $(p_1,p'_2)$ d\'epend de $X_2$.
Cette remarque reste vraie pour toutes les hypoth\`eses qui font intervenir des param\`etres expertisables \`a partir de statistiques ind\'ependantes.

\vfill\eject

\centerline{\soustitre ANNEXE I.}
\vskip 2cm
\bigskip

Soient $P_0$ et $P_1$, deux probabilit\'es d\'efinies sur $(\Omega,\cal A)$.
Elles admettent toujours des densit\'es $p_0$ et $p_1$ par rapport \`a une mesure $\mu$ sur $(\Omega,\cal A)$ (cf. [Leh.] p. 74).
En 2.2 nous avons d\'efini une statistique $K$ \`a valeurs dans $\overline{\IR^+}$ qui est \'egale au rapport $p_0/p_1$ quand il est d\'efini, c'est-\`a-dire dans le compl\'ementaire de $A=\{\omega\in\Omega\, ;\, p_0(\omega)=0\, et \,  p_1(\omega)=0\}$.
Rappelons que cette statistique v\'erifie :
\halign{$#$&#$\ =\ $&$#$\hfill\cr
\{K=k\}&&\{\omega\in\Omega\, ;\, p_0(\omega)=k.p_1(\omega)\}\bigcap A^c $ 
pour $k\in\IR^+,\cr
\{K=\infty\}&&\{\omega\in\Omega\, ;\, p_1(\omega)=0\}.\cr}

La forme ind\'etermin\'ee $0/0$ prend ici la valeur $+\,\infty$. 
Tout autre statistique \'egale au rapport $p_0/p_1$ quand il est d\'efini, est  $P_0$ et $P_1$ presque s\^urement \'egale \`a $K$. En fait,
cette statistique est unique, $P_0$ et $P_1$ presque s\^urement, au sens suivant : elle ne d\'epend pas de la mesure et des densit\'es choisies pour exprimer $P_0$ et $P_1$.
\medskip
{\leftskip=15mm \dli {\bf D\'emonstration}
\medskip
Soient $p'_0$ et $p'_1$ des densit\'es de $P_0$ et $P_1$ par rapport \`a une mesure $\mu'$ sur $(\Omega,\cal A)$ et $K'$ la statistique associ\'ee au rapport $p'_0/p'_1$. Nous devons montrer que $K$ et $K'$ sont \'egales $P_0$ et $P_1$ presque s\^urement.
C'est-\`a-dire que les \'ev\'enements $B=\{K>K'\}$ et $C=\{K<K'\}$ sont de probabilit\'es nulles pour $P_0$ et $P_1$.

D\'emontrons que $P_1(B)$ est nulle. Il est \'equivalent d'avoir $P_1(B')=0$ avec $B'=B\cap\{p_1>0\}$. Sur $B'$, la statistique $K$ est finie, on a $K'<K<+\,\infty$ et donc $p_0=K.p_1$, $p'_0=K'.p'_1$ ; ceci entra\^{\i}ne :
\dli $P_0(B')=\int_{B'}p_0\,d\mu=\int_{B'}K.p_1\,d\mu=\int_{B'}K\,dP_1$ et
\dli $P_0(B')=\int_{B'}p'_0\,d\mu'=\int_{B'}K'.p'_1\,d\mu'=\int_{B'}K'\,dP_1$ ;
\dli on en d\'eduit  $\int_{B'}(K-K')\,dP_1=0$ et donc  $P_1(B')=0$ puisque $(K-K')>0$ sur $B'$.

On obtient de m\^eme $P_0(B)=0$ en d\'emontrant que $B''=B\cap\{p'_0>0\}$ est de probabilit\'e nulle sous $P_0$. Sur $B''$, la statistique $K'$ est strictement positive, on a $0<K'<K\leq +\infty$ et donc $p'_1=p'_0/K'$, $p_1=p_0/K$ ; ceci entra\^{\i}ne :
\dli $P_1(B'')=\int_{B''}p_1\,d\mu=\int_{B''}(p_0/K)\,d\mu=\int_{B''}(1/K)\,dP_0$ et
\dli $P_1(B'')=\int_{B''}p'_1\,d\mu'=\int_{B''}(p'_0/K')\,d\mu'=\int_{B''}(1/K')\,dP_0$ ;
\dli on en d\'eduit  $\int_{B''}(1/K')-(1/K)\,dP_0=0$ et donc  $P_0(B'')=0$ puisque $(1/K')-(1/K)>0$ sur $B''$.

De fa\c con semblable on d\'emontre : $P_1(C)=0$ et $P_0(C)=0$. Il faut simplement prendre $C'=C\cap\{p'_1>0\}$ afin que les statistiques $K'$ et $K$ soient finies ; pour l'autre cas on pose $C''=C\cap\{p_0>0\}$, ce qui rend $K$ et $K'$ strictement positives.

\medskip\centerline{\hbox to 3cm{\bf \hrulefill}}\par}

\vfill\eject

\centerline{\soustitre ANNEXE II.}
\vskip 2cm
\bigskip

\dli Dans un probl\`eme de d\'ecision \`a hypoth\`eses stables 
$(\Omega ,{\cal A},(P_\theta=p_\theta.\mu)_{\theta\in\Theta_0\cup \Theta_1})$, 
on suppose l'existence d'une statistique r\'eelle $T$ 
et de fonctions croissantes 
$h''_{(\theta_0,\theta_1)} :\, \IR \rightarrow \overline{\IR^+}$
v\'erifiant
$p_{\theta_0}/p_{\theta_1}=h''_{(\theta_0,\theta_1)}(T)$
quand ce rapport n'est pas ind\'etermin\'e (voir la d\'efinition 4.1.2).
Son \'etude est facilit\'ee lorsque les fonctions $h''_{(\theta_0,\theta_1)}$ sont normalis\'ees (voir la d\'efinition 4.2.1.).
Nous allons montrer que ceci est toujours possible. Partant d'une famille
$\{h''_{(\theta_0,\theta_1)}\}_{(\theta_0,\theta_1)\in\Theta_0\times\Theta_1}$
nous allons en construire une normalis\'ee :
$\{h_{(\theta_0,\theta_1)}\}_{(\theta_0,\theta_1)\in\Theta_0\times\Theta_1}$.
\medskip
\dli \souli{1\up{\`ere} \'etape}.

Soit $(\theta_0,\theta_1)\in\Theta_0\times\Theta_1$. Nous allons commencer par construire une fonction $h'_{(\theta_0,\theta_1)}$ constante 
sur tout intervalle $I\!\subset\IR$ d\'efinissant un \'ev\'enement $T^{-1}(I)=\{\omega\in\Omega\, ;\,T(\omega)\in I\}$ 
sur lequel le rapport $p_{\theta_0}/p_{\theta_1}$ est ind\'etermin\'e :
$p_{\theta_0}(\omega)=p_{\theta_1}(\omega)=0$.
Ceci revient \`a trouver $h'_{(\theta_0,\theta_1)}$ constante sur chacun des intervalles maximaux de $N_{(\theta_0,\theta_1)}=\{t\in\IR\, ;\,\forall\omega\in T^{-1}(t)\  p_{\theta_0}(\omega)=p_{\theta_1}(\omega)=0\}$. Notons ${\cal I}_m$ l'ensemble de ces intervalles (certains pouvant \^etre non born\'es).
La fonction $h'_{(\theta_0,\theta_1)}$ est \'egale \`a $h''_{(\theta_0,\theta_1)}$
en dehors de $N_{(\theta_0,\theta_1)}$ et sur tout intervalle $I$ de ${\cal I}_m$ elle est \'egale \`a une constante $b_I$. Pour que $h'_{(\theta_0,\theta_1)}$ soit croissante cette constante $b_I$ doit v\'erifier : 
\dli $b_I\geq
a_I=sup\{h''_{(\theta_0,\theta_1)}(t)\, ;\, t<I\}\cup\{0\}$ et
$b_I\leq c_I=inf\{h''_{(\theta_0,\theta_1)}(t)\, ;\, t>I\}\cup\{+\infty\}$. \dli La fonction $\{h''_{(\theta_0,\theta_1)}\}$ n'ayant \'et\'e modifi\'ee que sur $N_{(\theta_0,\theta_1)}$, $h'_{(\theta_0,\theta_1)}(T)$ est encore \'egale au rapport des densit\'es $p_{\theta_0}/p_{\theta_1}$, sur le domaine de d\'efinition de ce rapport.
\medskip
\dli \souli{2\up{\`eme} \'etape}.

Afin d'obtenir une famille normalis\'ee 
$\{h_{(\theta_0,\theta_1)}\}_{(\theta_0,\theta_1)\in\Theta_0\times\Theta_1}$
nous allons modifier les fonctions $h'_{(\theta_0,\theta_1)}$ pr\'ec\'edentes sur $D_i$ et $D_s$. Rappelons que  $D_i=]-\infty,t)$ (resp. $D_s=(t,+\infty[$) d\'esigne la plus grande demi-droite ouverte ou ferm\'ee telle que $T^{-1}(D_i)\subseteq\cap_{\theta\in\Theta_0}\{p_\theta=0\}$
(resp. $T^{-1}(D_s)\subseteq\cap_{\theta\in\Theta_1}\{p_\theta=0\}$).
Pour avoir les propri\'et\'es ii) et iii) de la d\'efinition 4.2.1 nous posons :
$$h_{(\theta_0,\theta_1)}(t)\quad=\quad\left\{\matrix{
0\hfill & si & t\in D_i-D_s\hfill\cr
h'_{(\theta_0,\theta_1)}(t)\hfill & si & t\in\ \IR-(D_i\cup D_s)\hfill \cr
+\infty\hfill & si & t\in D_s\hfill \cr
}\right. $$
$h_{(\theta_0,\theta_1)}$ est bien une fonction croissante et $h_{(\theta_0,\theta_1)}(T)$ est toujours \'egale \`a $p_{\theta_0}/p_{\theta_1}$ sur le domaine de d\'efinition de ce rapport.
En effet, pour tout $\omega$ de $T^{-1}(D_i)$ (resp. $T^{-1}(D_s)$) on a
$p_{\theta_0}(\omega)=0$ (resp. $p_{\theta_1}(\omega)=0$) et le rapport
$p_{\theta_0}(\omega)/p_{\theta_1}(\omega)$ est soit ind\'etermin\'e soit \'egal \`a $0$ (resp. $+\infty$) ; bien entendu, s'il existe $\omega$ appartenant \`a $D_i\cap D_s$ le rapport des densit\'es est ind\'etermin\'e.
\dli Les fonctions $h_{(\theta_0,\theta_1)}$ v\'erifient aussi les propri\'et\'es de la d\'efinition 4.2.1, elles d\'efinissent donc une famille normalis\'ee.

Les propri\'et\'es ii) et iii) des familles normalis\'ees conduisent \`a deux lemmes utiles pour les d\'emonstrations des propositions 4.2.1 et 4.2.2.

\medskip
{\bf Lemme 1}
\nobreak
\medskip
\medskip
\moveleft 10.4pt\hbox{\vrule\kern 10pt\vbox{\defpro

Soit $(\Omega ,{\cal A},(p_\theta.\mu)_{\theta\in\Theta_0\cup\Theta_1})$ un 
probl\`eme de d\'ecision \`a hypoth\`eses stables. On note $\Delta_s$ l'ensemble des fonctions de test simples d\'efinies \`a partir d'une famille
$\{h_{(\theta_0,\theta_1)}\}_{(\theta_0,\theta_1)\in \Theta_0\times\Theta_1}$ normalis\'ee.
\dli Pour tout $\theta_0\in\Theta_0$ ; on a $sup\{\phi_{(\infty,0)}^{(\theta_0,\theta_1)}\}_{\theta_1\in\Theta_1}=f_{(t_1,u_1)}=sup\Delta_s$
}}\medskip

{\leftskip=15mm \dli {\bf D\'emonstration}
\medskip

Soit $\theta_0\in\Theta_0$. Posons $f_{(t,u)}= sup\{\phi_{(\infty,0)}^{(\theta_0,\theta_1)}\}_{\theta_1\in\Theta_1}$ et
notons $D_t$ la demi-droite $[t,+\infty[$ (resp. $]t,+\infty[$) si $u=0$ 
(resp. $u=1$). 
\dli Lorsque $(t,u)=(+\infty,0)$ on a \'evidemment $f_{(t,u)}=sup\Delta_s$.
Dans le cas contraire $D_t$ est non vide et est inclus dans $D_s$ ( voir la d\'efinition 4.2.1)
puisque pour tout  $\theta'\in\Theta_1$ et tout $t'\in D_t$, 
$\phi_{(\infty,0)}^{(\theta_0,\theta')}$ est nulle sur $\{T=t'\}$, ce qui implique $h_{(\theta_0,\theta')}(t')=+\infty$, donc $p_{\theta'}=0$ sur $\{T=t'\}$.

Par d\'efinition de $f_{(t_1,u_1)}=sup\Delta_s$ on a $f_{(t,u)}\leq f_{(t_1,u_1)}$. On doit montrer l'\'egalit\'e. 
\dli Si on avait $f_{(t,u)}<f_{(t_1,u_1)}$, il existerait
$(\theta'_0,\theta'_1)\in\Theta_0\times\Theta_1$ tel que :
$f_{(t,u)}<\phi_{(\infty,0)}^{(\theta'_0,\theta'_1)}\leq f_{(t_1,u_1)}$ ;
l'\'ev\'enement $A=\{\phi_{(\infty,0)}^{(\theta'_0,\theta'_1)}-f_{(t,u)}=1\}$
serait non vide et pour $\omega'\in A$ on aurait
$h_{(\theta'_0,\theta'_1)}(T(\omega'))<+\infty$ avec $t'=T(\omega')$ appartenant \`a $D_t\subseteq D_s$ ce qui est impossible puisque sur $D_s$ la 
fonction normalis\'ee $h_{(\theta'_0,\theta'_1)}$ vaut $+\infty$.

\medskip\centerline{\hbox to 3cm{\bf \hrulefill}}\par}

\medskip
{\bf Lemme 2}
\medskip
\medskip
\moveleft 10.4pt\hbox{\vrule\kern 10pt\vbox{\defpro

Soit $(\Omega ,{\cal A},(p_\theta.\mu)_{\theta\in\Theta_0\cup\Theta_1})$ un 
probl\`eme de d\'ecision \`a hypoth\`eses stables. On note $\Delta_s$ l'ensemble des fonctions de test simples d\'efinies \`a partir d'une famille
$\{h_{(\theta_0,\theta_1)}\}_{(\theta_0,\theta_1)\in \Theta_0\times\Theta_1}$ normalis\'ee.
\dli Pour tout $\theta_1\in\Theta_1$ ; on a $inf\{\phi_{(0,1)}^{(\theta_0,\theta_1)}\}_{\theta_0\in\Theta_0}=f_{(t_0,u_0)}=inf\Delta_s$
}}\medskip

{\leftskip=15mm \dli {\bf D\'emonstration}
\medskip
La d\'emonstration est semblable \`a celle du lemme 1.
\dli Soit $\theta_1\in\Theta_1$. On pose $f_{(t,u)}= inf\{\phi_{(0,1)}^{(\theta_0,\theta_1)}\}_{\theta_0\in\Theta_0}$ et on
note $D_t$ la demi-droite $]-\infty,t[$ (resp. $]-\infty,t]$) si $u=0$ 
(resp. $u=1$). 
\dli Lorsque $(t,u)=(-\infty,0)$ on a \'evidemment $f_{(t,u)}=inf\Delta_s$.
Dans le cas contraire $D_t$ est non vide et est inclus dans $D_i$ ( voir la d\'efinition 4.2.1)
car pour tout  $\theta\in\Theta_0$ et tout $t'\in D_t$ 
on a $h_{(\theta,\theta_1)}(t')=0$, donc $p_{\theta}=0$ sur $\{T=t'\}$.
On a m\^eme $D_t\subseteq D_i-D_s$ puisque 
$f_{(t,u)}\leq f_{(t_1,u_1)}$ et $\II_{D_s}(T)=1-f_{(t_1,u_1)}$ (voir la proposition 4.2.1).

Par d\'efinition de $f_{(t_0,u_0)}=inf\Delta_s$ on a $f_{(t,u)}\geq f_{(t_0,u_0)}$. On doit montrer l'\'egalit\'e. 
\dli Si on avait $f_{(t,u)}>f_{(t_0,u_0)}$, il existerait
$(\theta'_0,\theta'_1)\in\Theta_0\times\Theta_1$ tel que :
$f_{(t_0,u_0)}\leq\phi_{(0,1)}^{(\theta'_0,\theta'_1)}<f_{(t,u)}$ ;
l'\'ev\'enement $A=\{f_{(t,u)}-\phi_{(0,1)}^{(\theta'_0,\theta'_1)}=1\}$
serait non vide et pour $\omega'\in A$ on aurait
$h_{(\theta'_0,\theta'_1)}(T(\omega'))>0$ avec $t'=T(\omega')$ appartenant \`a $D_t\subseteq D_i-D_s$ ; ceci est impossible puisque sur $D_i - D_s$ la 
fonction normalis\'ee $h_{(\theta'_0,\theta'_1)}$ vaut $0$.
\nobreak
\medskip\centerline{\hbox to 3cm{\bf \hrulefill}}\par}

\vfill\eject

\centerline{\soustitre ANNEXE III.}
\vskip 2.0cm
\bigskip

\dli Soit $(\Omega ,{\cal A},(P_\theta=p_\theta.\mu)_{\theta\in\Theta_0\cup \Theta_1=\Theta})$,  un probl\`eme de d\'ecision \`a hypoth\`eses stables par rapport \`a la statistique r\'eelle $T$. D'apr\`es la d\'efinition 4.1.2, ceci implique  l'existence,  pour tout couple $(\theta_0,\theta_1)\in\Theta_0\times\Theta_1$, d'une fonction croissante 
$h_{(\theta_0,\theta_1)} :\, \IR \rightarrow \overline{\IR^+}$ v\'erifiant
$p_{\theta_0}/p_{\theta_1}=h_{(\theta_0,\theta_1)}(T)$ sur le domaine de d\'efinition de ce rapport.
\dli Les votes $\{Q_\theta\}_{\theta\in\Theta}$ (voir la proposition 4.3.1) ont \'et\'e construits \`a partir d'une famille normalis\'ee :
$\{h_{(\theta_0,\theta_1)}\}_{(\theta_0,\theta_1)\in\Theta_0\times\Theta_1}$
(voir la d\'efinition 4.2.1.). Il existe toujours une famille normalis\'ee 
(voir l'annexe II) mais elle n'est g\'en\'eralement pas unique. Il est l\'egitime de se demander si le choix de cette famille normalis\'ee influence les votes
$\{Q_\theta\}$.

Consid\'erons deux familles normalis\'ees : 
$$\{h_{(\theta_0,\theta_1)}\}_{(\theta_0,\theta_1)\in\Theta_0\times\Theta_1}\quad et \quad \{h'_{(\theta_0,\theta_1)}\}_{(\theta_0,\theta_1)\in\Theta_0\times\Theta_1}$$
\dli Notons respectivement $K(T)$ et $K'(T)$ les statistiques essentielles qu'elles permet\-tent de construire (voir la d\'efinition 4.3.1). Les votes
$\{Q_\theta\}_{\theta\in\Theta}$ (resp. $\{Q'_\theta\}_{\theta\in\Theta}$) ne d\'ependent que de la statistique $K(T)$ (resp. $K'(T)$)(voir la proposition 4.3.1). La fonction croissante $K :\, \IR \rightarrow \overline{\IR}$ (resp. $K'$) est construite \`a partir de 
$\Delta_s=\cup_{(\theta_0,\theta_1)\in\Theta_0\times\Theta_1}\Phi_{s}^{(\theta_0,\theta_1)}$ (resp. $\Delta'_s=\cup_{(\theta_0,\theta_1)\in\Theta_0\times\Theta_1}\Phi_{s'}^{(\theta_0,\theta_1)}$), $\Phi_{s}^{(\theta_0,\theta_1)}$ (resp. $\Phi_{s'}^{(\theta_0,\theta_1)}$) \'etant l'ensemble des fonctions de test simples bas\'ees sur le rapport des densit\'es, $p_{\theta_0}/p_{\theta_1}$, d\'efini par $K_{(\theta_0,\theta_1)}=h_{(\theta_0,\theta_1)}(T)$ (resp. 
$K'_{(\theta_0,\theta_1)}=h'_{(\theta_0,\theta_1)}(T)$) (voir la d\'efinition 2.2.1).

Nous allons commencer par caract\'eriser les r\'eels pour lesquels les valeurs des fonctions normalis\'ees $h_{(\theta_0,\theta_1)}$ et $h'_{(\theta_0,\theta_1)}$ peuvent \^etre diff\'erentes.
La d\'efinition 4.2.1  sur les familles normalis\'ees fait intervenir deux demi-droites $D_i$ et $D_s$.
$D_i=]-\infty,t)$ (resp. $D_s=(t,+\infty[$) d\'esigne la plus grande demi-droite ouverte ou ferm\'ee telle que $T^{-1}(D_i)\subseteq\cap_{\theta\in\Theta_0}\{p_\theta=0\}$
(resp. $T^{-1}(D_s)\subseteq\cap_{\theta\in\Theta_1}\{p_\theta=0\}$).

\medskip
{\bf Lemme 1}
\medskip
\medskip
\moveleft 10.4pt\hbox{\vrule\kern 10pt\vbox{\defpro

Soient $(\theta_0,\theta_1)\in\Theta_0\times\Theta_1$ et $t\in\IR$. 
Un intervalle $I\!\subset\IR$ est dit ind\'etermin\'e pour $(\theta_0,\theta_1)$ lorsque les densit\'es $p_{\theta_0}$ et $p_{\theta_1}$ sont nulles sur $T^{-1}(I)$.
\dli Si $t\in\IR\!-(D_i\cup D_s)$ et si $[t]$ est ind\'etermin\'e pour 
$(\theta_0,\theta_1)$, notons $I_t$ le plus grand intervalle de $\IR\!-(D_i\cup D_s)$ contenant $t$ et ind\'etermin\'e pour $(\theta_0,\theta_1)$ ; $h_{(\theta_0,\theta_1)}$ \'etant une fonction normalis\'ee posons :
$h_{(\theta_0,\theta_1)}(I_t^-)=sup\{h_{(\theta_0,\theta_1)}(x)\, ;\, x<I_t\}\cup\{0\}$ et 
$h_{(\theta_0,\theta_1)}(I_t^+)=inf\{h_{(\theta_0,\theta_1)}(x)\, ;\, x>I_t\}\cup\{+\infty\}$.

1) Il existe une fonction normalis\'ee diff\'erente de $h_{(\theta_0,\theta_1)}$ en $t$ si et seulement si : 
\dli $t\in\IR\!-(D_i\cup D_s)$, $t$ est ind\'etermin\'e pour $(\theta_0,\theta_1)$ et 
$h_{(\theta_0,\theta_1)}(I_t^-)<h_{(\theta_0,\theta_1)}(I_t^+)$
\dli (les limites $h_{(\theta_0,\theta_1)}(I_t^-)$ et $h_{(\theta_0,\theta_1)}(I_t^+)$ ne d\'ependent pas de la fonction normalis\'ee choisie)

2) Si $t\in\IR\!-(D_i\cup D_s)$ est ind\'etermin\'e pour $(\theta_0,\theta_1)$, 
il est totalement ind\'etermin\'e :
$\forall\omega\in T^{-1}(t)\quad \forall\theta\in\Theta\quad p_{\theta}(\omega)=0$, lorsque $h_{(\theta_0,\theta_1)}(I_t^-)>0$ ou $h_{(\theta_0,\theta_1)}(I_t^+)<+\infty$.
}}\medskip

\medskip
{\leftskip=15mm \dli {\bf D\'emonstration}
\medskip
{\parindent=-10mm I --- Condition n\'ecessaire de 1).}

La condition $t\in\IR\!-(D_i\cup D_s)$ est n\'ecessaire puisque pour $t\in D_s$ (resp. $t\in D_i-D_s$) toutes les fonctions normalis\'ees valent $+\infty$ (resp. $0$). C'est une cons\'equence directe des propri\'et\'es ii) et iii) de la d\'efinition 4.2.1.

La condition $t$ ind\'etermin\'e pour $(\theta_0,\theta_1)$, c'est-\`a-dire : 
\dli $\forall\omega\in T^{-1}(t)\quad p_{\theta_0}(\omega)=p_{\theta_1}(\omega)=0$, est aussi n\'ecessaire. En effet, lorsqu'il existe $\omega\in T^{-1}(t)$ tel que $p_{\theta_0}(\omega)>0$ ou $p_{\theta_1}(\omega)>0$ le rapport 
$p_{\theta_0}(\omega)/p_{\theta_1}(\omega)$ est d\'efini, il prend une valeur 
$k\in\overline{\IR^+}$ qui doit \^etre la valeur en $t$ de toute fonction normalis\'ee pour le couple $(\theta_0,\theta_1)$.

La derni\`ere condition a un sens lorsque les deux premi\`eres sont v\'erifi\'ees, on a alors $I_t\not=\emptyset$.
La propri\'et\'e i) de la d\'efinition 4.2.1 impose \`a la fonction normalis\'ee $h_{(\theta_0,\theta_1)}$ d'\^etre constante sur $I_t$.
Lorsque $h_{(\theta_0,\theta_1)}(I_t^-)=h_{(\theta_0,\theta_1)}(I_t^+)$, $h_{(\theta_0,\theta_1)}(t)$ ne peut prendre que cette valeur commune des deux limites pour \^etre croissante. Nous aurons d\'emontr\'e la n\'ecessit\'e de la condition $h_{(\theta_0,\theta_1)}(I_t^-)<h_{(\theta_0,\theta_1)}(I_t^+)$, si nous montrons que ces limites ne d\'ependent pas de la fonction normalis\'ee choisie pour les d\'efinir.

i) $h_{(\theta_0,\theta_1)}(I_t^-)$ ne d\'epend pas de la fonction normalis\'ee choisie pour la d\'efinir.
\dli Si $\{x<I_t\}=D_i-D_s$, $h_{(\theta_0,\theta_1)}(I_t^-)$ ne peut prendre que la valeur $0$, que 
$D_i-D_s$ soit vide ou pas. Dans le cas contraire il existe une suite 
$(x_n)_{n\in\IN}$ croissant vers $inf I_t$ et dont les \'el\'ements n'appartiennent ni \`a $I_t$ ni \`a $D_i-D_s$. Par d\'efinition de $I_t$ on peut m\^eme supposer  que les $x_n$ ne sont pas ind\'etermin\'es pour $(\theta_0,\theta_1)$. D'apr\`es ce qui pr\'ec\`ede, toutes les fonctions normalis\'ees pour $(\theta_0,\theta_1)$ prennent la m\^eme valeur $h_n$ en $x_n$. On a bien s\^ur $lim_{n\rightarrow +\infty}h_n=h_{(\theta_0,\theta_1)}(I_t^-)$.

ii) $h_{(\theta_0,\theta_1)}(I_t^+)$ ne d\'epend pas de la fonction normalis\'ee choisie pour la d\'efinir.
\dli La d\'emonstration est semblable \`a la pr\'ec\'edente. Lorsque $\{x>I_t\}=D_s$, $h_{(\theta_0,\theta_1)}(I_t^+)$ ne peut prendre que la valeur $+\infty$, que 
$D_s$ soit vide ou pas. Dans le cas contraire il existe une suite 
$(x_n)_{n\in\IN}$ d\'ecroissant vers $sup I_t$ et dont les \'el\'ements n'appartiennent ni \`a $I_t$ ni \`a $D_s$. Par d\'efinition de $I_t$ on peut m\^eme supposer  que les $x_n$ ne sont pas ind\'etermin\'es pour $(\theta_0,\theta_1)$. Toutes les fonctions normalis\'ees pour $(\theta_0,\theta_1)$ prennent alors la m\^eme valeur $h_n$ en $x_n$ et $lim_{n\rightarrow +\infty}h_n=h_{(\theta_0,\theta_1)}(I_t^+)$.
\medskip
{\parindent=-10mm II -- Condition suffisante de 1).}

Les trois conditions n\'ecessaires \'etant r\'ealis\'ees, nous devons trouver une fonction normalis\'ee $h'_{(\theta_0,\theta_1)}$ dont la valeur en $t$ est diff\'erente de $h_{(\theta_0,\theta_1)}(t)$. Consid\'erons :
$$\kern 13mm h'_{(\theta_0,\theta_1)}(x)\,=\,\left\{\matrix{
h_{(\theta_0,\theta_1)}(x)\hfill &si& x\notin I_t\hfill\cr
c\in[h_{(\theta_0,\theta_1)}(I_t^-),h_{(\theta_0,\theta_1)}(I_t^+)]-\{h_{(\theta_0,\theta_1)}(t)\}\hfill &si& x\in I_t\hfill \cr
}\right. $$
Comme $h_{(\theta_0,\theta_1)}(I_t^-)<h_{(\theta_0,\theta_1)}(I_t^+)$, on a $[h_{(\theta_0,\theta_1)}(I_t^-),h_{(\theta_0,\theta_1)}(I_t^+)]-\{h_{(\theta_0,\theta_1)}(t)\}\not=\emptyset$ et la fonction $h'_{(\theta_0,\theta_1)}$ est bien d\'efinie. Elle est croissante par d\'efinition de $h_{(\theta_0,\theta_1)}(I_t^-)$ et $h_{(\theta_0,\theta_1)}(I_t^+)$ ; elle est \'egale \`a $p_{\theta_0}/p_{\theta_1}$ sur le domaine de d\'efinition de ce rapport, mais diff\'erente de $h_{(\theta_0,\theta_1)}$ en $t$ et m\^eme sur $I_t$.
\dli Il reste \`a v\'erifier que $h'_{(\theta_0,\theta_1)}$ est normalis\'ee.
\dli Elle poss\`ede les propri\'et\'es ii) et iii) de la d\'efinition 4.2.1 puisque $h_{(\theta_0,\theta_1)}$ est normalis\'ee et que l'intervalle $I_t$ est d'intersection vide avec $D_i\cup D_s$.
\dli D'autre part, tout intervalle $I\subseteq\IR\!-(D_i\cup D_s)$ ind\'etermin\'e pour $(\theta_0,\theta_1)$ v\'erifie : $I\cap I_t=\emptyset$ ou $I\subseteq I_t$, par d\'efinition de $I_t$ ; 
$h'_{(\theta_0,\theta_1)}$ poss\`ede donc aussi la propri\'et\'e i).

\medskip
{\parindent=-10mm III - $t$ est totalement ind\'etermin\'e lorsque $h_{(\theta_0,\theta_1)}(I_t^-)>0$ ou \dli $h_{(\theta_0,\theta_1)}(I_t^+)<+\infty$.}

$t\in\IR\!-(D_i\cup D_s)$ est pris ind\'etermin\'e pour $(\theta_0,\theta_1)$ : 
\dli $\forall\omega\in T^{-1}(t)\quad p_{\theta_0}(\omega)=p_{\theta_1}(\omega)=0$. 

Nous allons d\'emontrer le r\'esultat recherch\'e, c'est-\`a-dire : \dli $\forall\omega\in T^{-1}(t)\quad \forall\theta\in\Theta\quad p_{\theta}(\omega)=0$, en raisonnant par l'absurde.
\dli Supposons qu'il existe $\omega'\in T^{-1}(t)$ et $\theta'\in\Theta$ tels que $p_{\theta'}(\omega')>0$. Nous allons montrer qu'on aurait $h_{(\theta_0,\theta_1)}(I_t^-)=0$ et $h_{(\theta_0,\theta_1)}(I_t^+)=+\infty$.
Etudions d'abord les cons\'equences de cette supposition suivant que $\theta'$ appartient \`a $\Theta_0$ ou $\Theta_1$.

1\up{er} cas : si $\theta'\in\Theta_0$ alors $h_{(\theta_0,\theta_1)}(I_t^+)=+\infty$.
\dli On a $p_{\theta'}(\omega')/p_{\theta_1}(\omega')=+\infty$ donc
$h_{(\theta',\theta_1)}(t)=+\infty$. $h_{(\theta',\theta_1)}$ est aussi infinie sur  $[t,+\infty[$, on a donc $p_{\theta_1}$ nulle sur 
$\Omega_t=T^{-1}([t,+\infty[)$. Quant au rapport $p_{\theta_0}/p_{\theta_1}$ il est soit ind\'etermin\'e soit \'egal \`a $+\infty$ sur $\Omega_t$.
\dli Montrons que $h_{(\theta_0,\theta_1)}(I_t^+)=inf\{h_{(\theta_0,\theta_1)}(x)\, ;\, x>I_t\}\cup\{+\infty\}$ est \'egal \`a $+\infty$. C'est \'evident lorsque $\{x>I_t\}=D_s$, que $D_s$ soit vide ou pas. Dans le cas contraire il existe une suite $(x_n)_{n\in\IN}$ d\'ecroissant vers $sup I_t$ et dont les \'el\'ements ne sont pas ind\'etermin\'es pour $(\theta_0,\theta_1)$, puisque $I_t$ est un intervalle maximum de $\IR\!-(D_i\cup D_s)$ ind\'etermin\'e pour $(\theta_0,\theta_1)$. Comme $T^{-1}(x_n)\subseteq\Omega_t$ on a
$h_{(\theta_0,\theta_1)}(x_n)=+\infty$ et donc $h_{(\theta_0,\theta_1)}(I_t^+)=+\infty$.

2\up{\`eme} cas : si $\theta'\in\Theta_1$ alors $h_{(\theta_0,\theta_1)}(I_t^-)=0$.
\dli On a $p_{\theta_0}(\omega')/p_{\theta'}(\omega')=0$ donc
$h_{(\theta_0,\theta')}(t)=0$. $h_{(\theta_0,\theta')}$ est aussi \'egale \`a $0$ sur  $]-\infty,t]$, on a donc $p_{\theta_0}$ nulle sur 
$T^{-1}(]-\infty,t])$. Tout \'el\'ement $x\in]-\infty,t]$ qui n'est pas ind\'etermin\'e pour $(\theta_0,\theta_1)$ v\'erifie $h_{(\theta_0,\theta_1)}(x)=0$.
\dli Montrons que $h_{(\theta_0,\theta_1)}(I_t^-)=sup\{h_{(\theta_0,\theta_1)}(x)\, ;\, x<I_t\}\cup\{0\}$ est \'egal \`a $0$. C'est \'evident lorsque $\{x<I_t\}=D_i-D_s$, que $D_i-D_s$ soit vide ou pas. Dans le cas contraire il existe une suite $(x_n)_{n\in\IN}$ croissant vers $inf I_t$ et dont les \'el\'ements ne sont pas ind\'etermin\'es pour $(\theta_0,\theta_1)$, puisque $I_t$ est un intervalle maximum de $\IR\!-(D_i\cup D_s)$ ind\'etermin\'e pour $(\theta_0,\theta_1)$. De plus $x_n\in]-\infty,t]$, on a donc
$h_{(\theta_0,\theta_1)}(x_n)=0$ et bien s\^ur $h_{(\theta_0,\theta_1)}(I_t^-)=0$.

Ces deux cas nous am\`enent \`a la conclusion recherch\'ee : 
\dli $h_{(\theta_0,\theta_1)}(I_t^-)=0$ et $h_{(\theta_0,\theta_1)}(I_t^+)=+\infty$, lorsqu'on a les deux conditions :
\dli 1) il existe $\omega'\in T^{-1}(t)$ et $\theta'\in\Theta_0$ tels que $p_{\theta'}(\omega')>0$
\dli 2) il existe $\omega''\in T^{-1}(t)$ et $\theta''\in\Theta_1$ tels que $p_{\theta''}(\omega'')>0$.
\dli Nous n'en avons suppos\'e qu'une seule vraie au d\'epart, mais nous allons montrer que l'on ne peut pas avoir l'une sans l'autre.
\dli Si on avait la premi\`ere condition sans la seconde, on aurait pour tout $\theta_1$ de $\Theta_1$ : $p_{\theta_1}(\omega')=0$ et $p_{\theta'}(\omega')>0$, donc $h_{(\theta',\theta_1)}(t)=+\infty$ ; ce qui implique : $\forall\theta_1\in\Theta_1$ et $\forall\omega\in T^{-1}([t,+\infty[)$, $p_{\theta_1}(\omega)=0$ ; ceci est impossible puisque 
$t\notin D_s$.
\dli De m\^eme lorsque la deuxi\`eme condition est v\'erifi\'ee, on ne peut pas avoir : $\forall\theta_0\in\Theta_0$ $\forall\omega\in T^{-1}(t)$ $p_{\theta_0}(\omega)=0$, car on aurait $h_{(\theta_0,\theta'')}(t)=0$ et donc : $\forall\theta_0\in\Theta_0$ et $\forall\omega\in T^{-1}(]-\infty,t])$, $p_{\theta_0}(\omega)=0$ ; ce qui est impossible puisque $t\notin D_i$.

\medskip\centerline{\hbox to 3cm{\bf \hrulefill}}\par}

D\'emontrons maintenant que le choix de la famille normalis\'ee n'influence pas les votes $\{Q_\theta\}_{\theta\in\Theta}$.
\medskip
{\bf Proposition 1}
\medskip
\medskip
\moveleft 10.4pt\hbox{\vrule\kern 10pt\vbox{\defpro

Soient $\{h_{(\theta_0,\theta_1)}\}_{(\theta_0,\theta_1)\in\Theta_0\times\Theta_1}$ et $\{h'_{(\theta_0,\theta_1)}\}_{(\theta_0,\theta_1)\in\Theta_0\times\Theta_1}$ deux familles normalis\'ees associ\'ees \`a un probl\`eme de d\'ecision \`a hypoth\`eses stables par rapport \`a la statistique r\'eelle $T$ d\'efinie sur 
$(\Omega ,{\cal A},(P_\theta=p_\theta.\mu)_{\theta\in\Theta})$. Elles d\'efinissent des votes $Q_\theta$ et $Q'_\theta$ identiques en dehors d'une partie $T^{-1}(M)$ de $\Omega$, sur laquelle les densit\'es $p_\theta$ sont toutes nulles
($\forall\omega\notin T^{-1}(M)\ \forall\theta\in\Theta\quad Q_{\theta}^\omega(\{1\})=Q_{\theta}'^\omega(\{1\})$ et
$\forall\omega\in T^{-1}(M)\ \forall\theta\in\Theta\quad p_{\theta}(\omega)=0$).
}}\medskip

\medskip
{\leftskip=15mm \dli {\bf D\'emonstration}
\medskip
Notons $G_\theta$ et $G'_\theta$ les fonctions de r\'epartition moyenne des statistiques essentielles $K(T)$ et $K'(T)$ d\'efinies par les deux familles normalis\'ees (voir la d\'efinition 4.3.1). D'apr\`es la proposition 4.3.1, 
pour chaque r\'ealisation $\omega$, les votes $Q_{\theta}^\omega$ sont d\'efinis par :
$$Q^\omega_\theta(\{1\})=\left\{\matrix{
1\hfill & si & T(\omega)\in D_i-D_s\hfill\cr
1-G_\theta(K(T(\omega)))\hfill & si & T(\omega)\in\IR\!-(D_i\cup D_s)\hfill \cr
0\hfill & si & T(\omega)\in D_s\hfill \cr
}\right. $$
Ces votes, comme les votes $Q_{\theta}'^\omega$ ne d\'ependent que de la valeur $T(\omega)=t$. 

Lorsque $t\in D_i\cup D_s$ on obtient les m\^emes  votes $Q_{\theta}^\omega$ et $Q_{\theta}'^\omega$ puisque $D_i$ et $D_s$ ne d\'ependent que des densit\'es
$p_\theta$.

Il nous reste \`a comparer les votes lorsque $t\in\IR\!-(D_i\cup D_s)\not=\emptyset$, ce qui peut s'\'ecrire :
$\II_{D_i}(T)\leq f_{(t,0)}< f_{(t,1)}\leq 1-\II_{D_s}(T)$.
\dli La valeur de $K(t)$ est alors construite \`a partir de deux experts (voir la d\'efinition 4.3.1) :
\dli $f_{(a_t,u_t)}=sup\{f\in\Delta_s\, ;\, f\leq f_{(t,0)}\}$ et 
$f_{(b_t,v_t)}=inf\{f\in\Delta_s\, ;\, f\geq f_{(t,1)}\}$.
\dli Ceux permettant d'obtenir $K'(t)$ sont not\'es :
\dli $f_{(a'_t,u'_t)}=sup\{f\in\Delta'_s\, ;\, f\leq f_{(t,0)}\}$ et 
$f_{(b'_t,v'_t)}=inf\{f\in\Delta'_s\, ;\, f\geq f_{(t,1)}\}$.
\dli Nous n'avons pas ajout\'e dans les d\'efinitions les \'el\'ements $f_{(-\infty,0)}$ et $f_{(+\infty,1)}$ car lorsque $t\notin D_i\cup D_s$ on a : $f_{(-\infty,1)}\leq inf\Delta_s=\II_{D_i}(T)\leq f_{(t,0)}< f_{(t,1)}\leq 1-\II_{D_s}(T)=sup\Delta_s\leq f_{(+\infty,0)}$ (voir la proposition 4.2.1).
Nous avons aussi remplac\'e $\overline{\Delta_s}$ par $\Delta_s$ (resp. $\overline{\Delta'_s}$ par $\Delta'_s$), ce qui ne change rien puisque $]f_{(t,0)},f_{(t,1)}[$ ne contient aucun \'el\'ement de $\Delta_s$ (resp. $\Delta'_s$).
\dli La valeur de $K(t)$ est alors donn\'e par :
$$K(t)\quad=\quad\left\{\matrix{
b_t-1\hfill & si & (a_t,u_t)=(-\infty,1)\hfill\cr
[a_t+b_t]/2\hfill & si & -\infty<a_t\leq b_t<+\infty\hfill \cr
a_t+1\hfill & si & (b_t,v_t)=(+\infty,0)\ et\ a_t>-\infty\hfill \cr
}\right. $$
\dli On obtient la d\'efinition de $K'(t)$ en rempla\c cant $(a_t,u_t)$ et $(b_t,v_t)$ par $(a'_t,u'_t)$ et $(b'_t,v'_t)$. Posons

\dli $\displaystyle\hfill N=\{t\in\IR\ ;\ \forall\omega\in T^{-1}(t)\quad \forall\theta\in\Theta\quad p_\theta(\omega)=0\}\hfill$

\dli Nous allons montrer que pour $t$ appartenant \`a $IR\!-(D_i\cup D_s\cup N)$, nous avons les deux \'egalit\'es suivantes :
\dli $\displaystyle\hfill f_{(a_t,u_t)}\relmont{=}{p.s.}f_{(a'_t,u'_t)}\quad et\quad f_{(b_t,v_t)}\relmont{=}{p.s.}f_{(b'_t,v'_t)}\hfill$
\dli Ceci conduit au r\'esultat recherch\'e avec $M\!\subseteq N$. En effet, on a
\dli $\{K(T)=K(t)\}=\{f_{(b_t,v_t)}-f_{(a_t,u_t)}=1\}$ et 
\dli $\{K'(T)=K'(t)\}=\{f_{(b'_t,v'_t)}-f_{(a'_t,u'_t)}=1\}$
\dli (voir le 2\up{\`eme} cas de la d\'emonstration de la proposition 4.3.1), donc pour tout $\theta$ de $\Theta$ :
$$\vbox{\cleartabs
\+ \kern 15mm $G_\theta(K(t))$&$=[E_\theta(f_{(a_t,u_t)})+E_\theta(f_{(b_t,v_t)})]/2$ \cr
\+&$=[E_\theta(f_{(a'_t,u'_t)})+E_\theta(f_{(b'_t,v'_t)})]/2=G'_\theta(K'(t))$ \cr
}$$

{\parindent=-10mm I --- $f_{(a_t,u_t)}\relmont{=}{p.s.}f_{(a'_t,u'_t)}$ pour $t\in\IR\!-(D_i\cup D_s\cup N)$.}

Pour tout couple $(\theta_0,\theta_1)\in\Theta_0\times\Theta_1$ on consid\`ere :
\dli $e_t^{(\theta_0,\theta_1)}=sup\{f\in\Phi_s^{(\theta_0,\theta_1)}\, ;\, f\leq f_{(t,0)}\}\cup\{f_{(-\infty,1)}\}$ et 
\dli $\varphi_t^{(\theta_0,\theta_1)}=sup\{f\in\Phi_{s'}^{(\theta_0,\theta_1)}\, ;\, f\leq f_{(t,0)}\}\cup\{f_{(-\infty,1)}\}$.
\dli Nous devons d\'emontrer :

\dli $\displaystyle\hfill sup\{e_t^{(\theta_0,\theta_1)}\}_{(\theta_0,\theta_1)\in\Theta_0\times\Theta_1}\relmont{=}{p.s.}sup\{\varphi_t^{(\theta_0,\theta_1)}\}_{(\theta_0,\theta_1)\in\Theta_0\times\Theta_1}\hfill$

\dli Une condition suffisante est de pouvoir associer \`a tout couple $(\theta_0,\theta_1)$ un couple $(\theta'_0,\theta'_1)$ tel que :
$e_t^{(\theta_0,\theta_1)}\relmont{\leq}{p.s.}\varphi_t^{(\theta'_0,\theta'_1)}$ et 
$\varphi_t^{(\theta_0,\theta_1)}\relmont{\leq}{p.s.}e_t^{(\theta'_0,\theta'_1)}$.
Ceci est en particulier r\'ealis\'e lorsque
$e_t^{(\theta_0,\theta_1)}\relmont{=}{p.s.}\varphi_t^{(\theta_0,\theta_1)}$.
Nous allons distinguer deux cas.

{\parindent=-5mm 1\up{er} cas : les densit\'es $p_{\theta_0}$ et $p_{\theta_1}$ ne sont pas s\'epar\'ees.}

Les densit\'es $p_{\theta_0}$ et $p_{\theta_1}$ sont dites non s\'epar\'ees lorsque : 
\dli $\IR\!-(D_i^{\theta_0}\cup D_s^{\theta_1})\not=\emptyset$ ; 
$D_i^{\theta_0}=]-\infty,t_i)$ (resp. $D_s^{\theta_1}=(t_s,+\infty[$) \'etant la plus grande demi-droite ouverte ou ferm\'ee pour laquelle la densit\'e 
$p_{\theta_0}$ (resp. $p_{\theta_1}$) est nulle sur $T^{-1}(D_i^{\theta_0})$ (resp. $T^{-1}(D_s^{\theta_1})$).
\dli $t\in\IR\!-(D_i\cup D_s\cup N)$ n'est pas totalement ind\'etermin\'e puisque $t\notin N$. Nous allons montrer qu'il n'est pas ind\'etermin\'e pour $(\theta_0,\theta_1)$.

a) Il existe $\omega'\in T^{-1}(t)$ tel que  $p_{\theta_0}(\omega')>0$ ou $p_{\theta_1}(\omega')>0$.
\dli Nous devons d\'emontrer que si $x\in\IR\!-(D_i\cup D_s)$ est ind\'etermin\'e pour $(\theta_0,\theta_1)$, les densit\'es $p_{\theta_0}$ et $p_{\theta_1}$ \'etant non s\'epar\'ees, $x$ est totalement ind\'etermin\'e. 
\dli Nous allons utiliser le r\'esultat 2) du lemme pr\'ec\'edent. Il suffit de montrer que l'on ne peut pas avoir $h_{(\theta_0,\theta_1)}(I_x^-)=0$ et $h_{(\theta_0,\theta_1)}(I_x^+)=+\infty$. Si cela \'etait, on aurait $h_{(\theta_0,\theta_1)}(y)=0$ (resp. $h_{(\theta_0,\theta_1)}(z)=+\infty$) sur $\{y<I_x\}$ (resp. $\{z>I_x\}$) ; $p_{\theta_0}$ (resp. $p_{\theta_1}$) serait nulle  sur $T^{-1}(\{y<I_x\})$ (resp. $T^{-1}(\{z>I_x\})$) ; 
$p_{\theta_0}$ et $p_{\theta_1}$ \'etant nulles sur $T^{-1}(I_x)$, on aurait 
$D_i^{\theta_0}\cup D_s^{\theta_1}=\IR$ ; ce qui est impossible puisque 
$p_{\theta_0}$ et $p_{\theta_1}$ ne sont pas s\'epar\'ees.

b) $e_t^{(\theta_0,\theta_1)}\relmont{=}{p.s.}\varphi_t^{(\theta_0,\theta_1)}$.
\dli D'apr\`es a) il existe $\omega'\in T^{-1}(t)$ tel que  $p_{\theta_0}(\omega')>0$ ou $p_{\theta_1}(\omega')>0$.
On a alors $h_{(\theta_0,\theta_1)}(t)=h'_{(\theta_0,\theta_1)}(t)=p_{\theta_0}(\omega')/p_{\theta_1}(\omega')=k\in\overline{\IR^+}$.
\dli Lorsque $k=0$ ceci implique 
$e_t^{(\theta_0,\theta_1)}=\varphi_t^{(\theta_0,\theta_1)}=f_{(-\infty,1)}$.
\dli Lorsque $k>0$ on a
$e_t^{(\theta_0,\theta_1)}=\phi_{(k,0)}^{(\theta_0,\theta_1)}$ et $\varphi_t^{(\theta_0,\theta_1)}=\phi_{(k,0)}^{'(\theta_0,\theta_1)}$ (voir la d\'efinition 2.2.1 avec respectivement $K=h_{(\theta_0,\theta_1)}(T)$ et $K=h'_{(\theta_0,\theta_1)}(T)$). Pour obtenir 
$e_t^{(\theta_0,\theta_1)}\relmont{=}{p.s.}\varphi_t^{(\theta_0,\theta_1)}$, on doit d\'emontrer que l'\'ev\'enement $E=\{e_t^{(\theta_0,\theta_1)}\not=\varphi_t^{(\theta_0,\theta_1)}\}$ est n\'egligeable pour tout $P_\theta$. C'est \'evident lorsque $E=\emptyset$. Dans le cas contraire $E=T^{-1}(I)$, $I$ \'etant l'intervalle non vide de $]-\infty,t[$ dont les \'el\'ements appartiennent  \`a un seul des intervalles $\{h_{(\theta_0,\theta_1)}(x)=k\}$ et $\{h'_{(\theta_0,\theta_1)}(x)=k\}$. Tout \'el\'ement $x$ de $I$ appartient \`a  $\IR\!-(D_i\cup D_s)$ et est ind\'etermin\'e pour $(\theta_0,\theta_1)$, sinon les fonctions $h_{(\theta_0,\theta_1)}$ et $h'_{(\theta_0,\theta_1)}$ seraient \'egales en $x$. Nous avons vu en a) que $x$, donc $I$, est totalement ind\'etermin\'e. $E=T^{-1}(I)$ est alors bien n\'egligeable.

{\parindent=-5mm 2\up{\`eme} cas : les densit\'es $p_{\theta_0}$ et $p_{\theta_1}$ sont s\'epar\'ees.}

Dans ce cas $D_i^{\theta_0}\cup D_s^{\theta_1}=\IR$. On a bien s\^ur :
$D_i\subseteq D_i^{\theta_0}$ et $D_s\subseteq D_s^{\theta_1}$.
De plus $D_i-D_s=D_i$ puisque $t\in\IR\!-(D_i\cup D_s\cup N)$, ce qui suppose $D_i\cup D_s\not=\IR$.
\dli Les diff\'erentes positions de  $D_s^{\theta_1}$ par rapport \`a $D_i$ et
de  $D_i^{\theta_0}$ par rapport \`a $D_s$ nous conduisent \`a d\'efinir une partition de $\IR$ en trois intervalles :
\dli $A=D_i\cup(\IR-D_s^{\theta_1})$, $C=D_s\cup(\IR-D_i^{\theta_0})$ et 
$B= \IR\!-(A+C)\subseteq D_i^{\theta_0}\cap D_s^{\theta_1}$ (on a l'\'egalit\'e lorsque $D_i\cap D_s^{\theta_1}=\emptyset$ et $D_s\cap D_i^{\theta_0}=\emptyset$).
\dli Ceci revient \`a poser $\II_A(T)=sup\{\II_{D_i}(T) , 1-\II_{D_s^{\theta_1}}(T)\}$ et
\dli $1-\II_C(T)=inf\{1-\II_{D_s}(T) , \II_{D_i^{\theta_0}}(T)\}$.

a) Etude de $\Phi_s^{(\theta_0,\theta_1)}$ et $\Phi_{s'}^{(\theta_0,\theta_1)}$.

Nous allons d\'emontrer que $\Phi_s^{(\theta_0,\theta_1)}$ et $\Phi_{s'}^{(\theta_0,\theta_1)}$ sont inclus dans $\{\II_A(T),1-\II_C(T)\}$. Pour cela nous allons montrer que les  fonctions $h_{(\theta_0,\theta_1)}$ et $h'_{(\theta_0,\theta_1)}$ sont \'egales \`a $0$ sur $A$, constantes sur $B$ et \'egales \`a $+\infty$ sur $C$. $\Phi_s^{(\theta_0,\theta_1)}$ (resp. $\Phi_{s'}^{(\theta_0,\theta_1)}$) est r\'eduit \`a $\{\II_A(T)\}$ ou $\{1-\II_C(T)\}$ lorsque la fonction normalis\'ee correspondante est \'egale \`a $+\infty$ ou $0$ sur l'intervalle $B$.

i) Si $x\in A$, $h_{(\theta_0,\theta_1)}(x)=h'_{(\theta_0,\theta_1)}(x)=0$.
\dli C'est \'evident lorsque $x\in D_i$. Dans le cas contraire $x$ n'appartient pas \`a $D_s^{\theta_1}$ ; par d\'efinition de $D_s^{\theta_1}$, il existe $y\in[x,+\infty[-D_s^{\theta_1}$ et $\omega_y\in T^{-1}(y)$ tels que 
$p_{\theta_1}(\omega_y)>0$ ; comme $y\in D_i^{\theta_0}$ on a en plus
$p_{\theta_0}(\omega_y)=0$, donc 
$h_{(\theta_0,\theta_1)}(y)=h'_{(\theta_0,\theta_1)}(y)=0$ ; la croissance des deux fonctions normalis\'ees entra\^{\i}ne la m\^eme propri\'et\'e en $x$.

ii) Sur $B$, $h_{(\theta_0,\theta_1)}$ et $h'_{(\theta_0,\theta_1)}$ sont constantes.
\dli $B$ \'etant inclus dans $D_i^{\theta_0}\cap D_s^{\theta_1}$, il est ind\'etermin\'e pour $(\theta_0,\theta_1)$ ; on a aussi $B\subseteq \IR\!-(D_i\cup D_s)$ ; $h_{(\theta_0,\theta_1)}$ et $h'_{(\theta_0,\theta_1)}$ \'etant normalis\'ees, elles sont constantes sur l'intervalle $B$ (voir la propri\'et\'e i) de la d\'efinition 4.2.1).

iii) Si $x\in C$, $h_{(\theta_0,\theta_1)}(x)=h'_{(\theta_0,\theta_1)}(x)=+\infty$.
\dli C'est \'evident lorsque $x\in D_s$. Dans le cas contraire $x$ n'appartient pas \`a $D_i^{\theta_0}$ ; par d\'efinition de $D_i^{\theta_0}$, il existe $y\in]-\infty,x]-D_i^{\theta_0}$ et $\omega_y\in T^{-1}(y)$ tels que 
$p_{\theta_0}(\omega_y)>0$ ; comme $y\in D_s^{\theta_1}$ on a en plus
$p_{\theta_1}(\omega_y)=0$, donc 
$h_{(\theta_0,\theta_1)}(y)=h'_{(\theta_0,\theta_1)}(y)=+\infty$ ; la croissance des deux fonctions normalis\'ees entra\^{\i}ne la m\^eme propri\'et\'e en $x$.

b) $e_t^{(\theta_0,\theta_1)}\relmont{\leq}{p.s.}\varphi_t^{(\theta'_0,\theta'_1)}$ et 
$\varphi_t^{(\theta_0,\theta_1)}\relmont{\leq}{p.s.}e_t^{(\theta'_0,\theta'_1)}$.

Soit $t\in\IR\!-(D_i\cup D_s\cup N)$. Il y a trois possibilit\'es par rapport \`a la partition $\IR=A+B+C$.

i) $t\in A$.
\dli D'apr\`es a-i) on a $h_{(\theta_0,\theta_1)}(t)=h'_{(\theta_0,\theta_1)}(t)=0$, donc \dli $e_t^{(\theta_0,\theta_1)}=\varphi_{t}^{(\theta_0,\theta_1)}=f_{(-\infty,1)}$ ; $(\theta'_0,\theta'_1)=(\theta_0,\theta_1)$ convient.

ii) $t\in B$.
\dli D'apr\`es a), $e_t^{(\theta_0,\theta_1)}$ (resp. $\varphi_{t}^{(\theta_0,\theta_1)}$) est \'egale \`a $\II_A(T)$  ou $f_{(-\infty,1)}$.
\dli Nous avons vu en a-ii) que $t$ appartient \`a $\IR\!-(D_i\cup D_s)$ et qu'il est ind\'etermin\'e pour $(\theta_0,\theta_1)$ ; mais $t$ n'est pas totalement ind\'etermin\'e puisque $t\notin N$ ; d'apr\`es la fin de la partie III de la d\'emonstration du lemme pr\'ec\'edent, il existe $\theta'\in\Theta_0$ et $\omega'\in T^{-1}(t)$ tels que $p_{\theta'}(\omega')>0$.
On a donc $h_{(\theta',\theta_1)}(t)=h'_{(\theta',\theta_1)}(t)=+\infty$, ce qui implique : $e_t^{(\theta',\theta_1)}=\phi_{(\infty,0)}^{(\theta',\theta_1)}$ et $\varphi_t^{(\theta',\theta_1)}=\phi_{(\infty,0)}^{'(\theta',\theta_1)}$ (voir la partie b) du 1\up{er} cas) ; mais pour $x\in A$ on a 
$h_{(\theta',\theta_1)}(x)=h'_{(\theta',\theta_1)}(x)<+\infty$ (suivre le raisonnement fait en a-i) avec $\theta_0=\theta'$ et $p_{\theta'}(\omega_y)\geq 0$), donc $e_t^{(\theta',\theta_1)}\geq\II_A(T)$ et
$\varphi_t^{(\theta',\theta_1)}\geq\II_A(T)$.
\dli Nous venons ainsi de trouver un couple $(\theta',\theta_1)$ tel que :
$e_t^{(\theta_0,\theta_1)}\leq\varphi_t^{(\theta',\theta_1)}$ et 
$\varphi_t^{(\theta_0,\theta_1)}\leq e_t^{(\theta',\theta_1)}$.

iii) $t\in C$.
\dli D'apr\`es a), $e_t^{(\theta_0,\theta_1)}$ (resp. $\varphi_{t}^{(\theta_0,\theta_1)}$) est \'egale \`a $\II_A(T)$ ou $1-\II_C(T)$.
Pour obtenir les in\'egalit\'es recherch\'ees on va trouver un couple $(\theta_0,\theta'_1)$ tel que : $e_t^{(\theta_0,\theta'_1)}\relmont{\geq}{p.s.}1-\II_C(T)$ et $\varphi_{t}^{(\theta_0,\theta'_1)}\relmont{\geq}{p.s.}1-\II_C(T)$.
\dli En fait $t\in C-(D_i\cup D_s\cup N)$, il appartient donc \`a $D_s^{\theta_1}$ sans appartenir \`a $D_s$ ; par d\'efinition de $D_s$, ceci permet de trouver $\theta'_1\in\Theta_1$ tel que $t\notin D_s^{\theta'_1}$.
Comme $t$ n'appartient pas \`a $D_i^{\theta_0}$, il existe $x\in ]-\infty,t]-D_i^{\theta_0}$ et $\omega\in T^{-1}(x)$ tels que $p_{\theta_0}(\omega)>0$ ; ce qui implique : 
$h_{(\theta_0,\theta'_1)}(x)>0$ et $h'_{(\theta_0,\theta'_1)}(x)>0$ ; ces deux fonctions \'etant croissantes on a aussi : 
$h_{(\theta_0,\theta'_1)}(t)=k>0$ et $h'_{(\theta_0,\theta'_1)}(t)=k'>0$, donc 
$e_t^{(\theta_0,\theta'_1)}=\phi_{(k,0)}^{(\theta_0,\theta'_1)}$ et $\varphi_t^{(\theta_0,\theta'_1)}=\phi_{(k',0)}^{'(\theta_0,\theta'_1)}$.
\dli Nous allons d\'emontrer l'in\'egalit\'e 
$e_t^{(\theta_0,\theta'_1)}\relmont{\geq}{p.s.}1-\II_C(T)$ 
\dli (resp. $\varphi_{t}^{(\theta_0,\theta'_1)}\relmont{\geq}{p.s.}1-\II_C(T)$)en montrant que l'intervalle 
\dli $I=\{h_{(\theta_0,\theta'_1)}(y)=k\}-C$ (resp. $I'=\{h'_{(\theta_0,\theta'_1)}(y)=k'\}-C$) est totalement ind\'etermin\'e.
Comme $I$ (resp. $I'$) est inclus dans $D_i^{\theta_0}$, la densit\'e $p_{\theta_0}$ est nulle sur cet intervalle ; $k$ (resp. $k'$) \'etant non nul, cet intervalle est ind\'etermin\'e pour $(\theta_0,\theta'_1)$. Les densit\'es $p_{\theta_0}$ et $p_{\theta'_1}$ n'\'etant pas s\'epar\'ees puisque 
$t\notin D_i^{\theta_0}\cup D_s^{\theta'_1}$, d'apr\`es la partie a) du 1\up{er} cas les intervalles $I$ et $I'$, qui sont inclus dans $\IR\!-(D_i\cup D_s)$ ($I<t<D_s$, $k>0$ et $I'<t<D_s$, $k'>0$), sont totalement ind\'etermin\'es. 

{\parindent=-10mm II -- $f_{(b_t,v_t)}\relmont{=}{p.s.}f_{(b'_t,v'_t)}$ pour $t\in\IR\!-(D_i\cup D_s\cup N)$.}

Cette d\'emonstration est semblable \`a celle de la partie I. Pour tout 
couple $(\theta_0,\theta_1)\in\Theta_0\times\Theta_1$ on consid\`ere :
\dli $g_t^{(\theta_0,\theta_1)}=inf\{f\in\Phi_s^{(\theta_0,\theta_1)}\, ;\, f\geq f_{(t,1)}\}\cup\{f_{(+\infty,0)}\}$ et 
\dli $\psi_t^{(\theta_0,\theta_1)}=inf\{f\in\Phi_{s'}^{(\theta_0,\theta_1)}\, ;\, f\geq f_{(t,1)}\}\cup\{f_{(+\infty,0)}\}$.
\dli Nous devons d\'emontrer :

\dli $\displaystyle\hfill inf\{g_t^{(\theta_0,\theta_1)}\}_{(\theta_0,\theta_1)\in\Theta_0\times\Theta_1}\relmont{=}{p.s.}inf\{\psi_t^{(\theta_0,\theta_1)}\}_{(\theta_0,\theta_1)\in\Theta_0\times\Theta_1}\hfill$

\dli Une condition suffisante est de pouvoir associer \`a tout couple $(\theta_0,\theta_1)$ un couple $(\theta'_0,\theta'_1)$ tel que :
$g_t^{(\theta_0,\theta_1)}\relmont{\geq}{p.s.}\psi_t^{(\theta'_0,\theta'_1)}$ et 
$\psi_t^{(\theta_0,\theta_1)}\relmont{\geq}{p.s.}g_t^{(\theta'_0,\theta'_1)}$.
Ceci est en particulier r\'ealis\'e lorsque
$g_t^{(\theta_0,\theta_1)}\relmont{=}{p.s.}\psi_t^{(\theta_0,\theta_1)}$.
Nous allons encore distinguer deux cas.

{\parindent=-5mm 1\up{er} cas : les densit\'es $p_{\theta_0}$ et $p_{\theta_1}$ ne sont pas s\'epar\'ees.}

On va d\'emontrer : $g_t^{(\theta_0,\theta_1)}\relmont{=}{p.s.}\psi_t^{(\theta_0,\theta_1)}$.
\dli La partie a) du 1\up{er} cas de I reste valable, il existe donc $\omega'\in T^{-1}(t)$ tel que  $p_{\theta_0}(\omega')>0$ ou $p_{\theta_1}(\omega')>0$ et on a encore : \dli $h_{(\theta_0,\theta_1)}(t)=h'_{(\theta_0,\theta_1)}(t)=p_{\theta_0}(\omega')/p_{\theta_1}(\omega')=k\in\overline{\IR^+}$.
\dli Lorsque $k=+\infty$ ceci implique 
$g_t^{(\theta_0,\theta_1)}=\psi_t^{(\theta_0,\theta_1)}=f_{(+\infty,0)}$.
\dli Lorsque $k<+\infty$ on a
$g_t^{(\theta_0,\theta_1)}=\phi_{(k,1)}^{(\theta_0,\theta_1)}$ et $\psi_t^{(\theta_0,\theta_1)}=\phi_{(k,1)}^{'(\theta_0,\theta_1)}$ (voir la d\'efinition 2.2.1 avec respectivement $K=h_{(\theta_0,\theta_1)}(T)$ et $K=h'_{(\theta_0,\theta_1)}(T)$). Pour obtenir 
$g_t^{(\theta_0,\theta_1)}\relmont{=}{p.s.}\psi_t^{(\theta_0,\theta_1)}$, on doit d\'emontrer que l'\'ev\'enement $E=\{g_t^{(\theta_0,\theta_1)}\not=\psi_t^{(\theta_0,\theta_1)}\}$ est n\'egligeable pour tout $P_\theta$. C'est \'evident lorsque $E=\emptyset$. Dans le cas contraire $E=T^{-1}(I)$, $I$ \'etant l'intervalle non vide de $]t,+\infty[$ dont les \'el\'ements appartiennent  \`a un seul des intervalles $\{h_{(\theta_0,\theta_1)}(x)=k\}$ et $\{h'_{(\theta_0,\theta_1)}(x)=k\}$. Tout \'el\'ement $x$ de $I$ appartient \`a  $\IR\!-(D_i\cup D_s)$ et est ind\'etermin\'e pour $(\theta_0,\theta_1)$, sinon les fonctions $h_{(\theta_0,\theta_1)}$ et $h'_{(\theta_0,\theta_1)}$ seraient \'egales en $x$. Nous avons vu que la partie a) du 1\up{er} cas de I reste valable, l'intervalle $I$ est donc totalement ind\'etermin\'e. $E=T^{-1}(I)$ est alors bien n\'egligeable.

{\parindent=-5mm 2\up{\`eme} cas : les densit\'es $p_{\theta_0}$ et $p_{\theta_1}$ sont s\'epar\'ees.}

Comme dans le 2\up{\`eme} cas de la partie I, $\Phi_s^{(\theta_0,\theta_1)}$ et $\Phi_{s'}^{(\theta_0,\theta_1)}$ sont inclus dans $\{\II_A(T),1-\II_C(T)\}$. Nous allons montrer qu'il existe un couple $(\theta'_0,\theta'_1)$ tel que : $g_t^{(\theta_0,\theta_1)}\relmont{\geq}{p.s.}\psi_t^{(\theta'_0,\theta'_1)}$ et $\psi_t^{(\theta_0,\theta_1)}\relmont{\geq}{p.s.} g_t^{(\theta'_0,\theta'_1)}$.

Soit $t\in\IR\!-(D_i\cup D_s\cup N)$. On consid\`ere les trois possibilit\'es de la partition $\IR=A+B+C$ d\'efinie, \`a partir de deux densit\'es s\'epar\'ees, dans le 2\up{\`eme} cas de la partie I.

i) $t\in A$.
\dli $g_t^{(\theta_0,\theta_1)}$ (resp. $\psi_{t}^{(\theta_0,\theta_1)}$) est \'egale \`a $\II_A(T)$ ou $1-\II_C(T)$.
Pour obtenir les in\'egalit\'es recherch\'ees on va trouver un couple $(\theta'_0,\theta_1)$ tel que : $g_t^{(\theta'_0,\theta_1)}\relmont{\leq}{p.s.}\II_A(T)$ et $\psi_{t}^{(\theta'_0,\theta_1)}\relmont{\leq}{p.s.}\II_A(T)$.
\dli En fait $t\in A-(D_i\cup D_s\cup N)$, il appartient donc \`a $D_i^{\theta_0}$ sans appartenir \`a $D_i$ ; par d\'efinition de $D_i$, ceci permet de trouver $\theta'_0\in\Theta_0$ tel que $t\notin D_i^{\theta'_0}$.
Comme $t$ n'appartient pas \`a $D_s^{\theta_1}$, il existe $x\in [t,+\infty,[-D_s^{\theta_1}$ et $\omega\in T^{-1}(x)$ tels que $p_{\theta_1}(\omega)>0$ ; ce qui implique : 
$h_{(\theta'_0,\theta_1)}(x)<+\infty$ et $h'_{(\theta'_0,\theta_1)}(x)<+\infty$ ; ces deux fonctions \'etant croissantes on a aussi : 
$h_{(\theta'_0,\theta_1)}(t)=k<+\infty$ et $h'_{(\theta'_0,\theta_1)}(t)=k'<+\infty$, donc 
$g_t^{(\theta'_0,\theta_1)}=\phi_{(k,1)}^{(\theta'_0,\theta_1)}$ et $\psi_t^{(\theta'_0,\theta_1)}=\phi_{(k',1)}^{'(\theta'_0,\theta_1)}$.
\dli Nous allons d\'emontrer l'in\'egalit\'e 
$g_t^{(\theta'_0,\theta_1)}\relmont{\leq}{p.s.}\II_A(T)$ 
\dli (resp. $\psi_{t}^{(\theta'_0,\theta_1)}\relmont{\leq}{p.s.}\II_A(T)$)en montrant que l'intervalle 
\dli $I=\{h_{(\theta'_0,\theta_1)}(y)=k\}-A$ (resp. $I'=\{h'_{(\theta'_0,\theta_1)}(y)=k'\}-A$) est totalement ind\'etermin\'e.
Comme $I$ (resp. $I'$) est inclus dans $D_s^{\theta_1}$, la densit\'e $p_{\theta_1}$ est nulle sur cet intervalle ; $k$ (resp. $k'$) n'\'etant pas infini, cet intervalle est ind\'etermin\'e pour $(\theta'_0,\theta_1)$. Les densit\'es $p_{\theta'_0}$ et $p_{\theta_1}$ \'etant non s\'epar\'ees puisque 
$t\notin D_i^{\theta'_0}\cup D_s^{\theta_1}$, d'apr\`es la propri\'et\'e a) du 1\up{er} cas de la partie I les intervalles $I$ et $I'$, qui sont inclus dans $\IR\!-(D_i\cup D_s)$ ($D_i<t<I$, $k<+\infty$ et $D_i<t<I'$, $k'<+\infty$), sont totalement ind\'etermin\'es. 

ii) $t\in B$.
\dli $g_t^{(\theta_0,\theta_1)}$ (resp. $\psi_{t}^{(\theta_0,\theta_1)}$) est \'egale \`a $1-\II_C(T)$  ou $f_{(+\infty,0)}$.
\dli Nous avons vu en a-ii) du 2\up{\`eme} cas de la partie I   que $t$ appartient \`a $\IR\!-(D_i\cup D_s)$ et qu'il est ind\'etermin\'e pour $(\theta_0,\theta_1)$ ; mais $t$ n'est pas totalement ind\'etermin\'e puisque $t\notin N$ ; d'apr\`es la fin de la partie III de la d\'emonstration du lemme pr\'ec\'edent, il existe $\theta''\in\Theta_1$ et $\omega''\in T^{-1}(t)$ tels que $p_{\theta''}(\omega'')>0$.
On a donc $h_{(\theta_0,\theta'')}(t)=h'_{(\theta_0,\theta'')}(t)=0$, ce qui implique : $g_t^{(\theta_0,\theta'')}=\phi_{(0,1)}^{(\theta_0,\theta'')}$ et $\psi_t^{(\theta_0,\theta'')}=\phi_{(0,1)}^{'(\theta_0,\theta'')}$ (voir le 1\up{er} cas) ; mais pour $x\in C$ on a 
$h_{(\theta_0,\theta'')}(x)=h'_{(\theta_0,\theta'')}(x)>0$ (suivre le raisonnement fait en a-iii) du 2\up{\`eme} cas de la partie I avec $\theta_1=\theta''$ et $p_{\theta''}(\omega_y)\geq 0$), donc $g_t^{(\theta_0,\theta'')}\leq 1-\II_C(T)$ et
$\psi_t^{(\theta_0,\theta'')}\leq 1-\II_C(T)$.
\dli Nous venons ainsi de trouver un couple $(\theta_0,\theta'')$ tel que :
$g_t^{(\theta_0,\theta_1)}\geq\psi_t^{(\theta_0,\theta'')}$ et 
$\psi_t^{(\theta_0,\theta_1)}\geq g_t^{(\theta_0,\theta'')}$.

iii) $t\in C$.
\dli On a $g_t^{(\theta_0,\theta_1)}=\psi_{t}^{(\theta_0,\theta_1)}=f_{(+\infty,0)}$ ; $(\theta'_0,\theta'_1)=(\theta_0,\theta_1)$ convient.
\nobreak

\smallskip\centerline{\hbox to 3cm{\bf \hrulefill}}\par}

\vfill\eject

\centerline{\soustitre ANNEXE IV.}
\vskip 2cm
\bigskip
Soit $(\Omega ,{\cal A},(p_\theta.\mu)_{\theta\in\Theta})$ un mod\`ele statistique \`a rapport de vraisemblance monotone 
par rapport \`a la statistique r\'eelle $T$, $\Theta$ \'etant muni de la relation d'ordre totale : $\preceq$. Ceci suppose l'existence, pour $\theta'\!\prec\theta''$, d'une fonction croissante $h_{(\theta'',\theta')} :\, \IR \rightarrow \overline{\IR^+}$ v\'erifiant
$p_{\theta''}/p_{\theta'}=h_{(\theta'',\theta')}(T)$
sur le domaine de d\'efinition de ce rapport c'est-\`a-dire en dehors de
$\{\omega\in\Omega\, ;\, p_{\theta'}(\omega)=p_{\theta''}(\omega)=0\}$.

Nous dirons que la fonction $h_{(\theta'',\theta')}$, $\theta'\!\prec\theta''$, est normalis\'ee si elle est constante sur tout intervalle $I$ ind\'etermin\'e pour $(\theta',\theta'')$ : 
\dli $\forall\omega\in T^{-1}(I)\quad p_{\theta'}(\omega)=p_{\theta''}(\omega)=0$.
\dli Nous allons travailler avec une famille $\{h_{(\theta'',\theta')}\}_{\theta'\prec\theta''}$ de fonctions normalis\'ees ; ceci est toujours possible d'apr\`es la premi\`ere \'etape de l'annexe II ($\theta_0=\theta''$ et $\theta_1=\theta'$). Une telle famille permet de construire facilement une famille normalis\'ee $\{h^{\Theta_1}_{(\theta_0,\theta_1)}\}_{(\theta_0,\theta_1)\in\Theta_0\times\Theta_1}$ pour les hypoth\`eses unilat\'erales $\{\Theta_1,\Theta_0\}$ (voir les d\'efinitions 4.2.1 et 5.1.2). D'apr\`es la deuxi\`eme \'etape de l'annexe II, il suffit de tronquer les  $h_{(\theta_0,\theta_1)}$ de la fa\c con suivante :
$$h^{\Theta_1}_{(\theta_0,\theta_1)}(t)\quad=\quad\left\{\matrix{
0\hfill & si & t\in D^{\Theta_0}_i-D^{\Theta_1}_s\hfill\cr
h_{(\theta_0,\theta_1)}(t)\hfill & si & t\in\ \IR-(D^{\Theta_0}_i\cup D^{\Theta_1}_s)\hfill \cr
+\infty\hfill & si & t\in D^{\Theta_1}_s\hfill \cr
}\right. $$
avec $D^{\Theta_0}_i=\cap_{\theta\in\Theta_0}D_i^\theta$ et $D^{\Theta_1}_s=\cap_{\theta\in\Theta_1}D_s^\theta$, $D^\theta_i=]-\infty,t)$ (resp. $D^\theta_s=(t,+\infty[$) d\'esignant la plus grande demi-droite ouverte ou ferm\'ee pour laquelle $p_\theta$ est nulle sur $T^{-1}(D^\theta_i)$ (resp. $T^{-1}(D^\theta_s)$).

Nous allons d\'emontrer quelques propri\'et\'es des probl\`emes de choix entre deux hypoth\`eses simples $\theta'\in\Theta$ et $\theta''\in\Theta$, $\theta'\!\prec\theta''$.
\dli $\Phi^{(\theta'',\theta')}_s=\Bigl\{\phi^{(\theta'',\theta')}_{(0,1)}\, ;\,\{\phi^{(\theta'',\theta')}_{(k,\beta)}\}_{k\in\IR_*^+,\,\beta\in \{0,1\}}\, ;\,\phi^{(\theta'',\theta')}_{(\infty,0)}\Bigr\}$ est l'ensemble des fonctions de test simples construites \`a partir du rapport des densit\'es $K=h_{(\theta'',\theta')}(T)$ (voir la d\'efinition 2.2.1).

\medskip
{\bf Lemme 1}
\medskip
\medskip
\moveleft 10.4pt\hbox{\vrule\kern 10pt\vbox{\defpro
Soient $\theta'\!\prec\theta''$, on a :
\dli $\phi^{(\theta'',\theta')}_{(0,1)}\leq\II_{D_i^{\theta''}}(T)$ et 
$\phi^{(\theta'',\theta')}_{(\infty,0)}\geq 1-\II_{D_s^{\theta'}}(T)\quad  ;
\quad D_i^{\theta'}\subseteq D_i^{\theta''}$ et 
$D_s^{\theta''}\subseteq D_s^{\theta'}$.
}}\medskip

\medskip
{\leftskip=15mm \dli {\bf D\'emonstration}
\medskip
{\parindent=-10mm a) --- $\phi^{(\theta'',\theta')}_{(0,1)}\leq\II_{D_i^{\theta''}}(T)$.}

Nous devons d\'emontrer que $\{h_{(\theta'',\theta')}=0\}$ est inclus dans $D_i^{\theta''}$, puisque $\{\phi^{(\theta'',\theta')}_{(0,1)}=1\}=\{h_{(\theta'',\theta')}(T)=0\}$.
\dli Soit $t$ tel que $h_{(\theta'',\theta')}(t)=0$, pour tout $\omega$ de $T^{-1}(t)$ on a $p_{\theta''}(\omega)=0$, sinon on aurait $h_{(\theta'',\theta')}(t)=p_{\theta''}(\omega)/p_{\theta'}(\omega)>0$.
La demi-droite inf\'erieure $\{h_{(\theta'',\theta')}=0\}$ est donc bien incluse dans $D_i^{\theta''}$.

{\parindent=-10mm b) --- 
$\phi^{(\theta'',\theta')}_{(\infty,0)}\geq 1-\II_{D_s^{\theta'}}(T)$.}

Ceci est \'equivalent \`a $\{h_{(\theta'',\theta')}=+\infty\}\subseteq D_s^{\theta'}$, puisque :  
\dli $\{\phi^{(\theta'',\theta')}_{(\infty,0)}=0\}=\{h_{(\theta'',\theta')}(T)=+\infty\}$.
\dli Soit $t$ tel que $h_{(\theta'',\theta')}(t)=+\infty$, pour tout $\omega$ de $T^{-1}(t)$ on a $p_{\theta'}(\omega)=0$, sinon on aurait $h_{(\theta'',\theta')}(t)=p_{\theta''}(\omega)/p_{\theta'}(\omega)<+\infty$. La demi-droite sup\'erieure $\{h_{(\theta'',\theta')}=+\infty\}$ est donc bien incluse dans $D_s^{\theta'}$.

{\parindent=-10mm c) --- $D_i^{\theta'}\subseteq D_i^{\theta''}$.} 

Si on avait $I=D_i^{\theta'}-D_i^{\theta''}\not=\emptyset$, il existerait  $\omega\in T^{-1}(t)$, $t\in I$, tel que $p_{\theta''}(\omega)>0$ et $p_{\theta'}(\omega)=0$ ; on aurait $h_{(\theta'',\theta')}(t)=+\infty$ et $[t,+\infty[$ serait inclus dans $D_s^{\theta'}$ (voir b)) ; ceci est impossible car on ne peut pas avoir $D_i^{\theta'}\cup D_s^{\theta'}=\IR$.

{\parindent=-10mm d) --- $D_s^{\theta''}\subseteq D_s^{\theta'}$.}

La d\'emonstration est semblable \`a la pr\'ec\'edente.
Si on avait $I=D_s^{\theta''}-D_s^{\theta'}\not=\emptyset$, il existerait  $\omega\in T^{-1}(t)$, $t\in I$, tel que $p_{\theta'}(\omega)>0$ et $p_{\theta''}(\omega)=0$ ; on aurait $h_{(\theta'',\theta')}(t)=0$ et $]-\infty,t]$ serait inclus dans $D_i^{\theta''}$ (voir a)) ; ceci est impossible car on ne peut pas avoir $D_i^{\theta''}\cup D_s^{\theta''}=\IR$.

\medskip\centerline{\hbox to 3cm{\bf \hrulefill}}\par}
\vfill\eject

\medskip
{\bf Lemme 2}
\medskip
\medskip
\moveleft 10.4pt\hbox{\vrule\kern 10pt\vbox{\defpro
Soient $\theta_a\prec\theta_b\prec\theta_c$ et $t\in\IR$.
Posons $h_{(\theta_b,\theta_a)}(t)=k_1$, 
\dli $h_{(\theta_c,\theta_a)}(t)=k_2$ et $h_{(\theta_c,\theta_b)}(t)=k_3$ (les fonctions $h$ \'etant normalis\'ees).

1) S'il existe $\omega_t\in T^{-1}(t)$ tel que : $p_{\theta_a}(\omega_t)>0$ et $p_{\theta_b}(\omega_t)>0$ 
\dli (resp. $p_{\theta_b}(\omega_t)>0$ et $p_{\theta_c}(\omega_t)>0$), alors $k_2=k_3.k_1$.

2) Dans le premier cas : $p_{\theta_a}(\omega_t)>0$ et $p_{\theta_b}(\omega_t)>0$, on a
\dli $\phi^{(\theta_c,\theta_b)}_{(k_3,0)}\relmont{\leq}{p.s.}sup\{\phi^{(\theta_b,\theta_a)}_{(k_1,0)},\phi^{(\theta_c,\theta_a)}_{(k_2,0)},\II_{D_i^{\{\theta\succ\theta_a\}}}(T)\}$ et 
\dli $\phi^{(\theta_c,\theta_b)}_{(k_3,1)}\relmont{\geq}{p.s.}inf\{\phi^{(\theta_b,\theta_a)}_{(k_1,1)},\phi^{(\theta_c,\theta_a)}_{(k_2,1)},1-\II_{D_s^{\{\theta\prec\theta_b\}}}(T)\}$.

3) Dans le second cas : $p_{\theta_b}(\omega_t)>0$ et $p_{\theta_c}(\omega_t)>0$, on a
\dli $\phi^{(\theta_b,\theta_a)}_{(k_1,0)}\relmont{\leq}{p.s.}sup\{\phi^{(\theta_c,\theta_a)}_{(k_2,0)},\phi^{(\theta_c,\theta_b)}_{(k_3,0)},\II_{D_i^{\{\theta\succ\theta_b\}}}(T)\}$ et 
\dli $\phi^{(\theta_b,\theta_a)}_{(k_1,1)}\relmont{\geq}{p.s.}inf\{\phi^{(\theta_c,\theta_a)}_{(k_2,1)},\phi^{(\theta_c,\theta_b)}_{(k_3,1)},1-\II_{D_s^{\{\theta\prec\theta_c\}}}(T)\}$.

Ces \'ecritures peuvent contenir des fonctions de test de la forme $\phi^{(\theta'',\theta')}_{(0,0)}$ (resp. $\phi^{(\theta'',\theta')}_{(\infty,1)}$), elles repr\'esentent $\II_\emptyset$ (resp. $\II_\Omega$).
}}\medskip

\medskip
{\leftskip=15mm \dli {\bf D\'emonstration}

On consid\`ere les trois intervalles de $\IR$ contenant $t$ d\'efinis par :
$I_1=\{h_{(\theta_b,\theta_a)}=k_1\}$, $I_2=\{h_{(\theta_c,\theta_a)}=k_2\}$ et 
$I_3=\{h_{(\theta_c,\theta_b)}=k_3\}$.

\medskip
{\parindent=-10mm 1\up{er} cas : $\exists\,\omega_t\in T^{-1}(t)$ tel que  $p_{\theta_a}(\omega_t)>0$ et $p_{\theta_b}(\omega_t)>0$.}

\dli $k_1=h_{(\theta_b,\theta_a)}(t)=h_{(\theta_b,\theta_a)}(T(\omega_t))=p_{\theta_b}(\omega_t)/p_{\theta_a}(\omega_t)\in ]0,+\infty[$.
\dli Que $p_{\theta_c}(\omega_t)$ soit nul ou pas, on a toujours \dli  $p_{\theta_c}(\omega_t)/p_{\theta_a}(\omega_t)=[p_{\theta_c}(\omega_t)/p_{\theta_b}(\omega_t)].[p_{\theta_b}(\omega_t)/p_{\theta_a}(\omega_t)]$, donc $k_2=k_3.k_1$ (ce qui d\'emontre le 1\up{er} cas de la propri\'et\'e 1) du lemmme).

\medskip
{\parindent=-5mm a) --- D\'emonstration de  $\phi^{(\theta_c,\theta_b)}_{(k_3,0)}\relmont{\leq}{p.s.}sup\{\phi^{(\theta_b,\theta_a)}_{(k_1,0)},\phi^{(\theta_c,\theta_a)}_{(k_2,0)},\II_{D_i^{\{\theta\succ\theta_a\}}}(T)\}$.}

\dli Lorsque $\{x<I_3\}\cap I_1=\emptyset$ ou $\{x<I_3\}\cap I_2=\emptyset$ ou $\{x<I_3\}\subseteq D_i^{\{\theta\succ\theta_a\}}$, l'in\'egalit\'e recherch\'ee est bien v\'erifi\'ee puisque l'on a alors respectivement 
$\phi^{(\theta_c,\theta_b)}_{(k_3,0)}\leq\phi^{(\theta_b,\theta_a)}_{(k_1,0)}$ ou  $\phi^{(\theta_c,\theta_b)}_{(k_3,0)}\leq\phi^{(\theta_c,\theta_a)}_{(k_2,0)}$ ou  $\phi^{(\theta_c,\theta_b)}_{(k_3,0)}\leq\II_{D_i^{\{\theta\succ\theta_a\}}}(T)$.

Dans le cas contraire, on a $\{x<I_3\}\cap I_1\not=\emptyset$, $\{x<I_3\}\cap I_2\not=\emptyset$ et $\{x<I_3\}\supset D_i^{\{\theta\succ\theta_a\}}$ ; ce qui est \'equivalent \`a 
$I=(\{x<I_3\}\cap I_1\cap I_2)-D_i^{\{\theta\succ\theta_a\}}\not=\emptyset$ puisque l'intervalle $I_1\cap I_2\cap I_3$ n'est pas vide (il contient $t$).

\dli L'intervalle $I$ est \'egal \`a $\{x<I_3\}\cap I_1$ ou $\{x<I_3\}\cap I_2$ ou 
\dli $\{x<I_3\}- D_i^{\{\theta\succ\theta_a\}}$, pour d\'emontrer l'in\'egalit\'e recherch\'ee il suffit de montrer que $I$ est totalement ind\'etermin\'e :
$\forall \omega\in T^{-1}(I)$ $\forall \theta\in \Theta$ $p_\theta(\omega)=0$.

Soit $x\in I$.
\dli i) Commen\c cons par montrer que $x$ est ind\'etermin\'e pour $(\theta_b,\theta_a)$.
\dli Supposons qu'il existe $\omega_x\in T^{-1}(x)$ tel que  $p_{\theta_b}(\omega_x)>0$, comme $x\in I_1$ et $k_1<+\infty$ on aurait aussi $p_{\theta_a}(\omega_x)>0$ donc $k_2=h_{(\theta_c,\theta_b)}(x).k_1$ (d'apr\`es le 1\up{er} cas de la propri\'et\'e 1) du lemme et le fait que $x$ appartient \`a $I_1\cap I_2$) ; ce qui est impossible car $k_2=k_3.k_1$ et $h_{(\theta_c,\theta_b)}(x)<k_3$ puisque $x<I_3$.
\dli La densit\'e $p_{\theta_b}$ est donc nulle sur $T^{-1}(x)$ ; comme 
$h_{(\theta_b,\theta_a)}(x)=k_1>0$ la densit\'e $p_{\theta_a}$ est aussi nulle sur $T^{-1}(x)$.

\dli ii) $x$ est totalement ind\'etermin\'e.
\dli Nous allons appliquer la propri\'et\'e 2) du lemme 1 de l'annexe III \`a $x$ ind\'etermin\'e pour $(\theta_b,\theta_a)$, les hypoth\`eses unilat\'erales \'etant d\'efinies par $\Theta_1=\{\theta\preceq\theta_a\}$ et $\Theta_0=\{\theta\succ\theta_a\}$.
\dli Montrons d'abord que $x$ est un \'el\'ement de $\IR\!-(D_i^{\Theta_0}\cup D_s^{\Theta_1})$ ; par d\'efinition il n'appartient pas \`a $D_i^{\{\theta\succ\theta_a\}}=D_i^{\Theta_0}$ ; d'apr\`es le lemme 1 de cette annexe on a $D_s^{\Theta_1}=D_s^{\theta_a}$, $x$ n'appartient donc pas non plus \`a $D_s^{\Theta_1}$ puisque $x<t\in I_3$ et $p_{\theta_a}(\omega_t)>0$.
\dli Notons $I_x$ le plus grand intervalle de $\IR\!-(D_i^{\Theta_0}\cup D_s^{\Theta_1})$ contenant $x$ et ind\'etermin\'e pour $(\theta_b,\theta_a)$. Nous allons montrer que $h^{\Theta_1}_{(\theta_b,\theta_a)}(I_x^+)=inf\{h^{\Theta_1}_{(\theta_b,\theta_a)}(y)\, ;\, y>I_x\}\cup\{+\infty\}$ n'est pas infini, la partie 2) du lemme 1 de l'annexe III implique alors la propri\'et\'e recherch\'ee. $t$ n'\'etant pas ind\'etermin\'e pour $(\theta_b,\theta_a)$, on a $t>I_x$ et $t\notin D_s^{\Theta_1}=D_s^{\theta_a}$ donc \dli
$h^{\Theta_1}_{(\theta_b,\theta_a)}(I_x^+)\leq h^{\Theta_1}_{(\theta_b,\theta_a)}(t)\leq h_{(\theta_b,\theta_a)}(t)=k_1<+\infty$.

\medskip
{\parindent=-5mm b) --- D\'emonstration de 
$\phi^{(\theta_c,\theta_b)}_{(k_3,1)}\relmont{\geq}{p.s.}inf\{\phi^{(\theta_b,\theta_a)}_{(k_1,1)},\phi^{(\theta_c,\theta_a)}_{(k_2,1)},1-\II_{D_s^{\{\theta\prec\theta_b\}}}(T)\}$.}

\dli Lorsque $\{x>I_3\}\cap I_1=\emptyset$ ou $\{x>I_3\}\cap I_2=\emptyset$ ou $\{x>I_3\}\subseteq D_s^{\{\theta\prec\theta_b\}}$, l'in\'egalit\'e recherch\'ee est bien v\'erifi\'ee puisque l'on a alors respectivement 
$\phi^{(\theta_c,\theta_b)}_{(k_3,1)}\geq\phi^{(\theta_b,\theta_a)}_{(k_1,1)}$ ou  $\phi^{(\theta_c,\theta_b)}_{(k_3,1)}\geq\phi^{(\theta_c,\theta_a)}_{(k_2,1)}$ ou  $\phi^{(\theta_c,\theta_b)}_{(k_3,1)}\geq 1-\II_{D_s^{\{\theta\prec\theta_b\}}}(T)$.

Dans le cas contraire, on a $\{x>I_3\}\cap I_1\not=\emptyset$, $\{x>I_3\}\cap I_2\not=\emptyset$ et $\{x>I_3\}\supset D_s^{\{\theta\prec\theta_b\}}$ ; ce qui est \'equivalent \`a 
$J=(\{x>I_3\}\cap I_1\cap I_2)-D_s^{\{\theta\prec\theta_b\}}\not=\emptyset$ puisque l'intervalle $I_1\cap I_2\cap I_3$ n'est pas vide (il contient $t$).

\dli L'intervalle $J$ est \'egal \`a $\{x>I_3\}\cap I_1$ ou $\{x>I_3\}\cap I_2$ ou 
\dli $\{x>I_3\}- D_s^{\{\theta\prec\theta_b\}}$, pour d\'emontrer l'in\'egalit\'e recherch\'ee il suffit de montrer que $J$ est totalement ind\'etermin\'e :
$\forall \omega\in T^{-1}(J)$ $\forall \theta\in \Theta$ $p_\theta(\omega)=0$.

Soit $x\in J$.
\dli i) Commen\c cons par montrer que $x$ est ind\'etermin\'e pour $(\theta_b,\theta_a)$.
\dli C'est le m\^eme raisonnement qu'en a)i), la contradiction venant du fait que $h_{(\theta_c,\theta_b)}(x)>k_3$, puisque $x>I_3$.

\dli ii) $x$ est totalement ind\'etermin\'e.
\dli Nous allons  appliquer la propri\'et\'e 2) du lemme 1 de l'annexe III \`a $x$ ind\'etermin\'e pour $(\theta_b,\theta_a)$, les hypoth\`eses unilat\'erales \'etant d\'efinies par $\Theta_1=\{\theta\prec\theta_b\}$ et $\Theta_0=\{\theta\succeq\theta_b\}$.
\dli Montrons d'abord que $x$ est un \'el\'ement de $\IR\!-(D_i^{\Theta_0}\cup D_s^{\Theta_1})$ ; par d\'efinition il n'appartient pas \`a $D_s^{\{\theta\prec\theta_b\}}=D_s^{\Theta_1}$ ; d'apr\`es le lemme 1 de cette annexe on a $D_i^{\Theta_0}=D_i^{\theta_b}$, $x$ n'appartient donc pas non plus \`a $D_i^{\Theta_0}$ puisque $x>t\in I_3$ et $p_{\theta_b}(\omega_t)>0$.
\dli Notons $I_x$ le plus grand intervalle de $\IR\!-(D_i^{\Theta_0}\cup D_s^{\Theta_1})$ contenant $x$ et ind\'etermin\'e pour $(\theta_b,\theta_a)$. Nous allons montrer que $h^{\Theta_1}_{(\theta_b,\theta_a)}(I_x^-)=sup\{h^{\Theta_1}_{(\theta_b,\theta_a)}(y)\, ;\, y<I_x\}\cup\{0\}$ n'est pas nul, la partie 2) du lemme 1 de l'annexe III implique alors la propri\'et\'e recherch\'ee. $t$ n'\'etant pas ind\'etermin\'e pour $(\theta_b,\theta_a)$, on a $t<I_x$ et $t\notin D_i^{\Theta_0}=D_i^{\theta_b}$ donc 
\dli $h^{\Theta_1}_{(\theta_b,\theta_a)}(I_x^-)\geq h^{\Theta_1}_{(\theta_b,\theta_a)}(t)\geq h_{(\theta_b,\theta_a)}(t)=k_1>0$.

\medskip
{\parindent=-10mm 2\up{\`eme} cas : $\exists\,\omega_t\in T^{-1}(t)$ tel que $p_{\theta_b}(\omega_t)>0$ et $p_{\theta_c}(\omega_t)>0$.}

\dli La d\'emonstration est semblable \`a celle du  1\up{er} cas.
\dli  $k_3=h_{(\theta_c,\theta_b)}(t)=h_{(\theta_c,\theta_b)}(T(\omega_t))=p_{\theta_c}(\omega_t)/p_{\theta_b}(\omega_t)\in ]0,+\infty[$.
\dli Que $p_{\theta_a}(\omega_t)$ soit nul ou pas, on a toujours 
\dli  $p_{\theta_c}(\omega_t)/p_{\theta_a}(\omega_t)=[p_{\theta_c}(\omega_t)/p_{\theta_b}(\omega_t)].[p_{\theta_b}(\omega_t)/p_{\theta_a}(\omega_t)]$, donc $k_2=k_3.k_1$ (ce qui d\'emontre le 2\up{\`eme} cas de la propri\'et\'e 1) du lemmme).

\medskip
{\parindent=-5mm a) --- D\'emonstration de $\phi^{(\theta_b,\theta_a)}_{(k_1,0)}\relmont{\leq}{p.s.}sup\{\phi^{(\theta_c,\theta_a)}_{(k_2,0)},\phi^{(\theta_c,\theta_b)}_{(k_3,0)},\II_{D_i^{\{\theta\succ\theta_b\}}}(T)\}$.}

\dli Lorsque $\{x<I_1\}\cap I_2=\emptyset$ ou $\{x<I_1\}\cap I_3=\emptyset$ ou $\{x<I_1\}\subseteq D_i^{\{\theta\succ\theta_b\}}$, l'in\'egalit\'e recherch\'ee est bien v\'erifi\'ee puisque l'on a alors respectivement 
$\phi^{(\theta_b,\theta_a)}_{(k_1,0)}\leq\phi^{(\theta_c,\theta_a)}_{(k_2,0)}$ ou  $\phi^{(\theta_b,\theta_a)}_{(k_1,0)}\leq\phi^{(\theta_c,\theta_b)}_{(k_3,0)}$ ou  $\phi^{(\theta_b,\theta_a)}_{(k_1,0)}\leq\II_{D_i^{\{\theta\succ\theta_b\}}}(T)$.

Dans le cas contraire, on a $\{x<I_1\}\cap I_2\not=\emptyset$, $\{x<I_1\}\cap I_3\not=\emptyset$ et $\{x<I_1\}\supset D_i^{\{\theta\succ\theta_b\}}$ ; ce qui est \'equivalent \`a 
$I=(\{x<I_1\}\cap I_2\cap I_3)-D_i^{\{\theta\succ\theta_b\}}\not=\emptyset$ puisque l'intervalle $I_1\cap I_2\cap I_3$ n'est pas vide (il contient $t$).

\dli L'intervalle $I$ est \'egal \`a $\{x<I_1\}\cap I_2$ ou $\{x<I_1\}\cap I_3$ ou 
\dli $\{x<I_1\}- D_i^{\{\theta\succ\theta_b\}}$, pour d\'emontrer l'in\'egalit\'e recherch\'ee il suffit de montrer que $I$ est totalement ind\'etermin\'e :
$\forall \omega\in T^{-1}(I)$ $\forall \theta\in \Theta$ $p_\theta(\omega)=0$.

Soit $x\in I$.
\dli i) Commen\c cons par montrer que $x$ est ind\'etermin\'e pour $(\theta_c,\theta_b)$.
\dli Supposons qu'il existe $\omega_x\in T^{-1}(x)$ tel que  $p_{\theta_b}(\omega_x)>0$, comme $x\in I_3$ et $k_3>0$ on aurait aussi $p_{\theta_c}(\omega_x)>0$ donc $k_2=k_3.h_{(\theta_b,\theta_a)}(x)$ (d'apr\`es le 2\up{\`eme} cas de la propri\'et\'e 1) du lemme et le fait que $x$ appartient \`a $I_2\cap I_3$) ; ce qui est impossible car $k_2=k_3.k_1$ et $h_{(\theta_b,\theta_a)}(x)<k_1$ puisque $x<I_1$.
\dli La densit\'e $p_{\theta_b}$ est donc nulle sur $T^{-1}(x)$ ; comme 
$h_{(\theta_c,\theta_b)}(x)=k_3<+\infty$ la densit\'e $p_{\theta_c}$ est aussi nulle sur $T^{-1}(x)$.

\dli ii) $x$ est totalement ind\'etermin\'e.
\dli La d\'emonstration est celle de la partie a)ii) du 1\up{er} cas en rempla\c cant $\theta_b$ par $\theta_c$, $\theta_a$ par $\theta_b$, $I_3$ par $I_1$ et $k_1$ par $k_3$. 

\medskip
{\parindent=-5mm b) --- D\'emonstration de $\phi^{(\theta_b,\theta_a)}_{(k_1,1)}\relmont{\geq}{p.s.}inf\{\phi^{(\theta_c,\theta_a)}_{(k_2,1)},\phi^{(\theta_c,\theta_b)}_{(k_3,1)},1-\II_{D_s^{\{\theta\prec\theta_c\}}}(T)\}$.}

\dli Lorsque $\{x>I_1\}\cap I_2=\emptyset$ ou $\{x>I_1\}\cap I_3=\emptyset$ ou $\{x>I_1\}\subseteq D_s^{\{\theta\prec\theta_c\}}$, l'in\'egalit\'e recherch\'ee est bien v\'erifi\'ee puisque l'on a alors respectivement 
$\phi^{(\theta_b,\theta_a)}_{(k_1,1)}\geq\phi^{(\theta_c,\theta_a)}_{(k_2,1)}$ ou  $\phi^{(\theta_b,\theta_a)}_{(k_1,1)}\geq\phi^{(\theta_c,\theta_b)}_{(k_3,1)}$ ou  $\phi^{(\theta_b,\theta_a)}_{(k_1,1)}\geq 1-\II_{D_s^{\{\theta\prec\theta_c\}}}(T)$.

Dans le cas contraire, on a $\{x>I_1\}\cap I_2\not=\emptyset$, $\{x>I_1\}\cap I_3\not=\emptyset$ et $\{x>I_1\}\supset D_s^{\{\theta\prec\theta_c\}}$ ; ce qui est \'equivalent \`a 
$J=(\{x>I_1\}\cap I_2\cap I_3)-D_s^{\{\theta\prec\theta_c\}}\not=\emptyset$ puisque l'intervalle $I_1\cap I_2\cap I_3$ n'est pas vide (il contient $t$).

\dli L'intervalle $J$ est \'egal \`a $\{x>I_1\}\cap I_2$ ou $\{x>I_1\}\cap I_3$ ou 
\dli $\{x>I_1\}- D_s^{\{\theta\prec\theta_c\}}$, pour d\'emontrer l'in\'egalit\'e recherch\'ee il suffit de montrer que $J$ est totalement ind\'etermin\'e :
$\forall \omega\in T^{-1}(J)$ $\forall \theta\in \Theta$ $p_\theta(\omega)=0$.

Soit $x\in J$.
\dli i) Commen\c cons par montrer que $x$ est ind\'etermin\'e pour $(\theta_c,\theta_b)$.
\dli C'est le m\^eme raisonnement que celui de la partie a)i) du 2\up{\`eme} cas, la contradiction venant du fait que $h_{(\theta_b,\theta_a)}(x)>k_1$, puisque $x>I_1$.

\dli ii) $x$ est totalement ind\'etermin\'e.
\dli La d\'emonstration est celle de la partie b)ii) du 1\up{er} cas en rempla\c cant $\theta_b$ par $\theta_c$, $\theta_a$ par $\theta_b$, $I_3$ par $I_1$ et $k_1$ par $k_3$. 

\medskip\centerline{\hbox to 3cm{\bf \hrulefill}}\par}

\vfill\eject

\null\vglue 3cm
\centerline{\soustitre BIBLIOGRAPHIE}
\vglue 2cm
\nom{Bar.}{J.R. Barra}
\refl{Notions fondamentales de statistique math\'ematique}
{Dunod, Paris, 1971}
\nom{Ber.}{J.O. Berger}
\refl{Statistical decision theory and Bayesian analysis (second edition)}{\dli Springer-Verlag, New York, 1985}
\nom{BerD}{J.O. Berger, M. Delampady}
\refa{Testing precise hypotheses}{Statist. Science}{2, p. 317-352, 1987}
\nom{Bor.}{A. Borovkov}
\refl{Statistique math\'ematique}{Mir, Moscou, 1987}
\nom{Bre.}{L. Breiman}
\refl{Probability}{Addison-Wesley, Reading, Massachusetts, 1968}
\nom{DacD}{D. Dacunha-castelle, M. Duflo}
\refl{Probabilit\'es et statistiques. Tome 1 : probl\`emes \`a temps fixe (2\up{\`eme} \'edition)}{Masson, Paris, 1994}
\nom{DhaM}{I.D. Dhariyal, N. Misra, R.K.S. Rathore}
\refa{Selecting the better of two binomial populations : optimal decision rules}{Calcutta Statist. Assoc. Bull.}{38, p. 157-167, 1989}
\nom{Dic.}{J.M. Dickey}
\refa{Three multidimensional-integral identities with bayesian applications}{Ann. Math. Statist.}{39, p. 1615-1628, 1968}
\nom{Gei.}{S. Geisser}
\refa{On prior distributions for binary trials}{American Statist.}{38, p. 244-251, 1984}
\nom{HenT}{P.L. Hennequin, A. Tortrat}
\refl{Th\'eorie des probabilit\'es et quelques applications}{Masson, Paris, 1965}
\nom{HwaC}{J.T. Hwang, G. Casella, C. Robert, M.T. Wells, R.H. Farrell}
\refa{Estimation of accuracy in testing}{Ann. Statist.}{20, p. 490-509, 1992}
\nom{Kar.}{S. Karlin}
\refa{Decision theory for P\'olya type distributions. Case of two actions, I}
{Proc. Third Berkeley Symposium on Math. Statist. and Prob.}{Vol 1, Univ. of Calif. Press, Berkeley, p. 115-128, 1955}
\nom{KarR}{S. Karlin, H. Rubin}
\refa{The theory of decision procedures for distributions with monotone likelihood ratio}
{Ann. Math. Statist.}{27, p. 272-299, 1956}
\nom{KroM}{A.H. Kroese, E.A. van der Meulen, K. Poortema, \-W. \-Schaafsma}
\refa{Distributional inference}{Statistica Neerlandica}{49, p. 63-82, 1995}
\nom{Leh.}{E.L. Lehmann}
\refl{Testing statistical hypotheses (second edition)}{Wiley, New York, 1986}
\nom{Mon.}{A. Monfort}
\refl{[1] Cours de probabilit\'es}{Economica, Paris, 1980}
\refl{[2] Cours de statistique math\'ematique}{Economica, Paris, 1982}
\nom{Mor.}{G. Morel}
\refl{[1] Proc\'edures statistiques pour espace de d\'ecisions totalement ordonn\'e et famille de lois \`a vraisemblance monotone (th\`ese)}
{Universit\'e de Rouen, France, 1987}
\refa{[2] D\'ecisions li\'ees aux intervalles d'une partition : le probl\`eme du choix dans la pratique de la recherche}
{Pub. Inst. Stat. Univ.}{XXXII, fasc. 1-2, p. 93-111, 1987}
\nom{Ney.}{J. Neyman, E.S. Pearson}
\refa{On the testing of statistical hypotheses in relation to probability a priori}{Proc. Cambridge Phil. Soc.}{29, p. 492-510, 1933}
\nom{Nev.}{J. Neveu}
\refl{Calcul des probabilit\'es}{Masson, Paris, 1970}
\nom{Nik.}{M.S. Nikulin.}
\refa{[1] On a result of L.N. Bol'shev from the theory of the statistical testing of hypotheses}
{Zap. Nauchn. Sem. Leningr. Otd. Mat. Inst.}{153, p. 129-137, 1986}
\refa{[2] Estimation of the efficiency of Bol'shev's decision rule in the problem of distinguishing of two hypotheses}
{J. Soviet Math.}{52, p.2955-2964, 1990}
\nom{Pfa.}{J. Pfanzagl}
\refa{A technical lemma for monotone likelihood ratio families}
{Ann. Math. Statist.}{38, p. 611-613, 1967}
\nom{R\'en.}{A. R\'enyi}
\refl{Calcul des probabilit\'es}{Dunod, Paris, 1966}
\nom{Reu.}{M. Reuchlin}
\refa{Epreuves d'hypoth\`eses nulles et inf\'erence fiduciaire en psychologie}{J. de Psychologie}{3, p. 277-292, 1977}
\nom{Rob.}{C. Robert}
\refl{L'analyse statistique bay\'esienne}{Economica, Paris, 1992}
\nom{Rou.}{R.D. Routledge}
\refa{Practicing safe statistics with the mid-{\it p*}}
{Canad. J. Statist.}{22, p. 103-110, 1994}
\nom{SchT}{W. Schaafsma, J. Tolboom, B. van der Meulen}
\refa{Discussing truth or falsity by computing a Q-value}{Statistical Data Analysis and Inference}{Dodge, North-Holland, p. 85-100, 1989}
\nom{Sch.}{H. Scheff\'e}
\refl{The analysis of variance(6\up{\`eme} \'edition)}{Wiley, New York, 1970}
\nom{Ste.}{W.L. Stevens}
\refa{Shorter intervals for the parameter of the binomial and poisson distributions}{Biometrika}{44, p. 436-440, 1957}
\nom{Wan.}{C. Wang}
\refl{Sense and nonsense of statistical inference}{Dekker, New York, 1993}
\vfill\eject

\null\vglue 0cm
\centerline{\soustitre TABLE des D\'EFINITIONS}
\vglue 1cm
\dli 1--INTRODUCTION.

\dli 2--CHOIX ENTRE DEUX PROBABILIT\'ES.
\dli\hskip 2cm D\'efinition 2.1.1 p. 11
\dli\hskip 2cm D\'efinition 2.2.1 p. 13
\dli\hskip 2cm D\'efinition 2.4.1 p. 21
\dli\hskip 2cm D\'efinition 2.5.1 p. 23

\dli 3--R\`EGLES DE D\'ECISION DE BOL'SHEV.

\dli\hskip 2cm D\'efinition 3.1.1 p. 25
\dli\hskip 2cm D\'efinition 3.1.2 p. 26
\dli\hskip 2cm D\'efinition 3.2.1 p. 31

\dli 4--CHOIX ENTRE DEUX HYPOTH\`ESES STABLES.

\dli\hskip 2cm D\'efinition 4.1.1 p. 38
\dli\hskip 2cm D\'efinition 4.1.2 p. 39
\dli\hskip 2cm D\'efinition 4.2.1 p. 42
\dli\hskip 2cm D\'efinition 4.3.1 p. 57
\dli\hskip 2cm D\'efinition 4.3.2 p. 64

\dli 5--MOD\`ELES \`A RAPPORT DE VRAISEMBLANCE MONOTONE.

\dli\hskip 2cm D\'efinition 5.1.1 p. 66
\dli\hskip 2cm D\'efinition 5.1.2 p. 67
\dli\hskip 2cm D\'efinition 5.2.1 p. 78
\dli\hskip 2cm D\'efinition 5.3.1 p. 106

\dli 6--HYPOTH\`ESES STABLES ET PARAM\`ETRES FANT\^OMES.

\dli\hskip 2cm D\'efinition 6.1.1 p. 131
\dli\hskip 2cm D\'efinition 6.1.2 p. 133

\vfill\eject

\null\vglue 0cm
\centerline{\soustitre TABLE des PROPOSITIONS}
\vglue 1cm
\dli 1--INTRODUCTION.

\dli 2--CHOIX ENTRE DEUX PROBABILIT\'ES.
\dli\hskip 2cm Proposition 2.3.1 p. 16

\dli 3--R\`EGLES DE D\'ECISION DE BOL'SHEV.

\dli\hskip 2cm Proposition 3.1.1 p. 26
\dli\hskip 2cm Proposition 3.1.2 p. 28
\dli\hskip 2cm Proposition 3.2.1 p. 32
\dli\hskip 2cm Proposition 3.2.2 p. 34

\dli 4--CHOIX ENTRE DEUX HYPOTH\`ESES STABLES.

\dli\hskip 2cm Proposition 4.2.1 p. 42
\dli\hskip 2cm Proposition 4.2.2 p. 47
\dli\hskip 2cm Proposition 4.3.1 p. 58
\dli\hskip 2cm Proposition 4.3.2 p. 62

\dli 5--MOD\`ELES \`A RAPPORT DE VRAISEMBLANCE MONOTONE.

\dli\hskip 2cm Proposition 5.1.1 p. 68
\dli\hskip 2cm Proposition 5.2.1 p. 79
\dli\hskip 2cm Proposition 5.2.2 p. 81
\dli\hskip 2cm Proposition 5.2.3 p. 91
\dli\hskip 2cm Proposition 5.2.4 p. 100
\dli\hskip 2cm Proposition 5.3.1 p. 107
\dli\hskip 2cm Proposition 5.3.2 p. 108

\dli 6--HYPOTH\`ESES STABLES ET PARAM\`ETRES FANT\^OMES.

\vfill\eject

\dli ANNEXE I.

\dli ANNEXE II.

\dli\hskip 2cm Lemme 1 p. 146
\dli\hskip 2cm Lemme 2 p. 147

\dli ANNEXE III.

\dli\hskip 2cm Lemme 1 p. 149
\dli\hskip 2cm Proposition 1 p. 153

\dli ANNEXE IV.
\dli\hskip 2cm Lemme 1 p. 162
\dli\hskip 2cm Lemme 2 p. 163

\bye